\documentclass[11pt]{ucthesis}
\usepackage{amssymb, amsmath, euscript, latexsym, amsfonts}
\usepackage{fancyhdr} 
\usepackage[all, knot]{xy}
\usepackage{diagram}

\author{Alissa S. Crans}
\title{Lie 2-Algebras}
\degreeyear{2004} \degree{Doctor of Philosophy} \chair{Professor
John Baez} \othermembers{Professor Vyjayanthi Chari \\Professor
Xiao-Song Lin}
\prevdegrees{B.S. (University of Redlands) 1999\\
             M.S. (University of California, Riverside) 2000}
\field{Mathematics}
\campus{Riverside}

\input{epsf}
\newcommand{\et}{\hspace{-0.08in}{\bf .}\hspace{0.1in}}
\newcommand{\BOX}{\hbox {$\sqcap$ \kern -1em $\sqcup$}}
\newcommand{\qed}{\hskip 2em \hbox{\BOX} \vskip 2ex}

\newcommand{\ttimes}{\! \times \!}
\newcommand{\To}{\Rightarrow}
\newcommand{\Set}{{\rm Set}}
\newcommand{\FinSet}{{\rm FinSet}}
\newcommand{\Cat}{{\rm Cat}}
\newcommand{\LieGrp}{{\rm LieGrp}}
\newcommand{\LieAlg}{{\rm LieAlg}}

\newcommand{\Grp}{{\rm Grp}}
\newcommand{\Quand}{{\rm Quand}}
\newcommand{\Diff}{{\rm Diff}}

\newcommand{\Vect}{{\rm Vect}}
\newcommand{\Term}{{\rm Term}}

\renewcommand{\to}{\rightarrow}
\newcommand{\tensor}{\otimes}
\newcommand{\maps}{\colon}
\newcommand{\inv}{{\rm inv}}

\newcommand{\iso}{\cong}

\newcommand{\tr}{{\rm tr}}
\newcommand{\ad}{{\rm ad}}
\renewcommand{\deg}{{\rm deg}}

\newcommand{\id}{{\rm id}}

\newcommand{\g}{{\mathfrak g}}

\newcommand{\Co}{{\mathcal C}}

\newcommand{\R}{{\mathbb R}}

\newcommand{\N}{{\mathbb N}}

        \newcommand{\be}{\begin{equation}}
        \newcommand{\ee}{\end{equation}}
        \newcommand{\ba}{\begin{eqnarray}}
        \newcommand{\ea}{\end{eqnarray}}
        \newcommand{\ban}{\begin{eqnarray*}}
        \newcommand{\ean}{\end{eqnarray*}}
        \newcommand{\barr}{\begin{array}}
        \newcommand{\earr}{\end{array}}

\renewcommand{\ltimes}{\! \otimes \!}

\newtheorem{thm}{Theorem}    
\newtheorem{cor}[thm]{Corollary}
\newtheorem{lem}[thm]{Lemma}

\newtheorem{prop}[thm]{Proposition}
\newtheorem{defn}[thm]{Definition}
\newtheorem{example}[thm]{Example}

\usepackage{makeidx}

\let\cal\mathcal

\def\N{\mathbb N}

\def\R{\mathbb R}

\def\cplus{\hbox{$\supset${\raise1.05pt\hbox{\kern -0.55em
${\scriptscriptstyle +}$}}\ }}

\def\dsp{\def\baselinestretch{2.0}\large\normalsize}
\dsp

\def\dsp{\def\baselinestretch{2.0}\large\normalsize}
\dsp

\fancypagestyle{plain}{ \lhead{} \chead{} \rhead{} \lfoot{}
\fancyfoot[C]{\vspace{-1pt}\thepage}

}

\begin{document}

\thispagestyle{empty} \

\ssp

\begin{center}\

UNIVERSITY OF CALIFORNIA\\ RIVERSIDE\\ \vspace{.7in} Lie $2$-Algebras\\
\vspace{.7in}

A Dissertation submitted in partial satisfaction\\
of the requirements for the degree
of\\\vspace{.2in} Doctor of Philosophy\\[12pt] in\\[12pt] Mathematics\\[12pt]
by\\[12pt]
Alissa Susan Crans\\[12pt] August 2004

\end{center}

\vspace{1.5in} \ssp

\noindent Dissertation Committee:\\
\indent Professor John Baez, Chairperson\\
\indent Professor Vyjayanthi Chari\\
\indent Professor Xiao-Song Lin

\newpage

\ssp \thispagestyle{empty}
\ \vspace{7.38in} \begin{center}Copyright by\\Alissa Susan Crans\\
2004\end{center}

\newpage

\renewcommand{\thepage}{\roman{page}}

\begin{acknowledgements}
I would like to express my sincere gratitude to my advisor, John
Baez, for his patience, guidance, enthusiasm, encouragement,
inspiration, and humor. Additionally, I am indebted to Vyjayanthi
Chari, James Dolan, and Xiao-Song Lin for their willingness to
share their knowledge with me numerous times during my graduate
studies.  I also thank my `mathematical brothers': Miguel
Carri\'on-\'Alvarez, Toby Bartels, Jeffrey Morton, and Derek Wise
for their friendship and engaging, educational conversations. I am
grateful for the assistance of Aaron Lauda in drawing various
braid diagrams, and thank Ronnie Brown, Andr\'ee Ehresmann, Thomas
Larsson, James Stasheff, J. \ Scott Carter, and Masahico Saito for
helpful discussions and correspondence. Finally, I am extremely
appreciative of the love and support of my family, friends, and
former professors during my time as a graduate student. I
certainly could not have accomplished all that I have without
them.
\end{acknowledgements}

\newpage

\thispagestyle{fancy}

\begin{center}

ABSTRACT OF THE DISSERTATION\\\vspace{.3in}

Lie $2$-Algebras\\[12pt]

by\\[12pt] Alissa Susan Crans\\\vspace{.3in} Doctor of Philosophy, Graduate
Program in Mathematics\\
University of California, Riverside, August 2004\\
Professor John Baez, Chairperson

\end{center}

\ssp \vspace{.7in}

\noindent We categorify the theory of Lie algebras beginning with
a new notion of categorified vector space, or `2-vector space',
which we define as an internal category in $\Vect$, the category
of vector spaces.  We then define a `semistrict Lie 2-algebra' to
be a 2-vector space $L$ equipped with a skew-symmetric bilinear
functor $[\cdot,\cdot] \maps L \times L \to L$ satisfying the
Jacobi identity up to a completely antisymmetric trilinear natural
transformation called the `Jacobiator', which in turn must satisfy
a certain law of its own. We construct a 2-category of semistrict
Lie 2-algebras and show that it is 2-equivalent to the 2-category
of `2-term $L_\infty$-algebras' in the sense of Stasheff.  We also
investigate strict and skeletal Lie 2-algebras.  We show how to
obtain the strict ones from Lie 2-groups and we use the skeletal
ones to classify Lie 2-algebras in terms of 3rd cohomology classes
in Lie algebra cohomology. This classification allows us to
construct for any finite-dimensional Lie algebra $\g$ a canonical
1-parameter family of Lie 2-algebras $\g_\hbar$ which reduces to
$\g$ at $\hbar = 0$.

We then explore the relationship between Lie algebras and
algebraic structures called `quandles'.  A quandle is a set $Q$
equipped with two binary operations \newline $\rhd \maps Q \times
Q \to Q$ and $\lhd \maps Q \times Q \to Q$ satisfying axioms that
capture the essential properties of the operations of conjugation
in a group and algebraically encode the three Reidemeister moves.
Indeed, we describe the relation to groups and show that quandles
give invariants of braids.  We further show that both Lie algebras
and quandles give solutions of the Yang--Baxter equation, and
explain how conjugation plays a prominent role in the both the
theories of Lie algebras and quandles.  Inspired by these
commonalities, we provide a novel, conceptual passage from Lie
groups to Lie algebras using the language of quandles.  Moreover,
we propose relationships between higher Lie theory and
higher-dimensional braid theory.  We conclude with evidence of
this connection by proving that any semistrict Lie 2-algebra gives
a solution of the Zamolodchikov tetrahedron equation, which is the
higher-dimensional analog of the Yang--Baxter equation.
\newpage

\thispagestyle{fancy} \tableofcontents

\newpage

\renewcommand{\thepage}{\arabic{page}}

\setcounter{page}{1}

\chapter{Introduction} \label{intro}
\pagenumbering{arabic} \setcounter{page}{1} \ssp
Higher-dimensional algebra is the study of generalizations of
algebraic concepts obtained through a process called
`categorification'.  In the mid-1990's, Crane \cite{Crane, CF}
coined this term to refer to the process of developing
category-theoretic analogs of set-theoretic concepts. In this
process we replace elements with objects, sets with categories,
and functions with functors.  We replace equations between
elements by isomorphisms between objects, and replace equations
between functions by natural isomorphisms between functors.
Finally, we require that these isomorphisms satisfy equations of
their own, called coherence laws. Finding the correct coherence
laws is often the most difficult aspect of this generalization
process. For example, the category $\FinSet$ of finite sets
together with functions is a categorification of the set of
natural numbers. The functions of sum and product in $\N$ are
replaced by the functors disjoint union and Cartesian product of
finite sets. The equational laws satisfied by addition and
multiplication in $\N$, such as commutativity, associativity, and
distributivity, now hold for disjoint union and Cartesian product,
but only up to natural isomorphism.  For instance, the associative
law $$(xy)z = x(yz)$$ is replaced by a natural isomorphism called
the \emph{associator}:
\[ a_{x,y,z} \maps \xymatrix@1{(xy)z \; \ar[r]^<<<<<{\sim} & \; x
(yz)},
\]
which we then require to satisfy a coherence law known as the
pentagon identity:
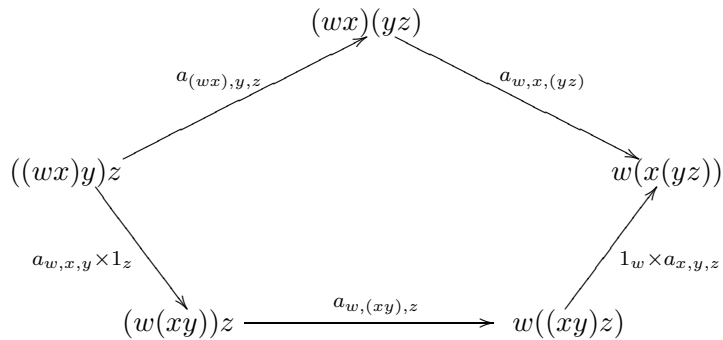
\[
\xy
 (0,20)*{(w  x)  (y  z)}="1";
 (40,0)*{w  (x  (y  z))}="2";
 (25,-20)*{ \quad w  ((x  y)  z)}="3";
 (-25,-20)*{(w  (x  y))  z}="4";
 (-40,0)*{((w  x)  y)  z}="5";
     {\ar^{a_{w,x,(yz)}}     "1";"2"}
     {\ar_{1_w \times  a _{x,y,z}}  "3";"2"}
     {\ar^{a _{w,(xy),z}}    "4";"3"}
     {\ar_{a _{w,x,y} \times  1_z}  "5";"4"}
     {\ar^{a _{(wx),y,z}}    "5";"1"}
\endxy
\\
\]
\vskip 1em
\noindent Ultimately, by iterating this process,
mathematicians wish to obtain and apply the $n$-categorical
generalizations of as many mathematical concepts as possible to
strengthen and simplify the connections between different
subfields of mathematics.

Perhaps the greatest strength of categorification is that it
allows us to refine our concept of `sameness' by enabling us to
distinguish between equality and isomorphism.  In a set, two
elements are either the same or different, while in a category,
two objects can be `the same in a way' while remaining different.
That is, they can be isomorphic, but not equal.  Better still, we
are able to explicitly keep track of how two objects are the same:
the isomorphism itself.  Moreover, two objects can be the same in
multiple ways, since there can be various different isomorphisms
between them.  An object can even be the same as itself in
multiple ways:  this leads to the concept of \emph{symmetry},
since the automorphisms of an object form a group, its `symmetry
group'.  This more careful consideration of the notion of sameness
is the reason that categorification plays an increasingly
important role not only in mathematics, but also in physics and
computer science, where symmetry plays a significant role and a
precise treatment of the notion of sameness is crucial.

As an example, consider the problem of determining whether two
sets are isomorphic.  We know that we can solve this problem
easily by simply counting the number of elements in each set.  By
comparing the resulting numbers, we will know that our sets are
isomorphic, but we will not know {\it how} they are isomorphic. In
examples, there will often be a specific isomorphism demonstrating
that the two sets are isomorphic.  Categorification is an attempt
to preserve this type of information.

At various stages in the process of categorifying, we are forced
to make certain decisions about the laws governing our objects and
morphisms.  On the one hand, we can insist that these laws hold
`on the nose', that is, that they hold {\it strictly} as
equational laws.  In this approach, the issue of coherence laws
does not arise.  On the other hand, we can impose these laws only
{\it up to isomorphism,} with our isomorphisms then satisfying
certain coherence laws.  We refer to this second method as
`weakening'. Categories constructed in this second way tend to be
more appealing because they typically arise naturally in
applications.

Typically when categorifying, we are attempting to find a
category-theoretic analog of some algebraic structure defined in
the category of sets, so we seek to create a hybrid of two
notions: a category and that algebraic structure. Here a technique
invented by Ehresmann \cite{E}, called `internalization', can
serve as a useful first step.

This process allows one to take concepts that were defined in the
category $\Set$ and transport them to other categories.  For
example, the concepts of `group' and `(small) category' live in
the world of sets, but they can actually live in other categories
as well! Internalizing a concept consists of first expressing it
completely in terms of commutative diagrams and then interpreting
those diagrams in some ambient category, $K$.  In order to define
an `internal group' or `group in $K$,' $K$ must have finite
products, while the definition of an `internal category' or
`category in $K$' requires that $K$ have finite limits.  The
notion of internalization tends to be somewhat confusing the first
time around, so it is a beneficial exercise to convince yourself
that when $K = \Set$, these reduce to the usual definitions of a
group and category.

A quick and easy way to categorify an algebraic structure is to
internalize it taking $K$ to be $\Cat$, the category whose objects
are (small) categories and whose morphisms are functors.
Unfortunately, this method does not weaken the concepts being
categorified: equations between functions are promoted to
equations between functors, not natural isomorphisms, and thus
coherence laws are not introduced. So, while we make extensive use
of internalization in this work, we still must then weaken various
concepts --- notably the concept of a Lie algebra --- in a
somewhat {\it ad hoc} way.



One should not be tempted to think that categorification is simply
generalization for its own sake!  For example, group theory serves
as a powerful tool in all branches of science where symmetry plays
a role.  But, in many situations where we are tempted to use
groups, it is actually more natural to use a richer sort of
structure, where in addition to group elements describing
symmetries, we also have isomorphisms between these, describing
{\it symmetries between symmetries}. This sort of structure is a
categorified group, or `$2$-group'.  The theory of $2$-groups, or
`higher-dimensional group theory', dates back to the mid $1940$'s
and has a complicated history, which is briefly sketched in the
introduction to Baez and Lauda's review article on $2$-groups
\cite{BLau}:  {\it Higher Dimensional Algebra V: $2$-groups},
henceforth denoted as HDA5.

The first sort of categorified group is a {\it strict} $2$-group
--- one in which all the laws hold on the nose as equations. As one
might expect, strict $2$-groups can be defined exactly as one
would hope:  using internalization! They are precisely groups in
$\Cat$, so they are often called `categorical groups'.  They have
numerous applications to homotopy theory \cite{Bro,BS},
topological quantum field theory \cite{Yetter}, nonabelian
cohomology \cite{Breen,Breen2,Giraud}, the theory of nonabelian
gerbes \cite{Breen2,BrM}, categorified gauge field theory
\cite{Attal,Baez,GP,Pfeiffer}, and even quantum gravity
\cite{CS,CY}.  However, the strict notion is not the best for all
applications, so a weaker concept of `coherent $2$-group' has also
been introduced --- originally by Sinh \cite{S}, who used the term
`$gr$-category'.  A careful introduction to coherent $2$-groups
can be found in HDA5 \cite{BLau}.


In recent years, thanks to the work of Ronald Brown \cite{Brown2},
$2$-groups have secured their place in the mathematical, and more
recently, physical, literature as useful and important objects. As
an example, the concept of `Lie $2$-group' becomes especially
important in applications of $2$-groups to geometry and physics,
in particular, to gauge theory.  A Lie $2$-group roughly amounts
to a $2$-group in which the set of objects and the set of
morphisms are manifolds, and all relevant maps are smooth.  Until
recently, only strict Lie $2$-groups had been defined \cite{Baez}.
However, in HDA5 \cite{BLau}, the concept of a coherent Lie
$2$-group was introduced, which we hope will contribute to the
study of generalized gauge theories.  Moreover, we suspect that
just as Lie groups and Lie algebras arise wherever differential
geometry and symmetry appear, the same will be true for Lie
$2$-groups and `Lie $2$-algebras'.  That is, Lie $2$-algebras will
ideally contribute to geometry and physics much in the same way as
Lie $2$-groups, and will hopefully have applications in situations
where the Jacobi identity need not hold on the nose.

The goal of this work is to develop and explore the theory of
categorified Lie algebras, or `Lie $2$-algebras'.  Just as every
Lie algebra has an underlying vector space, a Lie $2$-algebra will
be a `$2$-vector space' equipped with extra structure.  A
$2$-vector space blends together the notions of category and
vector space, so, roughly speaking, it is a category where
everything is linear.  More precisely, we use internalization to
define a $2$-vector space to be a category in $\Vect$, the
category of vector spaces.   Then, to obtain the notion of a Lie
$2$-algebra, we start with a $2$-vector space and equip it with a
weakened version of the structure of a Lie algebra.  That is,
instead of a bracket function, a Lie $2$-algebra has a bracket
{\it functor}, which we require to be `linear'.  We then weaken
the Jacobi identity so that now it holds only {\it up to a
`linear' natural isomorphism}, which we call the `Jacobiator'.
Following the Tao of weakening, we require that the Jacobiator
satisfy an identity of its own, which we call the `Jacobiator
identity'.


The theory of semistrict Lie $2$-algebras is investigated in
Chapter \ref{ch1} in a presentation parallel to the development of
the theory of $2$-groups in HDA5 \cite{BLau}. In particular, we
show that just as coherent $2$-groups were classified up to
equivalence using group cohomology, semistrict Lie $2$-algebras
can be classified up to equivalence using Lie algebra cohomology.
This classification allows us to construct a $1$-parameter family
of Lie $2$-algebras $\g_{\hbar}$ deforming any finite-dimensional
Lie algebra, which are especially interesting when $\g$ is
semisimple, and may be related to affine Lie algebras and quantum
groups.  Furthermore, we demonstrate the relationship between
semistrict Lie $2$-algebras and special $L_{\infty}$-algebras, or
sh-Lie algebras, which are generalizations of Lie algebras defined
by Stasheff \cite{SS} that blend together the notions of Lie
algebra and chain complex.  Much of the content of this chapter
has already appeared in a separate paper coauthored with John Baez
\cite{BC}, {\it Higher Dimensional Algebra VI: Lie $2$-algebras},
henceforth denoted as HDA6.

In Chapter \ref{ch2} we describe a modern approach to obtaining
the Lie algebra of a Lie group, a result that holds a prominent
place in differential geometry, using algebraic structures known
as `quandles'.  We begin with definitions and examples of quandles
and continue by illustrating their relationship to group theory
and topology, in particular, to the three Reidemeister moves.  In
fact, we show that quandles and Lie algebras give solutions to the
Yang--Baxter equation.  We also demonstrate how quandles give
invariants of braids.

Our road map for the description of this new, conceptual, passage
of Lie group to Lie algebra takes the following form:

$$
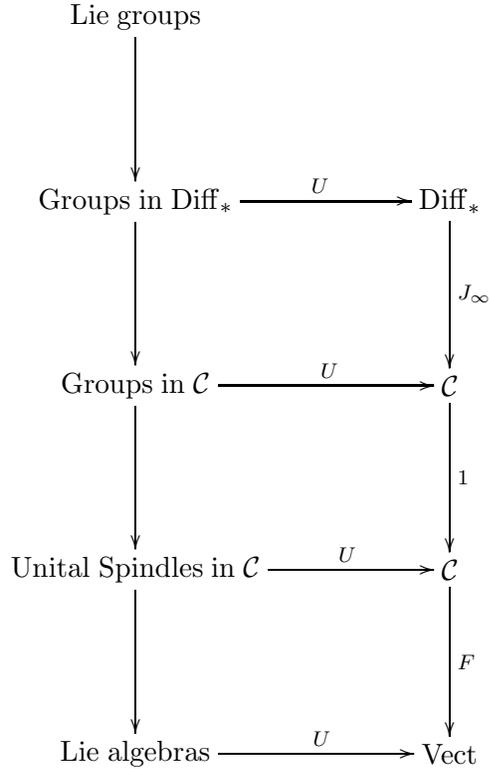
\xymatrix{
   {\rm Lie \; groups}
   \ar[dd] \\ \\
   {\rm Groups \; in \; \Diff_\ast}
   \ar[rr]^<<<<<<<<<<<{U}
   \ar[dd]
   && \Diff_{\ast}
   \ar[dd]^{J_{\infty}} \\ \\
   {\rm Groups \; in \; \mathcal{C}}
   \ar[rr]^<<<<<<<<<<<<<<<{U}
   \ar[dd]
   && \mathcal{C}
   \ar[dd]^{1} \\ \\
   {\rm Unital \; Spindles \; in \; \mathcal{C}}
   \ar[rr]^<<<<<<<<<<<{U}
   \ar[dd]
   && \mathcal{C}
   \ar[dd]^{F} \\ \\
   {\rm Lie \; algebras}
   \ar[rr]^<<<<<<<<<<<<<<{U}
   && \Vect
}$$ We explain this diagram thoroughly by recalling the notion of
a `jet', introducing the notion of a `special coalgebra', and
describing their roles in this process.  This chapter is a warmup
for a future attempt to obtain the Lie $2$-algebra of a Lie
$2$-group.  In the final section of this chapter, we outline what
we have done, and what remains to be done, in order to achieve
this `higher' result.

These two chapters are united by the fact that both Lie algebras
and quandles give solutions of the Yang--Baxter equation.
Furthermore, the theory of conjugation plays a prominent role in
both the theories of Lie algebras and quandles, which we describe
in detail in Chapter \ref{ch2}.  These observations inspired our
novel description of the passage from a Lie group to its Lie
algebra. In both chapters, we emphasize the relationship between
algebra and topology, which motivates the relationship between
higher-dimensional algebra and higher-dimensional topology.  While
we have only briefly sketched the contents of these chapters here,
more details can be found at the beginning of Chapter \ref{ch1}
and Chapter \ref{ch2}.

We conclude in Section \ref{conclusions} by outlining our future
goals related to this work, most importantly being the task of
showing that every Lie $2$-group has a Lie $2$-algebra.
Furthermore, we offer our guesses for how the relationship between
Lie theory and braids should generate a relationship between
higher Lie theory and higher-dimensional braid theory.  In
particular, as evidence of this relationship, we prove that just
as a Lie algebra gives a solution of the Yang--Baxter equation, or
third Reidemeister move, a Lie $2$-algebra gives a solution of the
Zamolodchikov tetrahedron equation, which is the
higher-dimensional analog of the Yang--Baxter equation.  This
result, while beautiful in its own right, also suggests that the
process of passing from Lie groups to Lie algebras given in
Section \ref{liealgofliegrp} may categorify, providing us with the
outcome we desire.

\chapter{Semistrict Lie 2-algebras} \label{ch1}
Just as the concept of Lie group can be categorified to obtain
various concepts of `Lie $2$-group', one can categorify the
concept of Lie algebra to obtain various concepts of `Lie
$2$-algebra'.  As mentioned in the Introduction, a Lie $2$-group,
roughly speaking, is a categorified group where everything in
sight is smooth.  Since we are going to categorify the notion of a
Lie algebra in approximately the same way, we begin by recalling
the manner in which we categorify a group.

We begin by using the technique of internalization, which we
recalled in the Introduction.  That is, given a group $G$, we
replace the underlying {\it set} with a {\it category},  and the
multiplication {\it function} $m \maps G \times G \to G$ with a
multiplication {\it functor} $m \maps G \otimes G \to G$.  And
this is where we stop if we want to describe the notion of a {\it
strict} $2$-group.  Recall from the Introduction that when we
define a strict concept, all laws hold strictly as equations.
Thus, the issue of coherence laws does not arise.  However, if we
choose to {\it weaken} our $2$-group, we impose the group laws
only {\it up to isomorphism}, and then require these isomorphisms
to satisfy certain laws of their own.  Thus, to weaken the concept
of $2$-group, we replace the {\it equation} expressing the
associative law by an {\it isomorphism} called the `associator':
\[ a_{x,y,z} \maps \xymatrix@1{(x \otimes y) \otimes z \; \ar[r]^<<<<<{\sim} & \;
x \otimes (y \otimes z)},
\]
which we then require to satisfy a coherence law known as the
`pentagon identity':
\[
\xy
 (0,20)*{(w \otimes x) \otimes (y \otimes z)}="1";
 (40,0)*{w \otimes (x \otimes (y \otimes z))}="2";
 (25,-20)*{ \quad w \otimes ((x \otimes y) \otimes z)}="3";
 (-25,-20)*{(w \otimes (x \otimes y)) \otimes  z}="4";
 (-40,0)*{((w \otimes  x)  \otimes y) \otimes  z}="5";
     {\ar^{a_{w,x,(y \otimes z)}}     "1";"2"}
     {\ar_{1_w \otimes  a _{x,y,z}}  "3";"2"}
     {\ar^{a _{w,(x \otimes y),z}}    "4";"3"}
     {\ar_{a _{w,x,y} \otimes  1_z}  "5";"4"}
     {\ar^{a _{(w \otimes x),y,z}}    "5";"1"}
\endxy
\\
\] \vskip 1em
\noindent Similarly, we replace the equations expressing the left
and right unit laws
\[        1 \tensor x = x, \qquad x \tensor 1 = x \]
by isomorphisms
\[
\ell_x \maps \xymatrix@1@=16pt{1 \tensor x \; \ar[r]^>>>>{\sim} &
\; x} , \qquad
  r_x \maps \xymatrix@1@=16pt{x \tensor 1\; \ar[r]^>>>>{\sim} &\; x}
\]
To obtain the notion of a `coherent' $2$-group, we weaken further
by replacing the equations
\[ x \tensor x^{-1} = 1, \qquad x^{-1} \tensor x = 1 \]
by isomorphisms called the `unit' and `counit'.  Thus, instead of
an inverse in the strict sense, the object $x$ only has a
specified `weak inverse'.

We now consider categorifying the notion of a Lie algebra.  In the
theory of Lie algebras, the analog to the associative law is the
Jacobi identity.  Therefore, in analogy to our discussion above,
to obtain a `Lie $2$-algebra, we will replace the Jacobi identity
by an isomorphism that we call the {\it Jacobiator}, which, we
will see, satisfies an interesting new law of its own.  Just as
the pentagon equation can be traced back to work done by Stasheff,
we will see that the coherence law for the Jacobiator is related
to his work on $L _\infty$-algebras, also known as strongly
homotopy Lie algebras \cite{LS,SS}. This demonstrates yet again
the close connection between categorification and homotopy theory.

Since our goal is to categorify a Lie algebra, which will be a
mixture of a Lie algebra and category, we first consider the
notion of categorified vector space.  Just as a Lie algebra has an
underlying vector space, a Lie $2$-algebra will have an underlying
`$2$-vector space'.  A $2$-vector space will blend together the
concept of a category with that of a vector space.  Thus, we begin
in Section \ref{intcats} by reviewing the theory of internal
categories, which we will use to create such a concept.  In
Section \ref{2vs} we focus specifically on categories in $\Vect$,
the category of vector spaces. We boldly call these `$2$-vector
spaces', despite the fact that this term is already used to refer
to a very different categorification of the concept of vector
space \cite{KV}, for it is our contention that our $2$-vector
spaces lead to a more interesting version of categorified linear
algebra than the traditional ones.  For example, the tangent space
at the identity of a Lie $2$-group is a $2$-vector space of our
sort, and this gives a canonical representation of the Lie
$2$-group: its `adjoint representation'. This is contrast to the
phenomenon observed by Barrett and Mackaay \cite{BM}, namely that
Lie $2$-groups have few interesting representations on the
traditional sort of $2$-vector space.  One reason for the
difference is that the traditional $2$-vector spaces do not have a
way to `subtract' objects, while ours do.  This will be especially
important for finding examples of Lie $2$-algebras, since we often
wish to set $[x,y] = xy - yx$.

At this point we should admit that our $2$-vector spaces are far
from novel entities!  In fact, a category in $\Vect$ is secretly
just the same as a $2$-term chain complex of vector spaces.  While
the idea behind this correspondence goes back to Grothendieck
\cite{G}, and is by now well-known to category-theorists, we
describe it carefully in Proposition \ref{1-1vs}, because it is
crucial for relating `categorified linear algebra' to more
familiar ideas from homological algebra.

In Section \ref{definitions} we introduce the key concept of
`semistrict Lie $2$-algebra'.  Roughly speaking, this is a
$2$-vector space $L$ equipped with a bilinear functor
\[          [\cdot,\cdot] \maps L \times L \to L ,\]
the Lie bracket, that is skew-symmetric and satisfies the Jacobi
identity up to a completely antisymmetric trilinear natural
isomorphism, the `Jacobiator' --- which in turn is required to
satisfy a law of its own, the `Jacobiator identity'. Since we do
not weaken the equation $[x,y] = -[y,x]$ to an isomorphism, we do
not reach the more general concept of `weak Lie $2$-algebra': this
remains a task for the future, which we describe in further detail
in Section \ref{conclusions} of the next chapter.

In Section \ref{Linftyalgs}, we recall the definition of an
$L_{\infty}$-algebra.  Briefly, this is a chain complex $V$ of
vector spaces equipped with a bilinear skew-symmetric operation
\newline $[\cdot,\cdot] \maps V \times V \to V$ which satisfies the
Jacobi identity up to an infinite tower of chain homotopies.  We
construct a $2$-category of `$2$-term' $L_\infty$-algebras, that
is, those with $V_i = \{0\}$ except for $i = 0,1$. Finally, we
show this $2$-category is equivalent to the previously defined
$2$-category of semistrict Lie $2$-algebras.

In the next two sections we study \emph{strict} and
\emph{skeletal} Lie $2$-algebras, the former being those where the
Jacobi identity holds `on the nose', while in the latter,
isomorphisms exist only between identical objects. Section
\ref{strictlie2algs} consists of an introduction to strict Lie
$2$-algebras and strict Lie $2$-groups, together with the process
for obtaining the strict Lie $2$-algebra of a strict Lie
$2$-group.  Section \ref{skeletallie2algs} begins with an
exposition of Lie algebra cohomology and its relationship to
skeletal Lie $2$-algebras.  We then show that Lie $2$-algebras can
be classified (up to equivalence) in terms of a Lie algebra
$\mathfrak{g}$, a representation of $\mathfrak{g}$ on a vector
space $V$, and an element of the Lie algebra cohomology group
$H^3(\mathfrak{g},V)$.  With the help of this result, we construct
from any finite-dimensional Lie algebra $\g$ a canonical
1-parameter family of Lie $2$-algebras $\g_\hbar$ which reduces to
$\g$ at $\hbar = 0$. This is a new way of deforming a Lie algebra,
in which the Jacobi identity is weakened in a manner that depends
on the parameter $\hbar$.  It is natural to speculate that this
deformation is somehow related to the theory of quantum groups and
affine Lie algebras.  However, we have only a little evidence for
this speculation at present.

{\bf Note:} In all that follows, we denote the composite of
morphisms $f\maps x \rightarrow y$ and $g\maps y \rightarrow z$ as
$fg\maps x \rightarrow z.$   All $2$-categories and $2$-functors
referred to are {\it strict}, though sometimes we include the word
`strict' to emphasize this fact.  We denote vertical composition
of $2$-morphisms by juxtaposition; we denote horizontal
composition and whiskering by the symbol $\circ$.


\section{Internal Categories} \label{intcats}

In order to create a hybrid of the notions of a vector space and a
category in the next section, we will use the technique of
internalization to blend together these concepts.  That is, we
need the concept of an `internal category', also called a
`category object', within some category. The idea is that given a
category $K$, we obtain the definition of a `category internal to
$K$', which we call `category in $K$' for short, by expressing the
definition of a usual (small) category completely in terms of
commutative diagrams and then interpreting those diagrams within
$K$.  The same idea allows us to define functors and natural
transformations in $K$, at least if $K$ has properties
sufficiently resembling those of the category of sets.

Internal categories were introduced by Ehresmann \cite{E} in the
1960s, and by now they are a standard part of category theory
\cite{Bo}. However, since not all readers may be familiar with
them, for the sake of a self-contained treatment we start with the
basic definitions.

\begin{defn} \et \label{co}  Let $K$ be a category.
A {\bf category internal to $K$} or {\bf category in $K$}, say
$X$, consists of:
\begin{itemize}
\item
an {\bf object of objects} $X_{0} \in K,$
\item
an {\bf object of morphisms} $X_{1} \in K,$
\end{itemize}
together with
\begin{itemize}
\item
{\bf source} and {\bf target} morphisms $s,t \maps X_{1}
\rightarrow X_{0},$
\item
a {\bf identity-assigning} morphism $i \maps X_{0} \rightarrow
X_{1},$
\item
a {\bf composition} morphism $\circ \maps X_{1} \times _{X_{0}}
X_{1} \rightarrow X_{1}$
\end{itemize}
such that the following diagrams commute, expressing the usual
category laws:
\begin{itemize}
\item laws specifying the source and target of identity morphisms:
\[
\xymatrix{
    X_{0}
      \ar[r]^{i}
      \ar[dr]_{1}
      & X_{1}
      \ar[d]^{s} \\
     & X_{0} }
 \hspace{.2in}
\xymatrix{
      X_{0}
      \ar[r]^{i}
      \ar[dr]_{1}
      & X_{1}
      \ar[d]^{t} \\
     & X_{0}}
\]
\item
laws specifying the source and target of composite morphisms:
\[
  \xymatrix{
   X_{1} \times _{X_{0}} X_{1}
     \ar[rr]^{\circ}
     \ar[dd]_{p_{1}}
     && X_{1}
     \ar[dd]^{s} \\ \\
   X_{1}
     \ar[rr]^{s}
     && X_{0} }
     \hspace{.2in}
\xymatrix{
  X_{1} \times_{X_{0}} X_{1}
     \ar[rr]^{\circ}
     \ar[dd]_{p_{2}}
      && X_{1}
     \ar[dd]^{t} \\ \\
      X_{1}
     \ar[rr]^{t}
      && X_{0} }
\]
\item the associative law for composition of morphisms:
\[
\xymatrix{
   X_{1} \times _{X_{0}} X_{1} \times _{X_{0}} X_{1}
     \ar[rr]^{\circ \times_{X_{0}} 1}
     \ar[dd]_{1 \times_{X_{0}} \circ}
      && X_{1} \times_{X_{0}} X_{1}
     \ar[dd]^{\circ} \\ \\
      X_{1} \times _{X_{0}} X_{1}
     \ar[rr]^{\circ}
      && X_{1} }
\]
\item the left and right unit laws for composition of morphisms:
\[
\xymatrix{
   X_{0} \times _{X_{0}} X_{1}
     \ar[r]^{i \times 1}
     \ar[ddr]_{p_2}
      & X_{1} \times _{X_{0}} X_{1}
     \ar[dd]^{\circ}
      & X_{1} \times_{X_{0}} X_{0}
     \ar[l]_{1 \times i}
     \ar[ddl]^{p_1} \\ \\
      & X_{1} }
\]
\end{itemize}
\end{defn}

The pullbacks referred to in the above definition should be clear
from the usual definition of category; for instance, composition
is defined on pairs of morphisms where the target of the first is
the source of the second, so the pullback $X_1 \times_{X_0} X_1$
is defined via the square
\[
\begin{diagram}[X_1 \times_{X_0} X_1]
\node{X_1 \times_{X_0} X_1} \arrow{e,t}{p_1} \arrow{s,l}{p_2}
\node{X_1} \arrow{s,r}{s} \\
\node{X_1} \arrow{e,t}{t} \node{X_0}
\end{diagram}
\]
Notice that inherent to this definition is the assumption that the
pullbacks involved actually exist.  This holds automatically when
the `ambient category' $K$ has finite limits, but there are some
important examples such as $K = \Diff$ where this is not the case.
Throughout this work, all of the categories considered have finite
limits:

\begin{itemize}
\item Set, \emph{the category whose objects are sets and whose
morphisms are functions.}

\item Vect, \emph{the category whose objects are vector spaces
over the field $k$ and whose morphisms are linear functions.}

\item Grp, \emph{the category
whose objects are groups and whose morphisms are homomorphisms.}

\item Cat, \emph{the category whose objects are small categories and whose
morphisms are functors.}

\item LieGrp, \emph{the category whose objects are Lie groups and whose
morphisms are Lie group homomorphisms.}

\item LieAlg, \emph{the category whose objects are Lie algebras over
the field $k$ and whose morphisms are Lie algebra homomorphisms.}
\end{itemize}

Having defined `categories in $K$', we can now internalize the
notions of functor and natural transformation in a similar manner
to obtain functors and natural transformations in $K$. We shall
use these to construct a $2$-category $K\Cat$ consisting of
categories, functors, and natural transformations in $K$.

\begin{defn} \et \label{cofunctor} Let $K$ be a category.
Given categories $X$ and $X'$ in $K$, an {\bf internal functor} or
{\bf functor in $K$} between them, say $F\maps X \rightarrow X',$
consists of:
\begin{itemize}
\item
a morphism $F_{0} \maps X_{0} \to X_{0}'$,
\item
a morphism $F_{1} \maps X_{1} \rightarrow X_{1}'$
\end{itemize}
such that the following diagrams commute, corresponding to the
usual laws satisfied by a functor:
\begin{itemize}
\item preservation of source and target:
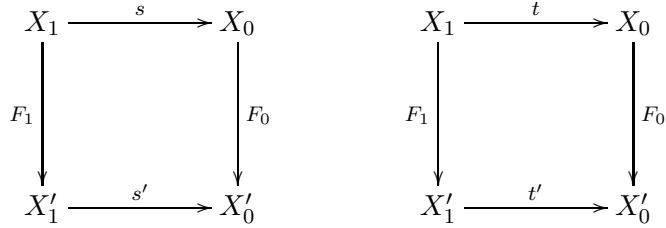
\[
\xymatrix{
  X_{1}
   \ar[rr]^{s}
   \ar[dd]_{F_{1}}
    && X_{0}
   \ar[dd]^{F_{0}} \\ \\
    X_{1}'
   \ar[rr]^{s'}
    && X_{0}' }
\qquad \qquad \xymatrix{
  X_{1}
   \ar[rr]^{t}
   \ar[dd]_{F_{1}}
    && X_{0}
   \ar[dd]^{F_{0}} \\ \\
    X_{1}'
   \ar[rr]^{t'}
    && X_{0}' }
\]
\item preservation of identity morphisms:
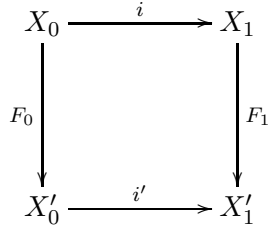
\[
\xymatrix{
    X_{0}
   \ar[rr]^{i}
   \ar[dd]_{F_{0}}
    && X_{1}
   \ar[dd]^{F_{1}} \\ \\
    X_{0}'
   \ar[rr]^{i'}
    && X_{1}' }
\]
\item preservation of composite morphisms:
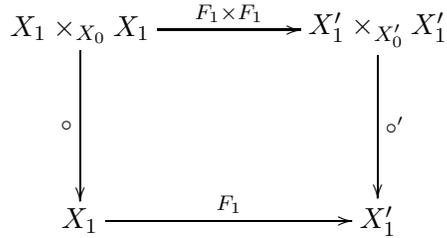
\[
\xymatrix{
   X_{1} \times _{X_{0}} X_{1}
    \ar[rr]^{F_{1} \times F_{1}}
    \ar[dd]_{\circ}
     && X_{1}' \times_{X_{0}'} X_{1}'
    \ar[dd]^{\circ '} \\ \\
     X_{1}
    \ar[rr]^{F_{1}}
     && X_{1}' }
\]
\end{itemize}
\end{defn}

Given two functors $F\maps X \rightarrow X'$ and $G\maps X'
\rightarrow X''$ in some category $K$, we define their composite
$FG\maps X \rightarrow X''$ by taking $(FG)_{0} = F_{0}G_{0}$ and
$(FG)_{1} = F_{1}G_{1}$. Similarly, we define the identity functor
in $K$, $1_X \maps X \rightarrow X$, by taking $(1_X)_0 = 1_{X_0}$
and $(1_X)_1 = 1_{X_1}$.  Showing that the composite $FG$
satisfies the diagrams of Definition \ref{cofunctor} is a
straightforward computation.  We now consider morphisms between
these functors in $K$:

\begin{defn} \et \label{conattrans} Let $K$ be a
category.  Given two functors $F,G\maps X \rightarrow X'$ in $K$,
an {\bf internal natural transformation} or {\bf natural
transformation in $K$} between them, say $\theta \maps F
\Rightarrow G$, is a morphism $\theta \maps X_0 \to X'_1$ for
which the following diagrams commute, expressing the usual laws
satisfied by a natural transformation:
\begin{itemize}
\item laws specifying the source and target of a natural transformation:
\[
\begin{diagram}
\node{X_0} \arrow{se,t}{F_{0}} \arrow{s,l}{\theta} \\
\node{X'_1}  \arrow{e,t}{s} \node{X'_0}
\end{diagram}
\qquad \qquad
\begin{diagram}
\node{X_0} \arrow{se,t}{G_{0}} \arrow{s,l}{\theta} \\
\node{X'_1}  \arrow{e,t}{t} \node{X'_0}
\end{diagram}
\]
\item the commutative square law:
\[  \xymatrix{
   X_1
    \ar[rr]^{\Delta (s\theta \times G)}
    \ar[dd]_{\Delta (F \times t\theta)}
     && X'_1 \times_{X_0'} X'_1
    \ar[dd]^{\circ'} \\ \\
     X'_1 \times_{X_0'} X'_1
    \ar[rr]^{\circ'}
     && X'_1
}
\]
\end{itemize}
\end{defn}

Just like ordinary natural transformations, natural
transformations in $K$ may be composed in two different, but
commuting, ways.  First, let $X$ and $X'$ be categories in $K$ and
let $F,G,H\maps X \rightarrow X'$ be functors in $K$.  If $\theta
\maps F \Rightarrow G$ and $\tau \maps G \Rightarrow H$ are
natural transformations in $K$, we define their {\bf vertical}
composite, $\theta \tau \maps F \Rightarrow H,$ by
$$ \theta \tau := \Delta(\theta \times \tau) \circ'.$$
The reader can check that when $K = \Cat$ this reduces to the
usual definition of vertical composition. We can represent this
composite pictorially as:
\[
\xymatrix{
    X
      \ar@/^2pc/[rr]_{\quad}^{F}="1"
      \ar@/_2pc/[rr]_{H}="2"
  && X'
    \ar@{}"1";"2"|(.2){\,}="7"
     \ar@{}"1";"2"|(.8){\,}="8"
    \ar@{=>}"7" ;"8"^{\theta \tau}
  }
\quad = \quad \xymatrix{
  X
    \ar@/^2pc/[rr]^{F}_{}="0"
    \ar[rr]^<<<<<<{G}^{}="1"
    \ar@/_2pc/[rr]_{H}_{}="2"
    \ar@{=>}"0";"1" ^{\theta}
    \ar@{=>}"1";"2" ^{\tau}
&&  X' }
\]

Next, let $X, X', X''$ be categories in $K$ and let $F,G\maps X
\rightarrow X'$ and $F', G'\maps X' \rightarrow X''$ be functors
in $K$.  If $\theta\maps F \Rightarrow G$ and $\theta'\maps F'
\Rightarrow G'$ are natural transformations in $K$, we define
their {\bf horizontal composite}, $\theta \circ \theta'\maps FF'
\Rightarrow GG'$, in either of two equivalent ways:
\begin{eqnarray*}
  \theta \circ \theta' &:=&
\Delta(F_{0} \times \theta)(\theta' \times G_{1}') \circ'  \\
&=& \Delta(\theta \times G_{0})(F_{1}' \times \theta') \circ' .
\end{eqnarray*}
Again, this reduces to the usual definition when $K = \Cat$. The
horizontal composite can be depicted as:
\[
\xymatrix{
    X
      \ar@/^2pc/[rr]_{\quad}^{FF'}="1"
      \ar@/_2pc/[rr]_{GG'}="2"
  && X''
    \ar@{}"1";"2"|(.2){\,}="7"
     \ar@{}"1";"2"|(.8){\,}="8"
    \ar@{=>}"7" ;"8"^{\theta \circ \theta '}
  }
\quad = \quad \xymatrix{
    X
      \ar@/^2pc/[rr]_{\quad}^{F}="1"
      \ar@/_2pc/[rr]_{G}="2"
  && X'
    \ar@{}"1";"2"|(.2){\,}="7"
     \ar@{}"1";"2"|(.8){\,}="8"
    \ar@{=>}"7" ;"8"^{\theta}
    \ar@/^2pc/[rr]_{\quad}^{F'}="1"
    \ar@/_2pc/[rr]_{G'}="2"
  && X''
    \ar@{}"1";"2"|(.2){\,}="7"
     \ar@{}"1";"2"|(.8){\,}="8"
    \ar@{=>}"7" ;"8"^{\theta '}
  }
\]

It is routine to check that these composites are again natural
transformations in $K$.  Finally, given a functor $F \maps X
\rightarrow X'$ in $K$, the identity natural transformation $1_F
\maps F \Rightarrow F$ in $K$ is given by $1_F = F_0 i$.

We now have all the ingredients of a $2$-category:

\begin{prop}\et Let $K$ be a category.  Then there
exists a strict $2$-category \textbf{\textit{K}{\bf Cat}} with
categories in $K$ as objects, functors in $K$ as morphisms, and
natural transformations in $K$ as $2$-morphisms, with composition
and identities defined as above.
\end{prop}

\noindent{\bf Proof. } It is straightforward to check that all the
axioms of a $2$-category hold; this result goes back to Ehresmann
\cite{E}. \qed

We now consider internal categories in Vect.


\section{$2$-Vector spaces} \label{2vs}

Since our goal is to categorify the concept of a Lie algebra, we
must first categorify the concept of a vector space. A
categorified vector space, or `$2$-vector space', should combine
the notions of category and vector space, and should therefore be
a category where everything in sight is linear.  That is, it will
be a category with structure analogous to that of a vector space,
with functors replacing the usual vector space operations.
Kapranov and Voevodsky \cite{KV} implemented this idea by taking a
finite-dimensional $2$-vector space to be a category of the form
$\Vect^n$, in analogy to how every finite-dimensional vector space
is of the form $k^n$. While this idea is useful in contexts such
as topological field theory \cite{Lawrence} and group
representation theory \cite{B2}, it has its limitations, which
arise from the fact that there is no notion of subtraction in
these $2$-vector spaces.

Here we instead define a $2$-vector space to be a category in
$\Vect$. Just as the main ingredient of a Lie algebra is a vector
space, a Lie $2$-algebra will have an underlying $2$-vector space
of this sort.  Thus, in this section we first define a
$2$-category of these $2$-vector spaces.  We then establish the
relationship between these $2$-vector spaces and $2$-term chain
complexes of vector spaces: that is, chain complexes having only
two nonzero vector spaces.  We conclude this section by developing
some `categorified linear algebra' --- the bare minimum necessary
for defining and working with Lie $2$-algebras in the next
section.

In the following we consider vector spaces over an arbitrary
field, $k$.

\begin{defn} \et A {\bf 2-vector space} is a category in {\rm Vect}.
\end{defn}

Thus, a $2$-vector space $V$ is a category with a vector space of
objects $V_0$ and a vector space of morphisms $V_1$, such that the
source and target maps $s,t \maps V_{1} \rightarrow V_{0}$, the
identity-assigning map $i \maps V_{0} \rightarrow V_{1}$, and the
composition map $\circ \maps V_{1} \times _{V_{0}} V_{1}
\rightarrow V_{1}$ are all {\it linear}.  As usual, we write a
morphism as $f \maps x \to y$ when $s(f) = x$ and $t(f) = y$, and
sometimes we write $i(x)$ as $1_x$.

In fact, the structure of a $2$-vector space is completely
determined by the vector spaces $V_{0}$ and $V_{1}$ together with
the source, target and identity-assigning maps. As the following
lemma demonstrates, composition can always be expressed in terms
of these, together with vector space addition:

\begin{lem} \et \label{watereddown} When $K = \Vect$, one can
omit all mention of composition in the definition of category in
$K$, without any effect on the concept being defined.
\end{lem}

\noindent{\bf Proof. }  First, given vector spaces $V_{0}$,
$V_{1}$ and maps $s,t \maps V_{1} \rightarrow V_{0}$ and $i\maps
V_{0} \rightarrow V_{1}$, we will define a composition operation
that satisfies the laws in Definition \ref{co}, obtaining a
$2$-vector space.

Given $f \in V_{1}$, we define the {\bf arrow part} of $f,$
denoted as $\vec{f}$, by
$$\vec{f} = f - i(s(f)).$$
Notice that $\vec{f}$ is in the kernel of the source map since
$$s(f - i(sf)) = s(f) - s(f) = 0.$$
While the source of $\vec{f}$ is always zero, its target may be
computed as follows:
$$t(\vec{f}) = t(f - i(s(f)) = t(f) - s(f).$$
The meaning of the arrow part becomes clearer if we write $f \maps
x \to y$ when $s(f) = x$ and $t(f) = y$. Then, given any morphism
$f \maps x \to y$, we have $\vec{f} \maps 0 \rightarrow y-x.$ In
short, taking the arrow part of $f$ has the effect of `translating
$f$ to the origin'.

We can always recover any morphism from its arrow part together
with its source, since $f = \vec{f} + i(s(f))$.  We shall take
advantage of this by identifying $f \maps x \to y$ with the
ordered pair $(x, \vec{f})$.  Note that with this notation we have
$$s(x,\vec{f}) = x , \qquad
t(x,\vec{f}) = x + t(\vec{f}).$$

Using this notation, given morphisms $f \maps x \to y$ and $g\maps
y \to z$, we define their composite by
$$f g := (x, \vec{f} + \vec{g}),$$
or equivalently,
$$(x,\vec{f}) (y,\vec{g}) := (x,\vec{f}+\vec{g}),$$
where the addition is the vector addition in $V_1$.   It remains
to show that with this composition, the diagrams of Definition
\ref{co} commute. The triangles specifying the source and target
of the identity-assigning morphism do not involve composition. The
second pair of diagrams commute since
$$s(f g) = x$$
and
$$t(f g) = x + t(\vec{f}) + t(\vec{g}) = x + (y-x) + (z-y) = z.$$
The associative law holds for composition because vector space
addition is associative.  Finally, the left unit law is satisfied
since given $f \maps x \to y$,
$$i(x) f = (x,0)  (x,\vec{f}) = (x,\vec{f}) = f$$
and similarly for the right unit law.  We thus have a $2$-vector
space.

Conversely, given a category $V$ in $\Vect$, we shall show that
its composition must be defined by the formula given above.
Suppose that $(f, g)= ((x, \vec{f}), (y, \vec{g}))$ and $(f', g')=
((x', \vec{f'}), (y', \vec{g'}))$ are composable pairs of
morphisms in $V_{1}$. Since the source and target maps are linear,
$(f + f', g + g')$ also forms a composable pair, and the linearity
of composition gives
$$(f + f') (g + g') = f g + f' g'.$$
If we set $g = 1_{y}$ and $f' = 1_{y'}$, the above equation
becomes
$$(f + 1_{y'}) (1_y + g') = f  1_y + 1_{y'} g' = f + g'.$$
Expanding out the left hand side we obtain
$$((x, \vec{f}) + (y', 0))
((y,0) + (y', \vec{g'})) = (x + y', \vec{f}) (y + y', \vec{g'}),$$
while the right hand side becomes
$$(x, \vec{f}) + (y, \vec{g'}) = (x+y', \vec{f} + \vec{g'}).$$
Thus we have $(x + y', \vec{f}) (y + y', \vec{g'}) = (x+y',
\vec{f} + \vec{g'})$, so the formula for composition in an
arbitrary $2$-vector space must be given by
$$f  g = (x, \vec{f})  (y, \vec{g}) = (x, \vec{f} + \vec{g})$$
whenever $(f,g)$ is a composable pair. This shows that we can
leave out all reference to composition in the definition of
`category in $K$' without any effect when $K = \Vect$. \qed

In order to simplify future arguments, we will often use only the
elements of the above lemma to describe a $2$-vector space.

We continue by defining the morphisms between $2$-vector spaces:

\begin{defn} \et Given $2$-vector spaces $V$ and
$W$, a {\bf linear functor} $F \maps V \to W$ is a functor in
$\Vect$ from $V$ to $W$.
\end{defn}

\noindent For now we let $2\Vect$ stand for the category of
$2$-vector spaces and linear functors between them; later we will
make $2\Vect$ into a $2$-category.

The reader may already have noticed that a $2$-vector space
resembles a {\bf 2-term chain complex} of vector spaces: that is,
a pair of vector spaces with a linear map between them, called the
`differential':
$$\xymatrix{
   C_{1} \ar[rr]^{d}
    && C_{0}.  }
$$
In fact, this analogy is very precise. Moreover, it continues at
the level of morphisms. A {\bf chain map} between $2$-term chain
complexes, say $\phi \maps C \to C'$, is simply a pair of linear
maps $\phi_0 \maps C_0 \to C'_0$ and $\phi_1 \maps C_1 \to C'_1$
that `preserves the differential', meaning that the following
square commutes:
$$\xymatrix{
    C_{1} \ar[rr]^{d}
    \ar[dd]_{\phi_{1}} &&
    C_{0} \ar[dd]^{\phi_{0}} \\ \\
    C_{1}' \ar[rr]^{d'} &&
    C_{0}'
}$$ There is a category $2\Term$ whose objects are $2$-term chain
complexes and whose morphisms are chain maps.  Moreover:

\begin{prop} \et \label{1-1vs}  The categories $2\Vect$
and $2\Term$ are equivalent.
\end{prop}

\noindent{\bf Proof. } We begin by introducing functors
$$S\maps {\rm 2Vect} \rightarrow {\rm 2Term}$$
and
$$T\maps {\rm 2Term} \rightarrow {\rm 2Vect}.$$
We first define $S$.  Given a $2$-vector space $V$, we define
$S(V)=C$ where $C$ is the $2$-term chain complex with
\begin{eqnarray*}
  C_{0} &=& V_{0}, \\
  C_{1} &=& ker(s) \subseteq V_{1}, \\
  d &=& t|_{C_{1}} ,
\end{eqnarray*}
and $s,t\maps V_{1} \rightarrow V_{0}$ are the source and target
maps associated with the $2$-vector space $V$.  It remains to
define $S$ on morphisms. Let $F \maps V \to V'$ be a linear
functor and let $S(V) = C$, $S(V') = C'.$ We define $S(F)= \phi$
where $\phi$ is the chain map with $\phi_{0} = F_{0}$ and
$\phi_{1} = F_{1}|_{C_{1}}$. Note that $\phi$ preserves the
differential because $F$ preserves the target map.

To show that $S$ is a functor, we consider composable internal
functors $F$ and $G$:
$$\xymatrix{
    V \ar[rr]^{F}
    &&
    V' \ar[rr]^{G}
    &&
    V''}$$
where $S(F) = \phi$ and $S(G) = \tau$.  Then, $S(FG) = \theta$
where $\theta$ is the chain map with
$$\theta_0 = (FG)_0 = F_0 G_0 = \phi_0 \tau_0$$ and
$$\theta_1 = (FG)_1 |_{C_{1}}  = F_1 G_1 |_{C_{1}} = \phi_1 \tau_1
|_{C_{1}}$$ where the second equality in each of the above follows
from the definition of the composite of linear functors.  We thus
have $\theta = \phi \tau$, so $S$ preserves composition.

To show $S$ preserves identities, consider a $2$-vector space $V$
with $1_V \maps V \to V$.   If we let $S(V)=C$, then $S(1_V) =
\phi$ where $\phi$ is the chain map with
$$\phi_0 = (1_V)_0 = 1_{V_0} = 1_{C_0} = 1_{S(V)_0}$$ and
$$\phi_1 = (1_V)_1 |_{C_1} = 1_{V_1} |_{C_1} = 1_{C_1} = 1_{S(V)_1}$$
where we have used the definition of the identity functor in
$\Vect$.  We thus have $\phi = 1_C$, so $S$ preserves identities.

We now turn to the second functor, $T$.  Given a $2$-term chain
complex $C$, we define $T(C) = V$ where $V$ is a $2$-vector space
with
\begin{eqnarray*}
  V_{0} & = & C_{0}, \\
  V_{1} & = & C_{0} \oplus C_{1}.
\end{eqnarray*}
To completely specify $V$ it suffices by Lemma \ref{watereddown}
to specify linear maps $s,t \maps V_1 \to V_0$ and $i \maps V_0
\to V_1$ and check that $s(i(x)) = t(i(x)) = x$ for all $x \in
V_0$. To define $s$ and $t$, we write any element $f \in V_1$ as a
pair $(x,\vec{f}) \in C_0 \oplus C_1$ and set
\[
\begin{array}{ccccl}
  s(f) &=& s(x, \vec{f}) &=& x, \\
  t(f) &=& t(x, \vec{f}) &=& x + d\vec{f}.
\end{array}
\]
For $i$, we use the same notation and set
\[  i(x) = (x,0)  \]
for all $x \in V_0$.  Clearly $s(i(x)) = t(i(x)) = x$. Note also
that with these definitions, the decomposition $V_1 = C_0 \oplus
C_1$ is precisely the decomposition of morphisms into their source
and `arrow part', as in the proof of Lemma \ref{watereddown}.
Moreover, given any morphism $f = (x,\vec{f}) \in V_1$, we have
$$t(f) - s(f) = d\vec{f}.$$

Next we define $T$ on morphisms.  Suppose $\phi \maps C
\rightarrow C'$ is a chain map between $2$-term chain complexes:
$$\xymatrix{
    C_{1} \ar[rr]^{d}
    \ar[dd]_{\phi_{1}} &&
    C_{0} \ar[dd]^{\phi_{0}} \\ \\
    C_{1}' \ar[rr]^{d'} &&
    C_{0}'
}$$ Let $T(C) = V$ and $T(C')=V'$. Then we define $T(\phi) = F$
where $F \maps V \to V'$ is the linear functor with $F_{0}
=\phi_{0}$ and $F_{1} = \phi_{0} \oplus \phi_{1}$.  To check that
$F$ really is a linear functor, note that it is linear on objects
and morphisms.  Moreover, it preserves the source and target,
identity-assigning and composition maps because all these are
defined in terms of addition and the differential in the chain
complexes $C$ and $C'$, and $\phi$ is linear and preserves the
differential.

To show that $T$ is a functor, we consider composable chain maps
$\phi$ and $\phi '$:
$$\xymatrix{
    C \ar[rr]^{\phi}
    &&
    C' \ar[rr]^{\phi '}
    &&
    C''}$$
where $T(\phi) = F$ and $T(\phi ') = G$.  Then, on one hand we
have $T(\phi \phi ') = H$ where $H$ is the linear functor with
$$H_0 = (\phi \phi ')_0 = \phi_0 \phi_0 '$$ and
$$H_1 = (\phi \phi ')_0 \oplus (\phi \phi ')_1 = \phi_0 \phi_0 ' \oplus \phi_1 \phi_1 '$$
On the other hand, $T(\phi)T(\phi ') = FG$ where $FG$ is the
linear functor with
$$(FG)_0 = F_0 G_0 = \phi _0 \phi _0 '$$ and
$$(FG)_1 = F_1 G_1 = (\phi_0 \oplus \phi_1)(\phi_0 ' \oplus \phi_1
') = \phi_0 \phi_0 ' \oplus \phi_1 \phi_1 '$$ which is equal to
the linear functor $H$, so that $T$ preserves composition.

To show $T$ preserves identities, let $C$ be a $2$-term chain
complex $C$ with $1_C \maps C \to C$.  Let $T(C)=V$.  Then $T(1_C)
= F$ where $F$ is the linear functor with
$$F_0 = (1_C)_0 = 1_{C_0} = 1_{V_0} = (1_V)_0 = 1_{T(C)_0}$$and
$$F_1 = 1_{C_0} \oplus 1_{C_1} = 1_{C_0 \oplus C_1} = 1_{V_1} =
(1_V)_1 = 1_{T(C)_1}$$ so that $T$ preserves identities.

To show that $S$ and $T$ form an equivalence, we construct natural
isomorphisms $\alpha \maps ST \To 1_{2\Vect}$ and $\beta \maps TS
\To 1_{2\Term}$.

To construct $\alpha$, consider a $2$-vector space $V$. Applying
$S$ to $V$ we obtain the $2$-term chain complex
$$\xymatrix{
    ker(s) \ar[rr]^<<<<<<<<<{t|_{ker(s)}}
   && V_{0}.
}$$ Applying $T$ to this result, we obtain a $2$-vector space $V'$
with the space $V_{0}$ of objects and the space $V_{0} \oplus
ker(s)$ of morphisms.  The source map for this $2$-vector space is
given by $s'(x,\vec{f}) = x$, the target map is given by
$t'(x,\vec{f}) = x + t(\vec{f})$, and the identity-assigning map
is given by $i'(x) = (x,0)$.  We thus can define an isomorphism
$\alpha_V \maps V' \to V$ by setting
\begin{eqnarray*}
(\alpha_V)_0(x) &=& x , \\
(\alpha_V)_1(x,\vec{f}) &=& i(x) + \vec{f} .
\end{eqnarray*}
It is easy to check that $\alpha_V$ is a linear functor. It is an
isomorphism thanks to the fact, shown in the proof of Lemma
\ref{watereddown}, that every morphism in $V$ can be uniquely
written as $i(x) + \vec{f}$ where $x$ is an object and $\vec{f}
\in ker(s)$.

To show that $\alpha$ is indeed a natural transformation, we must
show that given $2$-vector spaces $V,W$ and a linear functor $F
\maps V \to W$, the following diagram commutes:
$$\xymatrix{
    ST(V) \ar[rr]^{ST(F)}
    \ar[dd]_{\alpha_V} &&
    ST(W) \ar[dd]^{\alpha_W} \\ \\
    V \ar[rr]^{F} &&
    W
}$$ where $ST(F)$ is the linear functor with
$$ST(F)_0 = F_0 \; \; \; \textrm{and} \; \; \; ST(F)_1 = F_0
\oplus F_1 |_{ker(s)}.$$  This amounts to a straightforward
computation.

To construct $\beta$, consider a $2$-term chain complex, $C$,
given by
$$\xymatrix{
   C_{1} \ar[rr]^{d}
    && C_{0}.
}$$ Then $T(C)$ is the $2$-vector space with the space $C_{0}$ of
objects, the space $C_{0} \oplus C_{1}$ of morphisms, together
with the source and target maps $s \maps (x,\vec{f}) \mapsto x$,
$t \maps (x,\vec{f}) \mapsto x + d\vec{f}$ and the
identity-assigning map $i \maps x \mapsto (x,0)$. Applying the
functor $S$ to this $2$-vector space we obtain a $2$-term chain
complex $C'$ given by:
$$\xymatrix{
   ker(s) \ar[rr]^<<<<<<<<<<<{t|_{ker(s)}}
    && C_{0}.
}$$ Since $ker(s) = \{(x,\vec{f}) | x = 0 \} \subseteq C_0 \oplus
C_1$, there is an obvious isomorphism $ker(s) \cong C_{1}$.  Using
this we obtain an isomorphism $\beta_C \maps C' \to C$ given by:
$$\xymatrix{
    ker(s) \ar[rr]^{t|_{ker(s)}}
    \ar[dd]_{\sim} &&
    C_{0} \ar[dd]^{1} \\ \\
    C_{1} \ar[rr]^{d} &&
    C_{0}
}$$ where the square commutes because of how we have defined $t$.

To show that $\beta$ is indeed a natural transformation, we must
show that given $2$-term chain complexes $C,D$ and a chain map
$\phi \maps C \to D$, the following diagram commutes:
$$\xymatrix{
    TS(C) \ar[rr]^{TS(\phi)}
    \ar[dd]_{\beta_C} &&
   TS(D) \ar[dd]^{\beta_D} \\ \\
    C \ar[rr]^{\phi} &&
   D
}$$ where $TS(\phi)$ is the chain map with
$$TS(\phi)_0 = \phi_0 \; \; \; \textrm{and} \; \; \; TS(\phi)_1 = (\phi_0 \oplus \phi_1)|_{C_1} =
\phi_1 .$$  This amounts to a straightforward computation. \qed

As mentioned in the introduction to this chapter, the idea behind
Proposition \ref{1-1vs} goes back at least to Grothendieck
\cite{G}, who showed that groupoids in the category of abelian
groups are equivalent to $2$-term chain complexes of abelian
groups. There are many elaborations of this idea, some of which we
will mention later, but for now the only one we really {\it need}
involves making $2\Vect$ and $2\Term$ into $2$-categories and
showing that they are $2$-equivalent as $2$-categories.  To do
this, we require the notion of a `linear natural transformation'
between linear functors.  This will correspond to a chain homotopy
between chain maps.

\begin{defn} \et Given two linear functors $F,G \maps V \to W$
between $2$-vector spaces, a {\bf linear natural transformation}
$\alpha \maps F \To G$ is a natural transformation in $\Vect$.
\end{defn}

\begin{defn} \et \label{22Vect} We define {\bf 2Vect} to be
$\Vect\Cat$, or in other words, the $2$-category of $2$-vector
spaces, linear functors and linear natural transformations.
\end{defn}

Recall that in general, given two chain maps $\phi, \psi \maps C
\to C'$, a {\bf chain homotopy} $\tau \maps \phi \To \psi$ is a
family of linear maps $\tau \maps C_{p} \rightarrow C_{p+1}'$ such
that $\tau_{p} d_{p+1}' + d_{p} \tau_{p-1} = \psi_{p} - \phi_{p}$
for all $p$. In the case of $2$-term chain complexes, a chain
homotopy amounts to a map $\tau \maps C_{0} \rightarrow C_{1}'$
satisfying $\tau d' = \psi_{0} - \phi_{0}$ and $d \tau = \psi_{1}
- \phi_{1}$.

\begin{defn} \et \label{22Term}
We define {\bf 2Term} to be the $2$-category of $2$-term chain
complexes, chain maps, and chain homotopies.
\end{defn}

\noindent We will continue to sometimes use $2\Term$ and $2\Vect$
to stand for the underlying categories of these (strict)
$2$-categories. It will be clear by context whether we mean the
category or the $2$-category.

The next result strengthens Proposition \ref{1-1vs}.

\begin{thm} \et \label{equivof2vs} The $2$-category {\rm $2$Vect}
is $2$-equivalent to the $2$-category {\rm $2$Term}.
\end{thm}

\noindent{\bf Proof. } We begin by constructing $2$-functors
$$S\maps {\rm 2Vect} \rightarrow {\rm 2Term}$$
and
$$T\maps {\rm 2Term} \rightarrow {\rm 2Vect}.$$
By Proposition \ref{1-1vs}, we need only to define $S$ and $T$ on
$2$-morphisms. Let $V$ and $V'$ be $2$-vector spaces, $F,G \maps V
\to V'$ linear functors, and $\theta \maps F \To G$ a linear
natural transformation.  Then we define the chain homotopy
$S(\theta) \maps S(F) \To S(G)$ by
$$  S(\theta)(x) = \vec{\theta}_x ,$$
using the fact that a 0-chain $x$ of $S(V)$ is the same as an
object $x$ of $V$.

To show that $S(\theta)$ really is a chain homotopy we must show
that
$$S(\theta) t|_{ker(s')} = G_{0} - F_{0}$$ and
$$t|_{ker(s)} S(\theta) = G_{1}|_{ker(s)} - F_{1}|_{ker(s)}.$$
To demonstrate the first, we use the fact that $S(\theta)(x) =
\vec{\theta}_{x}$ is the arrow part of the linear transformation
$\theta$, so
$$\vec{\theta}_{x} \maps 0 \to t(\theta_{x}) - s(\theta_{x}) = G_0
(x) - F_0(x)$$ which shows that the target of $S(\theta)$ is $G_0
- F_0$, which is the first condition above. Next, consider $f
\maps 0 \to x$, so that $f \in ker(s)$. Then the commutative
square law for a linear natural transformation says
$$\vec{\theta}_0 G(f) = F(f) \vec{\theta} _x.$$
To obtain these composites, we add arrow parts, so that we have
$$G_1 (f) = F_1(f) + \vec{\theta}_x = F_1(f) + t| _{ker(s)} S(\theta)$$
which is precisely the second condition we sought.  We thus have
that $S(\theta)$ is a chain homotopy.

Conversely, let $C$ and $C'$ be $2$-term chain complexes,
$\phi,\psi \maps C \to C'$ chain maps and $\tau \maps \phi \To
\psi$ a chain homotopy.  Then we define the linear natural
transformation $T(\tau) \maps T(\phi) \To T(\psi)$ by
$$  T(\tau)(x) = (\phi_{0}(x),\tau(x)), $$
where we use the description of a morphism in $S(C')$ as a pair
consisting of its source and its arrow part, which is a 1-chain in
$C'$.

To show that $T(\tau)$ is really a linear natural transformation,
we must show that the three diagrams of Definition
\ref{conattrans} commute.  We show that the target is correct;
showing the source is correct is easier.  To show that this
diagram
\[
\begin{diagram}
\node{C_0} \arrow{se,t}{\psi_{0}} \arrow{s,l}{T(\tau)} \\
\node{C_0 ' \oplus C_1 '}  \arrow{e,t}{t} \node{C'_0}
\end{diagram} \]
commutes, we observe that $$t(T(\tau(x)) = t(\phi_0(x), \tau(x)) =
\phi_0 (x) + d'(\tau(x)) = \phi_0(x) + \psi_0 (x) - \phi_0 (x) =
\psi_0 (x)$$ where the third equality follows from the first of
the two properties of a chain homotopy following Definition
\ref{22Vect}. It remains to demonstrate that the commutative
square law holds by showing that the following diagram commutes:
\[  \xymatrix{
   C_0 \oplus C_1
    \ar[rr]^<<<<<<<<<<<<<<<{\Delta (s T(\tau) \times \psi)}
    \ar[dd]_{\Delta (\phi \times t T(\tau))}
     && (C_0 ' \oplus C_1 ') \otimes (C_0 ' \oplus C_1 ')
    \ar[dd]^{\circ'} \\ \\
   (C_0 ' \oplus C_1 ') \otimes (C_0 ' \oplus C_1 ')
    \ar[rr]^<<<<<<<<<<<<<<<{\circ'}
     && C_0 ' \oplus C_1 '
}
\]
Sending $(x, \vec{f}) \in C_0 \oplus C_1$ along the top-right path
of the above diagram produces
$$(\phi_0 (x), \tau(x)) \circ ' (\psi_0 (x), \psi_1(\vec{f})) =
(\phi_0 (x), \tau(x) + \psi_1(\vec{f}))$$ whereas sending it along
the left-bottom path produces
\begin{eqnarray*}
(\phi_0 (x),\phi_1(\vec{f})) \circ ' (\phi_0(x + d(\vec{f})),
\tau(x + d(\vec{f})) &=& (\phi_0 (x), \phi_1(\vec{f}) + \tau(x +
d(\vec{f}))) \\ & =& (\phi_0 (x), \phi_1(\vec{f}) + \tau(x) +
\tau(d(\vec{f}))) \\ &=& (\phi_0 (x), \phi_1(\vec{f}) + \tau(x) +
\psi_1(\vec{f}) - \phi_1(\vec{f}) )
\end{eqnarray*}
where the final equality follows from the second of the two
properties of a chain homotopy. Thus, we have that the commutative
square law holds.

We leave it to the reader to check that the natural isomorphisms
$\alpha \maps ST \To 1_{2\Vect}$ and $\beta \maps TS \To
1_{2\Term}$ defined in the proof of Proposition \ref{1-1vs} extend
to this $2$-categorical context. \qed

We conclude this section with a little categorified linear
algebra. We consider the direct sum and tensor product of
$2$-vector spaces.

\begin{prop} \et Given $2$-vector spaces $V=(V_{0}, V_{1},
s,t,i ,\circ)$ and $V'=(V_{0}', V_{1}',$ $s',t',i',\circ')$, there
is a $2$-vector space $V \oplus V'$ having:
\begin{itemize}
\item $V_0 \oplus V'_0$ as its vector space of objects,
\item $V_1 \oplus V'_1$ as its vector space of morphisms,
\item $s\oplus s'$ as its source map,
\item $t \oplus t'$ as its target map,
\item $i \oplus i'$ as its identity-assigning map, and
\item $\circ \oplus \circ'$ as its composition map.
\end{itemize}
\end{prop}

\noindent{\bf Proof. }  The proof amounts to a routine
verification that the diagrams in Definition \ref{co} commute.
\qed

\begin{prop} \label{tenprod}
\et Given $2$-vector spaces $V=(V_{0}, V_{1}, s,t,i ,\circ)$ and
$V'=(V_{0}', V_{1}',$ $s',t',i',\circ')$, there is a $2$-vector
space $V \otimes V'$ having:
\begin{itemize}
\item $V_0 \otimes V'_0$ as its vector space of objects,
\item $V_1 \otimes V'_1$ as its vector space of morphisms,
\item $s\otimes s'$ as its source map,
\item $t \otimes t'$ as its target map,
\item $i \otimes i'$ as its identity-assigning map, and
\item $\circ \otimes \circ'$ as its composition map.
\end{itemize}
\end{prop}

\noindent{\bf Proof. }  Again, the proof is a routine
verification. \qed

We now check the correctness of the above definitions by showing
the universal properties of the direct sum and tensor product.
These universal properties only require the category structure of
$2\Vect$, not its $2$-category structure, since the necessary
diagrams commute `on the nose' rather than merely up to a
$2$-isomorphism, and uniqueness holds up to isomorphism, not just
up to equivalence. The direct sum is what category theorists call
a `biproduct': both a product and coproduct, in a compatible way
\cite{Mac}:

\begin{prop} \et The direct sum $V \oplus V'$ is
the biproduct of the $2$-vector spaces $V$ and $V'$, with the
obvious inclusions
\[    i \maps V \to V \oplus V', \qquad i' \maps V' \to V \oplus V' \]
and projections
\[   p \maps V \oplus V' \to V, \qquad p' \maps V \oplus V' \to V' .\]
\end{prop}

\noindent{\bf Proof. } Showing that $V \oplus V'$ is a biproduct
requires verifying that the following hold
$$ ip=1_V \qquad i'p'=1_{V'} \;\;\; \textrm{and} \;\;\; pi + p'i' = 1_{V \oplus V'}$$
which is clear. \qed

Since the direct sum $V \oplus V'$ is a product in the categorical
sense, we may also denote it by $V \times V'$, as we do now in
defining a `bilinear functor', which is used in stating the
universal property of the tensor product:

\begin{defn} \et \label{bilinearfunct} Let $V, V',$ and $W$ be
$2$-vector spaces.  A {\bf bilinear functor} $F \maps V \times V'
\rightarrow W$ is a functor such that the underlying map on
objects
$$F_{0}\maps V_0 \times V'_0 \rightarrow W_{0}$$ and
the underlying map on morphisms
$$F_{1}\maps V_1 \times V'_1 \rightarrow W_{1}$$
are bilinear.
\end{defn}

\begin{prop} \et Let $V, V',$ and $W$ be $2$-vector spaces.  Given
a bilinear functor \newline $F \maps V \times V' \rightarrow W$
there exists a unique linear functor $\tilde{F} \maps V \otimes V'
\rightarrow W$ such that
$$\xymatrix{
   V \times V'
   \ar[rr]^{F}
   \ar[dd]_{i}
   &&
   W
   \\ \\
   V \otimes V'
   \ar[rruu]_{\tilde{F}}
}$$ commutes, where $i \maps V \times V' \rightarrow V \otimes V'$
is given by $(v,w) \mapsto v \otimes w$ for $(v,w) \in (V \times
V')_{0}$ and $(f,g) \mapsto f \otimes g$ for $(f,g) \in (V \times
V')_{1}$.
\end{prop}

\noindent{\bf Proof. } The existence and uniqueness of
$\tilde{F_0} \maps (V \otimes V')_0 \to W_0$ and $\tilde{F_1}
\maps (V \otimes V')_1 \to W_1$ follow from the universal property
of the tensor product of vector spaces, and it is then
straightforward to check that $\tilde{F}$ is a linear functor.
\qed

We can also form the tensor product of linear functors. Given
linear functors $F \maps V \to V'$ and $G \maps W \to W'$, we
define $F \otimes G \maps V \otimes V' \to W \otimes W'$ by
setting:
\[
\begin{array}{ccl}
(F \otimes G)_{0} &=& F_{0} \otimes G_{0}, \\
(F \otimes G)_{1} &=& F_{1} \otimes G_{1}.
\end{array}
\]
Furthermore, there is an `identity object' for the tensor product
of $2$-vector spaces.    In $\Vect$, the ground field $k$ acts as
the identity for tensor product: there are canonical isomorphisms
$k \otimes V \cong V$ and $V \otimes k \cong V$.  For $2$-vector
spaces, a categorified version of the ground field plays this
role:

\begin{prop} \et There exists a unique $2$-vector space $K$,
the {\bf categorified ground field}, with $K_{0} = K_{1} = k$ and
$s,t,i = 1_{k}.$
\end{prop}

\noindent{\bf Proof. }  Lemma \ref{watereddown} implies that there
is a unique way to define composition in $K$ making it into a
$2$-vector space.  In fact, every morphism in $K$ is an identity
morphism.  \qed

\begin{prop} \et Given any $2$-vector space $V$,
there is an isomorphism $\ell_V \maps K \otimes V \to V$, which is
defined on objects by $a \tensor v \mapsto av$ and on morphisms by
$a \tensor f \mapsto af$.   There is also an isomorphism $r_V
\maps V \otimes K \to V$, defined similarly.
\end{prop}

\noindent{\bf Proof. }  This is straightforward. \qed

The functors $\ell_V$ and $r_V$ are a categorified version of left
and right multiplication by scalars.  Our $2$-vector spaces also
have a categorified version of addition, namely a linear functor
\[     + \maps V \oplus V \to V \]
mapping any pair $(x,y)$ of objects or morphisms to $x+y$.
Combining this with scalar multiplication by the object $-1 \in
K$, we obtain another linear functor
\[      - \maps V \oplus V \to V \]
mapping $(x,y)$ to $x-y$.  This is the sense in which our
$2$-vector spaces are equipped with a categorified version of
subtraction. All the usual rules governing addition of vectors,
subtraction of vectors, and scalar multiplication hold `on the
nose' as equations.

One can show that with the above tensor product, the category
$2\Vect$ becomes a symmetric monoidal category. One can go further
and make the $2$-category version of $2\Vect$ into a symmetric
monoidal $2$-category \cite{DS}, but we will not need this here.
Now that we have a definition of $2$-vector space and some basic
tools of categorified linear algebra we may proceed to the main
focus of this chapter:  the definition of a categorified Lie
algebra.


\section{Semistrict Lie 2-algebras} \label{Lie2algs}
\subsection{Definitions} \label{definitions}

We now introduce the concept of a `Lie $2$-algebra', which blends
together the notion of a Lie algebra with that of a category.  As
mentioned previously, a Lie $2$-algebra is a $2$-vector space
equipped with a bracket {\it functor}, which satisfies the Jacobi
identity {\it up to a natural isomorphism}, the `Jacobiator'. Then
we require that the Jacobiator satisfy a new coherence law of its
own, the `Jacobiator identity'.  We shall assume the bracket is
bilinear in the sense of Definition \ref{bilinearfunct}, and also
skew-symmetric:

\begin{defn} \et Let $V$ and $W$ be $2$-vector spaces.
A bilinear functor $F\maps V \times V \to W$ is {\bf
skew-symmetric} if $F(x,y) = -F(y,x)$ whenever $(x,y)$ is an
object or morphism of $V \times V$.  If this is the case we also
say the corresponding linear functor $\tilde{F} \maps V \otimes V
\to W$ is skew-symmetric.
\end{defn}

\noindent We shall also assume that the Jacobiator is trilinear
and completely antisymmetric:

\begin{defn} \et Let $V$ and $W$ be $2$-vector spaces.
A functor $F \maps V^n \to W$ is \textbf{n}-{\bf linear} if
$F(x_1,\dots, x_n)$ is linear in each argument, where $(x_1,
\dots, x_n)$ is an object or morphism of $V^n$.  Given $n$-linear
functors $F,G \maps V^n \to W$, a natural transformation $\theta
\maps F \To G$ is \textbf{n}-{\bf linear} if
$\theta_{x_1,\dots,x_n}$ depends linearly on each object $x_i$,
and {\bf completely antisymmetric} if the arrow part of
$\theta_{x_1,\dots,x_n}$ is completely antisymmetric under
permutations of the objects.
\end{defn}

\noindent Since we do not weaken the bilinearity or skew-symmetry
of the bracket, we call the resulting sort of Lie $2$-algebra
`semistrict':

\begin{defn} \et \label{defnlie2alg}
A {\bf semistrict Lie $2$-algebra} consists of:

\begin{itemize}
\item a $2$-vector space $L$
\end{itemize}
equipped with
\begin{itemize}
\item a skew-symmetric
bilinear functor, the {\bf bracket}, $[\cdot, \cdot]\maps L \times
L \rightarrow L$
\item a completely antisymmetric
trilinear natural isomorphism, the {\bf Jacobiator},
$$J_{x,y,z} \maps [[x,y],z] \to [x,[y,z]] + [[x,z],y],$$
\end{itemize}
that is required to satisfy
\begin{itemize}
 \item the {\bf Jacobiator identity}:
    $$J_{[w,x],y,z} ([J_{w,x,z},y] + 1) (J_{w, [x,z], y} +
       J_{[w,z],x,y} + J_{w,x, [y,z]}) = $$
$$[J_{w,x,y},z] (J_{[w,y],x,z} + J_{w, [x,y],z}) ([J_{w,y,z},x] + 1)
([w, J_{x,y,z}] + 1)$$
\end{itemize}
for all $w,x,y,z \in L_{0}$, where the identity morphisms are
chosen to ensure that the target of each morphism being composed
is the source of the next one. (There is only one choice of
identity morphism that can be added to each term to make the
composite well-defined.)
\end{defn}

The Jacobiator identity looks quite intimidating at first. But if
we draw it as a commutative diagram, we see that it relates two
ways of using the Jacobiator to rebracket the expression
$[[[w,x],y],z]$:

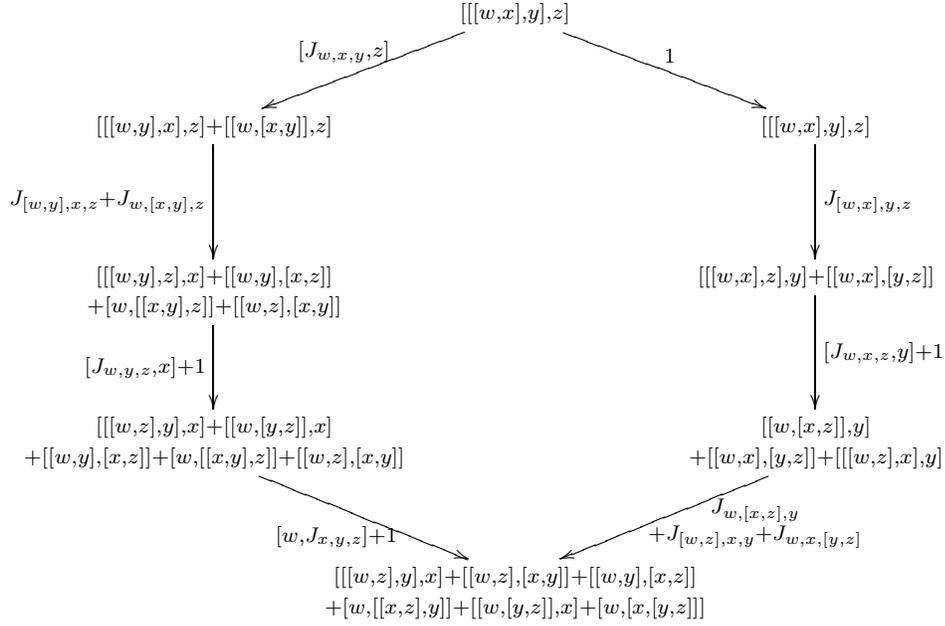
$$ \def\objectstyle{\scriptstyle}
  \def\labelstyle{\scriptstyle}
   \xy
   (0,35)*+{[[[w,x],y],z]}="1";
   (-40,20)*+{[[[w,y],x],z] + [[w,[x,y]],z]}="2";
   (40,20)*+{[[[w,x],y],z]}="3";
   (-40,0)*+{[[[w,y],z],x] + [[w,y],[x,z]]}="4'";
   (-40,-4)*+{+ [w,[[x,y],z]] + [[w,z],[x,y]]}="4";
   (40,0)*+{[[[w,x],z],y] + [[w,x],[y,z]]}="5'";
   (-40,-20)*+{[[[w,z],y],x] + [[w,[y,z]],x]}="6'";
   (-40,-24)*+{+ [[w,y], [x,z]] + [w,[[x,y],z]] + [[w,z],[x,y]]}="6";
   (40,-20)*+{[[w,[x,z]],y]}="7'";
   (40,-24)*+{+ [[w,x],[y,z]] + [[[w,z],x],y]}="7";
   (0,-40)*+{[[[w,z],y],x] + [[w,z],[x,y]]  + [[w,y],[x,z]]}="8'";
   (0,-44)*+{+ [w,[[x,z],y]]  + [[w,[y,z]],x] + [w,[x,[y,z]]]}="8";
            (32,-31)*{J_{w,[x,z],y} }; 
            (32,-34.5)*{+  J_{[w,z],x,y} + J_{w,x,[y,z]}};
        {\ar_{[J_{w,x,y},z]}                   "1";"2"};
        {\ar^{1}                               "1";"3"};
        {\ar_{J_{[w,y],x,z} + J_{w,[x,y],z}}   "2";"4'"};
        {\ar_{[J_{w,y,z},x]+1}                   "4";"6'"};
        {\ar^{J_{[w,x],y,z}}                   "3";"5'"};
        {\ar^{[J_{w,x,z},y]+1}                   "5'";"7'"};
        {\ar_{[w,J_{x,y,z}]+1}                   "6";"8'"};
        {\ar^{}                                "7";"8'"};
\endxy
\\ \\
$$

\vskip 1em

\noindent Here the identity morphisms come from terms on which we
are not performing any manipulation.

Typically, it is unusual to label an edge of a commutative diagram
solely by an identity morphism, as we do for the first arrow on
the right. We include it here because in Section \ref{conclusions}
we show that the Jacobiator identity is really just a disguised
version of the `Zamolodchikov tetrahedron equation', which plays
an important role in the theory of higher-dimensional knots and
braided monoidal $2$-categories \cite{BN,CS, Crans,DS,KV}. The
Zamolochikov tetrahedron equation says that two $2$-morphisms are
equal, each of which is the vertical composite of four factors.
However, when we translate this equation into the language of Lie
$2$-algebras, one of these factors is an identity $2$-morphism.

Until Section \ref{conclusions} of the following chapter, the term
`Lie $2$-algebra' will always refer to a semistrict one as defined
above.   We continue by setting up a $2$-category of these Lie
$2$-algebras. A homomorphism between Lie $2$-algebras should
preserve both the $2$-vector space structure and the bracket.
However, we shall require that it preserve the bracket only {\it
up to isomorphism}
--- or more precisely, up to a natural isomorphism satisfying a
suitable coherence law.  Thus, we make the following definition.

\begin{defn} \et \label{lie2algfunct} Given
Lie $2$-algebras $L$ and $L'$, a {\bf homomorphism} $F \maps L
\rightarrow L'$ consists of:

\begin{itemize}
  \item A linear functor $F$ from the underlying $2$-vector space of $L$
    to that of $L'$, and

  \item a skew-symmetric bilinear natural transformation
    $$F_{2}(x,y)\maps [F_{0}(x), F_{0}(y)] \rightarrow F_{0}[x,y]$$
\end{itemize}

such that the following diagram commutes:

$$\xymatrix{
    [F_{0}(x), [F_{0}(y), F_{0}(z)]]
      \ar[rrrr]^<<<<<<<<<<<<<<<<<<<<<<{J_{F_{0}(x), F_{0}(y), F_{0}(z)}}
      \ar[dd]_{[1, F_{2}]}
       &&&& [[F_{0}(x), F_{0}(y)], F_{0}(z)] + [F_{0}(y), [F_{0}(x),
F_{0}(z)]]
      \ar[dd]^{[F_{2}, 1] + [1, F_{2}]} \\ \\
       [F_{0}(x), F_{0}[y,z]]
      \ar[dd]_{F_{2}}
       &&&& [F_{0}[x,y], F_{0}(z)] + [F_{0}(y), F_{0}[x,z]]
      \ar[dd]^{F_{2} + F_{2}} \\ \\
       F_{0}[x,[y,z]]
      \ar[rrrr]^{F_{1}(J_{x,y,z})}
       &&&& F_{0}[[x,y],z] + F_{0}[y,[x,z]]}$$

\end{defn}

\noindent Here and elsewhere we omit the arguments of natural
transformations such as $F_2$ and $G_2$ when these are obvious
from context.

We also have `$2$-homomorphisms' between homomorphisms, which are
linear natural transformations with an extra property:

\begin{defn} \et \label{lie2algnattrans} Let $F,G \maps
L \to L'$ be Lie $2$-algebra homomorphisms.  A {\bf
2-homomorphism} $\theta \maps F \To G$ is a linear natural
transformation from $F$ to $G$ such that the following diagram
commutes:

$$\xymatrix{
    [F_{0}(x), F_{0}(y)]
     \ar[rr]^{F_{2}}
     \ar[dd]_{[\theta_{x}, \theta_{y}]}
      && F_{0}[x,y]
     \ar[dd]^{\theta_{[x,y]}} \\ \\
      [G_{0}(x), G_{0}(y)]
     \ar[rr]^{G_{2}}
      && G_{0}[x,y] }$$

\end{defn}

\noindent Definitions \ref{lie2algfunct} and \ref{lie2algnattrans}
are closely modelled after the usual definitions of `monoidal
functor' and `monoidal natural transformation' \cite{Mac}.

Next we introduce composition and identities for homomorphisms and
$2$-homomorphisms. The composite of a pair of Lie $2$-algebra
homomorphisms $F\maps L \rightarrow L'$ and $G\maps L' \rightarrow
L''$ is given by letting the functor $FG \maps L \to L''$ be the
usual composite of $F$ and $G$:

$$\xymatrix{
    L
     \ar[rr]^{F}
      &&  L'
     \ar[rr]^{G}
      && L''}$$

\noindent while letting $(FG)_{2}$ be defined as the following
composite:

$$\xymatrix{
    [(FG)_{0}(x), (FG)_{0}(y)]
     \ar[rr]^<<<<<<<<<{(FG)_{2}}
     \ar[dd]_{G_{2}}
      && (FG)_{0}[x,y]  \\ \\
      G_{0}[F_{0}(x), F_{0}(y)]
     \ar[uurr]^{F_{2}\circ G}}.
$$

\noindent where $F_2 \circ G$ is the result of whiskering the
functor $G$ by the natural transformation $F_2$. The identity
homomorphism $1_L \maps L \to L$ has the identity functor as its
underlying functor, together with an identity natural
transformation as $(1_L)_2$. Since $2$-homomorphisms are just
natural transformations with an extra property, we vertically and
horizontally compose these the usual way, and an identity
$2$-homomorphism is just an identity natural transformation.  We
obtain:

\begin{prop} \et There is a strict $2$-category {\bf Lie2Alg}
with semistrict Lie $2$-algebras as objects, homomorphisms between
these as morphisms, and $2$-homomorphisms between those as
$2$-morphisms, with composition and identities defined as above.
\end{prop}

\noindent{\bf Proof. } We leave it to the reader to check the
details, including that the composite of homomorphisms is a
homomorphism, this composition is associative, and the vertical
and horizontal composites of $2$-homomorphisms are again
$2$-homomorphisms.  \qed

Finally, note that there is a forgetful $2$-functor from {\rm
Lie$2$Alg} to {\rm $2$Vect}, which is analogous to the forgetful
functor from {\rm LieAlg} to {\rm Vect}.

We continue by exhibiting the correlation between our semistrict
Lie $2$-algebras and special versions of Stasheff's
$L_{\infty}$-algebras.

\subsection{$L_{\infty}$-algebras} \label{Linftyalgs}

An $L_\infty$-algebra is a generalization, or homotopy version, of
a Lie algebra.  More specifically, an $L_\infty$-algebra is a
chain complex equipped with a bilinear skew-symmetric bracket
operation that satisfies the Jacobi identity up to an infinite
tower of chain homotopies, thereby blending together the notion of
a Lie algebra with that of a chain complex. Such structures are
also called are `strongly homotopy Lie algebras' or `sh Lie
algebras' for short. Though they had existed in the literature
beforehand, they made their first notable appearance in a 1985
paper on deformation theory by Schlessinger and Stasheff
\cite{SS}. Since then, they have been systematically explored and
applied in a number of other contexts \cite{LM, LS, Mar, P}.

Since $2$-vector spaces are equivalent to $2$-term chain
complexes, as described in Section \ref{2vs}, it should not be
surprising that $L_\infty$-algebras are related to the
categorified Lie algebras we defined in the previous section.
Indeed, after recalling the definition of an $L_\infty$-algebra we
prove that the $2$-category of semistrict Lie $2$-algebras is
equivalent to a $2$-category of `$2$-term' $L_\infty$-algebras:
that is, those having a zero-dimensional space of $j$-chains for
$j > 1$.

Henceforth, all algebraic objects mentioned are considered over a
fixed field $k$ of characteristic other than 2.  We make
consistent use of the usual sign convention when dealing with
graded objects.  That is, whenever we interchange something of
degree $p$ with something of degree $q,$ we introduce a sign of
$(-1)^{pq}.$ The following conventions regarding graded vector
spaces, permutations, unshuffles, etc., follow those of Lada and
Markl \cite{LM}.

For graded indeterminates $x_{1}, \ldots, x_{n}$ and a permutation
$\sigma \in S_{n}$ we define the {\bf Koszul sign} $\epsilon
(\sigma) = \epsilon(\sigma; x_{1}, \dots, x_{n})$ by
$$x_{1} \wedge \cdots \wedge x_{n} = \epsilon(\sigma; x_{1},
\ldots, x_{n}) \cdot x_{\sigma(1)} \wedge \cdots \wedge
x_{\sigma(n)},$$ which must be satisfied in the free
graded-commutative algebra on $x_{1}, \ldots, x_{n}.$  This is
nothing more than a formalization of what has already been said
above.  Furthermore, we define
$$\chi(\sigma) = \chi(\sigma;
x_{1}, \dots, x_{n}) := \textrm{sgn} (\sigma) \cdot
\epsilon(\sigma; x_{1}, \dots, x_{n}).$$ Thus, $\chi(\sigma)$
takes into account the sign of the permutation in $S_{n}$ and the
sign obtained from iteration of the basic convention.

If $n$ is a natural number and $1 \leq j \leq n-1$ we say that
$\sigma \in S_{n}$ is an $(j,n-j)${\bf -unshuffle} if
$$ \sigma(1) \leq\sigma(2) \leq \cdots \leq \sigma(j)
\hspace{.2in} \textrm{and} \hspace{.2in} \sigma(j+1) \leq
\sigma(j+2) \leq \cdots \leq \sigma(n).$$ Readers familiar with
shuffles will recognize unshuffles as their inverses. A {\it
shuffle} of two ordered sets (such as a deck of cards) is a
permutation of the ordered union preserving the order of each of
the given subsets. An {\it unshuffle} reverses this process.  A
simple example should clear up any confusion:

\begin{example} \et When $n=3$, the $(1,2)$-unshuffles in $S_{3}$ are:
$$\emph{id} =  \begin{pmatrix}
  1 & 2 & 3 \\
  1 & 2 & 3
\end{pmatrix}, \hspace{.2in}
(132) = \begin{pmatrix}
  1 & 2 & 3 \\
  3 & 1 & 2
\end{pmatrix}, \hspace{.1in} \emph{and} \hspace{.1in}
(12) = \begin{pmatrix}
  1 & 2 & 3 \\
  2 & 1 & 3
\end{pmatrix}.$$
\end{example}

The following definition of an $L_{\infty}$-structure was
formulated by Stasheff in 1985, see \cite{SS}.  This definition
will play an important role in what will follow.

\begin{defn} \et \label{L-alg} An
{\bf $\mathbf{L_{\infty}}$-algebra} is a graded vector space $V$
equipped with a system $\{l_{k}| 1 \leq k < \infty\}$ of linear
maps $l_{k} \maps V^{\otimes k} \rightarrow V$ with $\deg(l_{k}) =
k-2$ which are totally antisymmetric in the sense that
\begin{eqnarray}
   l_{k}(x_{\sigma(1)}, \dots,x_{\sigma(k)}) =
   \chi(\sigma)l_{k}(x_{1}, \dots, x_{n})
\label{antisymmetry}
\end{eqnarray}
for all $\sigma \in S_{n}$ and $x_{1}, \dots, x_{n} \in V,$ and,
moreover, the following generalized form of the Jacobi identity
holds for $0 \le n < \infty :$
\begin{eqnarray}
   \displaystyle{\sum_{i+j = n+1}
   \sum_{\sigma}
   \chi(\sigma)(-1)^{i(j-1)} l_{j}
   (l_{i}(x_{\sigma(1)}, \dots, x_{\sigma(i)}), x_{\sigma(i+1)},
   \ldots, x_{\sigma(n)}) =0,}
\label{megajacobi}
\end{eqnarray}
where the summation is taken over all $(i,n-i)$-unshuffles with $i
\geq 1.$
\end{defn}

While it appears complicated at first, this definition truly does
combine the important aspects of Lie algebras and chain complexes.
The map $l_1$ makes $V$ into a chain complex, since this map has
degree $-1$ and equation (\ref{megajacobi}) says its square is
zero. Moreover, the map $l_{2}$ resembles a Lie bracket, since it
is skew-symmetric in the graded sense by equation
(\ref{antisymmetry}). In what follows, we usually denote $l_1(x)$
as $dx$ and $l_2(x,y)$ as $[x,y]$.  The higher $l_k$ maps are
related to the Jacobiator, the Jacobiator identity, and the higher
coherence laws that would appear upon further categorification of
the Lie algebra concept.

To make this more precise, let us refer to an $L_{\infty}$-algebra
$V$ with $V_{n} = 0$ for $n \geq k$ as a \textbf{\textit{k}-term
$\mathbf{L_{\infty}}$-algebra.} Note that a $1$-term
$L_{\infty}$-algebra is simply an ordinary Lie algebra, where
$l_{3} =0$ gives the Jacobi identity.  However, in a $2$-term
$L_{\infty}$-algebra, we no longer have a trivial $l_{3}$ map.
Instead, equation (\ref{megajacobi}) says that the Jacobi identity
for the 0-chains $x,y,z$ holds up to a term of the form
$dl_3(x,y,z)$.  We do, however, have $l_{4} = 0$, which provides
us with the coherence law that $l_{3}$ must satisfy.

Since we will be making frequent use of these $2$-term
$L_{\infty}$-algebras, it will be advantageous to keep track of
their ingredients.

\begin{lem} \et \label{rmk}
A $2$-term $L_{\infty}$-algebra, $V,$ consists of the following
data:

\begin{itemize}
  \item two vector spaces $V_{0}$ and
   $V_{1}$ together with a linear map $d\maps V_{1} \rightarrow V_{0},$

  \item a bilinear map $l_{2}\maps V_{i} \times V_{j}
   \rightarrow V_{i+j},$ where $0 \leq i + j \leq 1$,
   \newline which we denote more suggestively as $[\cdot, \cdot],$

  \item a trilinear map $l_{3}\maps V_{0} \times V_{0} \times
   V_{} \rightarrow V_{1}.$

\end{itemize}

These maps satisfy:

\begin{itemize}
  \item[(a)] $[x,y] = -[y,x]$,
  \item[(b)] $[x,h] = -[h,x]$,
  \item[(c)] $[h,k]=0$,
  \item[(d)] $l_{3}(x,y,z)$ is totally antisymmetric in the
    arguments $x,y,z$,
  \item[(e)] $d([x,h]) = [x,dh]$,
  \item[(f)] $[dh,k] = [h,dk]$,
  \item[(g)] $d(l_{3}(x,y,z)) = -[[x,y],z] + [[x,z],y] + [x,[y,z]]$,
  \item[(h)] $l_{3}(dh,x,y) = - [[x,y],h] + [[x,h],y] + [x,[y,h]]$,
  \item[(i)] $[l_{3}(w,x,y),z]
    + [l_{3}(w,y,z),x] + l_{3}([w,y],x,z) + l_{3}([x,z],w,y) =$
    $$[l_{3}(w,x,z),y] + [l_{3}(x,y,z),w] + l_{3}([w,x],y,z) +
    l_{3}([w,z],x,y)+ l_{3}([x,y],w,z)+ l_{3}([y,z],w,x),$$
\end{itemize}
for all $w,x,y,z \in V_{0}$ and $h, k \in V_{1}.$
\end{lem}

\noindent{\rm Proof. } Note that $(a)-(d)$ hold by equation
(\ref{antisymmetry}) of Definition \ref{L-alg} while $(e)-(i)$
follow from (\ref{megajacobi}). \qed

We notice that $(a)$ and $(b)$ are the usual skew-symmetric
properties satisfied by the bracket in a Lie algebra; $(c)$ arises
simply because there are no $2$-chains.  Properties $(e)$ and
$(f)$ tell us how the differential and bracket interact, while
condition $(g)$ says that the Jacobi identity no longer holds on
the nose, but up to chain homotopy.  We will use $(g)$ to define
the Jacobiator in the Lie $2$-algebra corresponding to a $2$-term
$L_\infty$-algebra.  Equation $(h)$ will give the naturality of
the Jacobiator.   Similarly, $(i)$ will give the Jacobiator
identity.

We continue by defining homomorphisms between $2$-term
$L_{\infty}$-algebras:

\begin{defn} \et \label{Linftyhomo}
Let $V$ and $V'$ be $2$-term $L_{\infty}$-algebras. An
\textbf{$\mathbf{L_{\infty}}$-homomorphism} \newline $\phi \maps V
\rightarrow V'$ consists of:

\begin{itemize}
  \item a chain map $\phi \maps V \to V'$ (which consists of
  linear maps $\phi_0 \maps V_0 \to V'_0$ and
              $\phi_1 \maps V_1 \to V'_1$ preserving the differential),
  \item a skew-symmetric bilinear map $\phi_{2} \maps V_{0} \times V_{0} \to
  V_{1}'$,
\end{itemize}

such that the following equations hold for all $x,y,z \in V_0$,
$h,k \in V_{1}:$

\begin{itemize}
\item $d \phi_{2}(x,y) = \phi_{0}[x,y] - [\phi_{0}(x),
\phi_{0}(y)]$
\item $\phi_{2}(d(h),d(k)) = \phi_{1}[h,k] - [\phi_{1}(h),
\phi_{1}(k)]$
\item
$ \phi_1(l_3(x,y,z)) + \phi_2(x,[y,z]) + [\phi_0(x), \phi_2(y,z)]
= $ $ l_3(\phi_0(x),\phi_0(y), \phi_0(z)) + \phi_2([x,y],z) +
\phi_2(y,[x,z]) + [\phi_2(x,y),\phi_0(z)] + [\phi_0(y),
\phi_2(x,z)]  $
\end{itemize}
\end{defn}

\noindent The reader should note the similarity between this
definition and that of homomorphisms between Lie $2$-algebras
(Definition \ref{lie2algfunct}).  In particular, the first two
equations say that $\phi_{2}$ defines a chain homotopy from
$[\phi(\cdot), \phi(\cdot)]$ to $\phi[\cdot, \cdot]$, where these
are regarded as chain maps from $V \otimes V$ to $V'$.  The third
equation in the above definition is just a chain complex version
of the commutative square in Definition \ref{lie2algfunct}.

To make $2$-term $L_\infty$-algebras and $L_\infty$-homomorphisms
between them into a category, we must describe composition and
identities. We compose a pair of $L_{\infty}$-homomorphisms $\phi
\maps V \rightarrow V'$ and $\psi \maps V' \rightarrow V''$ by
letting the chain map $\phi \psi \maps V \to V''$ be the usual
composite:
$$
\xymatrix{
    V
     \ar[rr]^{\phi}
      &&  V'
     \ar[rr]^{\psi}
      && V''}
$$

\noindent while defining $(\phi \psi)_{2}$ as follows:
$$ (\phi \psi)_2(x,y) = \psi_2(\phi_0(x),\phi_0(y)) +
                              \psi_0(\phi_2(x,y)).$$
This is just a chain complex version of how we compose
homomorphisms between Lie $2$-algebras.  The identity homomorphism
$1_V \maps V \to V$ has the identity chain map as its underlying
map, together with $(1_V)_2 =0$.

With these definitions, we obtain:

\begin{prop} \et
There is a category {\bf 2TermL$_\mathbf\infty$} with $2$-term
$L_{\infty}$-algebras as objects and $L_\infty$-homomorphisms as
morphisms.
\end{prop}

\noindent{\bf Proof. } We must show that the identity homomorphism
behaves as it should and that composition of
$L_\infty$-homomorphisms is associative. Obviously, we need only
to check these for the skew-symmetric bilinear maps, since the
chain maps satisfy these conditions.  To check the left unit law,
we consider
$$
\xymatrix{
    V
     \ar[rr]^{1_V}
      &&  V
     \ar[rr]^{\phi}
      && V'}
$$
Then $$(1_V \phi)_2 = \phi_2(1_{V_0} (x,y)) + \phi_0 (0) =
\phi_2(x,y)$$ as desired.   The proof of the right unit law is
similar.  To demonstrate associativity, consider
$$
\xymatrix{
    V
     \ar[rr]^{\tau}
      &&  V'
     \ar[rr]^{\phi}
      && V''
      \ar[rr]^{\psi}
      && V'''}
$$
On one hand,
\begin{eqnarray*}
(\tau (\phi \psi))_2 (x,y) &=& (\phi \psi)_2 (\tau_0 (x), \tau_0
(y)) + (\phi \psi)_0 (\tau _2 (x,y)) \\ &=& \psi_2 (\phi_0 (\tau
_0 (x)), \phi_0 (\tau _0 (y))) + \psi_0 (\phi _2 (\tau _0 (x),
\tau_0 (y))) + (\phi \psi)_0 (\tau _2 (x,y))
\end{eqnarray*}
while on the other hand,
\begin{eqnarray*}
((\tau \phi) \psi)_2 (x,y) &=& \psi_2 ( (\tau \phi)_0 (x), (\tau
\phi)_0 (y)) + \psi _0 ((\tau \phi)_2 (x,y)) \\
&=& \psi_2 ( (\tau \phi)_0 (x), (\tau \phi)_0 (y)) + \psi_0(
\phi_2 (\tau_0 (x), \tau_0 (y)) + \phi_0 (\tau_2 (x,y))).
\end{eqnarray*}
Recalling our conventions regarding composition, these two
expressions are identical. \qed

Next we establish the equivalence between the category of Lie
$2$-algebras and that of $2$-term $L_{\infty}$-algebras.  This
result is based on the equivalence between $2$-vector spaces and
$2$-term chain complexes described in Proposition \ref{1-1vs}.

\begin{thm} \et \label{1-1}The categories {\rm Lie$2$Alg}
and {\rm 2TermL$_{\infty}$} are equivalent.
\end{thm}

\noindent{\bf Proof. }  First we sketch how to construct a functor
$T\maps {\rm 2TermL_{\infty}} \rightarrow {\rm Lie2Alg}$.  Given a
$2$-term $L_{\infty}$-algebra $V$ we construct the Lie $2$-algebra
$L=T(V)$ as follows.

We construct the underlying $2$-vector space of $L$ as in the
proof of Proposition \ref{1-1vs}.  Thus $L$ has vector spaces of
objects and morphisms
\begin{eqnarray*}
  L_{0} & = & V_{0}, \\
  L_{1} & = & V_{0} \oplus V_{1},
\end{eqnarray*}
and we denote a morphism $f\maps x \rightarrow y$ in $L_{1}$ by
$f=(x, \vec{f})$ where $x \in V_{0}$ is the source of $f$ and
$\vec{f} \in V_{1}$ is its arrow part.  The source, target, and
identity-assigning maps in $L$ are given by
\begin{eqnarray*}
  s(f) &=& s(x, \vec{f}) = x ,\\
  t(f) &=& t(x, \vec{f}) = x + d\vec{f}, \\
  i(x) &=& (x, 0),
\end{eqnarray*}
and we have $t(f) - s(f) = d\vec{f}$.  The composite of two
morphisms in $L_{1}$ is given as in the proof of Lemma
\ref{watereddown}.  That is, given $f = (x, \vec{f})\maps x
\rightarrow y,$ and $g = (y, \vec{g})\maps y \rightarrow z$, we
have
$$f g:= (x, \vec{f} + \vec{g}).$$

We continue by equipping $L=T(V)$ with the additional structure
which makes it a Lie $2$-algebra.  First, we use the degree-zero
chain map $l_2 \maps V \otimes V \to V $ to define the bracket
functor $[\cdot, \cdot]\maps L \times L \rightarrow L.$  For a
pair of objects $x,y \in L_0$ we define $[x,y] = l_2(x,y)$, where
we use the `$l_2$' notation in the $L_\infty$-algebra $V$ to avoid
confusion with the bracket in $L$. The bracket functor is
skew-symmetric and bilinear on objects since $l_{2}$ is.  This is
not sufficient, however.  It remains to define the bracket functor
on pairs of morphisms.

We begin by defining the bracket on pairs of morphisms where one
morphism is an identity.  We do this as follows: given a morphism
$f= (x, \vec{f})\maps x \rightarrow y $ in $L_{1}$ and an object
$z \in L_{0}$, we define
$$[1_z,f]:= (l_2(z,x), l_2(z, \vec{f})),$$
$$[f,1_z]:= (l_2(x,z), l_2(\vec{f}, z)). $$
Clearly these morphisms have the desired sources; we now verify
that they also have the desired targets.  Using the fact that
$t(f) = s(f) + d\vec{f}$ for any morphism $f \in L_{1}$, we see
that:
\begin{eqnarray*}
t[1_{z},f] &=& s[1_{z},f] + dl_2(z,\vec{f}) \cr
       &=& l_2(z,x) + l_2(z, d\vec{f}) \; \; \; \; \textrm{by $(e)$ of
       Lemma \ref{rmk}} \cr
       &=& l_2(z,x) + l_2(z, y-x) \cr
       &=& l_2(z, y)
\end{eqnarray*}
as desired, using the bilinearity of $l_2$. Similarly we have:
\begin{eqnarray*}
t[f, 1_{z}] &=& s[f,1_{z}] + dl_2(\vec{f},z) \cr
       &=& l_2(x,z) + l_2(d\vec{f}, z) \; \; \; \; \textrm{by $(e)$ of
       Lemma \ref{rmk}} \cr
       &=& l_2(x,z) + l_2(y-x, z) \cr
       &=& l_2(y, z)
\end{eqnarray*}
again using the bilinearity of $l_2$.

These definitions together with the desired functoriality of the
bracket force us to define the bracket of an arbitrary pair of
morphisms $f\maps x \rightarrow y$, $g\maps a \rightarrow b$ as
follows:
\begin{eqnarray*}
[f,g] &=& [f 1_y, 1_a  g] \cr
       &:=& [f,1_a]  [1_y,g] \cr
       &=& (l_2(x,a), l_2(\vec{f}, a))  (l_2(y,a), l_2(y, \vec{g})) \cr
       &=& (l_2(x,a), l_2(\vec{f}, a) + l_2(y, \vec{g})).
\end{eqnarray*}
On the other hand, they also force us to define it as:
\begin{eqnarray*}
[f,g] &=& [1_x f, g  1_b] \\
      &:=& [1_x,g]  [f,1_b]   \\
      &=& (l_2(x,a), l_2(x,\vec{g})  (l_2(x,b), l_2(\vec{f},b)) \\
      &=& (l_2(x,a), l_2(x,\vec{g}) + l_2(\vec{f},b))  .
\end{eqnarray*}
Luckily these are compatible!  We have
\begin{equation}
  l_2(\vec{f}, a) + l_2(y, \vec{g}) =
    l_2(x,\vec{g}) + l_2(\vec{f},b)
\label{magic.equation}
\end{equation}
because the left-hand side minus the right-hand side equals
$l_2(d\vec{f},\vec{g}) - l_2(\vec{f}, d\vec{g})$, which vanishes
by $(f)$ of Lemma \ref{rmk}.

At this point we relax and use $[\cdot,\cdot]$ to stand both for
the bracket in $L$ and the $L_\infty$-algebra $V$. In this new
relaxed notation the bracket of morphisms $f \maps x \to y$, $g
\maps a \to b$ in $L$ is given by
\begin{eqnarray*}
   [f,g] &=& ([x,a], [\vec{f}, a] + [y, \vec{g}])   \\
         &=& ([x,a], [x,\vec{g}] + [\vec{f},b])  .
\end{eqnarray*}

The bracket $[\cdot,\cdot] \maps L \times L \to L$ is clearly
bilinear on objects.  Either of the above formulas shows it is
also bilinear on morphisms, since the source, target and arrow
part of a morphism depend linearly on the morphism, and the
bracket in $V$ is bilinear. The bracket is also skew-symmetric:
this is clear for objects, and can be seen for morphisms if we use
{\it both} the above formulas.

To show that $[\cdot, \cdot]\maps L \times L \rightarrow L$ is a
functor, we must check that it preserves identities and
composition. We first show that $[1_{x}, 1_{y}] = 1_{[x,y]}$,
where $x, y \in L_{0}$.  For this we use the fact that identity
morphisms are precisely those with vanishing arrow part.  Either
formula for the bracket of morphisms gives
\begin{eqnarray*}
  [1_{x}, 1_{y}] &=& ([x,y],0) \cr
                 &=& 1_{[x,y]}.
\end{eqnarray*}

\noindent To show that the bracket preserves composition, consider
the morphisms $f=(x,\vec{f}), f'= (y, \vec{f'}), g= (a, \vec{g}),$
and $g'=(b,\vec{g'})$ in $L_{1},$ where $f\maps x \rightarrow y,$
$f'\maps y \rightarrow z,$ $g\maps a \rightarrow b,$ and $g'\maps
b \rightarrow c.$ We must show
$$[f  f', g  g'] = [f,g]  [f',g'].$$
On the one hand, the definitions give
$$
  [f f', g g'] =
        ([x,a], [\vec{f}, a] +[\vec{f'},a] + [z, \vec{g}] + [z, \vec{g'}]),
$$
while on the other, they give
$$
  [f,g]  [f',g'] =
       ([x,a], [\vec{f},a] + [y, \vec{g}] + [\vec{f'},b] + [z, \vec{g'}])
$$
Therefore, it suffices to show that
$$[\vec{f'},a] + [z,\vec{g}] = [y,\vec{g}] + [\vec{f'},b].$$
After a relabelling of variables, this is just equation
(\ref{magic.equation}).

Next we define the Jacobiator for $L$ and check its properties. We
set
$$J_{x,y,z} := ([[x,y],z], l_{3}(x,y,z)).  $$
Clearly the source of $J_{x,y,z}$ is $[[x,y],z]$ as desired, while
its target is $[x, [y,z]] + [[x,z],y]$ by condition $(g)$ of Lemma
\ref{rmk}. To show $J_{x,y,z}$ is natural one can check that is
natural in each argument.  We check naturality in the third
variable, leaving the other two as exercises for the reader.  Let
$f \maps z \rightarrow z'.$ Then, $J_{x,y,z}$ is natural in $z$ if
the following diagram commutes:
$$\xymatrix{
    [[x,y],z]
     \ar[rrr]^{[[1_x,1_y],f]}
     \ar[dd]_{J_{x,y,z}}
      &&& [[x,y],z']
     \ar[dd]^{J_{x,y,z'}} \\ \\
      [[x,z],y]+ [x,[y,z]]
     \ar[rrr]^{[[1_x,f],1_y]+ [1_x,[1_y,f]]}
      &&& [[x,z'],y]+[x,[y,z']] }$$
\\
\noindent Using the formula for the composition and bracket in $L$
this means that we need
$$([[x,y],z], \vec{J}_{x,y,z'} + [[x,y],\vec{f}]) =
([[x,y],z], [[x, \vec{f}], y] + [x,[y,\vec{f}]] +
\vec{J}_{x,y,z}).$$ Thus, it suffices to show that
$$\vec{J}_{x,y,z'} + [[x,y],\vec{f}] =
[[x,\vec{f}],y] + [x,[y,\vec{f}]] + \vec{J}_{x,y,z}.$$ But
$\vec{J}_{x,y,z}$ has been defined as $l_{3}(x,y,z)$ (and
similarly for $\vec{J}_{x,y,z'}$), so now we are required to show
that:
$$l_{3}(x,y,z') + [[x,y],\vec{f}] =
l_{3}(x,y,z) + [[x,\vec{f}],y] + [x, [y,\vec{f}],$$ or in other
words,
$$ [[x,y],\vec{f}] +
l_{3}(x,y,d\vec{f}) = [[x,\vec{f}],y] + [x, [y,\vec{f}]].$$ This
holds by condition $(h)$ in Lemma \ref{rmk} together with the
complete antisymmetry of $l_3$.

The trilinearity and complete antisymmetry of the Jacobiator
follow from the corresponding properties of $l_3$. Finally,
condition $(i)$ in Lemma \ref{rmk} gives the Jacobiator identity:
 $$J_{[w,x],y,z} ([J_{w,x,z},y]+1) (J_{w, [x,z], y} +
       J_{[w,z],x,y} + J_{w,x, [y,z]}) = $$
$$[J_{w,x,y},z] (J_{[w,y],x,z} + J_{w, [x,y],z}) ([J_{w,y,z},x]+1)
([w, J_{x,y,z}]+1).$$

This completes the construction of a Lie $2$-algebra $T(V)$ from
any $2$-term $L_\infty$-algebra $V$.  Next we construct a Lie
$2$-algebra homomorphism $T(\phi) \maps T(V) \to T(V')$ from any
$L_\infty$-homomorphism $\phi \maps V \to V'$ between $2$-term
$L_\infty$-algebras.

Let $T(V) = L$ and $T(V') = L'$.  We define the underlying linear
functor of $T(\phi)=F$ as in Proposition \ref{1-1vs}.  To make $F$
into a Lie $2$-algebra homomorphism we must equip it with a
skew-symmetric bilinear natural transformation $F_{2}$ satisfying
the conditions in Definition \ref{lie2algfunct}. We do this using
the skew-symmetric bilinear map $\phi_{2} \maps V_{0} \times V_{0}
\to V_{1}'$.  In terms of its source and arrow parts, we let
$$F_{2}(x,y) = ([\phi_{0}(x), \phi_{0}(y)], \phi_{2}(x,y)).$$
Computing the target of $F_{2}(x,y)$ we have:
\begin{eqnarray*}
tF_{2}(x,y) &=& sF_{2}(x,y) + d\vec{F_{2}}(x,y) \\
            &=& [\phi_{0}(x), \phi_{0}(y)] + d\phi_{2}(x,y) \\
            &=& [\phi_{0}(x), \phi_{0}(y)] + \phi_{0}[x,y] - [\phi_{0}(x),
            \phi_{0}(y)]\\
            &=& \phi_{0}[x,y] \\
            &=& F_0[x,y]
\end{eqnarray*}
by the first equation in Definition \ref{Linftyhomo} and the fact
that $F_0 = \phi_0$. Thus we have \newline $F_{2}(x,y) \maps
[F_{0}(x), F_{0}(y)] \to F_{0}[x,y]$. Notice that $F_{2}(x,y)$ is
bilinear and skew-symmetric since $\phi_{2}$ and the bracket are.
$F_{2}$ is a natural transformation by Theorem \ref{equivof2vs}
and the fact that $\phi_2$ is a chain homotopy from $[\phi(\cdot),
\phi(\cdot)]$ to $\phi([\cdot,\cdot])$, thought of as chain maps
from $V \tensor V$ to $V'$. Finally, the equation in Definition
\ref{Linftyhomo} gives the commutative diagram in Definition
\ref{lie2algfunct}, since the composition of morphisms corresponds
to addition of their arrow parts.

We leave it to the reader to check that $T$ is indeed a functor.
Next, we describe how to construct a functor $S\maps {\rm Lie2Alg}
\to {\rm 2TermL_{\infty}}$.

Given a Lie $2$-algebra $L$ we construct the $2$-term
$L_\infty$-algebra $V = S(L)$ as follows.  We define:
\begin{eqnarray*}
  V_{0} &=& L_{0} \\
  V_{1} &=& ker(s) \subseteq L_{1}.
\end{eqnarray*}
In addition, we define the maps $l_{k}$ as follows:
\begin{itemize}
\item $l_1h = t(h)$ for $h \in V_1 \subseteq L_1$.
\item $l_{2}(x,y) = [x,y]$ for $x,y \in V_0 = L_0$.
\item $l_{2}(x,h) = -l_{2}(h,x) = [1_x, h]$
for $x \in V_0 = L_0$ and $h \in V_1 \subseteq L_1$.
\item $l_2(h,k) = 0$ for $h,k \in V_1 \subseteq L_1$.
\item $l_{3}(x,y,z) = \vec{J}_{x,y,z}$ for $x,y,z \in V_0 = L_0$.
\end{itemize}
As usual, we abbreviate $l_1$ as $d$ and $l_2$ as $[\cdot,\cdot]$.

With these definitions, conditions $(a)$ and $(b)$ of Lemma
\ref{rmk} follow from the antisymmetry of the bracket functor.
Condition $(c)$ is automatic.  Condition $(d)$ follows from the
complete antisymmetry of the Jacobiator.

For $h \in V_{1}$ and $x \in V_{0}$, the functoriality of $[\cdot,
\cdot]$ gives
\begin{eqnarray*}
   d([x,h]) &=& t([1_x,h]) \cr
            &=& [t(1_x), t(h)] \cr
            &=& [x, dh],
\end{eqnarray*}
which is $(e)$ of Lemma \ref{rmk}. To obtain $(f)$, \ref{rmk}, we
let $h\maps 0 \rightarrow x$ and $k\maps 0 \rightarrow y$ be
elements of $V_{1}$.  We then consider the following square in $L
\times L$,
$$\xymatrix{
   &
     0
    \ar[rr]^{h}
     && x \\
     0
    \ar[dd]_{k}
     & (0,0)
    \ar[rr]^{(h,1_0)}
    \ar[dd]_{(1_0,k)}
     && (x,0)
    \ar[dd]^{(1_x,k)} \\ \\
     y
     & (0,y)
    \ar[rr]^{(h,1_y)}
     && (x,y)}$$
which commutes by definition of a product category. Since $[\cdot,
\cdot]$ is a functor, it preserves such commutative squares, so
that
$$\xymatrix{
    [0,0]
     \ar[rr]^{[h,1_0]}
     \ar[dd]_{[1_0,k]}
      && [x,0]
     \ar[dd]^{[1_x,k]} \\ \\
      [0,y]
     \ar[rr]^{[h,1_y]}
      && [x,y]}$$
commutes.  Since $[h,1_0]$ and $[1_0,k]$ are easily seen to be
identity morphisms, this implies $[h,1_y] =[1_x,k]$.  This means
that in $V$ we have $[h,y] = [x,k]$, or, since $y$ is the target
of $k$ and $x$ is the target of $h$, simply $[h,dk]=[dh,k],$ which
is $(f)$ of Lemma \ref{rmk}.

Since $J_{x,y,z} \maps [[x,y],z] \to [x,[y,z]] + [[x,z],y]$, we
have
\begin{eqnarray*}
   d(l_{3}(x,y,z)) &=& t(\vec{J}_{x,y,z}) \cr
                   &=& (t-s)(J_{x,y,z}) \cr
                   &=& [x,[y,z]] + [[x,z],y] - [[x,y],z],
\end{eqnarray*}
which gives $(g)$ of Lemma \ref{rmk}.  The naturality of
$J_{x,y,z}$ implies that for any $f \maps z \to z'$, we must have
$$[[1_x,1_y],f] \circ J_{x,y,z'} = J_{x,y,z} \circ
([[1_x,f],1_y] +[1_x,[1_y,f]]).$$ This implies that in $V$ we have
$$[[x,y],\vec{f}] + l_{3}(x,y, z'-z) = [[x,\vec{f}],y] + [x, [y,\vec{f}]],$$
for $x,y \in V_{0}$ and $\vec{f} \in V_{1},$ which is $(h)$ of
Lemma \ref{rmk}.

Finally, the Jacobiator identity dictates that the arrow part of
the Jacobiator, $l_3$, satisfies the following equation:
$$[l_{3}(w,x,y),z] + [l_{3}(w,y,z),x] + l_{3}([w,y],x,z) +
l_{3}([x,z],w,y) =$$
$$[l_{3}(w,x,z),y] + [l_{3}(x,y,z),w] + l_{3}([w,x],y,z)
+ l_{3}([w,z],x,y)+ l_{3}([x,y],w,z)+ l_{3}([y,z],w,x).$$ This is
$(i)$ of Lemma \ref{rmk}.

This completes the construction of a $2$-term $L_\infty$-algebra
$S(L)$ from any Lie $2$-algebra $L$.  Next we construct an
$L_\infty$-homomorphism $S(F) \maps S(L) \to S(L')$ from any Lie
$2$-algebra homomorphism $F \maps L \to L'$.

Let $S(L) = V$ and $S(L') = V'$.  We define the underlying chain
map of $S(F) = \phi$ as in Proposition \ref{1-1vs}. To make $\phi$
into an $L_\infty$-homomorphism we must equip it with a
skew-symmetric bilinear map $\phi_{2} \maps V_{0} \times V_{0} \to
V_{1}'$ satisfying the conditions in Definition \ref{Linftyhomo}.
To do this we set
$$
\phi_{2}(x,y) = \vec{F}_{2}(x,y).
$$
The bilinearity and skew-symmetry of $\phi_{2}$ follow from that
of $F_{2}$.  Then, since $\phi_{2}$ is the arrow part of $F_{2}$,
\begin{eqnarray*}
d\phi_{2}(x,y) &=& (t-s)F_{2}(x,y) \\
&=& F_{0}[x,y] - [F_{0}(x), F_{0}(y)] \\
&=& \phi_{0}[x,y] - [\phi_{0}(x), \phi_{0}(y)],
\end{eqnarray*}
by definition of the chain map $\phi$.  The naturality of $F_{2}$
gives the second equation in Definition \ref{Linftyhomo}. Finally,
since composition of morphisms corresponds to addition of arrow
parts, the diagram in Definition \ref{lie2algfunct} gives:
$$ \phi_1(l_3(x,y,z)) + \phi_2(x,[y,z]) + [\phi_0(x), \phi_2(y,z)] =
 l_3(\phi_0(x),\phi_0(y), \phi_0(z)) +$$ $$ \phi_2([x,y],z) +
\phi_2(y,[x,z]) + [\phi_2(x,y),\phi_0(z)] + [\phi_0(y),
\phi_2(x,z)],  $$ since $\phi_{0} = F_{0}$, $\phi_{1} = F_{1}$ on
elements of $V_{1}$, and the arrow parts of $J$ and $F_{2}$ are
$l_{3}$ and $\phi_{2}$, respectively.

We leave it to the reader to check that $S$ is indeed a functor,
and to construct natural isomorphisms $\alpha \maps ST \To 1_{\rm
Lie2Alg}$ and $\beta \maps TS \To 1_{{\rm 2TermL}_\infty}$. \qed


The above theorem also has a $2$-categorical version. We have
defined a $2$-category of Lie $2$-algebras, but not yet defined a
$2$-category of $2$-term $L_\infty$-algebras. For this, we need
the following:

\begin{defn} \et \label{Linfty2homo} Let $V$ and $V'$ be
$2$-term $L_\infty$-algebras and let $\phi, \psi \maps V \to V'$
be $L_{\infty}$-homomorphisms.  An {\bf $\mathbf
L_{\infty}$-2-homomorphism} $\tau \maps \phi \To \psi$ is a chain
homotopy such that the following equation holds for all $x,y \in
V_0$:
\begin{itemize}
\item $\phi_2(x,y) + \tau_{[x,y]} = [\tau_x, \tau_y] + \psi_2(x,y)$
\end{itemize}
\end{defn}
\noindent
 We now define vertical and horizontal
composition for these $2$-homomorphisms.  First let $V$ and $V'$
be $2$-term $L_\infty$-algebras and let $\phi, \psi, \gamma \maps
V \to V'$ be $L_{\infty}$-homomorphisms.  If $\theta \maps \phi
\To \psi$ and $\tau \maps \psi \To \gamma$ are
$L_{\infty}$-$2$-homomorphisms, we define their {\bf vertical}
composite, $\theta \tau \maps \phi \to \gamma,$ by
$$\theta \tau(x) : = \theta(x) + \tau(x).$$

Next, let $V, V', V''$ be $2$-term $L_\infty$-algebras and let
$\phi, \psi \maps V \to V'$ and \newline $\phi ' , \psi ' \maps V'
\to V''$ be $L_{\infty}$-homomorphisms.  If $\tau \maps \phi \To
\psi$ and $\tau ' \maps \phi ' \To \psi '$ are
$L_{\infty}$-$2$-homomorphisms, we define their {\bf horizontal
composite}, $\tau \circ \tau ' \maps \phi \phi ' \To \psi \psi ',$
in either of two equivalent ways:
\begin{eqnarray*}
\tau \circ \tau ' (x)&:=& \tau '(\phi_0 (x)) + \psi_1 ' (\tau(x))
\\
&=& \phi_1 ' (\tau(x)) + \tau ' (\psi_0 (x)).
\end{eqnarray*}

Finally, given a $L_{\infty}$-homomorphism $\phi \maps V
\rightarrow V'$,  the identity $L_{\infty}$-$2$-homomorphism
$1_\phi \maps \phi \Rightarrow \phi$ is given by $1_\phi (x) = 1 _
{\phi_0 (x)}$.

With these definitions, it becomes a straightforward exercise to
prove the following:

\begin{prop} \et There is a strict $2$-category
{\bf 2TermL$_\mathbf\infty$} with $2$-term $L_{\infty}$-algebras
as objects, homomorphisms between these as morphisms, and
$2$-homomorphisms between those as $2$-morphisms.
\end{prop}

With these definitions one can strengthen Theorem \ref{1-1} as
follows:

\begin{thm} \label{1-1'} \et The $2$-categories {\rm Lie$2$Alg}
and {\rm 2TermL$_{\infty}$} are $2$-equivalent.
\end{thm}

The main benefit of the results in this section is that they
provide us with another method to create examples of Lie
$2$-algebras.  Instead of thinking of a Lie $2$-algebra as a
category equipped with extra structure, we may work with a
$2$-term chain complex endowed with the structure described in
Lemma \ref{rmk}.  In the next two sections we investigate the
results of trivializing various aspects of a Lie $2$-algebra, or
equivalently of the corresponding $2$-term $L_\infty$-algebra.


\section{Strict Lie 2-algebras} \label{strictlie2algs}

As the reader may expect, a `strict' Lie $2$-algebra is a mixture
of a Lie algebra and category in which all laws hold on the nose
as equations, not just up to isomorphism.  In one of his papers
\cite{Baez}, Baez showed how to construct these starting from
`strict Lie $2$-groups'. Here we describe this process in a
somewhat more highbrow manner, and explain how these `strict'
notions are special cases of the semistrict ones described here.

Since we only weakened the Jacobi identity in our definition of
`semistrict' Lie $2$-algebra, we need only require that the
Jacobiator be the identity to recover the `strict' notion:

\begin{defn} \et A semistrict Lie $2$-algebra is {\bf strict}
if its Jacobiator is the identity.
\end{defn}

\noindent Similarly, requiring that the bracket be strictly
preserved gives the notion of `strict' homomorphism between Lie
$2$-algebras:

\begin{defn} \et Given semistrict Lie $2$-algebras $L$ and $L'$, a
homomorphism $F \maps L \to L'$ is {\bf strict} if $F_2$ is the
identity.
\end{defn}

\begin{prop} \et {\rm Lie2Alg}
contains a sub-$2$-category {\bf SLie2Alg} with strict Lie
$2$-algebras as objects, strict homomorphisms between these as
morphisms, and arbitrary $2$-homomorphisms between those as
$2$-morphisms.
\end{prop}

\noindent {\bf Proof. }  One need only check that if the
homomorphisms $F \maps L \to L'$ and $G \maps L' \to L''$ have
$F_2 = 1$, $G_2 = 1$, then their composite has $(FG)_2 = 1$. \qed

The following proposition shows that strict Lie $2$-algebras as
defined here agree with those as defined previously \cite{Baez}:

\begin{prop} \et The $2$-category {\rm SLie2Alg}
is isomorphic to the $2$-category $\LieAlg\Cat$ consisting of
categories, functors and natural transformations in $\LieAlg$.
\end{prop}

\noindent {\bf Proof. }  This is just a matter of unravelling the
definitions.  \qed

Just as Lie groups have Lie algebras, `strict Lie $2$-groups' have
strict Lie $2$-algebras.  Before we can state this result
precisely, we must recall the concept of a strict Lie $2$-group,
which was treated in greater detail in HDA5 \cite{BLau}:

\begin{defn} \et \label{SLie2Grp} We define {\bf SLie2Grp} to be
the strict $2$-category $\LieGrp\Cat$ consisting of categories,
functors and natural transformations in $\LieGrp$.  We call the
objects in this $2$-category {\bf strict Lie 2-groups}; we call
the morphisms between these {\bf strict homomorphisms}, and we
call the $2$-morphisms between those {\bf 2-homomorphisms}.
\end{defn}

\begin{prop} \et  There exists a unique $2$-functor
$$d \maps {\rm SLie2Grp} \to {\rm SLie2Alg}$$
such that:
\begin{enumerate}

\item $d$ maps any strict Lie $2$-group $C$ to the strict Lie $2$-algebra
$dC = \mathfrak{c}$ for which $\mathfrak{c}_0$ is the Lie algebra
of the Lie group of objects $C_0$, $\mathfrak{c}_1$ is the Lie
algebra of the Lie group of morphisms $C_1$, and the maps $s,t
\maps \mathfrak{c}_1 \rightarrow \mathfrak{c}_0$, $i \maps
\mathfrak{c}_0 \rightarrow \mathfrak{c}_1$ and $\circ \maps
\mathfrak{c}_1 \times_{\mathfrak{c}_0} \mathfrak{c}_1 \rightarrow
\mathfrak{c}_1$ are the differentials of those for $C$.

\item $d$ maps any strict Lie $2$-group homomorphism
$F \maps C \rightarrow C'$ to the strict Lie $2$-algebra
homomorphism $dF\maps \mathfrak{c} \rightarrow \mathfrak{c}'$ for
which $(dF)_0$ is the differential of $F_0$ and $(dF)_1$ is the
differential of $F_1$.

\item $d$ maps any strict Lie $2$-group $2$-homomorphism $\alpha \maps F \To G$
where $F,G \maps C \to C'$ to the strict Lie $2$-algebra
$2$-homomorphism $d\alpha \maps dF \To dG$ for which the map
$d\alpha \maps \mathfrak{c}_0 \to \mathfrak{c}_1$ is the
differential of $\alpha \maps C_0 \to C_1$.
\end{enumerate}
\end{prop}

\noindent{Proof. } The proof of this long-winded proposition is a
quick exercise in internal category theory: the well-known functor
from $\LieGrp$ to $\LieAlg$ preserves pullbacks, so it maps
categories, functors and natural transformations in $\LieGrp$ to
those in $\LieAlg$, defining a $2$-functor $d \maps {\rm SLie2Grp}
\to {\rm SLie2Alg}$, which is given explicitly as above.  \qed

We would like to generalize this theorem and define the Lie
$2$-algebra not just of a strict Lie $2$-group, but of a general
Lie $2$-group as defined in HDA5 \cite{BLau}.  However, this may
require a weaker concept of Lie $2$-algebra than that studied
here.


\section{Skeletal Lie 2-algebras} \label{skeletallie2algs}

A semistrict Lie $2$-algebra is {\it strict} when the Jacobiator
is the identity, which means that the map $l_3$ vanishes in the
corresponding $L_\infty$-algebra.  We now investigate the
consequences of assuming the differential $d$ vanishes in the
corresponding $L_\infty$-algebra. Thanks to the formula
$$   d\vec{f} =  t(f) - s(f)   , $$
this implies that the source of any morphism in the Lie
$2$-algebra equals its target.  In other words, the Lie
$2$-algebra is {\it skeletal:}

\begin{defn} \et
A category is {\bf skeletal} if isomorphic objects are equal.
\end{defn}

Skeletal categories are useful in category theory because every
category is equivalent to a skeletal one formed by choosing one
representative of each isomorphism class of objects \cite{Mac}.
The same sort of thing is true in the context of $2$-vector
spaces:

\begin{lem}\label{skeletal2vs}\et
Any $2$-vector space is equivalent, as an object of $2\Vect$, to a
skeletal one.
\end{lem}

\noindent {\bf Proof.} Using the result of Theorem
\ref{equivof2vs} we may treat our $2$-vector spaces as $2$-term
chain complexes.  In particular, a $2$-vector space is skeletal if
the corresponding $2$-term chain complex has vanishing
differential, and two $2$-vector spaces are equivalent if the
corresponding $2$-term chain complexes are chain homotopy
equivalent. So, it suffices to show that any $2$-term chain
complex is chain homotopy equivalent to one with vanishing
differential. This is well-known, but the basic idea is as
follows. Given a $2$-term chain complex
$$\xymatrix{
   C_{1}
   \ar[rr]^{d}
    && C_{0}}
$$
we express the vector spaces $C_{0}$ and $C_{1}$ as $C_{0} = im(d)
\oplus C_{0}'$ and $C_{1} = ker(d) \oplus X$ where $X$ is a vector
space complement to $ker(d)$ in $C_{1}$. This allows us to define
a $2$-term chain complex $C'$ with vanishing differential:
$$\xymatrix{
   C_{1}'=ker(d)
   \ar[rr]^{0}
   && C_{0}'}
.$$ The inclusion of $C'$ in $C$ can easily be extended to a chain
homotopy equivalence. {\hbox{\hskip 30em} \qed}

Using this fact we obtain a result that will ultimately allow us
to classify Lie $2$-algebras:

\begin{prop} \et \label{skeletal}
Every Lie $2$-algebra is equivalent, as an object of {\rm
Lie2Alg}, to a skeletal one.
\end{prop}

\noindent {\bf Proof.} Given a Lie $2$-algebra $L$ we may use
Lemma \ref{skeletal2vs} to find an equivalence between the
underlying $2$-vector space of $L$ and a skeletal $2$-vector space
$L'$.  We may then use this to transport the Lie $2$-algebra
structure from $L$ to $L'$, and obtain an equivalence of Lie
$2$-algebras between $L$ and $L'$. \qed

It is interesting to observe that a skeletal Lie $2$-algebra that
is also strict amounts to nothing but a Lie algebra $L_{0}$
together with a representation of $L_{0}$ on a vector space
$L_{1}$. This is the infinitesimal analogue of how a strict
skeletal $2$-group $G$ consists of a group $G_0$ together with an
action of $G_0$ as automorphisms of an abelian group $G_1$.  Thus,
the representation theory of groups and Lie algebras is
automatically subsumed in the theory of $2$-groups and Lie
$2$-algebras!

To generalize this observation to other skeletal Lie $2$-algebras,
we recall some basic definitions concerning Lie algebra
cohomology:


\begin{defn} \et
Let $\mathfrak{g}$ be a Lie algebra and $\rho$ a representation of
$\mathfrak{g}$ on the vector space $V$.  Then a \textbf{V-{\bf
valued} n-{\bf cochain}} $\mathbf{\omega}$ on $\mathfrak{g}$ is a
totally antisymmetric map
$$\omega \maps \mathfrak{g}^{\otimes n} \to V.$$
The vector space of all $n$-cochains is denoted by
$C^{n}(\mathfrak{g},V)$. The {\bf coboundary operator}
$\delta\maps C^{n}(\mathfrak{g},V) \rightarrow
C^{n+1}(\mathfrak{g},V)$ is defined by:
\begin{eqnarray*}
 (\delta \omega) (v_{1}, v_{2}, \dots, v_{n+1}) &:=&
    \sum_{i=1} ^{n+1} (-1)^{i+1} \rho (v_{i})
    \omega_{n} (v_{1}, \dots, \hat{v_{i}}, \dots, v_{n+1}) \cr
 &+& \sum _{1 \leq j < k \leq n+1} (-1)^{j+k} \omega _{n}
    ([v_{j}, v_{k}], v_{1}, \dots, \hat{v_{j}}, \dots,
    \hat{v_{k}}, \dots, v_{n+1}),
\end{eqnarray*}
\end{defn}

\begin{prop} \et The Lie algebra coboundary operator $\delta$
satisfies $\delta ^{2} = 0$.
\end{prop}

\begin{defn} \et
A $V$-valued $n$-cochain $\omega$ on $\mathfrak{g}$ is an {\bf
$n$-cocycle} when $\delta \omega = 0$ and an {\bf $n$-coboundary}
if there exists an $(n-1)$-cochain $\theta$ such that $\omega =
\delta \theta.$  We denote the groups of $n$-cocycles and
$n$-coboundaries by $Z^{n}(\mathfrak{g},V)$ and
$B^{n}(\mathfrak{g},V)$ respectively. The $n$th {\bf Lie algebra
cohomology group} $H^{n}(\mathfrak{g},V)$ is defined by
$$H^{n}(\mathfrak{g},V) = Z^{n}(\mathfrak{g},V)/B^{n}(\mathfrak{g},V).$$
\end{defn}

The following result illuminates the relationship between Lie
algebra cohomology and $L_{\infty}$-algebras.

\begin{thm} \et \label{trivd}
There is a one-to-one correspondence between $L_{\infty}$-algebras
consisting of only two nonzero terms $V_{0}$ and $V_{n}$, with
$d=0,$ and quadruples $(\mathfrak{g}, V, \rho, l_{n+2})$ where
$\mathfrak{g}$ is a Lie algebra, $V$ is a vector space, $\rho$ is
a representation of $\mathfrak{g}$ on $V$, and $l_{n+2}$ is a
$(n+2)$-cocycle on $\mathfrak{g}$ with values in $V$.
\end{thm}

\noindent {\bf Proof. }

\noindent Given such an $L_{\infty}$-algebra $V$ we set
$\mathfrak{g}= V_{0}$. $V_0$ comes equipped with a bracket as part
of the $L_{\infty}$-structure, and since $d$ is trivial, this
bracket satisfies the Jacobi identity on the nose, making
$\mathfrak{g}$ into a Lie algebra. We define $V = V_{n},$ and note
that the bracket also gives a map $\rho\maps \mathfrak{g} \otimes
V \rightarrow V$, defined by $\rho(x)f = [x,f]$ for $x \in
\mathfrak{g}, f \in V$. We have
\begin{eqnarray*}
  \rho ([x,y])f &=& [[x,y],f] \cr
                &=& [[x,f],y] + [x,[y,f]] \; \; \; \;
                    \textrm{by $(2)$ of Definition \ref{L-alg}} \cr
                &=& [\rho(x)f, y] + [x, \rho(y) f]
\end{eqnarray*}

\noindent for all $x,y \in \mathfrak{g}$ and $f \in V$, so that
$\rho$ is a representation.  Finally, the $L_{\infty}$ structure
gives a map $l_{n+2}\maps \mathfrak{g}^{\otimes(n+2)} \rightarrow
V$ which is in fact a $(n+2)$-cocycle. To see this, note that
$$0 = \sum_{i+j = n+4} \sum_{\sigma}
l_{j}(l_{i}(x_{\sigma(1)}, \ldots, x_{\sigma(i)}),
x_{\sigma(i+1)}, \ldots, x_{\sigma(n+2)})$$ where we sum over $(i,
(n+3)-i)$-unshuffles $\sigma \in S_{n+3}$. However, the only
choices for $i$ and $j$ that lead to nonzero $l_{i}$ and $l_{j}$
are $i=n+2, j=2$ and $i=2, j=n+2.$  In addition, notice that in
this situation, $\chi(\sigma)$ will consist solely of the sign of
the permutation because all of our $x_{i}$'s have degree zero.
Thus, the above becomes, with $\sigma$ a $(n+2, 1)$-unshuffle and
$\tau$ a $(2, n+1)$-unshuffle:

\begin{eqnarray*}
  0 &=& \sum_{\sigma} \chi(\sigma) (-1)^{n+2}
     [l_{n+2}(x_{\sigma(1)}, \dots, x_{\sigma(n+2)}), x_{\sigma(n+3)}] \cr \\
     && + \sum_{\tau} \chi(\tau) l_{n+2}([x_{\tau(1)}, x_{\tau(2)}],
     x_{\tau(3)}, \dots, x_{\tau(n+3)}) \cr  \\
    &=& \sum_{i=1}^{n+3} (-1)^{n+3-i}(-1)^{n+2}
     [l_{n+2}(x_{1}, \dots, x_{i-1}, x_{i+1},\dots, x_{n+3}), x_{i}] \cr  \\
    & & + \sum_{1 \leq i < j \leq n+3} (-1)^{i+j+1}
     l_{n+2}([x_{i}, x_{j}], x_{1}, \dots,\hat{x_{i}},
     \dots, \hat{x_{j}}, \dots, x_{n+3})
     \qquad \qquad (\dag) \cr \\
    &=& \sum_{i=1}^{n+3} (-1)^{i+1} [l_{n+2}(x_{1}, \dots, x_{i-1},
     x_{i+1}, \dots, x_{n+3}), x_{i}] \cr  \\
    & & + \sum_{1 \leq i < j \leq n+3} (-1)^{i+j+1}
     l_{n+2}([x_{i}, x_{j}], x_{1}, \dots, \hat{x_{i}},
    \dots, \hat{x_{j}}, \dots, x_{n+3}) \cr \\
    &=& - \sum_{i=1}^{n+3} (-1)^{i+1} [x_{i},
     l_{n+2}(x_{1}, \dots, x_{i-1}, x_{i+1}, \dots, x_{n+3})]  \cr  \\
    & &  - \sum_{1 \leq i < j \leq n+3} (-1)^{i+j}
     l_{n+2}([x_{i}, x_{j}], x_{1}, \dots, \hat{x_{i}}, \dots,
      \hat{x_{j}}, \dots, x_{n+3})  \cr \\
    &=& -\delta l_{n+2}(x_{1}, x_{2}, \dots, x_{n+3}).
\end{eqnarray*}

\noindent The first sum in $(\dag)$ follows because we have
$(n+3)$ $(n+2,1)$-unshuffles and the sign of any such unshuffle is
$(-1)^{n+3-i}$. The second sum follows similarly because we have
$(n+3)$ $(2, n+1)$-unshuffles and the sign of a
$(2,n+1)$-unshuffle is $(-1)^{i+j+1}.$  Therefore, $l_{n+2}$ is a
$(n+2)$-cocycle.

\vskip 1em  Conversely, given a Lie algebra $\mathfrak{g}$, a
representation $\rho$ of $\mathfrak{g}$ on a vector space $V$, and
an $(n+2)$-cocycle $l_{n+2}$ on $\mathfrak{g}$ with values in $V$,
we define our $L_{\infty}$-algebra $V$ by setting $V_{0} =
\mathfrak{g}$, $V_{n} = V$, $V_{i}=\{0\}$ for $i \neq 0, n$ and
$d=0.$ It remains to define the system of linear maps $l_{k}$,
which we do as follows: Since $\mathfrak{g}$ is a Lie algebra, we
have a bracket defined on $V_{0}$. We extend this bracket to
define the map $l_2$, denoted by $[\cdot, \cdot] \maps V_{i}
\otimes V_{j} \rightarrow V_{i+j}$ where $i,j=0,n,$ as follows:
$$[x,f] = \rho (x) f,$$
$$[f,y] = - \rho(y)f,$$
$$[f,g] =0$$
for $x,y \in V_0$ and $f,g \in V_n$. With this definition, the map
$[\cdot, \cdot]$ satisfies condition $(1)$ of Definition
\ref{L-alg} . We define $l_{k}=0$ for $3 \leq k \leq n+1$ and $k>
n+2$, and take $l_{n+2}$ to be the given $(n+2)$ cocycle, which
satisfies conditions $(1)$ and $(2)$ of Definition \ref{L-alg} by
the cocycle condition. \qed

We can classify skeletal Lie $2$-algebras using the above
construction with $n=1$:

\begin{cor} \et \label{class} There is a one-to-one correspondence
between isomorphism classes of skeletal Lie $2$-algebras and
isomorphism classes of quadruples consisting of a Lie algebra
$\mathfrak{g}$, a vector space $V$, a representation $\rho$ of
$\mathfrak{g}$ on $V$, and a 3-cocycle on $\mathfrak{g}$ with
values in $V$.
\end{cor}

\noindent {\bf Proof. }  This is immediate from Theorem \ref{1-1}
and Theorem \ref{trivd}.  \qed

Since every Lie $2$-algebra is equivalent as an object of {\rm
Lie2Alg} to a skeletal one, this in turn lets us classify {\it
all} Lie $2$-algebras, though only up to equivalence:

\begin{thm} \et \label{class2} There is a one-to-one correspondence
between equivalence classes of Lie $2$-algebras (where equivalence
is as objects of the $2$-category {\rm Lie2Alg}) and isomorphism
classes of quadruples consisting of a Lie algebra $\mathfrak{g}$,
a vector space $V$, a representation $\rho$ of $\mathfrak{g}$ on
$V$, and an element of $H^3(\mathfrak{g},V)$.
\end{thm}

\noindent {\bf Proof. }  This follows from Theorem \ref{skeletal}
and Corollary \ref{class}; we leave it to the reader to verify
that equivalent skeletal Lie $2$-algebras give cohomologous
3-cocycles and conversely.  \qed

We conclude with perhaps the most interesting examples of
finite-dimensional Lie $2$-algebras coming from Theorem
\ref{class}. These make use of the following identities involving
the Killing form $\langle x,y\rangle := \tr(\ad(x)\ad(y)) $ of a
finite-dimensional Lie algebra:
$$\langle x, y \rangle = \langle y, x \rangle,$$
and
$$\langle [x,y], z \rangle = \langle x, [y,z] \rangle.$$

\begin{example} \label{ghbar} \et There is a skeletal Lie $2$-algebra
built using Theorem \ref{class} by taking $V_{0} = \mathfrak{g}$
to be a finite-dimensional Lie algebra over the field $k$, $V_{1}$
to be $k$, $\rho$ the trivial representation, and $l_3(x,y,z) =
\langle x,[y,z]\rangle $. We see that $l_3$ is a $3$-cocycle using
the above identities as follows:
\begin{eqnarray*}
  (\delta l_3)(w,x,y,z) &=& \rho(w) l_3(x,y,z) - \rho(x)l_3(w,y,z) +
                          \rho(y)l_3(w,x,z) - \rho(z) l_3(w,x,y) \cr
           & & -l_3([w,x],y,z) + l_3([w,y],x,z) -l_3([w,z],x,y) \cr
           & & -l_3([x,y],w,z) + l_3([x,z],w,y) -l_3([y,z],w,x) \cr
                      &=& -\langle [w,x],[y,z] \rangle
                          +\langle [w,y],[x,z] \rangle
                          -\langle [w,z],[x,y] \rangle \cr
                      & & -\langle [x,y],[w,z] \rangle
                          +\langle [x,z],[w,y] \rangle
                          -\langle [y,z],[w,x] \rangle
\end{eqnarray*}

\noindent This second step above follows because we have a trivial
representation. Continuing on, we have
\begin{eqnarray*}
   (\delta l_3)(w,x,y,z)  &=& -2\langle [w,x],[y,z] \rangle
                            +2\langle [w,y],[x,z] \rangle
                            -2\langle [w,z],[x,y] \rangle \cr
                        &=& -2\langle w,[x,[y,z]] \rangle
                            +2\langle w,[y,[x,z]] \rangle
                            -2\langle w,[z,[x,y]] \rangle \cr
                        &=& -2\langle w,[x,[y,z]] + [y,[z,x]]
                                +[z,[x,y]]\rangle  \cr
                        &=& -2\langle w, 0 \rangle   \cr
                        &=& 0. \; \; \;
\end{eqnarray*}
More generally, we obtain a Lie $2$-algebra this way taking
$l_3(x,y,z) = \hbar \langle x,[y,z]\rangle $ where $\hbar$ is any
element of $k$. We call this Lie $2$-algebra $\mathfrak{g}_\hbar$.
\end{example}

It is well known that the Killing form of $\mathfrak{g}$ is
nondegenerate if and only if $\mathfrak{g}$ is semisimple. In this
case the 3-cocycle described above represents a nontrivial
cohomology class when $\hbar \ne 0$, so by Theorem \ref{class2}
the Lie $2$-algebra $\mathfrak{g}_\hbar$ is not equivalent to a
skeletal one with vanishing Jacobiator. In other words, we obtain
a Lie $2$-algebra that is not equivalent to a skeletal strict one.

Suppose the field $k$ has characteristic zero, the Lie algebra
$\mathfrak{g}$ is finite dimensional and semisimple, and $V$ is
finite dimensional.  Then a version of Whitehead's Lemma
\cite{AIP} says that $H^3(\mathfrak{g},V) = \{0\}$ whenever the
representation of $\mathfrak{g}$ on $V$ is nontrivial and
irreducible. This places some limitations on finding interesting
examples of nonstrict Lie $2$-algebras other than those of the
form $\mathfrak{g}_\hbar$.

We expect the Lie $2$-algebras $\mathfrak{g}_\hbar$ to be related
to quantum groups, affine Lie algebras and other constructions
that rely crucially on the 3-cocycle $\langle x, [y,z] \rangle$ or
the closely related $2$-cocycle on $\mathfrak{g}[z,z^{-1}]$. The
smallest example comes from  $\mathfrak{su}(2)$. Since
$\mathfrak{su}(2)$ is isomorphic to $\R^3$ with its usual vector
cross product, and its Killing form is proportional to the dot
product, this Lie $2$-algebra relies solely on familiar properties
of the dot product and cross product: $$x \times y = - y \times
x,$$
$$x \cdot y = y \cdot x,$$
$$x \cdot (y \times z) = (x \times y) \cdot z,$$
$$x \times (y \times z) + y \times (z \times x) + z \times (x \times y) =0.$$
It will be interesting to see if this Lie $2$-algebra, where the
Jacobiator comes from the triple product, has any applications to
physics. Just for fun, we work out the details again in this case:

\begin{example} \et There is a skeletal Lie $2$-algebra
built using Theorem \ref{class} by taking $V_{0} = \R^3$ equipped
with the cross product, $V_{1} = \R$, $\rho$ the trivial
representation, and $l_3(x,y,z) = x \cdot (y \times z)$. We see
that $l_3$ is a 3-cocycle as follows:

\begin{eqnarray*}
   (\delta l_3)(w,x,y,z) &=&
      -l_3([w,x],y,z) + l_3([w,y],x,z) -l_3([w,z],x,y) \cr
                & & -l_3([x,y],w,z) + l_3([x,z],w,y) -l_3([y,z],w,x) \cr
                       &=& -(w \times x) \cdot (y \times z)
                           +(w \times y) \cdot (x \times z)
                           -(w \times z) \cdot (x \times y) \cr
                       & & -(x \times y) \cdot (w \times z)
                           +(x \times z) \cdot (w \times y)
                           -(y \times z) \cdot (w \times x) \cr
                       &=& -2(w \times x) \cdot (y \times z)
                           +2(w \times y) \cdot (x \times z)
                           -2(w \times z) \cdot (x \times y) \cr
                       &=& -2w \cdot (x \times (y \times z))
                           +2w \cdot (y \times (x \times z))
                           -2w \cdot (z \times (x \times y)) \cr
                       &=& -2w \cdot (x \times (y \times z)
                     +y \times (z \times x) + z \times (x \times y)) \cr
                      &=& 0.
\end{eqnarray*}

\end{example}

We now shift gears by moving away from the categorification
process and considering the passage from Lie groups to Lie
algebras. In the next chapter, we describe a new method of
obtaining the Lie algebra of a Lie group using algebraic
structures called `quandles'.

\chapter{Lie Theory, Quandles and Braids} \label{ch2}

Now that we have a theory of Lie $2$-algebras as well as that of
Lie $2$-groups, it becomes natural to wonder if the standard
results of Lie theory have higher-dimensional analogues.  The
obvious first question is whether Lie $2$-algebras arise from Lie
$2$-groups. Of course, as mentioned in the previous chapter, this
may require a weaker version of our semistrict Lie $2$-algebras.
It is our hope that the correct definition of a weak Lie
$2$-algebra will present itself as we unravel the passage of Lie
$2$-group to Lie $2$-algebra. Once we have established that Lie
$2$-groups give rise to Lie $2$-algebras, we would then like to
categorify the exponential map, and determine whether or not a
homomorphism between Lie $2$-groups induces a homomorphism between
the corresponding weak Lie $2$-algebras.

In this chapter we begin to explore the first question above.  In
this direction, we describe a novel method of obtaining Lie
algebras from Lie groups, which we intend to categorify to show
that every Lie $2$-algebra arises from a Lie $2$-group.  The
standard procedure for obtaining a Lie algebra from a Lie group
consists of taking the tangent space of the Lie group at the
identity element. Moreover, one often then shows that the space of
left-invariant vector fields on the Lie group is isomorphic as a
vector space to this tangent space at the identity.  As we do not,
yet, have notions of categorified vector fields, we seek an
alternate route.

Since the Lie algebra is a tangent space of the Lie group, it is
natural to obtain operations on the Lie algebra from those on the
Lie group by differentiation.  Considering that negatives in the
Lie algebra correspond to inverses in the Lie group, one may be
tempted to think that the bracket in the Lie algebra arises from
differentiating the multiplication.  However, it turns out that
the differential of the product map is simply the addition
operation in the Lie algebra.  The bracket results from
differentiating conjugation {\it twice}, which can be seen in the
following computation: Let $G$ be a Lie group and $\mathfrak{g}$
its Lie algebra.  Consider curves $exp(sx)$ and $exp(ty)$ in $G$
where $x, y \in \mathfrak{g}$.  Then, if $G$ is a matrix Lie
group, we have:
$$[x,y] = \frac{d^2}{dsdt} (e^{sx} e^{ty} e^{-sx}) |_{s=t=0}$$
Therefore, it is the operation of {\it conjugation}, and not
multiplication, in the Lie group which gives the bracket
operation. Our new description of the passage from a Lie group to
its Lie algebra captures this key aspect.

The first step in this process involves focusing our attention on
the conjugation operation.  Consider two elements $g$ and $h$ in
some group $G$.  If we denote left conjugation by $g$, $ghg^{-1}$,
as $g \rhd h$ and right conjugation by $g$, $g^{-1}hg$, as $h \lhd
g$ then Joyce has shown that all equational laws involving only
left and right conjugation can be derived from these three
\cite{J}:
\begin{itemize}
\item (idempotence) $g \rhd g = g$
\item (inverses) $(g \rhd h) \lhd g = h = g \rhd (h \lhd g)$
\item (self-distributivity) $g \rhd (h \rhd k) = (g \rhd h) \rhd (g \rhd k)$
\end{itemize}
for $g, h, k$ in some group G.  That is, conjugating an element by
itself acts as the identity, left and right conjugation are
inverses, and conjugation satisfies a self-distributive law. These
observations led to the creation of a new algebraic structure: a
`quandle'.  A quandle is a set $Q$ together with two binary
operations $\rhd \maps Q \times Q \to Q$ and $\lhd \maps Q \times
Q \to Q$ satisfying the three identities above. As anticipated,
the primordial example of a quandle is a group with the two
operations of left and right conjugation.

Since every law involving conjugation can be derived from the
three axioms above, the theory of quandles can be regarded as the
theory of conjugation.  Thus, since conjugation in a Lie group
gives rise to the bracket operation in its Lie algebra, we desire
a way to treat our Lie groups as though they were quandles.  In
fact, we can get away with less than a quandle; a structure we
call a `spindle' will do.

We begin this chapter in Section \ref{defns} by recalling the
definitions of quandles and racks and introducing two related
notions: `spindles' and `shelves'. We demonstrate the relationship
to groups and define the categories of Shelf, Rack, Spind, and
Quand.

In Section \ref{topology2}, we illustrate the connection between
these four algebraic concepts and topology. In particular, we
demonstrate a relationship to the Reidemeister moves, and show
that both quandles and Lie algebras give solutions to the
Yang--Baxter equation.  This suggests that a Lie algebra is
perhaps a quandle!  We investigate this suggestion later on.
Furthermore, we show that in some sense, the inverse properties
are more essential than the idempotence and self-distributivity
axioms.  This is because the latter two have completely
symmetrical counterparts that are implied by the three axioms
given above, whereas \emph{neither} inverse property implies the
other.

We continue in Section \ref{braids} by describing the connection
between these four algebraic concepts and various braid and framed
braid monoids and groups.  It turns out that quandles, racks,
spindles, and shelves, respectively, give an action of the braid
group, framed braid group, positive braid monoid, and positive
framed braid monoid! This is a consequence of the fact that each
of the three quandle axioms is equivalent, in a suitable context,
to one of the three Reidemeister moves.

Clearly, a Lie group may be thought of as quandle in the same way
that an ordinary group can, via conjugation.  However, this
description lacks knowledge about the manifold structure of our
Lie group. This deficiency suggests that we require a means by
which we can describe Lie groups in terms of special sorts of
quandles --- ones that are also manifolds!  Thus, in Section
\ref{internalization} we internalize these four concepts, which
will allow us to treat Lie groups as `quandles in $\Diff_{\ast}$',
the category of pointed smooth manifolds.  Just as a group gives a
quandle, we will use internalization to show that a group in any
category $K$ with products is a special sort of `quandle in $K$'.
Using the language of internalization, Lie groups, which are
groups in $\Diff_{\ast}$, will become quandles in $\Diff_{\ast}$.

Now that we have a way to think of our Lie groups as quandles,
which allows us to focus on conjugation, we must bring
differentiation into the game in order to obtain Lie algebras.
Since the bracket arises from differentiating conjugation twice,
it will not suffice to simply consider the tangent space of our
Lie group at the identity.  So instead we consider the symmetric
algebra on $T_{\ast}G$. This has two desirable qualities, namely
that the process of taking the symmetric algebra of a space
preserves products, and the symmetric algebra contains information
about the space $S^{2}(T_{\ast}G)$, which can be thought of as the
space of second-order Taylor expansions of polynomials in the
elements of $T_{\ast}G$. However, like the tangent space, the
symmetric algebra has its flaws as well.  The problem is that a
map $f \maps M \to N$ gives a map from $S(T_pM)$ to $S(T_{f(p)}
N)$ that only depends on $df$ and \emph{not} on the higher
derivatives of $f$, which is what we need in order to obtain the
bracket.

Therefore, in Section \ref{cojets}, we recall the concept of
`$k$-jets'.  Roughly speaking, a $k$-jet of a smooth function $f$
is a gadget which keeps track of the value of the function at a
point together with all derivatives of $f$ up to the $k$th order
at that point. That is, this is simply a way of describing the
Taylor expansion of $f$ at a given point up to order $k$.  The
algebra, $J^{k}(M,p)$, of $k$-jets of real-valued functions on $M$
at the point $p$, has the desirable properties of keeping track of
the necessary derivatives, and being isomorphic to an algebra that
we are familiar with: $\displaystyle{\bigoplus _{i=0} ^{k}
S^{i}(T^{\ast} _{p}M)}$. Nonetheless, this space of $k$-jets is
still not quite what we need, since it defines a {\it
contravariant} functor from $\Diff_{\ast}$, the category of smooth
manifolds, to $\Vect$, the category of vector spaces, whereas we
want a {\it covariant} functor.  Therefore, we recall the notion
of the coalgebra of `$k$-cojets', $J_{k}(M,p)$, which is defined
to be the dual of $J^{k}(M,p)$. While $J^{k} \maps \Diff_{\ast}
\to \Vect$ does define a covariant functor, it does not preserve
products, so $J^{k}$ will not send groups and quandles in
$\Diff_{\ast}$ to groups and quandles in $\Vect$, which we
require.

Thus, we continue in this section by recalling the concept of the
space of `cojets', $J_{\infty}(M,p)$, which is simply the union
over $k$ of the spaces of $k$-cojets.  $J_{\infty}(M,p)$ forms a
cocommutative coalgebra that is isomorphic as a coalgebra to
$S(T_{p}M)$, which allows us in Section \ref{specialcoalgs} to
define a covariant functor $J_{\infty} \maps \Diff_{\ast} \to
\mathcal{C}$ where $\mathcal{C}$ is the category of cocommutative
coalgebras that are isomorphic as coalgebras to $S(V)$ for some
vector space $V$. That is, the objects of this category are
triples $(C, V, \alpha)$, where $C$ is a coalgebra, $V$ is a
vector space, and $\alpha$ is a coalgebra isomorphism from $C$ to
$S(V)$.  We call the objects of $\mathcal{C}$ `special
coalgebras'. We show that $\mathcal{C}$ has products, which
enables us to define a group in $\mathcal{C}$, and then
demonstrate that $J_{\infty}$ preserves products which implies
that it sends groups in $\Diff_{\ast}$ to groups in $\mathcal{C}$.

Finally, in Section \ref{punchline} we use our special coalgebra
$(C, V, \alpha)$ to obtain the Lie algebra of the Lie group we
started with.  This requires defining a functor \newline $F \maps
\mathcal{C} \to \Vect$ that `picks off' the vector space. Then,
given a Lie group $G$, we use the fact that $(J_{\infty}(G,1),
T_{1}G, \alpha)$ is a group in $\mathcal{C}$, and hence quandle in
$\mathcal{C}$ to define the bracket on $T_{1}G$ in terms of the
quandle operation $\rhd$ on $J_{\infty}(G,1)$.  The Jacobi
identity for the bracket follows from the self-distributive law
for the quandle operation, while the antisymmetry of the bracket
comes from the idempotence law satisfied by the quandle operation.
We again remark that we do not require a full-fledged quandle for
this process, but can get by with a `spindle' since we use only
the self-distributivity and idempotence conditions.

We summarize our process in the following diagram:

$$
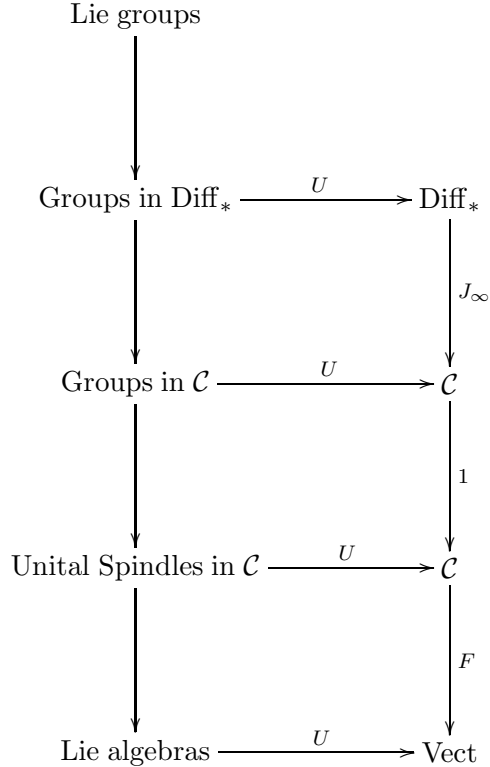
\xymatrix{
   {\rm Lie \; groups}
   \ar[dd] \\ \\
   {\rm Groups \; in \; \Diff_\ast}
   \ar[rr]^<<<<<<<<<<<{U}
   \ar[dd]
   && \Diff_{\ast}
   \ar[dd]^{J_{\infty}} \\ \\
   {\rm Groups \; in \; \mathcal{C}}
   \ar[rr]^<<<<<<<<<<<<<<<{U}
   \ar[dd]
   && \mathcal{C}
   \ar[dd]^{1} \\ \\
   {\rm Unital \; Spindles \; in \; \mathcal{C}}
   \ar[rr]^<<<<<<<<<<<{U}
   \ar[dd]
   && \mathcal{C}
   \ar[dd]^{F} \\ \\
   {\rm Lie \; algebras}
   \ar[rr]^<<<<<<<<<<<<<<{U}
   && \Vect
}$$ In the final section of this chapter, we conclude by
describing the categorified version of this process and present an
outline of our future work.

We remind the reader that in all that follows, we denote the
composition of two morphisms $f \maps x \to y$ and $g \maps y \to
z$ as $f  g \maps x \to z.$ In addition, $\pi_{1}$ and $\pi_{2}$
are the canonical projections of $Q \times Q$ onto the first and
second components, respectively.

\section{Shelves, Racks, Spindles and Quandles} \label{shelves}

\subsection{Definitions and Relation to Groups} \label{defns}
Quandles and racks are algebraic structures that capture the
essential properties of the operations of conjugation in a group.
While both provide interesting invariants of braids, quandles also
give a conceptual explanation of the passage from a Lie group to
its Lie algebra.  This is a result of the fact that the bracket in
the Lie algebra arises from differentiating conjugation in the Lie
group. Our goal is to explain this passage in a way which might
make it easier to categorify.

A comprehensive history of racks, beginning with the unpublished
correspondence between Conway and Wraith roughly $45$ years ago,
is provided in the work of Fenn and Rourke \cite{FR}. Quandles,
which are special cases of racks, were first introduced as a
source of knot invariants by David Joyce in $1982$ and have been
studied and explored since then in various papers \cite{Bri, FR,
J, K}. Fenn and Rourke demonstrate that racks provide an elegant
algebraic framework which gives complete algebraic invariants for
both framed links in $3$-manifolds and for the $3$-manifolds
themselves.  Likewise, Joyce uses the structure provided by
quandles to formulate additional invariants.  In fact, Fenn and
Rourke suggest that the theory of racks presents the opportunity
to find a complete sequence of computable invariants for framed
links and $3$-manifolds.

Both quandles and racks illustrate the connection between algebra
and topology since their axioms algebraically encode the three
Reidemeister moves.  In particular, we will show in the next
section that the `self-distributivity axiom' is equivalent to the
Yang--Baxter equation, or third Reidemeister move, in a suitable
context.

We begin with the simplest of these four structures, a `shelf':

\begin{defn} \et \label{lshelf}
A {\bf left shelf} $(S, \rhd)$ is a nonempty set $S$ equipped with
a binary operation $\rhd \maps S \times S \to S$ called {\bf left
conjugation}, which satisfies
\begin{description}
    \item[(i)] {\bf (left distributive law)} $x \rhd (y \rhd z) = (x \rhd y) \rhd (x \rhd z)$ for all
               $x, y, z \in S$.
\end{description}
\end{defn}

A typical example of a left shelf is a group $G$ where the shelf
operation is defined to be left conjugation, $x \rhd y :=
xyx^{-1}$ for all $x,y \in G.$  A trivial computation shows that
the left distributive law holds for this operation and that it
amounts to the distributive law satisfied by conjugation.  As we
have alluded to in the introduction to this chapter, it is
precisely this self-distributive law that corresponds to the
Jacobi identity in a Lie algebra once we differentiate it.

Since a group can also act on itself by right conjugation,
$x^{-1}yx$, which satisfies an analogous distributive law, we
obtain the symmetrical notion of a `right shelf':

\begin{defn} \et \label{rshelf}
A {\bf right shelf} $(S, \lhd)$ is a nonempty set $S$ equipped
with a binary operation $\lhd \maps S \times S \to S$ called {\bf
right conjugation}, which satisfies
\begin{description}
    \item[(ii)] {\bf (right distributive law)} $(x \lhd y) \lhd z = (x \lhd z) \lhd (y \lhd z)$
    for all $x, y, z \in S$.
\end{description}
\end{defn}

As mentioned above, an example of a right shelf is a group $G$
where the shelf operation is defined to be right conjugation, $y
\lhd x := x^{-1}yx$ for all $x,y \in G$.

Henceforth, the term {\bf shelf} will always refer to a left
shelf.  We remark that, in general, the shelf operation is neither
associative nor commutative.  Furthermore, we should think of
these shelf operations as \emph{actions}, where the element
pointed to is the one being acted upon.  The location of the
remaining element determines whether we have left or right
conjugation.  For instance, `$x \lhd y$' says that $x$ is being
conjugated on the right by $y$.

In the next section, we prove that we can define a Yang--Baxter
operator on a shelf, thus demonstrating that the left distributive
law is equivalent to the Yang--Baxter equation, or third
Reidemeister move.  We now turn to the concept of a `rack', which
is a set equipped with two operations satisfying properties
corresponding to the second and third Reidemeister moves.  Roughly
speaking, a rack results when a left shelf and a right shelf fit
together nicely.

\begin{defn} \et \label{rack}
A {\bf rack} $(Q, \rhd, \lhd)$ is a nonempty set $Q$ together with
two binary operations $\rhd \maps Q \times Q \to Q$ and $\lhd
\maps Q \times Q \to Q$, called {\bf left conjugation} and {\bf
right conjugation}, such that $(Q, \rhd)$ is a left shelf, $(Q,
\lhd)$ is a right shelf, and the following hold:
\begin{description}
    \item[(iii)] {\bf (right inverse property)}  $(x \rhd y) \lhd x = y$ for all $x, y \in Q$, and
    \item[(iv)] {\bf (left inverse property)} $x \rhd (y \lhd x) = y$ for all $x, y \in Q$.
\end{description}
\end{defn}
That is, $(iii)$ and $(iv)$ simply say that $x \rhd -$ and $- \lhd
x$ are inverse operations.  Once again, groups serve as the
primary examples of racks, where the rack operations are left and
right conjugation.  Again, a routine computation gives the two
inverse properties, which result from the fact that left and right
conjugation are inverses.

Just as a rack is comprised of two shelves which fit together
nicely, a quandle is made up of two `spindles':

\begin{defn} \et \label{lspindle} A {\bf left spindle} $(S, \rhd)$ is a
left shelf such that
\begin{description}
    \item[(v)] {\bf (left idempotence)} $x \rhd x = x$ for all $x \in S$.
\end{description}
\end{defn}

\begin{defn} \et \label{rspindle} A {\bf right spindle} $(S,
\lhd)$ is a right shelf such that
\begin{description}
    \item[(vi)]{\bf (right idempotence)} $x \lhd x = x$ for all $x \in S$.
\end{description}
\end{defn}
As in the case of shelves, we will use the term {\bf spindle} to
refer to left spindles.

We now arrive at the notion of a `quandle', which is a rack
satisfying an additional axiom that is related to the first
Reidemeister move:

\begin{defn} \et \label{quandle}
A {\bf quandle} $(Q, \rhd, \lhd)$ is a rack such that $(Q, \rhd)$
is a left spindle and $(Q, \lhd)$ is a right spindle.
\end{defn}

The motivating example of a quandle is a group with the operations
of left and right conjugation. In their paper \cite{FR}, Fenn and
Rourke mention that Conway and Wraith regarded a quandle as the
wreckage left behind when we start with a group and discard the
group operation so that only the notion of conjugacy remains.
Indeed, this is why Fenn and Rourke use the term `rack' to
describe a structure still weaker than that of a quandle, as in
the phrase ``rack and ruin".

Another example of a quandle is the {\bf reflection quandle} whose
underlying set is $\R^2$ and has quandle operations defined as:
$$a \rhd b = b \lhd a = \textrm{``reflect b about a"} = 2a - b.$$
Notice that left and right conjugation are the same in this case
because reflection in the plane is its own inverse.  Numerous
examples of quandles and racks can be found in the work of Fenn
and Rourke \cite{FR}.

In order to define the categories of these concepts, we must first
define shelf, rack, spindle and quandle homomorphisms. Such
morphisms must preserve the left and/or right conjugation.

\begin{defn} \et \label{quandlemorphism} Given quandles $Q, Q'$, a
function $f \maps Q \to Q'$ is a {\bf quandle morphism} if $f(x
\rhd y) = f(x) \rhd f(y)$ and $f(x \lhd y) = f(x) \lhd f(y)$.
\end{defn}
Shelf, rack and spindle morphisms are defined similarly.  With
these morphisms, we obtain the categories:  {\bf Shelf, Rack,
Spind,} and {\bf Quand}.

Due to the relationship between groups and quandles, we can
describe a functor $Conj \maps \Grp \to \Quand$ from the category
of groups to the category of quandles, which sends any group $G$
to $(G, \rhd, \lhd)$ where the operations are simply left and
right conjugation.  In addition, we can consider $n$-fold
conjugation: $x \rhd y = x^{n}yx^{-n}$ and $x \lhd y =
y^{-n}xy^{n}$ and obtain another functor $Conj_{n} \maps \Grp \to
\Quand$ that is defined in the obvious way.

It is important to notice that not every quandle is of the form
$Conj(G)$ for some group $G$.  In addition, since every quandle is
a rack, we have natural inclusion of the category of quandles in
the category of racks.  However, the converse is not necessarily
true. For example, the {\bf cyclic rack of order} $\mathbf{n}$ is
a rack but not a quandle. This rack consists of a set $C_{n} =
\{0, 1, \ldots, n-1\}$ with the rack operations defined by $x \rhd
y = (y + 1)$mod $n$ and $y \lhd x = (y-1)$mod $n$.

Readers familiar with the notions of `rack' and `quandle' may
recall these definitions and notation in another guise. Fenn and
Rourke used the following, equivalent, definition:

\begin{defn} \et \label{altdefnquandle}  A {\bf quandle} $(Q,
\lhd)$ is a set Q equipped with a binary operation \newline $\lhd
\maps Q \times Q \to Q$ defined as $(a,b) \mapsto a \lhd b$ which
satisfies:
\begin{itemize}
\item[(a)] $x \lhd x = x$ for all $x \in Q$
\item[(b)] For all $x,y \in Q$, there exists a unique $z \in Q$
such that $z \lhd x = y$.
\item[(c)] $(x \lhd y) \lhd z = (x \lhd z) \lhd (y \lhd z)$ for all $x,y,z \in Q$.
\end{itemize}
\end{defn}
The second condition in the above definition is equivalent to the
existence of a binary operation $\rhd \maps Q \times Q \to Q$ such
that the two operations $\rhd$ and $\lhd$ satisfy the two inverse
properties of Definition \ref{rack}.  While they used the above
definition, Fenn and Rourke chose to use exponential notation for
their operation, since the operation is not symmetric.  That is,
whereas we used $x \lhd y$ to represent the result of $y$ acting
on $x$ from the right, Fenn and Rourke notated this same quantity
as $x^{y}$. While the exponential notation certainly has its
advantages, we have chosen to use the notation $x \rhd y$ and $x
\lhd y$ because this language makes it easier to describe
generators and relations and will prove easier to categorify in
the future.  That is, when we define $2$-quandles, which will be
categories equipped with two conjugation \emph{functors}, we can
continue to use the notation $\rhd$ and $\lhd$ for these functors.

We continue by exhibiting the correlation between shelves, racks,
spindles, and quandles and the theory of braids.  In particular,
we illustrate the relationship between the Reidemeister moves and
the axioms of distributivity, inverses, and idempotence.

\subsection{Relation to Topology} \label{topology2}
The Yang--Baxter equation arises in many contexts in mathematics
and physics.  All these concepts are related by the fact that this
equation is an algebraic distillation of the `third Reidemeister
move' in knot theory:
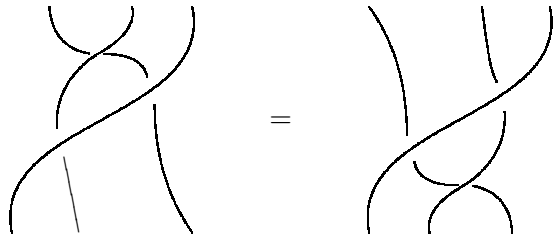
\[
\def\objectstyle{\scriptstyle}
\def\labelstyle{\scriptstyle}
  \xy
   (12,15)*{}="C";
   (4,15)*{}="B";
   (-7,15)*{}="A";
   (12,-15)*{}="3";
   (-3,-15)*{}="2";
   (-12,-15)*{}="1";
       "C";"1" **\crv{(15,0)& (-15,0)};
       (-5,-5)*{}="2'";
       (7,2)*{}="3'";
     "2'";"2" **\crv{};
     "3'";"3" **\crv{(7,-8)};
       \vtwist~{"A"}{"B"}{(-6,-1)}{(6,5.7)};
\endxy
 \qquad = \qquad
   \xy
   (-12,-15)*{}="C";
   (-4,-15)*{}="B";
   (7,-15)*{}="A";
   (-12,15)*{}="3";
   (3,15)*{}="2";
   (12,15)*{}="1";
       "C";"1" **\crv{(-15,0)& (15,0)};
       (5,5)*{}="2'";
       (-7,-2)*{}="3'";
     "2'";"2" **\crv{(4,6)};
     "3'";"3" **\crv{(-7,8)};
       \vtwist~{"A"}{"B"}{(6,1)}{(-6,-5.7)};
\endxy
\]
The Yang--Baxter equation makes sense in any monoidal category.
Originally mathematical physicists concentrated on solutions in
the category of {\it vector spaces} with its usual tensor product,
obtaining solutions from quantum groups. More recently there has
been interest in solutions of the Yang--Baxter equation in the
category of {\it sets}, with its Cartesian product \cite{G-I, GV,
ESS, EGS, H, LYZ, O, R, S, WX}. In this section we prove that any
shelf (and thus any rack, spindle, or quandle) gives a solution to
the Yang--Baxter equation in the category of sets. In particular,
this implies that any group gives a solution to the Yang--Baxter
equation in which the crossing:
$$
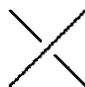
\begin{xy} \xoverv~{(-5,5)}{(5,5)}{(-5,-5)}{(5,-5)};
    (-5,6.4)*{};
    (7,6)*{};
    (-15,0)*{};
    (-6,-7)*{};
    (6,-7)*{};
\end{xy}$$
is related to the {\it action of conjugating one group element by
another}.  Similarly, we prove that any Lie algebra gives a
solution to the Yang--Baxter equation in the category of vector
spaces, in which the crossing is related to the {\it bracket
operation}.

These two results illustrate one of the connections between Lie
algebras and shelves, by demonstrating that the Jacobi identity
and self-distributivity axiom are, in a suitable context,
equivalent. Not surprisingly, our results for shelves and Lie
algebras are closely related!  We will make this relation much
more explicit in Section \ref{punchline}, where we show how
spindles in a certain category ${\mathcal C}$ come from Lie
groups, and give Lie algebras.  Furthermore, we show that the
first and second Reidemeister moves in knot theory are related to
the idempotence and inverse properties.

We begin by recalling the Yang--Baxter equation:

\begin{defn} \et Given a set $Q$ and an
isomorphism $B \maps Q \times Q \to Q \times Q$, we say $B$ is a
{\bf Yang--Baxter operator} if it satisfies the {\bf Yang--Baxter
equation}, which says that:
$$(B \times 1)(1 \times B)(B \times 1)
= (1 \times B)(B \times 1)(1 \times B),$$ or in other words, that
this diagram commutes:
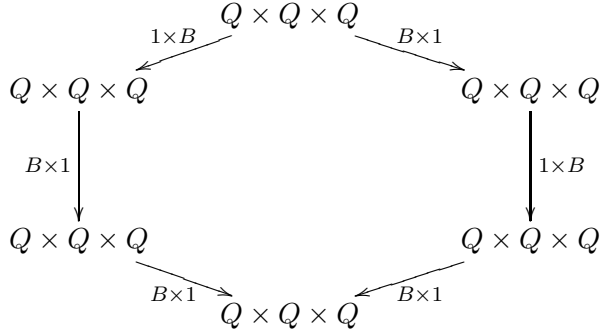
\[
\begin{xy}
    (-2, 20)*+{Q \times Q \times Q}="1";
    (-30,10)*+{Q \times Q \times Q}="2";
    (-30,-10)*+{Q \times Q \times Q}="3";
    (-2,-20)*+{Q \times Q \times Q}="4";
    (30,-10)*+{Q \times Q \times Q}="5";
    (30, 10)*+{Q \times Q \times Q}="6";
        {\ar^{B \times 1} "1";"6"};
        {\ar_{1 \times B} "1";"2"};
        {\ar_{B \times 1} "2";"3"};
        {\ar_{B \times 1} "3";"4"};
        {\ar^{B \times 1} "5";"4"};
        {\ar^{1 \times B} "6";"5"};
\end{xy}
\]
\end{defn}

If we draw $B \maps Q \times Q \to Q \times Q $ as a braiding:
$$
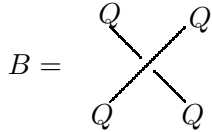
\begin{xy} \xoverv~{(-5,5)}{(5,5)}{(-5,-5)}{(5,-5)};
    (-5,6.4)*{Q};
    (7,6)*{Q};
    (-15,0)*{B=};
    (-6,-7)*{Q};
    (6,-7)*{Q};
\end{xy}$$
the Yang--Baxter equation says that:
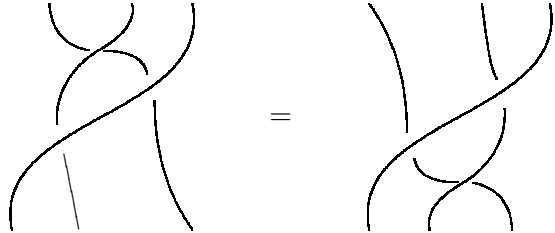
\[
\def\objectstyle{\scriptstyle}
\def\labelstyle{\scriptstyle}
  \xy
   (12,15)*{}="C";
   (4,15)*{}="B";
   (-7,15)*{}="A";
   (12,-15)*{}="3";
   (-3,-15)*{}="2";
   (-12,-15)*{}="1";
       "C";"1" **\crv{(15,0)& (-15,0)};
       (-5,-5)*{}="2'";
       (7,2)*{}="3'";
     "2'";"2" **\crv{};
     "3'";"3" **\crv{(7,-8)};
       \vtwist~{"A"}{"B"}{(-6,-1)}{(6,5.7)};
\endxy
 \qquad = \qquad
   \xy
   (-12,-15)*{}="C";
   (-4,-15)*{}="B";
   (7,-15)*{}="A";
   (-12,15)*{}="3";
   (3,15)*{}="2";
   (12,15)*{}="1";
       "C";"1" **\crv{(-15,0)& (15,0)};
       (5,5)*{}="2'";
       (-7,-2)*{}="3'";
     "2'";"2" **\crv{(4,6)};
     "3'";"3" **\crv{(-7,8)};
       \vtwist~{"A"}{"B"}{(6,1)}{(-6,-5.7)};
\endxy
\]

\noindent This is called the `third Reidemeister move' in knot
theory \cite{BZ}, and it gives the most important relations in
Artin's presentation of the braid group \cite{Birman}. As a
result, any solution of the Yang--Baxter equation gives an
invariant of braids.

In general, a `braiding' operation provides a diagrammatic
description of the process of switching the order of two things.
This idea is formalized in the concept of a braided monoidal
category, where the braiding is an isomorphism
\[    B_{x,y} \maps x \tensor y \to y \tensor x .\]
To show that the distributive law in the definition of a shelf is
equivalent to the Yang--Baxter equation, we define a braiding
operation in terms of the conjugation operation.  Given a set $Q$
and binary operation $\rhd \maps Q \times Q \to Q,$ we can define
a braiding map $B \maps Q \times Q \to Q \times Q$ as $B(x,y) =
(y, y \rhd x)$. That is, since we have switched the order of $x$
and $y$, we introduce a `correction' term involving the
conjugation operation.  We draw $B$ as a positive crossing:
$$
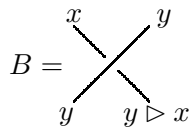
\begin{xy} \xoverv~{(-5,5)}{(5,5)}{(-5,-5)}{(5,-5)};
    (-5,6.4)*{x};
    (7,6)*{y};
    (-10,0)*{B=};
    (-6,-7)*{y};
    (6,-7)*{y \rhd x};
\end{xy}$$
Since the third Reidemeister move involves three strands and the
shelf distributive law involves three shelf elements, it should
not come as a surprise that the Yang--Baxter equation is actually
{\it equivalent} to the shelf distributive law is actually
equivalent in a suitable context:

\begin{lem} \et \label{lemma1} Let $Q$ be a set equipped with
a binary operation $\rhd \maps Q \times Q \to Q$.  The braiding
operation, $B$, defined above satisfies the Yang--Baxter equation:
$$(B \times 1)(1 \times B)(B \times 1) = (1 \times B)(B \times
1)(1 \times B)$$ if and only if $(Q,\rhd)$ is a left shelf.
\end{lem}

\noindent {\bf Proof.}  The braiding $B$ satisfies the Yang-Baxter
equation if and only if the following braid equation holds:
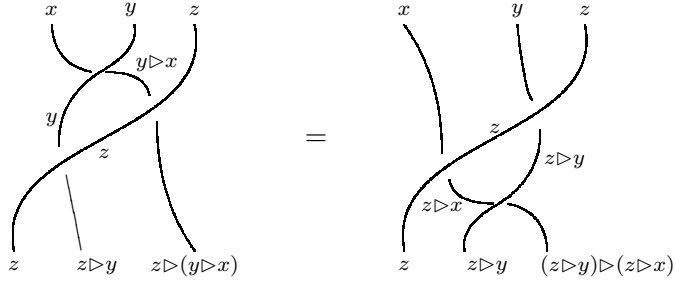
\[
\def\objectstyle{\scriptstyle}
\def\labelstyle{\scriptstyle}
  \xy
   (12,15)*{}="C";
   (4,15)*{}="B";
   (-7,15)*{}="A";
   (12,-15)*{}="3";
   (-3,-15)*{}="2";
   (-12,-15)*{}="1";
   (12,17)*{z};
   (3.5, 17)*{y};
   (-7, 17)*{x};
   (12,-17)*{z \rhd (y \rhd x)};
   (-1, -17)*{z \rhd y};
   (-12, -17)*{z};
   (0,-2)*{z};
   (-7,2.5)*{y};
   (7,10)*{y \rhd x};
       "C";"1" **\crv{(15,0)& (-15,0)};
       (-5,-5)*{}="2'";
       (7,2)*{}="3'";
     "2'";"2" **\crv{};
     "3'";"3" **\crv{(7,-8)};
       \vtwist~{"A"}{"B"}{(-6,-1)}{(6,5.7)};
\endxy
 \qquad = \qquad
   \xy
   (-12,-15)*{}="C";
   (-4,-15)*{}="B";
   (7,-15)*{}="A";
   (-12,15)*{}="3";
   (3,15)*{}="2";
   (12,15)*{}="1";
   (12,17)*{z};
   (3, 17)*{y};
   (-12, 17)*{x};
   (15,-17)*{(z \rhd y) \rhd (z \rhd x)};
   (-1, -17)*{z \rhd y};
   (-12, -17)*{z};
   (0,1.3)*{z};
   (9.3,-3)*{z \rhd y};
   (-7,-9)*{z \rhd x};
       "C";"1" **\crv{(-15,0)& (15,0)};
       (5,5)*{}="2'";
       (-7,-2)*{}="3'";
     "2'";"2" **\crv{(4,6)};
     "3'";"3" **\crv{(-7,8)};
       \vtwist~{"A"}{"B"}{(6,1)}{(-6,-5.7)};
\endxy
\]
That is, $B$ satisfies the Yang--Baxter equation if and only if $z
\rhd (y \rhd x) = (z \rhd y) \rhd (z \rhd x),$ which is the left
distributive law, $(i)$. Hence, the braiding gives a solution of
the Yang--Baxter equation if and only if $(Q, \rhd)$ is a left
shelf. \qed

Thus, left shelves, and hence left conjugation, correspond to
`positive' or right-handed crossings.  While this may seem
perverse, it has the effect that our `favorite' conjugation $x
\rhd y = xyx^{-1}$ corresponds to the topologist's `favorite'
crossing, namely the right-handed one.  Recall that our `favorite'
conjugation is the one that differentiates to give $[x,y]$ in a
Lie algebra. Analogously, right shelves and right conjugation
correspond to `negative' or left-handed crossings.  We now show
that we can define a Yang--Baxter operator on right shelves as
well.

Given a set $Q$ and binary operation $\lhd \maps Q \times Q \to
Q$, we may define a new braiding map $\bar{B}$ as $\bar{B}(x,y) =
(y \lhd x, x)$.  We draw $\bar{B}$ as a negative crossing:
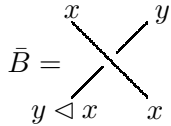
\[
\begin{xy}
\xunderv~{(-5,5)}{(5,5)}{(-5,-5)}{(5,-5)};
    (-5,6.4)*{x};
    (7,6)*{y};
    (-10,0)*{\bar{B}=};
    (-6,-7)*{y \lhd x};
    (6,-7)*{x};
\end{xy}
\]
Notice that in the illustrations of both $B$ and $\bar{B}$, the
strand above always acts on the strand below.

Analogous to the comment made above, since the right distributive
law involves three right shelf elements, it, too, is equivalent to
the Yang--Baxter equation in an appropriate context:

\begin{lem} \et \label{lemma4} Let $Q$ be a set equipped with a binary operation
$\lhd \maps Q \times Q \to Q $.  The braiding map $\bar{B}$
satisfies the Yang--Baxter equation if and only if $(Q, \lhd)$ is
a right shelf.
\end{lem}

\noindent {\bf Proof.} The braiding $\bar{B}$ satisfies the
Yang-Baxter equation if and only if the following braid equation
holds:
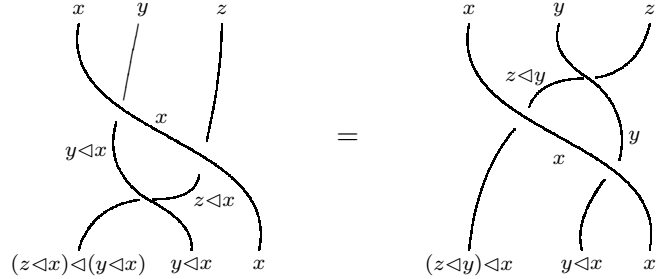
\[
\def\objectstyle{\scriptstyle}
\def\labelstyle{\scriptstyle}
    \xy
    (7,15)*{}="C";
    (-4,15)*{}="B";
    (-12,15)*{}="A";
    (12,-15)*{}="3";
    (3,-15)*{}="2";
    (-12,-15)*{}="1";
        "A";"3" **\crv{(-15,0)& (15,0)};
        \vcross~{(-7,2)}{(4,-5)}{"1"}{"2"};
    (-6,5)*{}="B'";
    (5,-.5)*{}="C'";
      "B'";"B" **\crv{};
      "C'";"C" **\crv{(7,8)};
         (-12,17)*{x};
         (-3.5,17)*{y};
         (7,17)*{z};
         (-12,-17)*{(z \lhd x) \lhd (y \lhd x)};
         (3,-17)*{y \lhd x};
         (12,-17)*{x};
         (-11,-2)*{y \lhd x};
         (6,-8.5)*{z \lhd x};
         (-1,2.5)*{x};
 \endxy
  \qquad = \qquad
  \xy
    (12,15)*{}="A";
    (0,15)*{}="B";
    (-12,15)*{}="C";
    (12,-15)*{}="1";
    (3,-15)*{}="2";
    (-12,-15)*{}="3";
    (6,-5.7)*{}="2'";
    (-6,1)*{}="3'";
        \vcross~{"B"}{"A"}{(-4,4)}{(8,-3)};
        "C";"1" **\crv{(-15,0)& (15,0)};
        "2'";"2" **\crv{(2,-11)};
        "3'";"3" **\crv{(-11,-5)};
           (-12,17)*{x};
           (0,17)*{y};
           (12,17)*{z};
           (-12,-17)*{(z \lhd y) \lhd x};
           (3,-17)*{y \lhd x};
           (12,-17)*{x};
             (0,-3)*{x};
             (-4.5,8)*{z \lhd y};
             (10,0)*{y};
 \endxy
\]
Thus, $\bar{B}$ satisfies the Yang--Baxter equation if and only if
$((z \lhd y) \lhd x, y \lhd x, x) = ((z \lhd x) \lhd (y \lhd x), y
\lhd x, x),$ which is the right distributive law, $(ii)$. Hence,
the braiding $\bar{B}$ gives a solution of the Yang--Baxter
equation if and only if $(Q, \rhd)$ is a right shelf. \qed

Therefore, shelves give set-theoretic solutions of the
Yang--Baxter equation.  Of course, since racks, spindles, and
quandles are shelves equipped with some combination of extra
structure and extra properties, they will also give solutions to
the Yang--Baxter equation.   We proceed by showing that Lie
algebras also give solutions of the Yang--Baxter equation, a
result due to James Dolan \cite{D}.  Now, we consider vector
spaces:

\begin{defn} \et Given a vector space $V$ and an
isomorphism $B \maps V \otimes V \to V \otimes V$, we say $B$ is a
{\bf Yang--Baxter operator} if it satisfies the {\bf Yang--Baxter
equation}, which says that:
$$(B \otimes 1)(1 \otimes B)(B \otimes 1)
= (1 \otimes B)(B \otimes 1)(1 \otimes B),$$ or in other words,
that this diagram commutes:
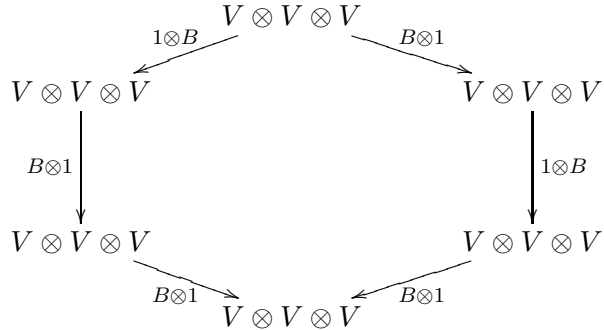
\[
\begin{xy}
    (-2, 20)*+{V \otimes V \otimes V}="1";
    (-30,10)*+{V \otimes V \otimes V}="2";
    (-30,-10)*+{V \otimes V \otimes V}="3";
    (-2,-20)*+{V \otimes V \otimes V}="4";
    (30,-10)*+{V \otimes V \otimes V}="5";
    (30, 10)*+{V \otimes V \otimes V}="6";
        {\ar^{B \otimes 1} "1";"6"};
        {\ar_{1 \otimes B} "1";"2"};
        {\ar_{B \otimes 1} "2";"3"};
        {\ar_{B \otimes 1} "3";"4"};
        {\ar^{B \otimes 1} "5";"4"};
        {\ar^{1 \otimes B} "6";"5"};
\end{xy}
\]
\end{defn}
As in the set-theoretic definition, if we draw $B$ as a braiding,
the Yang--Baxter equation is equivalent to the third Reidemeister
move.

We have seen that the braiding defined on a shelf $Q$ formalized
the failure of $x$ and $y$ to commute by introducing a correction
term involving the conjugation operation.  Since the bracket
$[x,y]$ in a Lie algebra measures the difference between $xy$ and
$yx$, it should not be too surprising that we can get a
Yang--Baxter operator from any Lie algebra. Here, when we switch
the order of two Lie algebra elements $x$ and $y$, we add on a
correction term involving the bracket: $[x,y]$. And, just as we
have already observed for shelves, since the third Reidemeister
move involves three strands, while the Jacobi identity involves
three Lie algebra elements, it should also not be surprising that
the Yang--Baxter equation is actually {\it equivalent} to the
Jacobi identity in a suitable context:

\begin{prop} \et \label{YBEJacobi} Let $L$ be a vector space equipped
with a skew-symmetric bilinear operation $[\cdot,\cdot] \maps L
\times L \to L$.  Let $L' = k \oplus L$ and define the isomorphism
$B \maps L' \tensor L' \to L' \tensor L'$ by
$$B((a,x) \otimes (b,y))
= (b,y) \otimes (a,x) + (1,0) \otimes (0, [x,y]).$$ Then $B$ is a
solution of the Yang--Baxter equation if and only if
$[\cdot,\cdot]$ satisfies the Jacobi identity.
\end{prop}

\noindent {\bf Proof. } Applying the left-hand side of the
Yang--Baxter equation to an object $(a,x) \otimes (b,y) \otimes
(c,z)$ of $L' \otimes L' \otimes L'$ yields:
\begin{eqnarray*}
(a,x) \otimes (b,y) \otimes (c,z)  & \longmapsto & (b,y) \otimes
(a,x) \otimes (c,z) + (1,0) \otimes (0,[x,y]) \otimes (c,z) \\
& \longmapsto & (b,y) \otimes (c,z) \otimes (a,x) + (b,y) \otimes
(1,0) \otimes (0,[x,z]) \\ & & + \; (1,0) \otimes (c,z) \otimes
(0,[x,y]) + (1,0) \otimes (1,0) \otimes (0,[[x,y],z]) \\
& \longmapsto & (c,z) \otimes (b,y) \otimes (a,x) + (1,0) \otimes
(0,[y,z]) \otimes (a,x) \\ & & + \; (1,0) \otimes (b,y) \otimes
(0, [x,z]) + (c,z) \otimes (1,0) \otimes (0,[x,y]) \\ & & + \;
(1,0) \otimes (1,0) \otimes (0,[[x,y],z])
\end{eqnarray*}
while applying the right-hand side produces:
\begin{eqnarray*}
(a,x) \otimes (b,y) \otimes (c,z) & \longmapsto & (a,x) \otimes
(c,z) \otimes (b,y) + (a,x) \otimes (1,0) \otimes (0,[y,z]) \\
& \longmapsto & (c,z) \otimes (a,x) \otimes (b,y) + (1,0) \otimes
(0,[x,z]) \otimes (b,y) \\ & & + \; (1,0) \otimes (a,x) \otimes
(0, [y,z]) \\
& \longmapsto & (c,z) \otimes (b,y) \otimes (a,x) + (c,z) \otimes
(1,0) \otimes (0, [x,y]) \\ & & + \; (1,0) \otimes (b,y) \otimes
(0,[x,z]) + (1,0) \otimes (1,0) \otimes (0,[[x,z],y]) \\ & & + \;
(1,0) \otimes (0, [y,z]) \otimes (a,x) + (1,0) \otimes (1,0)
\otimes (0,[x,[y,z]])
\end{eqnarray*}
Notice that both sides consist of the same four uninteresting
terms.  The remaining terms are equal if and only if the Jacobi
identity is satisfied in $L$. \qed

Now that we have seen that the shelf distributive law and Jacobi
identity are equivalent to the third Reidemeister move, it is
natural to wonder whether Lie algebras actually {\it are} shelves.
It turns out that they are!  Actually, they are a special sort of
spindle in a certain category.  Of course, it has to be a category
with products for the notion of spindle to make sense.
Furthermore, while forming the space $k \oplus L$ allowed us to
get a solution of the Yang--Baxter equation, this space is not
quite large enough to enable us to define a spindle structure on
it.  Therefore, we need something a bit larger and fancier to get
this special sort of spindle from a Lie algebra.  We will see that
we need a space more like the symmetric algebra, $SL$, when we
return to this idea and make it more precise in Section
\ref{punchline}. For now, we continue by showing the relationship
between the second Reidemeister move and the two inverse
properties.

The definitions given in Section \ref{defns} show that the six
quandle axioms: left and right self-distributivity, left and right
inverse laws, and left and right idempotence are completely
symmetrical. In fact, it turns out that two of these six axioms
are implied by the other four!  We will show that while we may
omit one of the two self-distributive laws and one of the two
idempotence laws in the definition of a quandle, we may \emph{not}
eliminate either of the two inverse properties.  We will use the
braidings $B$ and $\bar{B}$ described above to demonstrate that we
may omit one of the two self-distributive laws from the definition
of a rack, and therefore from the definition of a quandle, without
harm. Specifically we will prove that a shelf $(Q, \rhd)$ equipped
with an additional operation $\lhd \maps Q \times Q \to Q$
satisfying the two inverse properties $(iii)$ and $(iv)$ is
actually a rack. We begin by showing that $B$ and $\bar{B}$ are
inverses when $(Q, \rhd, \lhd)$ satisfies the two inverse
properties.

\begin{lem} \et \label{lemma2} Let $Q$ be a set equipped with two
binary operations $\rhd \maps Q \times Q \to Q$ and
\newline $\lhd \maps Q \times Q \to Q$.  Then the braiding
operations $B$ and $\bar{B}$,
$$\begin{xy} \xoverv~{(-5,5)}{(5,5)}{(-5,-5)}{(5,-5)};
    (-5,6.4)*{x};
    (7,6)*{y};
    (-10,0)*{B=};
    (-6,-7)*{y};
    (6,-7)*{y \rhd x};
\end{xy} \qquad \qquad
\begin{xy}
\xunderv~{(-5,5)}{(5,5)}{(-5,-5)}{(5,-5)};
    (-5,6.4)*{x};
    (7,6)*{y};
    (-10,0)*{\bar{B}=};
    (-6,-7)*{y \lhd x};
    (6,-7)*{x};
\end{xy}
$$
are inverses if and only if $(Q, \rhd, \lhd)$ satisfies the two
inverse properties.
\end{lem}

\noindent{\bf Proof.}  The braidings $B$ and $\bar{B}$ are
inverses if and only if the following braid equation holds:
\[
\def\objectstyle{\scriptstyle}
\def\labelstyle{\scriptstyle}
 \begin{xy}
\vcross~{(-5,10)}{(5,10)}{(-5,0)}{(5,0)};
\vtwist~{(-5,0)}{(5,0)}{(-5,-10)}{(5,-10)};
        (-5,12)*{x};
        (5,12)*{y};
        (-5,-12)*{x};
        (5,-12)*{x \rhd (y \lhd x)};
        (-8.3,0)*{y \lhd x};
        (6.3,0)*{x};
 \end{xy}
\qquad = \qquad
  \xy
    (-2.5,10)*{}="A";
    (2.5,10)*{}="B";
    (-2.5,-10)*{}="A'";
    (2.5,-10)*{}="B'";
      {\ar@{-} "A";"A'"};
      {\ar@{-} "B";"B'"};
        (-5,12)*{x};
        (5,12)*{y};
        (-5,-12)*{x};
        (5,-12)*{y};
\endxy
\qquad = \qquad
\begin{xy}
\vtwist~{(-5,10)}{(5,10)}{(-5,0)}{(5,0)};
\vcross~{(-5,0)}{(5,0)}{(-5,-10)}{(5,-10)};
        (-5,12)*{x};
        (5,12)*{y};
        (-5,-12)*{(y \rhd x) \lhd y};
        (5,-12)*{y};
        (-6,0)*{y};
        (8.3,0)*{y \rhd x};
 \end{xy}\]
Thus, $B$ and $\bar{B}$ are inverses if and only if $ x \rhd (y
\lhd x) = y$ and $(y \rhd x) \lhd y = x$, which are the two
inverse properties, $(iii)$ and $(iv)$.  Therefore, $B$ and
$\bar{B}$ are inverses if and only if $(Q, \rhd, \lhd)$ satisifies
these two properties. \qed When $(Q, \rhd)$ is a left shelf and
$(Q,\lhd)$ is a right shelf, we have an immediate corollary:

\begin{cor} \et \label{lemma2a}  Let $Q$ be a set equipped
with two binary operations $\rhd \maps Q \times Q \to Q$ and
\newline $\lhd \maps Q \times Q \to Q$ such that $(Q, \rhd)$ is a left
shelf and $(Q, \lhd)$ is a right shelf. Then the braiding
operations $B$ and $\bar{B}$ are inverses if and only if $(Q,
\rhd, \lhd)$ is a rack.
\end{cor}
Henceforth we will write $B^{-1}$ for $\bar{B}$. It is easy to see
that if $B$ satisfies the Yang--Baxter equation, then so does
$B^{-1}$:

\begin{lem} \et \label{lemma3} If $B \maps Q \times Q \to Q \times Q$ is invertible
and satisfies the Yang--Baxter equation, then $B^{-1}$ does as
well.
\end{lem}

\noindent {\bf Proof.}  We begin with the Yang--Baxter equation
satisfied by $B$ and apply inverses to both sides to see that
$B^{-1}$ also satisfies this equation: \begin{center} $(1 \times
B^{-1})(B^{-1} \times 1)(1 \times B^{-1}) = (B^{-1} \times 1)(1
\times B^{-1})(B^{-1} \times 1)$ \quad \qed \end{center}

As a corollary to Lemma \ref{lemma2}, note that if $(Q, \rhd)$ is
a left shelf and the braiding $B$ defined on it is invertible,
then $B^{-1}$ satisfies the Yang--Baxter equation by Lemma
\ref{lemma3}.  Therefore we can use $B^{-1}$ to {\it define} an
operation $\lhd \maps Q \times Q \to Q$ as $\lhd (x,y) = x \lhd y
= B^{-1} \pi_{1} (x,y)$, where we have expressed the composite in
the non-traditional sense.  It is not difficult to show that with
this definition,  $(Q, \lhd)$ is a right shelf, and thus $(Q,
\rhd, \lhd)$ is a rack.

We finally arrive at our desired conclusion:  A shelf together
with an additional operation satisfying the two inverse
properties, $x \rhd (y \lhd x) = y = (x \rhd y) \lhd x$, is a
rack. That is, the right distributive law is implied by the left
distributive law together with the two inverse properties.

\begin{prop} \et \label{rackwo(ii)} A left shelf $(Q, \rhd)$ together with an operation
$\lhd \maps Q \times Q \to Q$ that satisfies the inverse
properties, $(iii)$ and $(iv)$, also satisfies the right
distributive law, $(ii)$. That is, $(Q, \rhd, \lhd)$ is a rack.
\end{prop}

\noindent {\bf Proof.}  Define $B \maps Q \times Q \to Q \times Q$
in terms of $\rhd$ as before.  The braiding $B$ satisfies the
Yang--Baxter equation by Lemma \ref{lemma1} since $(Q, \rhd)$ is a
left shelf. By hypothesis, $(Q, \rhd, \lhd)$ satisfies the two
inverse properties, so Lemma \ref{lemma2} establishes that $B$ is
invertible.  Since $B$ is invertible and satisfies the
Yang--Baxter equation, Lemma \ref{lemma3} demonstrates that
$B^{-1}$ also satisfies the Yang--Baxter equation.  Finally, Lemma
\ref{lemma4} shows that $(Q, \rhd)$ is a right shelf, so that $(Q,
\rhd, \lhd)$ is a rack. \qed

Notice that we could have similarly shown that a right shelf $(Q,
\lhd)$ together with an operation $\rhd \maps Q \times Q \to Q$
that satisfies the inverse properties also satisfies the left
distributive law, $(i)$.  Therefore, we see that either the right
distributive law, $(ii)$, or the left distributive law, $(i)$, is
superfluous when defining a rack. Furthermore, we can omit either
the right idempotence law, $(vi)$, or the left idempotence law,
$(v)$, in the definition of a quandle because either one implies
the other with the help of the inverse properties.  For instance,
assuming $x \rhd x = x$, we have:
$$x \lhd x = (x \rhd x) \lhd x = x,$$ where the first equality
holds by the idempotence law $(v)$ and the second follows from the
inverse property $(iii)$.  We use the other inverse property to
show that $x \lhd x = x$ implies $x \rhd x = x$.  This result
together with Lemma \ref{rackwo(ii)} illustrates the power and
necessity of \emph{both} inverse properties.  That is, we have
taken advantage of the connection between algebra and topology to
show that we are unable to eliminate either of the two inverse
axioms from the definition of a rack.  This should not be terribly
surprising since one of the inverse laws says that $x \rhd -$ is a
left inverse of $- \lhd x$, while the other says that it is a
right inverse.

We conclude this section by illustrating the relation between the
first Reidemeister move and the idempotence conditions.  However,
because the first Reidemeister move involves `framed braids',
which are more general than braids, we are unable to apply our
braiding and inverse braiding maps to this move. Therefore, we
make the following conventions:
\[
 \xy 
 (0,10)*{}="T";
 (0,12)*{x};
 (0,-10)*{}="B";
 (0,-12)*{x \lhd x};
 (-10,0)*{x};
 (0,5)*{}="T'";
 (0,-5)*{}="B'";
 "T";"T'" **\dir{-};
 "B";"B'" **\dir{-};
 (-3,0)*{}="MB";
 (-7,0)*{}="LB";
 (-15,0)*{};
    "T'";"LB" **\crv{(-1,-4) & (-7,-4)}; \POS?(.25)*{\hole}="2z";
    "LB"; "2z" **\crv{(-8,6) & (-2,6)};
    "2z"; "B'"  **\crv{(0,-3)};
 \endxy
\quad = \quad \xy (0,12)*{x};
 (0,-12)*{x};
(0,10)*{}="T";
 (0,-10)*{}="B";
 "T";"B" **\dir{-};
\endxy
\quad = \quad
\xy 
 (0,10)*{}="T";
 (0,12)*{x};
 (10,0)*{x};
 (0,-10)*{}="B";
 (0,-12)*{x \rhd x};
 (0,5)*{}="T'";
 (0,-5)*{}="B'";
 "T";"T'" **\dir{-};
 "B";"B'" **\dir{-};
 (3,0)*{}="MB";
 (7,0)*{}="LB";
 (-7,0)*{};
    "T'";"LB" **\crv{(1,-4) & (7,-4)}; \POS?(.25)*{\hole}="2z";
    "LB"; "2z" **\crv{(8,6) & (2,6)};
    "2z"; "B'"  **\crv{(0,-3)};
 \endxy
\]

We continue to illustrate the relationship of shelves, racks,
spindles, and quandles to topology by showing how these structures
give invariants of different sorts of braids.

\subsection{Braid and Framed Braid Groups and Monoids} \label{braids}

Thus far we have demonstrated the relationship between the quandle
axioms and the Reidemeister moves.  We would now like to use that
relation to illustrate a stronger connection between these
algebraic structures and certain collections of braids.  For
reasons that are still somewhat mysterious, we will actually be
more interested in `quasi-idempotent' shelves, which are shelves
satisfying an additional property:

\begin{defn} \et \label{quasiidempotent}
A shelf $(Q, \rhd)$ is {\bf left quasi-idempotent} if it satisfies
\begin{description}
    \item[(vii)] {\bf (left quasi-idempotence)} $(x \rhd x) \rhd y = x \rhd y$ for all $x, y \in Q$.
\end{description}
\end{defn}
Of course, right shelves may also be quasi-idempotent and thus
satisfy $y \lhd (x \lhd x) = y \lhd x$.  We remark that quandles
are trivially left and right quasi-idempotent, and will show that
racks are quasi-idempotent as well.  Moreover, the rack axioms
imply other useful properties such as mixed distributive laws:

\begin{lem} \et \label{extrarackaxioms} Given a rack $(Q, \rhd, \lhd)$, the following
additional axioms hold:
\begin{itemize}
\item[(1)] {\bf (left quasi-idempotence)} $(x \rhd x) \rhd y = x \rhd y$
\item[(2)]{\bf (right quasi-idempotence)} $y \lhd (x \lhd x) = y \lhd x$
\item[(3)]{\bf (mixed distributive law)} $x \rhd (y \lhd z) = (x \rhd y) \lhd (x \rhd z)$
\item[(4)]{\bf (mixed distributive law)} $(x \rhd y) \lhd z = (x \lhd z) \rhd (y \lhd z)$
\end{itemize}
\end{lem}

\noindent{\bf Proof.} The mixed distributive laws imply the
quasi-idempotence laws, so we begin by verifying \emph{(3)} and
\emph{(4)}. In order to simplify the proofs, we will use the fact
that
$$x \rhd (y \lhd y) = A \quad \Leftrightarrow \quad [x \rhd (y
\lhd y)] \lhd x = A \lhd x \quad \Leftrightarrow \quad y \lhd y =
A \lhd x,$$ which holds by the left inverse property of racks.

Using this fact, demonstrating \emph{(3)} is equivalent to showing
$$[x \rhd (y \lhd z)] \lhd x = [(x \rhd y) \lhd (x \rhd z)] \lhd x\quad
\Leftrightarrow \quad y \lhd z = [(x \rhd y) \lhd (x \rhd z)] \lhd
x.$$  But, $$[(x \rhd y) \lhd (x \rhd z)] \lhd x  =  [(x \rhd y)
\lhd x] \lhd [(x \rhd z) \lhd x] = (y \lhd z)$$ where the first
equality holds by the self distributivity law, $(ii)$, and the
second equality follows from two applications of the right inverse
law, $(iii)$. Property \emph{(4)} is proved in a completely
analogous way using the other self distributivity law, $(i)$, and
the left inverse law, $(iv)$.

Again, using the above fact, property \emph{(1)} is equivalent to
$$[(x \rhd x) \rhd y] \lhd x = (x \rhd y) \lhd x \quad
\Leftrightarrow \quad [(x \rhd x) \rhd y] \lhd x = y$$ using the
right inverse property. But,
$$[(x \rhd x) \rhd y] \lhd x = [(x \rhd x) \lhd x] \rhd (y \lhd x)
= x \rhd (y \lhd x) = y$$ where the first equality results from
property \emph{(4)} above, the second equality holds by applying
the right inverse property to the first term, and the final
equality follows from the left inverse property.  Thus, \emph{(1)}
holds. Property \emph{(2)} is proved analogously, using \emph{(3)}
above. \qed

As suggested by the name, idempotence is stronger than
quasi-idempotence.  This follows from the fact that there exist
quandles which are not racks, such as the cyclic rack described in
the Section \ref{defns}.  We now turn to the task of showing how
the four categories of these algebraic structures are related to
certain braid groups.

In $1925$ E. Artin presented an algebraic description of braids.
We begin by reminding the reader of both the geometric and
algebraic descriptions of the braid group.

\begin{defn} \et Let $R$ be a rectangle in $3$-space.  On opposite
sides of $R$, label equidistant points $P_{i}$ and $Q_{i}$, where
$1 \leq i \leq n$.  Then, let $f_{i},$ $1 \leq i \leq n$ be $n$
simple disjoint smooth arcs in $\mathbb{R}^{3}$ which begin at
$P_{i}$ and end at $Q_{\tau(i)}$ where $i \mapsto \tau(i)$ is a
permutation on $\{1, 2, \ldots , n\}$.  We require that each
$f_{i}$ head downward as we move along any one of these $f_{i}$
from the top of $R$ to the bottom.  That is, each $f_{i}$
intersects any horizontal plane between the top and bottom of $R$
exactly once.  We say that the set of $n$ arcs $f_{i}$ constitute
an {\bf n-braid}.
\end{defn}
As an example:
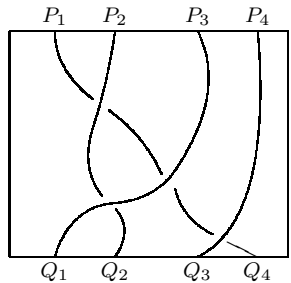
\[
\def\objectstyle{\scriptstyle}
\def\labelstyle{\scriptstyle}
    \xy
    (19,15)*{}="X";
    (-18,15)*{}="Y";
    (19,-15)*{}="X'";
    (-18,-15)*{}="Y'";
    (15,17)*{P_{4}}="P4";
    (7,17)*{P_{3}}="P3";
    (-4,17)*{P_{2}}="P2";
    (-12,17)*{P_{1}}="P1";
    (15,15)*{}="D";
    (7,15)*{}="C";
    (-4,15)*{}="B";
    (-12,15)*{}="A";
    (15,-15)*{}="3";
    (7,-15)*{}="4";
    (-4,-15)*{}="2";
    (-12,-15)*{}="1";
    (15,-17)*{Q_{4}}="Q4";
    (7,-17)*{Q_{3}}="Q3";
    (-4,-17)*{Q_{2}}="Q2";
    (-12,-17)*{Q_{1}}="Q1";
    (-7,6)*{}="A'";
    (-4.8,4.5)*{}="B'";
    (2.2,-4)*{}="B''";
    (4,-6)*{}="C'";
    (9,-12)*{}="B'''";
    (11,-13)*{}="D'";
    (-7,2)*{}="b";
    (4,-3.5)*{}="c";
    \vtwist~{(-7,2)}{(4,-3.5)}{"1"}{"2"};
      "A";"A'" **\crv{(-12,10)};
      "B'";"B''" **\crv{(0,1)};
      "C'";"B'''" **\crv{(5,-10)};
      "D'";"3" **\crv{};
      "B";"b" **\crv{(-5,8)};
      "C";"c" **\crv{(11,7)};
      "D";"4" **\crv{(17,-10)};
      "X";"Y" **\dir{-};
      "X'";"Y'" **\dir{-};
      "X";"X'" **\dir{-};
      "Y";"Y'" **\dir{-};
\endxy
\]
is a $4$-braid. Two braids are said to be \textbf{equivalent} if
they are ambient isotopic via an isotopy that fixes the boundary
points $P_{i}$ and $Q_{i}$.

Given two $n$-strand braids $b$ and $b'$ we may multiply, or
compose, them by identifying the terminal points $Q_{i}$ of $b$
with the initial points $P_{i}'$ of $b'$.  That is, we stack $b$
on top of $b'$:
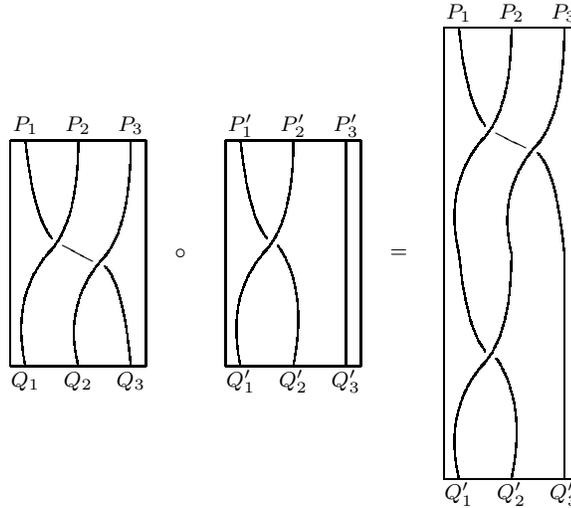
\[
\def\objectstyle{\scriptstyle}
\def\labelstyle{\scriptstyle}
    \xy
    (9,0)*{}="X";
    (-9,0)*{}="Y";
    (9,-30)*{}="X'";
    (-9,-30)*{}="Y'";
    (7,2)*{P_{3}}="P3";
    (0,2)*{P_{2}}="P2";
    (-7,2)*{P_{1}}="P1";
    (7,0)*{}="C";
    (0,0)*{}="B";
    (-7,0)*{}="A";
    (7,-30)*{}="3";
    (0,-30)*{}="2";
    (-7,-30)*{}="1";
    (7,-32)*{Q_{3}}="Q3";
    (0,-32)*{Q_{2}}="Q2";
    (-7,-32)*{Q_{1}}="Q1";
    (-3.5,-13)*{}="a";
    (-2,-14)*{}="b";
    (2,-16)*{}="c";
    (3.5,-16.7)*{}="d";
    (-4,-15)*{}="e";
    (3,-16)*{}="f";
      "X";"Y" **\dir{-};
      "X'";"Y'" **\dir{-};
      "X";"X'" **\dir{-};
      "Y";"Y'" **\dir{-};
      "B";"e" **\crv{(0,-11)};
      "C";"f" **\crv{(7,-12)};
      "e";"1" **\crv{(-9,-21)};
      "f";"2" **\crv{(-2,-22)};
      "A";"a" **\crv{(-6,-10)};
      "b";"c" **\crv{};
      "d";"3" **\crv{(6,-19)};
\endxy
\quad \xy
    (9,0)*{}="X";
    (-9,0)*{}="Y";
    (9,-30)*{}="X'";
    (-9,-30)*{}="Y'";
    (7,2)*{P_{3}'}="P3'";
    (0,2)*{P_{2}'}="P2'";
    (-7,2)*{P_{1}'}="P1'";
    (7,0)*{}="C";
    (0,0)*{}="B";
    (-7,0)*{}="A";
    (7,-30)*{}="3";
    (0,-30)*{}="2";
    (-7,-30)*{}="1";
    (7,-32)*{Q_{3}'}="Q3'";
    (0,-32)*{Q_{2}'}="Q2'";
    (-7,-32)*{Q_{1}'}="Q1'";
    (-3.5,-13)*{}="a";
    (-2,-14)*{}="b";
    (2,-16)*{}="c";
    (3.5,-16.7)*{}="d";
    (-4,-15)*{}="e";
    (3,-16)*{}="f";
    (-15,-15)*{\circ}="com";
      "X";"Y" **\dir{-};
      "X'";"Y'" **\dir{-};
      "X";"X'" **\dir{-};
      "Y";"Y'" **\dir{-};
      "C";"3" **\dir{-};
      "B";"e" **\crv{(0,-11)};
      "e";"1" **\crv{(-9,-21)};
      "A";"a" **\crv{(-6,-10)};
      "b";"2" **\crv{(2,-20)};
\endxy
\quad \xy
    (9,15)*{}="X";
    (-9,15)*{}="Y";
    (9,-15)*{}="X'";
    (-9,-15)*{}="Y'";
    (7,17)*{P_{3}}="P3";
    (0,17)*{P_{2}}="P2";
    (-7,17)*{P_{1}}="P1";
    (7,15)*{}="C";
    (0,15)*{}="B";
    (-7,15)*{}="A";
    (7,-15)*{}="3";
    (0,-15)*{}="2";
    (-7,-15)*{}="1";
    (-3.5,2)*{}="a";
    (-2,1)*{}="b";
    (2,-1)*{}="c";
    (3.5,-1.7)*{}="d";
    (-4,0)*{}="e";
    (3,-1)*{}="f";
    (-15, -15)*{=}="eq";
      "X";"Y" **\dir{-};
      "X";"X'" **\dir{-};
      "Y";"Y'" **\dir{-};
      "B";"e" **\crv{(0,4)};
      "C";"f" **\crv{(7,3)};
      "e";"1" **\crv{(-9,-6)};
      "f";"2" **\crv{(-2,-7)};
      "A";"a" **\crv{(-6,5)};
      "b";"c" **\crv{};
      "d";"3" **\crv{(6,-4)};
    (9,-45)*{}="X'a";
    (-9,-45)*{}="Y'a";
    (7,-15)*{}="Ca";
    (0,-15)*{}="Ba";
    (-7,-15)*{}="Aa";
    (7,-45)*{}="3a";
    (0,-45)*{}="2a";
    (-7,-45)*{}="1a";
    (7,-47)*{Q_{3}'}="Q3'a";
    (0,-47)*{Q_{2}'}="Q2'a";
    (-7,-47)*{Q_{1}'}="Q1'a";
    (-3.5,-28)*{}="aa";
    (-2,-29)*{}="ba";
    (2,-31)*{}="ca";
    (3.5,-31.7)*{}="da";
    (-4,-30)*{}="ea";
    (3,-31)*{}="fa";
      "Y'";"Y'a" **\dir{-};
      "X'";"X'a" **\dir{-};
      "X'a";"Y'a" **\dir{-};
      "Ca";"3a" **\dir{-};
      "Ba";"ea" **\crv{(0,-26)};
      "ea";"1a" **\crv{(-9,-36)};
      "Aa";"aa" **\crv{(-6,-25)};
      "ba";"2a" **\crv{(2,-35)};
\endxy
\]
to obtain the product.  There exists an obvious identity with
respect to this multiplication, namely the braid with no
crossings:
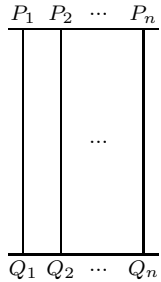
\[
\def\objectstyle{\scriptstyle}
\def\labelstyle{\scriptstyle}
    \xy
    (-23,17)*{P_{1}}="P1";
    (-18,17)*{P_{2}}="P2";
    (-13,17)*{\ldots}="P2'";
    (-7,17)*{P_{n}}="P3";
    (-25,15)*{}="X";
    (-23,15)*{}="A";
    (-18,15)*{}="B";
    (-7,15)*{}="C";
    (-5,15)*{}="Y";
    (-25,-15)*{}="X'";
    (-23,-17)*{Q_{1}}="Q1";
    (-18,-17)*{Q_{2}}="Q2";
    (-13,-17)*{\ldots}="Q2'";
    (-7,-17)*{Q_{n}}="Q3";
    (-23,-15)*{}="1";
    (-18,-15)*{}="2";
    (-7,-15)*{}="3";
    (-5,-15)*{}="Y'";
    (-13,0)*{\ldots}="d";
    "X";"Y" **\dir{-};
    "X'";"Y'" **\dir{-};
    "A";"1" **\dir{-};
    "B";"2" **\dir{-};
    "C";"3" **\dir{-};
\endxy
\]
In addition, each braid $b$
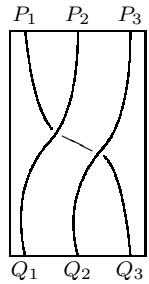
\[
\def\objectstyle{\scriptstyle}
\def\labelstyle{\scriptstyle}
\xy
    (9,15)*{}="X";
    (-9,15)*{}="Y";
    (9,-15)*{}="X'";
    (-9,-15)*{}="Y'";
    (7,17)*{P_{3}}="P3";
    (0,17)*{P_{2}}="P2";
    (-7,17)*{P_{1}}="P1";
    (7,15)*{}="C";
    (0,15)*{}="B";
    (-7,15)*{}="A";
    (7,-15)*{}="3";
    (0,-15)*{}="2";
    (-7,-15)*{}="1";
    (7,-17)*{Q_{3}}="Q3";
    (0,-17)*{Q_{2}}="Q2";
    (-7,-17)*{Q_{1}}="Q1";
    (-3.5,2)*{}="a";
    (-2,1)*{}="b";
    (2,-1)*{}="c";
    (3.5,-1.7)*{}="d";
    (-4,0)*{}="e";
    (3,-1)*{}="f";
      "X";"Y" **\dir{-};
      "X'";"Y'" **\dir{-};
      "X";"X'" **\dir{-};
      "Y";"Y'" **\dir{-};
      "B";"e" **\crv{(0,4)};
      "C";"f" **\crv{(7,3)};
      "e";"1" **\crv{(-9,-6)};
      "f";"2" **\crv{(-2,-7)};
      "A";"a" **\crv{(-6,5)};
      "b";"c" **\crv{};
      "d";"3" **\crv{(6,-4)};
\endxy
\]
has an inverse with respect to this multiplication, $b^{-1}$
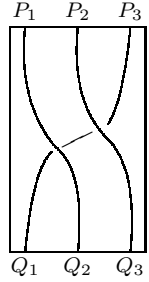
\[
\def\objectstyle{\scriptstyle}
\def\labelstyle{\scriptstyle}
\xy
    (9,15)*{}="X";
    (-9,15)*{}="Y";
    (9,-15)*{}="X'";
    (-9,-15)*{}="Y'";
    (7,17)*{P_{3}}="P3";
    (0,17)*{P_{2}}="P2";
    (-7,17)*{P_{1}}="P1";
    (7,15)*{}="C";
    (0,15)*{}="B";
    (-7,15)*{}="A";
    (7,-15)*{}="3";
    (0,-15)*{}="2";
    (-7,-15)*{}="1";
    (7,-17)*{Q_{3}}="Q3";
    (0,-17)*{Q_{2}}="Q2";
    (-7,-17)*{Q_{1}}="Q1";
    (4,2)*{}="a";
    (2,1)*{}="b";
    (-2,-1)*{}="c";
    (-3.5,-1.7)*{}="d";
    (4,0)*{}="e";
    (-3,-1)*{}="f";
      "X";"Y" **\dir{-};
      "X'";"Y'" **\dir{-};
      "X";"X'" **\dir{-};
      "Y";"Y'" **\dir{-};
      "B";"e" **\crv{(-1,5)};
      "e";"3" **\crv{(8,-4)};
      "f";"2" **\crv{(1,-4)};
      "A";"f" **\crv{(-8,5)};
      "C";"a" **\crv{(6,5)};
      "b";"c" **\crv{};
      "d";"1" **\crv{(-6,-4)};
\endxy
\]
obtained by reflecting $b$ in a plane perpendicular to it.  The
collection of $n$-braids together with the multiplication defined
above describes a group known as the {\bf braid group}, $B_{n}$.
To be precise, an element of this group is a braid $b$ together
with all the braids equivalent to it.  That is, an element of
$B_{n}$ is an equivalence class of braids, $[b]$.

Thus far we have described the braid group from a purely
topological viewpoint.  We can also give an algebraic description
of braids by listing which of the strands cross over and under one
another as we move down the braid.  We may always arrange our
braid so that no two crossing occur at exactly the same height. We
will use $s_{i}$ to denote the $(i+1)^{st}$ strand crossing over
the $i^{th}$ strand:
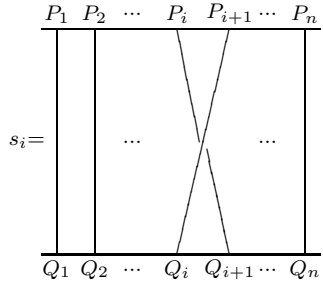
\[
\def\objectstyle{\scriptstyle}
\def\labelstyle{\scriptstyle}
    \xy
    (-23,17)*{P_{1}}="P1";
    (-18,17)*{P_{2}}="P2";
    (-13,17)*{\ldots}="P2'";
    (-7,17)*{P_{i}}="P3";
    (0,17)*{P_{i+1}}="P4";
    (5,17)*{\ldots}="P4'";
    (10,17)*{P_{n}}="P5";
    (-25,15)*{}="X";
    (-23,15)*{}="A";
    (-18,15)*{}="B";
    (-7,15)*{}="C";
    (0,15)*{}="D";
    (10,15)*{}="E";
    (12,15)*{}="Y";
    (-25,-15)*{}="X'";
    (-23,-17)*{Q_{1}}="Q1";
    (-18,-17)*{Q_{2}}="Q2";
    (-13,-17)*{\ldots}="Q2'";
    (-7,-17)*{Q_{i}}="Q3";
    (0,-17)*{Q_{i+1}}="Q4";
    (5,-17)*{\ldots}="Q4'";
    (10,-17)*{Q_{n}}="Q5";
    (-23,-15)*{}="1";
    (-18,-15)*{}="2";
    (-7,-15)*{}="3";
    (0,-15)*{}="4";
    (10,-15)*{}="5";
    (12,-15)*{}="Y'";
    (-13,0)*{\ldots}="d";
    (5,0)*{\ldots}="e";
    (-4,0)*{}="a";
    (-3,-1)*{}="b";
    (-27,0)*{s_{i}=}="c";
    "X";"Y" **\dir{-};
    "X'";"Y'" **\dir{-};
    "A";"1" **\dir{-};
    "B";"2" **\dir{-};
    "E";"5" **\dir{-};
    "D";"3" **\dir{-};
    "C";"a" **\dir{-};
    "b";"4" **\dir{-};
\endxy
\]
and $s_{i}^{-1}$ to denote the $i^{th}$ strand crossing over the
$(i+1)^{st}$ strand:
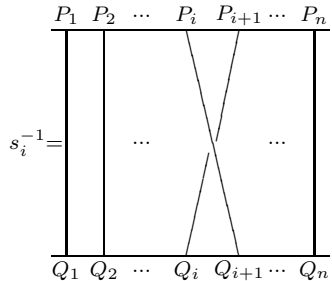
\[
\def\objectstyle{\scriptstyle}
\def\labelstyle{\scriptstyle}
    \xy
    (-23,17)*{P_{1}}="P1";
    (-18,17)*{P_{2}}="P2";
    (-13,17)*{\ldots}="P2'";
    (-7,17)*{P_{i}}="P3";
    (0,17)*{P_{i+1}}="P4";
    (5,17)*{\ldots}="P4'";
    (10,17)*{P_{n}}="P5";
    (-25,15)*{}="X";
    (-23,15)*{}="A";
    (-18,15)*{}="B";
    (-7,15)*{}="C";
    (0,15)*{}="D";
    (10,15)*{}="E";
    (12,15)*{}="Y";
    (-25,-15)*{}="X'";
    (-23,-17)*{Q_{1}}="Q1";
    (-18,-17)*{Q_{2}}="Q2";
    (-13,-17)*{\ldots}="Q2'";
    (-7,-17)*{Q_{i}}="Q3";
    (0,-17)*{Q_{i+1}}="Q4";
    (5,-17)*{\ldots}="Q4'";
    (10,-17)*{Q_{n}}="Q5";
    (-23,-15)*{}="1";
    (-18,-15)*{}="2";
    (-7,-15)*{}="3";
    (0,-15)*{}="4";
    (10,-15)*{}="5";
    (12,-15)*{}="Y'";
    (-3,.3)*{}="a";
    (-4,-1.3)*{}="b";
    (5,0)*{\ldots}="e";
    (-13,0)*{\ldots}="d";
    (-27,0)*{s_{i}^{-1}=}="c";
    "X";"Y" **\dir{-};
    "X'";"Y'" **\dir{-};
    "A";"1" **\dir{-};
    "B";"2" **\dir{-};
    "E";"5" **\dir{-};
    "C";"4" **\dir{-};
    "D";"a" **\dir{-};
    "b";"3" **\dir{-};
\endxy
\]
Using these two types of crossings, we may describe the braid
group in terms of generators and relations.  This presentation of
the braid group takes the following form:

\begin{prop} \et \label{artin} The braid group, $B_{n}$, may be presented as
$$B_{n} = \langle s_{1}, \ldots, s_{n-1} | s_{i} s_{i+1} s_{i}
= s_{i+1} s_{i} s_{i+1} \; \; (1 \leq i \leq n-2), \;\; s_{i}
s_{j} = s_{j} s_{i} \; \; (|i-j| > 1) \rangle.$$
\end{prop}
\noindent{\bf Proof.} This was proved by Artin \cite{A}; a
convenient reference is the book by Birman \cite{Birman}. \qed
Notice that the first relation, when drawn, is the third
Reidemeister move.  The second relation indicates that when two
pairs of strands are far enough apart, $|i - j| > 1$, the order in
which we braid them is irrelevant.

When we allow our strands to rotate, we obtain `framed braids':

\begin{defn} \et Let $D^{2}$ denote the closed unit ball in
$\R^{2}$.  Let $e_{i} \maps D^{2} \to \R^{2}$, $1 \leq i \leq n,$
be disjoint oriented balls embedded in $\R^{2}$.  Then an {\bf
n-strand framed braid} is an oriented embedding $F$ of the
disjoint union of $n$ solid cylinders $[0,1] \times D^{2}$ in
$[0,1] \times \R^{2}$, such that $F_{i}(0, \cdot) = F_{i}(1,
\cdot) = e_{\sigma(i)}$, for some $\sigma \in S_{n}$, where
$F_{i}$ denotes the embedding of the $i$th cylinder, and such that
$F_{i}(t,x) = (t, F_{i,t}(x))$ for some function $F_{i,t} \maps
D^{2} \to \R^{2}$.
\end{defn}
The set of homotopy classes of n-strand framed braids, where the
homotopy is required to preserve the above conditions on $F$, is
denoted by $FB_{n}$.  This becomes a group in an obvious way,
called the {\bf framed braid group}. The framed braid group
$FB_{n}$ has the braid group $B_{n}$ as a quotient group in an
obvious way, so there is a quotient map
$$\pi \maps FB_{n} \to B_{n}.$$  The framed braid group keeps
track of both the interchange of strands as well as their
rotations.  The rotations in the framed braid group correspond to
$2\pi$ rotations of the $i$th strand which will be drawn as
tangles:
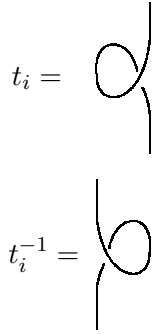
\[
 \xy 
 (0,10)*{}="T";
 (0,-10)*{}="B";
 (0,5)*{}="T'";
 (0,-5)*{}="B'";
 "T";"T'" **\dir{-};
 "B";"B'" **\dir{-};
 (-3,0)*{}="MB";
 (-7,0)*{}="LB";
 (-15,0)*{t_{i} =};
    "T'";"LB" **\crv{(-1,-4) & (-7,-4)}; \POS?(.25)*{\hole}="2z";
    "LB"; "2z" **\crv{(-8,6) & (-2,6)};
    "2z"; "B'"  **\crv{(0,-3)};
 \endxy
\]

\[\xy 
 (0,10)*{}="T";
 (0,-10)*{}="B";
 (0,5)*{}="T'";
 (0,-5)*{}="B'";
 "T";"T'" **\dir{-};
 "B";"B'" **\dir{-};
 (3,0)*{}="MB";
 (7,0)*{}="LB";
 (-7,0)*{t_{i}^{-1} =};
    "T'";"LB" **\crv{(1,-4) & (7,-4)}; \POS?(.25)*{\hole}="2z";
    "LB"; "2z" **\crv{(8,6) & (2,6)};
    "2z"; "B'"  **\crv{(0,-3)};
 \endxy
\]
As in the case of the braid group, the framed braid group can be
described in terms of generators and relations.  In addition to
the generators $s_{i}$ and $s_{i}^{-1}$ and the relations they
satisfy in the braid group, the framed braid group has tangles
$t_{i}$ and $t_{i}^{-1}$, $1 \leq i \leq n$, and relations:
\begin{eqnarray*}
t_i t_j &=& t_j t_i \\
t_i s_j &=& s_j t_i \; \; \; \textrm{if} \; \; i \geq j+2 \; \; \textrm{or} \; \; i < j \\
t_{i+1} s_i &=& s_i t_i \\
t_i s_i &=& s_i t_{i+1},
\end{eqnarray*}
These relations tell us how the braidings and twists interact with
one another.  We clearly have a natural inclusion: $$\iota \maps
B_{n} \hookrightarrow FB_{n}$$ such that $\iota \pi$ is the
identity on $B_{n}$, where $\pi$ is the quotient map mentioned
above that sends the generators $t_{i}$ to 1.

We can also consider the monoids generated by $s_{i}$ and $t_{i}$
without their inverses. The {\bf positive braid monoid}
$B_{n}^{+}$ is the monoid generated by the right-handed braidings,
$s_{i}$.  The {\bf positive framed braid monoid} $FB_{n}^{+}$ is
the monoid generated by the right-handed braidings, $s_{i}$, and
the right-handed twists, $t_{i}$. We have the following
commutative diagram of monoid homomorphisms:
\[
\xy (-15,0)*+{FB_{n}^{+}}="X1";
    (-15,-3.5)*{}="a";
    (15,0)*+{B_{n}^{+}}="X2";
    (15,-3.5)*{}="b";
    (-15,-30)*+{FB_{n}}="X3";
    (15,-30)*+{B_{n}}="X4";
             {\ar@{->>}"X1";"X2"};
             {\ar@{_{(}->} "a";"X3"};
             {\ar@{->>}"X3";"X4"};
             {\ar@{_{(}->} "b";"X4"};
\endxy
\]

We have a similar square illustrating the forgetful functors
between the categories of the algebraic structures described in
Section \ref{defns}:
$$\xymatrix{
   \textrm{Shelf}
    && \textrm{Spind}
     \ar[ll] \\ \\
   \textrm{Rack}
     \ar[uu]
     && \textrm{Quand}
     \ar[uu]
     \ar[ll] }$$
When considered together, these squares demonstrate a
contravariant relationship between the braid and framed braid
groups and monoids and our four categories:
\begin{thm} \et \label{squares} Given a set $Q$ equipped with operations $\rhd \maps Q \times Q \to Q$ and
$\lhd \maps Q \times Q \to Q$, and defining maps $s_{i}, t_{i},
s_{i}^{-1},$ and $t_{i}^{-1}$ on $Q^{n}$ as
\begin{itemize}
\item $s_{i} \maps (q_{1}, q_{2}, \ldots, q_{i}, \ldots, q_{n}) \mapsto (q_{1}, q_{2},
\ldots, q_{i-1}, q_{i+1}, q_{i+1} \rhd q_{i}, q_{i+2}, \ldots,
q_{n})$
\item $t_{i} \maps (q_{1}, q_{2}, \ldots, q_{i}, \ldots, q_{n}) \mapsto (q_{1},
q_{2}, \ldots, q_{i-1}, q_{i} \rhd q_{i}, q_{i+1}, \ldots q_{n})$
\item $s_{i}^{-1} \maps (q_{1}, q_{2}, \ldots, q_{i}, \ldots, q_{n}) \mapsto (q_{1}, q_{2},
\ldots, q_{i-1}, q_{i+1} \lhd q_{i}, q_{i}, q_{i+2}, \ldots,
q_{n})$
\item $t_{i}^{-1} \maps (q_{1}, q_{2}, \ldots, q_{i}, \ldots, q_{n}) \mapsto (q_{1},
q_{2}, \ldots, q_{i-1}, q_{i} \lhd q_{i}, q_{i+1}, \ldots q_{n})$
\end{itemize}
these maps give an action of $\begin{cases}
FB_{n}^{+} \\
B_{n}^{+} \\
FB_{n} \\
B_{n}
\end{cases}$ on $Q^{n}$ if and only if
$\begin{cases}
    \text{$(Q, \rhd)$ is a quasi-idempotent shelf} \\
    \text{$(Q, \rhd)$ is a spindle} \\
    \text{$(Q, \rhd, \lhd)$ is a rack} \\
    \text{$(Q, \rhd, \lhd)$ is a quandle}
  \end{cases}$
\end{thm}

\noindent{\bf Proof.}  We begin by proving that if

\begin{center}
$\begin{cases}
    \text{$(Q, \rhd)$ is a quasi-idempotent shelf} \\
    \text{$(Q, \rhd)$ is a spindle} \\
    \text{$(Q, \rhd, \lhd)$ is a rack} \\
    \text{$(Q, \rhd, \lhd)$ is a quandle}
  \end{cases}$
\\ \end{center}

\noindent then we obtain an action of \\
\begin{center}
$\begin{cases}
FB_{n}^{+} \\
B_{n}^{+} \\
FB_{n} \\
B_{n}
\end{cases}$
\end{center} on $Q^{n}$ where the generators $s_{i}, t_{i}, s_{i}^{-1},$ and
$t_{i}^{-1}$ act as above. It suffices for each case to show that
the relations in the appropriate collection of braids are
satisfied.

We first consider the case when $(Q, \rhd)$ is a quasi-idempotent
shelf.  Therefore, we must show that the relations in the positive
framed braid monoid hold.  We begin by examining the relations in
the braid group:
\begin{itemize}
\item[$(a)$] Let $q_{i} \in Q$. Assume $i < j$ and $|i-j|
> 1$. Then,
\begin{eqnarray*}
s_i s_j(q_1, \ldots, q_n) &=& s_i(q_1, q_2, \ldots, q_{j+1},
q_{j+1} \rhd q_{j}, q_{j+2}, \ldots, q_n) \\
&=& (q_1, q_2, \ldots, q_{i+1}, q_{i+1} \rhd q_{i}, q_{i+2},
\ldots, q_{j+1}, q_{j+1} \rhd q_{j}, q_{j+2}, \ldots, q_n) \\
&=& s_j s_i(q_1, q_2, \ldots, q_n).
\end{eqnarray*}

\item[$(b)$] Assume $1 \leq i \leq n-2$.  Then,
\begin{eqnarray*}
s_i s_{i+1}s_i(q_1, \ldots, q_n) &=& s_i s_{i+1}(q_1, q_2,
\ldots, q_{i+1}, q_{i+1} \rhd q_{i}, q_{i+2}, \ldots, q_n) \\
&=&s_i(q_1, q_2, \ldots, q_{i+1}, q_{i+2}, q_{i+2} \rhd (q_{i+1}
\rhd q_{i}), q_{i+3}, \ldots, q_n) \\
&=& (q_1, q_2, \ldots, q_{i+2}, q_{i+2} \rhd q_{i+1}, q_{i+2} \rhd
(q_{i+1} \rhd q_{i}), q_{i+3}, \ldots, q_n),
\end{eqnarray*}
while on the other hand,
\begin{eqnarray*}
s_{i+1} s_{i}s_{i+1}(q_1, \ldots, q_n) &=& s_{i+1} s_{i}(q_1,
\ldots, q_{i}, q_{i+2}, q_{i+2} \rhd q_{i+1}, \ldots, q_n) \\
&=& s_{i+1}(q_1, \ldots, q_{i+2}, q_{i+2} \rhd q_i, q_{i+2}
\rhd q_{i+1}, \ldots, q_n) \\
&=& (q_1, \ldots, q_{i+2}, q_{i+2} \rhd q_{i+1}, (q_{i+2} \rhd
q_{i+1}) \rhd (q_{i+2} \rhd q_{i}), \ldots, q_n).
\end{eqnarray*}
The final two expressions are equal by the left distributive law,
$(i)$, in the definition of a shelf. We continue by verifying the
relations involving the generators $t_{i}$:

\item[$(c)$]  Assume without loss of generality that $i < j$. We have
\begin{eqnarray*}
t_i t_j(q_1, q_2, \ldots, q_n) &=& t_i(q_1, q_2, \ldots, q_{j-1},
q_j \rhd q_j, q_{j+1}, \ldots, q_n) \\
&=& (q_1, q_2, \ldots, q_{i-1}, q_{i} \rhd q_i, q_{i+1}, \ldots,
q_{j-1}, q_j \rhd q_j, q_{j+1}, \ldots, q_n) \\
&=& t_jt_i(q_1, q_2, \ldots, q_n).
\end{eqnarray*}
We now consider the relations which describe how the braiding and
twisting may interact.

\item[$(d)$]  Let $i \geq j + 2$. Then,
\begin{eqnarray*}
t_{i} s_j(q_1, q_2, \ldots, q_n) &=& t_{i}(q_1, q_2, \ldots,
q_{j+1}, q_{j+1} \rhd q_{j}, q_{j+2}, \ldots, q_{i}, \ldots, q_n)
\\ &=& (q_1, q_2, \ldots,
q_{j+1}, q_{j+1} \rhd q_{j}, q_{j+2}, \ldots, q_{i} \rhd q_{i},
\ldots, q_n)
\end{eqnarray*}
while
\begin{eqnarray*}
s_jt_i(q_1, q_2, \ldots, q_n) &=& s_j(q_1, q_2, \ldots, q_{j},
\ldots, q_{i} \rhd q_{i}, \ldots, q_n)
\\ &=& (q_1, q_2, \ldots,
q_{j+1}, q_{j+1} \rhd q_{j}, q_{j+2}, \ldots, q_{i} \rhd q_{i},
\ldots, q_n)
\end{eqnarray*}
We arrive at the same conclusion when $i < j$, since $i$ and $j$
are sufficiently far enough apart so as to not interfere with one
another.

\item[$(e)$] Continuing on, we have
\begin{eqnarray*}
t_{i+1} s_i(q_1, q_2, \ldots, q_n) &=& t_{i+1}(q_1, q_2, \ldots,
q_{i+1}, q_{i+1} \rhd q_{i}, q_{i+2}, \ldots, q_n) \\
&=& (q_1, q_2, \ldots, q_{i+1}, (q_{i+1} \rhd q_{i}) \rhd (q_{i+1}
\rhd q_{i}), q_{i+2}, \ldots, q_n),
\end{eqnarray*}
while on the other hand,
\begin{eqnarray*}
s_i t_i (q_1, q_2, \ldots, q_n) &=& s_i(q_1, q_2,\ldots, q_{i-1},
q_i \rhd q_i, q_{i+1}, \ldots, q_n) \\
&=& (q_1, q_2, \ldots, q_{i-1}, q_{i+1}, q_{i+1} \rhd (q_{i} \rhd
q_{i}), q_{i+2}, \ldots, q_n).
\end{eqnarray*}
These statements are equal by the right-distributive law, $(i)$
satisfied by a shelf.

\item[$(f)$] Finally,
\begin{eqnarray*}
t_i s_i(q_1, q_2, \ldots, q_n) &=& t_i(q_1, q_2, \ldots, q_{i+1},
q_{i+1} \rhd q_{i}, q_{i+2}, \ldots, q_n) \\
&=& (q_1, q_2, \ldots, q_{i+1} \rhd q_{i+1}, q_{i+1} \rhd q_{i},
q_{i+2}, \ldots, q_n)
\end{eqnarray*}
while
\begin{eqnarray*}
s_i t_{i+1}(q_1, q_2, \ldots, q_n) &=& s_i(q_1, q_2,\ldots, q_{i},
q_{i+1} \rhd q_{i+1}, q_{i+2}, \ldots, q_n) \\
&=& (q_1, q_2,\ldots, q_{i+1} \rhd q_{i+1}, (q_{i+1} \rhd q_{i+1})
\rhd q_{i} , q_{i+2}, \ldots, q_n).
\end{eqnarray*}
These quantities are equal because of the quasi-idempotence of our
shelf.
\end{itemize}

We next consider the case when $Q$ is a spindle.  We need to
verify the relations of the positive braid monoid, and that the
map $t_{i}$ acts as the identity.  The two relations of the
positive braid monoid have already been demonstrated in $(a)$ and
$(b)$ above.  Therefore, we need only to examine the action of
$t_{i}$:
\begin{itemize}
\item[$(g)$]
\begin{eqnarray*} t_{i}(q_{1}, q_{2}, \ldots, q_{i}, \ldots, q_{n}) &=& (q_{1},
q_{2}, \ldots, q_{i-1}, q_{i} \rhd q_{i}, q_{i+1}, \ldots q_{n})
\\
&=& (q_{1}, q_{2}, \ldots, q_{i}, \ldots, q_{n}).
\end{eqnarray*}
\end{itemize}
Thus, because of the idempotence of a spindle, $t_{i}$ acts as the
identity, as desired.

We now turn to the case when $Q$ is a rack.  Since a rack is both
a left and right shelf, $(a)-(f)$ remain true.  All relations
involving inverses follow algebraically from these six.  It
remains to show that $s_{i}$ and $s_{i}^{-1}$ and $t_{i}$ and
$t_{i}^{-1}$ are indeed inverses.  We have:
\begin{eqnarray*}
s_{i}s_{i}^{-1}(q_1, q_2, \ldots, q_n) &=& s_{i}(q_1, q_2, \ldots
q_{i-1}, q_{i+1} \lhd q_{i}, q_{i},
q_{i+2}, \ldots, q_n) \\
&=&(q_1, q_2, \ldots, q_{i-1}, q_i, q_i \rhd (q_{i+1} \lhd q_{i}),
q_{i+2}, \ldots, q_n) \\
&=& (q_1, q_2, \ldots, q_i, q_{i+1}, q_{i+2}, \ldots, q_n),
\end{eqnarray*}
where the final equality is due to the right inverse property
applied to the $(i+1)$st entry.  Furthermore,
\begin{eqnarray*}
s_{i}^{-1}s_{i}(q_1, q_2, \ldots, q_n) &=& s_{i}^{-1}(q_1, q_2,
\ldots, q_{i-1}, q_{i+1}, q_{i+1} \rhd
q_{i}, q_{i+2}, \ldots, q_n) \\
&=& (q_1, q_2, \ldots, q_{i-1}, (q_{i+1} \rhd q_{i}) \lhd q_{i+1},
q_{i+1}, q_{i+2}, \ldots, q_n) \\
&=& (q_1, q_2, \ldots, q_{i-1}, q_i, q_{i+1}, \ldots, q_n),
\end{eqnarray*}
where the final equality follows from the left inverse property
applied to the $i$th entry.

Considering the twists, we have:
\begin{eqnarray*}
t_{i}^{-1}t_{i}(q_1, q_2, \ldots, q_n) &=& t_{i}^{-1}(q_1, q_2,
\ldots, q_{i-1}, q_i \rhd q_i, q_{i+1}, \ldots, q_n) \\
&=& (q_1, q_2, \ldots, q_{i-1}, (q_i \rhd q_i) \lhd (q_i \rhd
q_i), q_{i+1}, \ldots, q_n) \\
&=& (q_1, q_2, \ldots, q_i \rhd (q_i \lhd q_i), \ldots, q_n) \\
&=& (q_1, q_2, \ldots, q_i, \ldots, q_n)
\end{eqnarray*}
where the second to last equality is the mixed distributive law
\emph{(3)} from Lemma \ref{extrarackaxioms}, and the final
equality is due to the left inverse property of a rack.  Moreover,
\begin{eqnarray*}
t_i t_{i}^{-1}(q_1, q_2, \ldots, q_n) &=& t_i(q_1, q_2, \ldots,
q_{i-1}, q_i \lhd q_i, q_{i+1}, \ldots, q_n) \\
&=& (q_1, q_2, \ldots, q_{i-1}, (q_i \lhd q_i) \rhd (q_i \lhd
q_i), q_{i+1}, \ldots, q_n) \\
&=& (q_1, q_2, \ldots, (q_i \rhd q_i) \lhd q_i, \ldots, q_n) \\
&=& (q_1, q_2, \ldots, q_i, \ldots q_n),
\end{eqnarray*}
where the second to last equality is the mixed distributive law
\emph{(4)} from Lemma \ref{extrarackaxioms}, and the final
equality is due to the right inverse property of a rack.  Thus,
$s_{i}s_{i} ^{-1} = s_{i}^{-1}s_{i} = 1_{Q^{n}}$ and
$t_{i}t_{i}^{-1} = t_{i}^{-1}t_{i} = 1_{Q^{n}}$, as desired.

Finally, we consider the case when $Q$ is a quandle. Since a
quandle is a spindle, $(a)$ and $(b)$ above still hold. The
relations involving inverses follow algebraically from these two.
In the previous case we demonstrated that $s_i$ and $s_{i}^{-1}$
are indeed inverses.  It remains to show that $t_{i}$ acts as the
identity.  This again follows from the fact that a quandle is a
spindle.

Conversely, we show that given a set $Q$ equipped with operations
$\rhd \maps Q \times Q \to Q$ and $\lhd \maps Q \times Q \to Q$,
and defining maps $s_{i}, t_{i}, s_{i}^{-1},$ and $t_{i}^{-1}$ on
$Q^{n}$ as above, these maps give an action of

\begin{center}
$\begin{cases}
FB_{n}^{+} \\
B_{n}^{+} \\
FB_{n} \\
B_{n}
\end{cases}$
\\
\end{center}

\noindent on $Q^{n}$ if \\

\begin{center}
$\begin{cases}
    \text{$(Q, \rhd)$ is a quasi-idempotent shelf} \\
    \text{$(Q, \rhd)$ is a spindle} \\
    \text{$(Q, \rhd, \lhd)$ is a rack} \\
    \text{$(Q, \rhd, \lhd)$ is a quandle.}
  \end{cases}$
\end{center}

Consider a set $Q$ equipped with an operation $\rhd \maps Q \times
Q \to Q$ together with maps $s_{i}$ and $t_{i}$ defined as above
which give an action of $FB_{n}^{+}$.  The relations $s_{i}
s_{i+1} s_{i} = s_{i+1} s_{i} s_{i+1} $ and $t_{i}s_{i} =
s_{i}t_{i+1}$ imply the distributive and quasi-idempotence laws,
so that $(Q, \rhd)$ is a quasi-idempotent shelf.

Next, consider a set $Q$ equipped with an operation $\rhd \maps Q
\times Q \to Q$ and maps $s_{i}$ and $t_{i}$ defined as above
which give an action of $B_{n}^{+}$. Again, the relation $s_{i}
s_{i+1} s_{i} = s_{i+1} s_{i} s_{i+1}$ gives the distributive law.
Since $t_{i}$ acts as the identity in $B_{n}^{+},$ we have
$$(q_{1}, q_{2}, \ldots, q_{n}) = t_{i}(q_{1}, q_{2}, \ldots, q_{i}, \ldots, q_{n}) = (q_{1},
q_{2}, \ldots, q_{i-1}, q_{i} \rhd q_{i}, q_{i+1}, \ldots q_{n})$$
so that $q_{i} \rhd q_{i} = q_{i}$, which is the idempotence
condition.  Thus, $(Q, \rhd)$ is a spindle.

Consider a set $Q$ equipped with operations $\rhd \maps Q \times Q
\to Q$ and $\lhd \maps Q \times Q \to Q$, together with maps
$s_{i}, t_{i}, s_{i}^{-1}$ and $t_{i}^{-1}$ defined as above which
give an action of $FB_{n}$.  As in the previous two cases, the
relation $s_{i} s_{i+1} s_{i} = s_{i+1} s_{i} s_{i+1}$ implies the
distributive law.  Since $s_{i}$ and $s_{i}^{-1}$ are inverses, we
have:
\begin{eqnarray*}
(q_1, q_2, \ldots, q_n) = s_{i}s_{i}^{-1}(q_1, q_2, \ldots, q_n)
&=& s_{i}(q_1, q_2, \ldots q_{i-1}, q_{i+1} \lhd q_{i}, q_{i},
q_{i+2}, \ldots, q_n) \\
&=&(q_1, q_2, \ldots, q_{i-1}, q_i, q_i \rhd (q_{i+1} \lhd q_{i}),
q_{i+2}, \ldots, q_n),
\end{eqnarray*}
and
\begin{eqnarray*}
(q_1, q_2, \ldots, q_n) =
s_{i}^{-1}s_{i}(q_1, q_2, \ldots, q_n) &=& s_{i}^{-1}(q_1, q_2,
\ldots, q_{i-1}, q_{i+1}, q_{i+1} \rhd
q_{i}, q_{i+2}, \ldots, q_n) \\
&=& (q_1, q_2, \ldots, q_{i-1}, (q_{i+1} \rhd q_{i}) \lhd q_{i+1},
q_{i+1}, q_{i+2}, \ldots, q_n),
\end{eqnarray*}
so that $q_i \rhd (q_{i+1} \lhd q_{i}) = q_{i+1}$ and $(q_{i+1}
\rhd q_{i}) \lhd q_{i+1} = q_{i}$, which are the two inverse
properties.  Thus, $(Q, \rhd, \lhd)$ is a rack.

Finally, consider a set $Q$ equipped with operations $\rhd \maps Q
\times Q \to Q$ and $\lhd \maps Q \times Q \to Q$, together with
maps $s_{i}, t_{i}, s_{i}^{-1}$ and $t_{i}^{-1}$ defined as above
which give an action of $B_{n}$.  The distributive law follows
from the relation $s_{i} s_{i+1} s_{i} = s_{i+1} s_{i} s_{i+1}$,
and the two inverse properties follow from the fact that $s_{i}$
and $s_{i}^{-1}$ are inverses.  As in the spindle case, the
idempotence law is a result of the fact that $t_{i}$ acts as the
identity.  Thus, $(Q, \rhd, \lhd)$ is a quandle. \qed

The moral of this theorem is that the relations in these four
braid and framed braid groups and monoids {\it force} the shelf,
rack, spindle, and quandle laws to hold!  Therefore, if we
desired, we could actually {\it define} a shelf, rack, spindle, or
quandle to be a set equipped with binary operations, $(Q, \rhd,
\lhd)$, such that if we defined maps on $Q^{n}$ by the above
formulas, we would obtain an action of the appropriate braid or
framed braid group or monoid.  Doing so would enable us to easily
categorify these concepts, since the categorified versions of the
braid and framed braid groups and monoids are known.

Now that we have exhibited the connection between the categories
of shelves, racks, spindles and quandles, we will use them to
describe a new means by which we can obtain the Lie algebra of a
Lie group.  Specifically, we will show that we can think of our
Lie group as a spindle in $\Diff_{\ast}$, just as any group gives
a spindle.  Therefore, we continue in the next section by
internalizing the concepts of shelf, rack, spindle and quandle.

\subsection{Internalization} \label{internalization}
Our goal is to present a conceptual explanation of the passage
from a Lie group to its Lie algebra using the language of
spindles. Seeing as how the bracket arises from differentiating
conjugation, and spindles possess essential properties of
conjugation, we desire a way to treat our Lie groups as though
they were spindles. This idea should not be so surprising since we
saw in Section \ref{topology2} that the self-distributive law of a
shelf and the Jacobi identity of a Lie algebra are each equivalent
to the Yang--Baxter equation!  As mentioned previously, however,
in addition to the self-distributive law of a shelf, we need the
idempotence law of a spindle to obtain the antisymmetry of the
bracket.

Using the language of internalization, as was done in Section
\ref{2vs} of Chapter \ref{ch1} to define the notion of a
$2$-vector space, a Lie group can be thought of as a `group in
$\Diff_{\ast}$'.  Then, we can internalize the relationship
between groups and spindles to show that a Lie group is also a
`spindle in $\Diff_{\ast}$'.  Though we will only use the
internalized version of a spindle, we introduce the internalized
versions of all of the concepts from Section \ref{defns} for
completeness. In the following, we consider a category $K$ with
finite products.

\begin{defn} \et \label{lshelfobj}
A {\bf left shelf in $K$} consists of:
\begin{itemize}
\item an object $Q \in K$,
\end{itemize}
equipped with:
\begin{itemize}
\item a {\bf left conjugation} morphism $\rhd \maps Q \times Q \to Q$,
\end{itemize}
such that the following diagram commutes, expressing the usual
left distributive law:
\begin{itemize}
\item the {\bf left distributive law}, ($i^{ \prime}$):
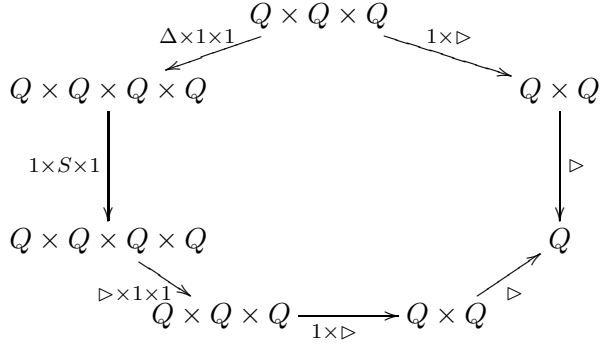
\[
\begin{xy}
    (-2, 20)*+{Q \times Q \times Q}="1";
    (-30,10)*+{Q \times Q \times Q \times Q}="2";
    (-30,-10)*+{Q \times Q \times Q \times Q}="3";
    (-15,-20)*+{Q \times Q \times Q}="4";
    (15,-20)*+{Q \times Q}="4a";
    (30,-10)*+{Q}="5";
    (30, 10)*+{Q \times Q}="6";
        {\ar^{1 \times \rhd} "1";"6"};
        {\ar_{\Delta \times 1 \times 1} "1";"2"};
        {\ar_{1 \times S \times 1} "2";"3"};
        {\ar_{\rhd \times 1 \times 1} "3";"4"};
        {\ar_{1 \times \rhd} "4";"4a"};
        {\ar_{\rhd} "4a";"5"};
        {\ar^{\rhd} "6";"5"};
\end{xy}
\]
\end{itemize}
where $\Delta \maps Q \to Q \times Q$ is the diagonal morphism in
$K$ and $S \maps Q \times Q \to Q \times Q$ is a morphism which
switches its inputs.
\end{defn}

\begin{defn} \et \label{rshelfobj}
A {\bf right shelf in $K$} consists of:
\begin{itemize}
\item an object $Q \in K$,
\end{itemize}
equipped with:
\begin{itemize}
\item a {\bf right conjugation} morphism $\lhd \maps Q \times Q \to Q$,
\end{itemize}
such that the following diagram commutes, expressing the usual
right distributive law:
\begin{itemize}
\item the {\bf right distributive law}, ($ii^{ \prime}$):
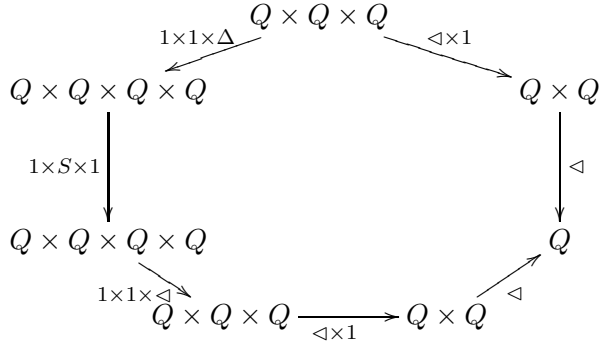
\[\begin{xy}
    (-2, 20)*+{Q \times Q \times Q}="1";
    (-30,10)*+{Q \times Q \times Q \times Q}="2";
    (-30,-10)*+{Q \times Q \times Q \times Q}="3";
    (-15,-20)*+{Q \times Q \times Q}="4";
    (15,-20)*+{Q \times Q}="4a";
    (30,-10)*+{Q}="5";
    (30, 10)*+{Q \times Q}="6";
        {\ar^{\lhd \times 1} "1";"6"};
        {\ar_{1 \times 1 \times \Delta} "1";"2"};
        {\ar_{1 \times S \times 1} "2";"3"};
        {\ar_{1 \times 1 \times \lhd} "3";"4"};
        {\ar_{\lhd \times 1} "4";"4a"};
        {\ar_{\lhd} "4a";"5"};
        {\ar^{\lhd} "6";"5"};
\end{xy}
\]
\end{itemize}
where $\Delta \maps Q \to Q \times Q$ is the diagonal morphism in
$K$ and $S \maps Q \times Q \to Q \times Q$ is a morphism which
switches its inputs.
\end{defn}

\begin{defn} \et \label{rackobj}
A {\bf rack in $K$} consists of:
\begin{itemize}
\item an object $Q \in K$,
\end{itemize}
equipped with:
\begin{itemize}
\item a {\bf left conjugation} morphism $\rhd \maps Q \times Q \to Q$,
\item a {\bf right conjugation} morphism $\lhd \maps Q \times Q \to Q$,
\end{itemize}
such that $(Q, \rhd)$ is a left shelf in $K$, $(Q, \lhd)$ is a
right shelf in $K$, and the following diagrams commute, expressing
the usual left and right inverse laws:
\begin{itemize}
\item the {\bf left inverse law}, ($iii^{ \prime}$):
\[
  \vcenter{
  \xymatrix@!C{
    Q \times Q
    \ar[ddrr]_{\pi_2}
    \ar[r]^{\Delta \times 1}
  & Q \times Q \times Q \ar[r]^{1 \times S}
  & Q \times Q \times Q \ar[d]^{\rhd \times 1} \\
  && Q \times Q \ar[d]^{\lhd} \\
  && Q
     }}
\]
\item the {\bf right inverse law}, ($iv^{ \prime}$):
\[
  \vcenter{
  \xymatrix@!C{
    Q \times Q
    \ar[ddrr]_{\pi_2}
    \ar[r]^{\Delta \times 1}
  & Q \times Q \times Q \ar[r]^{1 \times S}
  & Q \times Q \times Q \ar[d]^{1 \times \lhd} \\
  && Q \times Q \ar[d]^{\rhd} \\
  && Q
     }}
\]
\end{itemize}
where $\Delta \maps Q \to Q \times Q$ is the diagonal morphism in
$K$ and $S \maps Q \times Q \to Q \times Q$ is a morphism which
switches its inputs.
\end{defn}

\begin{defn} \et \label{spindleobj}
A {\bf left spindle in $K$} consists of:
\begin{itemize}
\item an object $Q \in K$,
\end{itemize}
equipped with:
\begin{itemize}
\item a {\bf left conjugation} morphism $\rhd \maps Q \times Q \to Q$,
\end{itemize}
such that $(Q, \rhd)$ is a left shelf in $K$ and the following
diagram commutes, expressing the usual left idempotence law:
\begin{itemize}
\item the {\bf left idempotence law}, ($v^{ \prime}$):
\[
  \vcenter{
  \xymatrix@!C{
  & Q \times Q \ar[d]^{\rhd} \\
    Q
    \ar[ur]^{\Delta}
    \ar[r]_{1}
  & Q
     }}
\]
\end{itemize}
where $\Delta \maps Q \to Q \times Q$ is the diagonal morphism in
$K$.
\end{defn}

\begin{defn} \et \label{quandleobj}
A {\bf quandle in $K$} consists of:
\begin{itemize}
\item an object $Q \in K$,
\end{itemize}
equipped with:
\begin{itemize}
\item a {\bf left conjugation} morphism $\rhd \maps Q \times Q \to Q$,
\item a {\bf right conjugation} morphism $\lhd \maps Q \times Q \to Q$,
\end{itemize}
such that $(Q, \rhd, \lhd)$ is a rack in $K$ and the following
diagrams commute, expressing the usual left and right idempotence
laws:
\begin{itemize}
\item  the {\bf left idempotence law}, ($v^{ \prime}$):
\[
  \vcenter{
  \xymatrix@!C{
  & Q \times Q \ar[d]^{\rhd} \\
    Q
    \ar[ur]^{\Delta}
    \ar[r]_{1}
  & Q
     }}
\]
\item the {\bf right idempotence law}, ($vi^{ \prime}$):
\[
  \vcenter{
  \xymatrix@!C{
  & Q \times Q \ar[d]^{\lhd} \\
    Q
    \ar[ur]^{\Delta}
    \ar[r]_{1}
  & Q
     }}
\]
\end{itemize}
where $\Delta \maps Q \to Q \times Q$ is the diagonal morphism in
$K$.
\end{defn}
In the special case where $K = \Set$, these concepts in $K$ reduce
to the ordinary notions of shelf, rack, spindle and quandle.
Therefore, we should expect that we can internalize the close
relationship between quandles and groups, and show that a group in
$K$, a category with finite products, gives a quandle in $K$. We
first recall the definition of an internalized group:

\begin{defn} \et \label{grpobj}
A {\bf group in $K$} consists of:
\begin{itemize}
\item an object $G \in K$,
\end{itemize}
equipped with:
\begin{itemize}
\item a {\bf multiplication} morphism $m \maps G \times G \to G$
\item an {\bf identity} morphism $\id \maps I \to G$, where $I$ is the
terminal object in $K$
\item an {\bf inverse} morphism $\inv \maps G \to G$,
\end{itemize}
such that the following diagrams commute, expressing the usual
group laws:
\begin{itemize}
\item the {\bf associative law}:
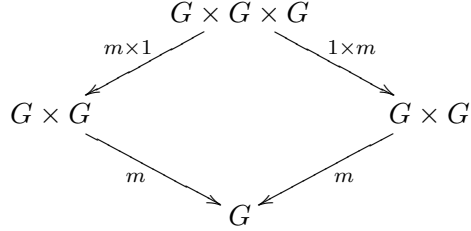
\[ \vcenter{
\xymatrix{ &   G \times G \times G \ar[dr]^{1 \times m}
   \ar[dl]_{m \times 1} \\
 G \times G \ar[dr]_{m}
&&  G \times G \ar[dl]^{m}  \\
&  G }}
\]
\item the {\bf right and left unit laws}:
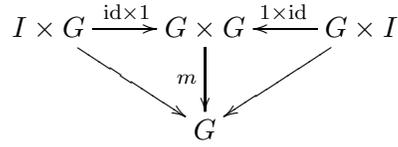
\[ \vcenter{
\xymatrix{
 I \times G \ar[r]^{\id \times 1} \ar[dr]
& G \times G \ar[d]_{m}
& G \times I \ar[l]_{1 \times \id} \ar[dl] \\
& G }}
\]
\item the {\bf right and left inverse laws}:
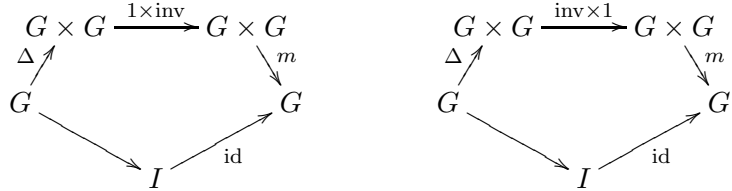
\[
\xy (-12,10)*+{G \times G}="TL"; (12,10)*+{G \times G}="TR";
(-18,0)*+{G}="ML"; (18,0)*+{G}="MR"; (0,-10)*+{I}="B";
     {\ar_{} "ML";"B"};
     {\ar^{\Delta} "ML";"TL"};
     {\ar_{\id} "B";"MR"};
     {\ar^{m} "TR";"MR"};
     {\ar^{1 \times \inv } "TL";"TR"};
\endxy
\qquad \qquad \xy (-12,10)*+{G \times G}="TL"; (12,10)*+{G \times
G}="TR"; (-18,0)*+{G}="ML"; (18,0)*+{G}="MR"; (0,-10)*+{I}="B";
     {\ar_{} "ML";"B"};
     {\ar^{\Delta} "ML";"TL"};
     {\ar_{\id} "B";"MR"};
     {\ar^{m} "TR";"MR"};
     {\ar^{\inv \times 1} "TL";"TR"};
\endxy
\]
\end{itemize}
where $\Delta \maps G \to G \times G$ is the diagonal morphism in
$K$.
\end{defn}
Indeed, for any category $K$ with finite products, there is a
category $K\Grp$ consisting of groups in $K$ and homomorphisms
between these, where a {\bf homomorphism} $f \maps G \to G'$ is a
morphism in $K$ that preserves multiplication, meaning that this
diagram commutes:
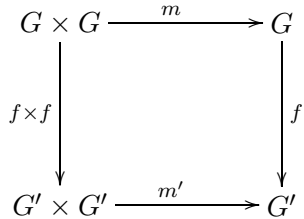
\[
\xymatrix{ G \times G
 \ar[rr]^{m}
 \ar[dd]_{f \times f}
  && G
 \ar[dd]^{f} \\ \\
  G' \times G'
 \ar[rr]^{m'}
  && G'}
\]

Given a group in $K$, say $G$, we can construct a quandle in $K$
by defining the left and right conjugation morphisms in terms of
conjugation in $G$, just as we did when $K = \Set$. That is, we
define left conjugation $\rhd \maps G \times G \to G$ by the
composite:
\[
 \xy
     (-70,0)*+{(x,y)}="1";
     (-50,0)*+{(x,x,y)}="2";
     (-30,0)*+{(x,y,x)}="3a";
     (-0,0)*+{(x,y,x^{-1})}="3";
     (25,0)*+{(x,yx^{-1})}="4";
     (50,0)*+{xyx^{-1}}="5";
            {\ar^{\Delta \times 1} "1";"2"};
            {\ar^{1 \times S}      "2";"3a"};
            {\ar^<<<<<<<{1 \times 1 \times \inv} "3a";"3"};
            {\ar^{1 \times m} "3";"4"};
            {\ar^{m} "4";"5"};
 \endxy
\]
while right conjugation $\lhd \maps G \times G \to G$ is defined
by the composite:
\[
 \xy
     (-70,0)*+{(x,y)}="1";
     (-50,0)*+{(x,y,y)}="2";
     (-30,0)*+{(y,x,y)}="3a";
     (-0,0)*+{(y^{-1},x,y)}="3";
     (25,0)*+{(y^{-1}x,y)}="4";
     (50,0)*+{y^{-1}xy}="5";
            {\ar^{1 \times \Delta} "1";"2"};
            {\ar^{S \times 1}      "2";"3a"};
            {\ar^<<<<<<<{\inv \times 1 \times 1} "3a";"3"};
            {\ar^{m \times 1} "3";"4"};
            {\ar^{m} "4";"5"};
 \endxy
\]
We show that these operations satisfy the quandle axioms by simply
internalizing the calculations which showed that left and right
conjugation in groups satisfy the quandle laws.  For example, the
left idempotence law for left conjugation holds since
$$x \rhd x = xxx^{-1} = x1 = x,$$
so that we expect the internalized version to use the internalized
right inverse and left unit laws.  The left idempotence law for a
quandle in $K$ follows from this diagram:
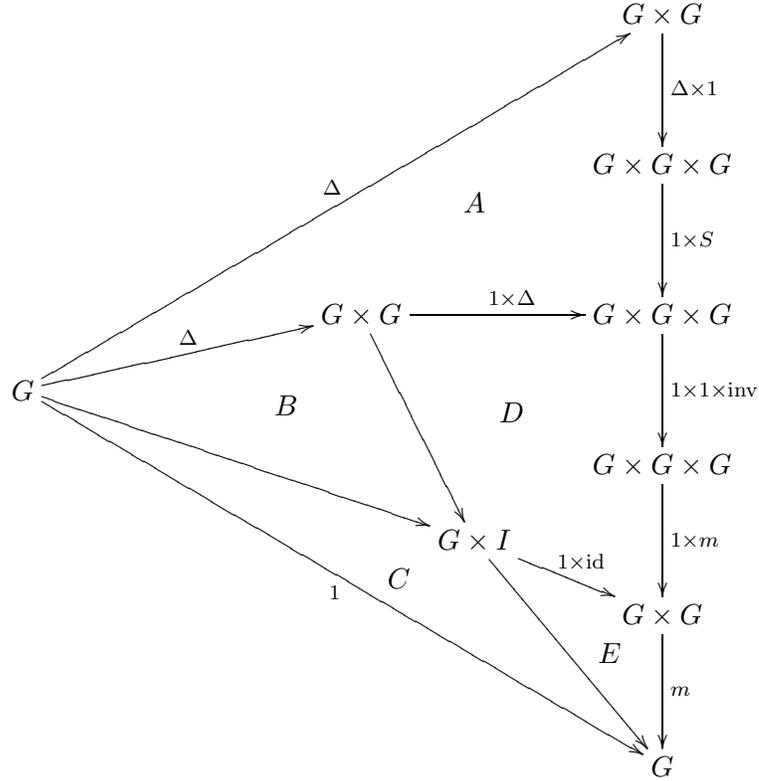
\[
\xy
    (0,50)*+{G \times G}="1";
    (0,30)*+{G \times G \times G}="2";
    (0,10)*+{G \times G \times G}="3";
    (0,-10)*+{G \times G \times G}="4";
    (0,-30)*+{G \times G}="5";
    (0,-50)*+{G}="6";
    (-85,0)*+{G}="7";
    (-40,10)*+{G \times G}="8";
    (-25,-20)*+{G \times I}="9";
    (-25,25)*{A}="10";
    (-50,-2)*{B}="11";
    (-35,-25)*{C}="12";
    (-20,-3)*{D}="13";
    (-7,-35)*{E}="14";
        {\ar^{\Delta}   "7";"1"};
        {\ar^{\Delta \times 1}  "1";"2"};
        {\ar^{1 \times S}   "2";"3"};
        {\ar^{1 \times 1 \times \inv}    "3";"4"};
        {\ar^{1 \times m}   "4";"5"};
        {\ar^{m}    "5";"6"};
        {\ar_{1}    "7";"6"};
        {\ar^{\Delta}   "7";"8"};
        {\ar^{1 \times \id}  "9";"5"};
        {\ar^{}             "7";"9"};
        {\ar^{}             "8";"9"};
        {\ar^{}             "9";"6"};
        {\ar^{1 \times \Delta}  "8";"3"};
\endxy
\]
Regions $A$, $B$, and $C$ commute trivially.  Region $D$ is the
right inverse law multiplied on the left by $1$, and region $E$ is
the left unit law.

The remaining internalized quandle axioms follow from similar
sorts of diagrams.  Thus, we obtain quandles in $K$ from groups in
$K$ just as we get quandles from groups.  We continue by showing
that there is an analogue of Proposition \ref{rackwo(ii)} for rack
objects. As in Section \ref{topology2}, we begin by illustrating
the relationship between the left shelf axiom and the Yang--Baxter
equation:

\begin{lem} \et The braiding $B \maps Q \times Q \to Q \times Q$
defined as $B = S(\Delta \times 1)(1 \times \rhd)$ satisfies the
Yang-Baxter equation if and only if $(Q, \rhd)$ is a left shelf
object.
\end{lem}

\noindent {\bf Proof.} The proof of this statement amounts to a
rather large, though mostly trivial, diagram that is nothing more
than an internalization of the proof of Lemma \ref{lemma1}. \qed
Moreover, internalizing the proofs of Lemma \ref{lemma4},
Corollary \ref{lemma2a}, and Lemma \ref{lemma3} results in the
following:

\begin{lem} \et The inverse braiding, $B^{-1}=S(1 \times \Delta)(\lhd
\times 1)$, satisfies the Yang--Baxter equation if and only if
$(Q, \lhd)$ is a right shelf object.
\end{lem}

\begin{lem} \et Let $Q$ be an object of a category $K$ with finite products
equipped with two morphisms $\rhd \maps Q \times Q \to Q$ and
$\lhd \maps Q \times Q \to Q$ such that $(Q, \rhd)$ is a left
shelf in $K$ and $(Q, \lhd)$ is a right shelf in $K$.  Then, the
braidings $B$ and $B^{-1}$ are inverses if and only if $(Q, \rhd,
\lhd)$ is a rack in $K$.
\end{lem}

\begin{lem} \et If $B$ is invertible and satisfies the Yang--Baxter equation,
then $B^{-1}$ does as well.
\end{lem}

Thus, a shelf in $K$ together with an additional morphism which
satisfies the two internalized inverse properties,
$(iii^{\prime})$ and $(iv^{\prime})$, is a rack in $K$.  This
statement follows from the previous lemmas as in the case of
racks.

\begin{prop} \et A left shelf $(Q, \rhd)$ in $K$ together with a
morphism $\lhd \maps Q \times Q \to Q$ satisfying the internalized
inverse properties $(iii^{\prime})$ and $(iv^{\prime})$ also
satisfies the internalized right distributive law,
$(ii^{\prime})$. That is, $(Q, \rhd, \lhd)$ is a rack in $K$.
\end{prop}

Just as we demonstrated that we could obtain one of the idempotent
laws from the other for quandles, the corresponding result holds
for quandles in $K$.  Now, however, instead of equations, we have
commutative diagrams.

\begin{prop} \label{internalidemp} \et In the definition of a quandle in $K$, the left
idempotence law, $(v^ {\prime})$, is satisfied if and only if the
right idempotence law, $(vi^{ \prime})$, is satisfied.
\end{prop}

\noindent {\bf Proof.}  We prove that the left idempotence law
implies the right.  That is, we will show that diagram $(vi^{
\prime})$ commutes using the two inverse properties, $(iii^{
\prime})$ and $(v^{ \prime})$.  The other implication is proved
analogously.

In the figure below, regions $A$, $B$, and $E$ clearly commute.
Region $C$ commutes because it is the internalized left inverse
property, axiom $(iii^{ \prime})$, and region $D$ results when
multiplying the internalized left idempotence property, axiom
$(v^{ \prime})$, on the right by $1$ and therefore commutes:
\[ \vcenter{
\begin{xy}
    (-60, 0)*+{Q}="1";
    (60,50)*+{Q \times Q}="2";
    (60,-50)*+{Q}="3";
    (8,10)*+{Q \times Q \times Q}="4";
    (30,20)*+{Q \times Q \times Q}="5";
    (45,-10)*+{Q \times Q}="6";
    (10,-10)*+{Q \times Q}="7";
    (2,-18)*{A}="8";
    (-13,8)*{B}="9";
    (28,-7)*{C}="10";
    (55,-5)*{E}="11";
    (45,10)*{D}="12";
        {\ar^{\Delta} "1";"2"};
        {\ar_{1} "1";"3"};
        {\ar^{\lhd} "2";"3"};
        {\ar^{\Delta} "1";"7"};
        {\ar^{\pi_{2}} "7";"3"};
        {\ar_{\lhd} "6";"3"};
        {\ar^{\Delta \times 1} "7";"4"};
        {\ar^{1 \times S} "4";"5"};
        {\ar_{\rhd \times 1} "5";"6"};
        {\ar_{\Delta \times 1} "2";"5"};
        {\ar_{1 \times 1} "2";"6"};
\end{xy}
} \] Thus, since each piece of the diagram commutes, the outer
edges do as well, which gives $(vi^{ \prime})$. \qed

In Theorem \ref{squares} from Section \ref{braids}, we
demonstrated how shelves, racks, spindles, and quandles give
actions of the various braid and framed braid groups and monoids
in the category of $\Set.$  In fact, this result generalizes in an
obvious way to shelves, racks, spindles, and quandles in a
category with finite products.  The generalized proof is just like
the proof given for sets where now everything is internalized,
just like the proof of the previous proposition.

Now that we have the ability to view our Lie groups as spindles in
$\Diff_{\ast}$, we continue by developing the additional necessary
machinery we need to obtain the corresponding Lie algebras.

\section{From Lie Groups to Lie Algebras} \label{liealgofliegrp}

Recall that our goal is to describe a novel, conceptual
explanation of the passage from a Lie group to its Lie algebra
using the language of spindles.  More specifically, we will use
the internalized concepts of the previous section to think of a
Lie group as a group, and therefore, spindle in $\Diff_{\ast}$,
and then show how to use this spindle in $\Diff_{\ast}$ to obtain
the Lie algebra of the given Lie group.  Recall from the
introduction to this chapter that our desired process takes the
following diagrammatic form:
$$
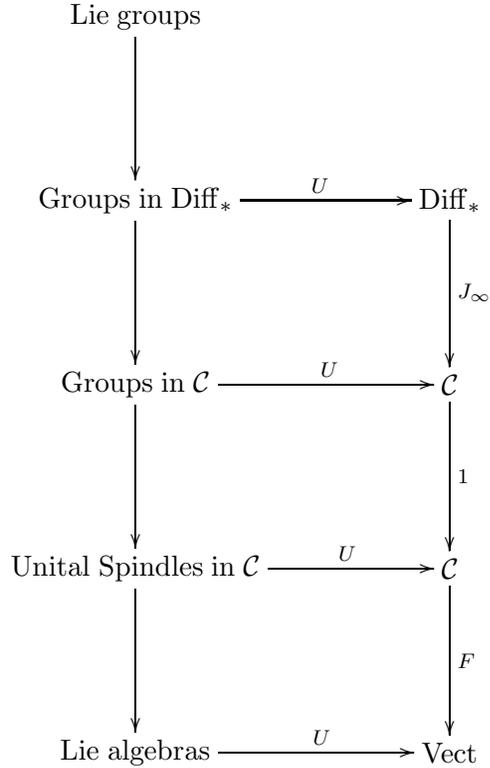
\xymatrix{
   {\rm Lie \; groups}
   \ar[dd] \\ \\
   {\rm Groups \; in \; \Diff_\ast}
   \ar[rr]^<<<<<<<<<<<{U}
   \ar[dd]
   && \Diff_{\ast}
   \ar[dd]^{J_{\infty}} \\ \\
   {\rm Groups \; in \; \mathcal{C}}
   \ar[rr]^<<<<<<<<<<<<<<<{U}
   \ar[dd]
   && \mathcal{C}
   \ar[dd]^{1} \\ \\
   {\rm Unital \; Spindles \; in \; \mathcal{C}}
   \ar[rr]^<<<<<<<<<<<{U}
   \ar[dd]
   && \mathcal{C}
   \ar[dd]^{F} \\ \\
   {\rm Lie \; algebras}
   \ar[rr]^<<<<<<<<<<<<<<{U}
   && \Vect
}$$ The first tasks consist of defining the `cojets' functor
$J_{\infty}$ and category of `special coalgebras', $\mathcal{C}$.

\subsection{Cojets} \label{cojets}
In the next two sections, we focus on this aspect of our diagram:
$$
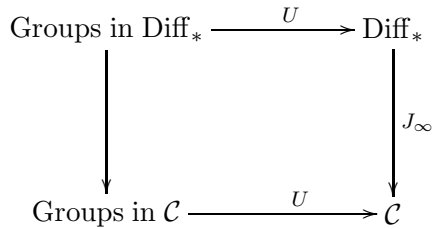
\xymatrix{
{\rm Groups \; in \; \Diff_\ast}
   \ar[rr]^<<<<<<<<<<<{U}
   \ar[dd]
   && \Diff_{\ast}
   \ar[dd]^{J_{\infty}} \\ \\
   {\rm Groups \; in \; \mathcal{C}}
   \ar[rr]^<<<<<<<<<<<<<<<{U}
   && \mathcal{C}}$$
by defining the `cojets' functor, $J_{\infty}$, and the category
$\mathcal{C}$ of `special coalgebras'.

Since we are interested in obtaining the Lie algebra of a Lie
group, and since the bracket in the Lie algebra arises from
differentiating conjugation in the Lie group {\it twice}, we need
a way to keep track of higher--order derivatives.  Roughly
speaking, a `$k$-jet' of a smooth function $f$ is a gadget which
keeps track of the value of the $f$ at a point together with all
its derivatives up to the $k$th order.  In other words, it is
simply a way of describing the Taylor expansion of $f$ at a given
point up to order $k$.

Let $M$ be a smooth manifold and let $C^\infty(M)$ be the algebra
of smooth real-valued functions on $M$. Suppose $p$ is a point of
$M$ and let $f \in C^\infty(M)$.  Recall that the first partial
derivatives of $f$ define a smooth map $df \maps TM \to \R$.

\begin{defn} \et \label{zerothorder} $f$ {\bf vanishes
to zeroth order at $p$} if $f(p) = 0$
\end{defn}

\begin{defn} \et \label{firstorder} $f$ {\bf vanishes
to first order at $p$} if it vanishes to zeroth order at $p$ and
$df \maps TM \to \R$ vanishes to zeroth order at every point of
$T_pM$.
\end{defn}

\begin{defn} \et $f$ {\bf vanishes to $k$th order at $p$}
if it vanishes to $(k-1)$st order at $p$ and $df \maps TM \to \R$
vanishes to $(k-1)$st order at every point of $T_{p}M$.
\end{defn}
Saying that $f$ vanishes to $k$th order simply means that $f$
vanishes along with all its partial derivatives up to the $k$th
order at $p$; the fancier definition given above simply shows that
this concept is independent of any choice of coordinates.

\begin{defn} \et Let $I^k(M,p)$ be the set of all smooth
functions $f \maps M \to \R$ that vanish to $k$th order at $p$.
\end{defn}
The product rule implies that $I^k(M,p)$ is an ideal in
$C^\infty(M)$.  This lets us define a finite-dimensional quotient
algebra, which we call the algebra of `$k$-jets':

\begin{defn} \et Let $J^{k}(M,p)$ denote the algebra
$C^\infty(M)/I^{k}(M,p)$.  We call an element of $J^k(M,p)$ a
\textbf{\textit{k}-jet at $p$}. Given a function $f \in
C^\infty(M)$, we denote its equivalence class in $J^k(M,p)$ by
$[f]_k$, and call this the \textbf{\textit{k}-jet of $f$ at $p$}.
\end{defn}
Thus, $[f]_{k}$ consists of all functions whose derivatives up to
order $k$ agree at $p$.  Moreover, the Taylor expansions of degree
$k$ at $p$ of functions in $[f]_{k}$ agree.

We now define a contravariant functor $J^{k} \maps \Diff_{\ast}
\to \Vect$ which sends \newline $(M, p) \in \Diff_{\ast}$ to
$J^{k}(M,p)$, the algebra of $k$-jets of real-valued functions on
$M$ at the point $p$. To define $J^k$ on morphisms, we consider a
smooth map $f \maps (M,p) \to (N,q)$, and define $J^{k}(f) \maps
J^{k}(N,q) \to J^{k}(M,p)$ by setting $J^{k}(f)[\phi]_{k} = [f
\phi]_{k}$, where $[\phi]_{k}$ is the $k$-jet at some point $q$ of
the function $\phi \maps N \to \mathbb{R}.$ This definition is
independent of the choice of representative of the equivalence
class $[\phi]_{k}$ since, by the chain rule, the first $k$
derivatives of $f \phi$ are determined by the first $k$
derivatives of $\phi$ so long as we know $f$. The contravariance
and functorality of $J^{k}$ follow from routine computations.

We have defined $J^{k}(M,p)$ in a completely
coordinate-independent way, making it seem somewhat mysterious. We
know that $J^{k}(M,p)$ contains all the information about the
derivatives of real valued functions on $M$ at $p$ up to order
$k$.  Therefore, it should not come as a surprise that the space
of $k$-jets is isomorphic to the space of Taylor polynomials of
degree $k$, since $k$th order Taylor expansions keep track of all
partial derivatives up to order $k$.  Furthermore, given a real
vector space $V$, the symmetric algebra, $SV$, is isomorphic to
the polynomial algebra in $n$ variables, so we can exhibit a
relationship between the space of $k$-jets and the sum of the
spaces of symmetric tensors of degree less than or equal to $k$.
This association will allow us to understand what the space of
$k$-jets `looks like'.

\begin{lem} \et \label{kjetstaylor}
The space of $k$-jets of $f$ at $p$, $J^{k}(M,p)$, is isomorphic
as an algebra to $\displaystyle{\bigoplus _{i=0} ^{k}
S^{i}(T^{\ast} _{p}M)}.$
\end{lem}

\noindent {\bf Proof.}  If $dim(M) = n$, then $S^{k}(T^{\ast}
_{p}M)$ is isomorphic to the algebra of polynomials in $n$
variables of degrees less than or equal to $k$. Choose coordinates
$x_{1}, x_{2}, \ldots, x_{n}$ on a neighborhood of $p \in M$.  To
simplify notation, we will use $D_{i}f$ to denote $\frac{\partial
f}{\partial x_{i}}$. Define a map $$\phi \maps J^{k}(M,p) \to
\bigoplus _{i=0} ^{k} S^{i}(T^{\ast} _{p}M)$$ which sends a
$k$-jet, $[f]_{k}$, to its Taylor polynomial of degree $k$:
$$\phi([f]_k) = \sum \frac{(D_{1}^{i_{1}} \cdots D_{n}^{i_{n}} f)(p)}{{i_{1}!
\cdots i_{n}!}} x_{1}^{i_{1}} \cdots x_{n}^{i_{n}},$$ where this
summation extends over all ordered $n$-tuples $(i_1, \ldots, i_n)$
such that each $i_j$ is a nonnegative integer and $i_1 + \cdots +
i_n \leq k$. This map depends on our choice of coordinates and
becomes well-defined once we have made our choice.  The map $\phi$
is injective since all functions within the same equivalence class
have the same derivatives, and therefore Taylor expansions, up to
order $k$ at $p$.

To show that $\phi$ is surjective, we first note that any
polynomial of degree $k$ is the Taylor expansion of degree $k$ of
some compactly supported smooth function on $\R^{n}$.  Given a
polynomial of degree $k$, we can multiply it by a smooth `bump'
function that is defined to be 1 in a neighborhood of the origin
and vanishes outside of some large ball.  By construction, the
product of our polynomial with this bump function has the same
derivatives at the origin as our original polynomial. Moreover,
given a manifold $M$ of dimension $n$ and coordinate chart $U \ni
p$, any polynomial $q(x)$ of degree $k$ is the Taylor expansion of
degree $k$ at $p$ of some smooth function on $M$.  To see this, we
can use our first observation to construct a compactly supported
smooth function $f$ on $U$ such that $q(x)$ is the $k$th degree
Taylor polynomial of $f$.  We can then extend $f$ to a smooth
function $\tilde{f}$ on all of $M$ by defining $\tilde{f}(x)= 0$
for all $x \in M \smallsetminus U.$ Thus, $q(x)$ is the Taylor
expansion of degree $k$ of $\tilde{f}$ on $M$, so that $\phi$ is
surjective.

Furthermore, these spaces are isomorphic as algebras.  The
multiplication on $J^{k}(M,p)$ is defined as $[f]_k \cdot [g]_k :=
[fg]_k$, while the multiplication on $\displaystyle{\bigoplus
_{i=0} ^{k} S^{i}(T^{\ast} _{p}M)}$ amounts to first multiplying
two polynomials, $q \cdot r$, and then forming the Taylor
expansion of degree $k$ of the product, $T_k(qr)$, at $p$.  That
is, we multiply the polynomials $q$ and $r$ and then throw out the
terms having degree higher than $k$. Thus, the isomorphism $\phi$
is a morphism of algebras since first multiplying two $k$-jets and
then forming the Taylor polynomial of degree $k$ is the same as
first forming the Taylor polynomials of the two $k$-jets and then
generating the $k$th degree Taylor polynomial of their product.
\qed

Recall that our goal is to construct a product preserving functor
$J_{\infty} \maps \Diff_{\ast} \to \mathcal{C}$, where
$\mathcal{C}$ will be the category of `special coalgebras', that
sends a pointed manifold $(M,p)$ to the coalgebra of `$k$-cojets
at $p$'.  In order to preserve products, this functor must be
covariant. Before we continue toward our goal, we recall the
definition of a coalgebra, which is the dual notion of an algebra.

\begin{defn} \et A {\bf coalgebra} is a vector space $C$ together
with a {\bf comultiplication} \newline $\Delta \maps C \to C
\otimes C$ which is bilinear and coassociative.  A {\bf coalgebra
with counit} is a coalgebra with a {\bf counit} $\epsilon \maps C
\to k$ such that
\[ \vcenter{
\xymatrix{
 k \times C \ar[r]^{\epsilon \times 1} \ar[dr]
& C \times C
& C \times k \ar[l]_{1 \times \epsilon} \ar[dl] \\
& C \ar[u]_{\Delta}}}
\] commutes.
A coalgebra is {\bf cocommutative} if $\Delta(c) = S(\Delta(c))$
for $c \in C$ where \newline $S \maps C \otimes C \to C \otimes C$
is the map that switches its inputs.
\end{defn}

\begin{defn} \et Given two coalgebras $C$ and $D$, $f \maps C \to
D$ is a {\bf coalgebra homomorphism} if $f$ is linear and the
following diagrams commute:
\[
  \xymatrix{
   C
     \ar[rr]^{\Delta_C}
     \ar[dd]_{f}
     && C \otimes C
     \ar[dd]^{f \otimes f} \\ \\
   D
     \ar[rr]^{\Delta_D}
     && D \otimes D }
\]
\[
\xymatrix{
     C
     \ar[r]^{\epsilon_C}
      \ar[dr]_{f}
      & k \\
     & D \ar[u]_{\epsilon_D}}
\]
\end{defn}

We remark that the dual of an algebra is a coalgebra so long as
the algebra is finite-dimensional.  To see this, consider an
algebra $A$ with multiplication $m \maps A \otimes A \to A$.
Then, $A^{\ast}$ becomes a coalgebra with comultiplication $\Delta
\maps A^{\ast} \to A^{\ast} \otimes A^{\ast}$ defined as the
composite:
\[
\xymatrix{
     A^{\ast}
     \ar[r]^<<<<<{m^{\ast}}
      \ar[dr]_{\Delta}
      & (A \otimes A)^{\ast} \ar[d]^{\wr} \\
     & A^{\ast} \otimes A^{\ast} }
\]
where the isomorphism from $(A \otimes A)^{\ast}$ to $A^{\ast}
\otimes A^{\ast}$ exists only when $A$ is finite-dimensional.

Therefore, seeing as how we have a finite-dimensional
\emph{algebra}, $J^{k}(M,p)$, and \emph{contravariant} functor,
$J^{k}$, it is natural to consider the duals of these as a first
step toward achieving our goal. Therefore, we consider
$$J_k(M,p) = J^{k}(M,p)^{\ast},$$ which is the space of
\textbf{\emph{k}-cojets at \emph{p}}. Since $J_{k}(M,p)$ is the
dual of a finite-dimensional algebra, it is a coalgebra.
Furthermore, we define the functor $J_{k} \maps \Diff_{\ast} \to
\Vect$ which sends any object $(M,p) \in \Diff_{\ast}$ to
$J_{k}(M,p)$ and any smooth map $f \maps (M,p) \to (N,q)$ to
$J_{k}(f) \maps J_{k}(M,p) \to J_{k}(N,q)$, which is the adjoint
of the linear map $J^{k}(f) \maps J^{k}(N,q) \to J^{k}(M,p)$.  The
functorality and covariance of $J_{k}$ hold since it is the
adjoint of a contravariant functor.

Lemma \ref{kjetstaylor} established an algebra isomorphism from
the space of $k$-jets to the space $\bigoplus _{i=0} ^{k}
S^{i}(T^{\ast} _{p}M)$.  The dual statement gives us an
isomorphism of coalgebras:

\begin{lem} \et \label{cojetsisotoSV}
The space of $k$-cojets of $f$ at $p$, $J_{k}(M,p)$ is isomorphic
as a coalgebra to $\displaystyle{\bigoplus _{i=0} ^{k}
S^{i}(T_{p}M)}.$
\end{lem}

\noindent {\bf Proof.} \begin{eqnarray*} J_{k}(M,p) =
J^{k}(M,p)^{\ast} & \cong & \left( \bigoplus _{i=0}
^{k} S^{i}(T^{\ast} _{p}M) \right)^{\ast} \\
& \cong & \bigoplus _{i=0} ^{k} (S^{i}(T^{\ast}_{p}M))^{\ast} \\
& \cong & \bigoplus _{i=0} ^{k} S^{i}(T^{\ast \ast}_{p}M) \\
& \cong & \bigoplus _{i=0} ^{k} S^{i}(T_{p}M)
\end{eqnarray*}
\hskip 30 em \qed

Thus far, we have achieved precisely what we desired:  a covariant
functor $J_{k}$ that sends any pointed smooth manifold $(M,p)$ to
the coalgebra of $k$-cojets of real-valued functions on $M$ at the
point $p$.  Unfortunately, however, this functor does not preserve
products since
\begin{eqnarray*}
J_{k}(M \times N, (p,q))
& \cong & \bigoplus_{i=0}^{k} S^{i}(T_{(p,q)}(M \times N)) \\
& \cong & \bigoplus_{i=0}^{k} S^{i}(T_{p}M \oplus T_{q}(N)) \\
& \cong & \bigoplus_{i=0}^{k} \bigoplus_{j = 0}^i S^{j}(T_{p}M)
\otimes
S^{i-j}(T_{q}N) \\
& \ncong & \bigoplus_{i=0}^{k} S^{i}(T_{p}M)
\otimes \bigoplus_{j=0}^{k} S^{j}(T_{q}N) \\
& \cong & J_{k}(M,p) \otimes J_{k}(N,q)
\end{eqnarray*}

Failing to preserve products implies that $J_{k}$ will not
preserve groups, and hence spindles, which is crucial. In order to
overcome this shortcoming, we consider a space that more closely
resembles the symmetric algebra, since it preserves products.
Thus, we turn to the task of constructing a product-preserving
functor that will send a pointed smooth manifold $(M,p)$ to the
coalgebra of `cojets': the union of the coalgebras of $k$-cojets.

For each $k$, we have a surjective map from $J^{k+1}(M,p)$ to
$J^{k}(M,p)$ which simply forgets the $(k+1)$st partial
derivatives.  More precisely, this map arises as a result of the
fact that we mod out by more when we form the quotient algebra
$J^{k}(M,p)$ than when we form the algebra $J^{k+1}(M,p)$.  Thus
we obtain a diagram:
\begin{eqnarray*} \label{jetdiagram}
\ldots \to J^{3}(M,p) \to J^{2}(M,p) \to J^{1}(M,p) \to
J^{0}(M,p)= \R
\end{eqnarray*}
where each arrow is surjective.

Taking the dual, or adjoint, of this diagram produces a
corresponding inclusion of coalgebras of $k$-cojets.  Forming the
dual amounts to turning the arrows around and taking the duals of
the algebras of $k$-jets. We obtain:
\begin{eqnarray*}
\R = J_{0}(M,p) \to J_{1}(M,p) \to J_{2}(M,p) \to J_{3}(M,p) \to
\ldots
\end{eqnarray*}
where now each arrow is injective.  Thus, we can think of
$J_{k}(M,p)$ as sitting inside of $J_{k+1}(M,p)$ and therefore we
may form the category-theoretic colimit, or union, of these spaces
$$J_{\infty}(M,p) := \bigcup_{k \geq 0} J_{k}(M,p),$$ which we call the
\textbf{space of cojets at p}.

We remark that we could have formed the corresponding space of
jets, $J^{\infty}(M,p)$, by taking the limit of the diagram of
$k$-jets above.  However, unlike in the case of $k$-jets and
$k$-cojets, the dual of the space of jets is not isomorphic to the
space of cojets.  That is, $J_{\infty}(M,p) \neq
J^{\infty}(M,p)^{\ast}$.

The space $J_{\infty}(M,p)$ of cojets clearly is a vector space,
as it is the union of $J_{k}(M,p)$ which are vector spaces. In
fact, the categorically inclined reader will notice that
$J_{\infty}(M,p)$ is just the colimit of the coalgebras
$J_{k}(M,p)$!  That is, $J_{\infty}(M,p)$ has a much richer
structure:  that of a cocommutative coalgebra.  In the next
section, we describe the coalgebra structures of $J_{\infty}(M,p)$
and the symmetric algebra $S(T_{p}M)$ and show that they are
compatible.

\subsection{Special Coalgebras} \label{specialcoalgs}
In this section, we describe the coalgebra structures of the space
of cojets and the symmetric algebra, and then provide a definition
of the category $\mathcal{C}$ of `special coalgebras', which is
the target category of our functor $J_{\infty}$.

Since we have already recalled the algebraic definition of a
coalgebra, we now offer a diagrammatic description of both
algebras and coalgebras. Since multiplication is a process that
takes two things and fuses them together into one new thing of the
same sort, we can represent this process diagrammatically as:
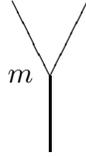
\[
\xy
(-5,15)*{}="1";
(5,15)*{}="2";
(0,5)*{}="3";
(-4,5)*{m}="a";
(0,-5)*{}="4";
"1";"3" **\dir{-};
"2";"3" **\dir{-};
"3";"4" **\dir{-};
\endxy
\]
Then, since comultiplication is the dual process, we `cothink' and turn everything upside down!
Thus, comultiplication takes the form:
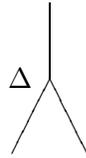
\[
\xy (0,15)*{}="1";
(0,5)*{}="2";
(-4,5)*{\Delta}="a";
(-5,-5)*{}="3";
(5,-5)*{}="4";
"1";"2" **\dir{-};
"2";"3" **\dir{-};
"2";"4" **\dir{-};
\endxy
\]
In terms of these diagrams, associativity and coassociativity are expressed as:
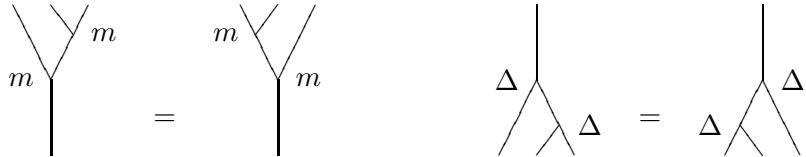
\[
\xy
(-5,15)*{}="1";
(5,15)*{}="2";
(0,15)*{}="5";
(3,11)*{}="6";
(7,11)*{m}="b";
(0,5)*{}="3";
(-4,5)*{m}="a";
(0,-5)*{}="4";
"1";"3" **\dir{-};
"2";"3" **\dir{-};
"3";"4" **\dir{-};
"6";"5" **\dir{-};
\endxy
\quad = \quad
\xy
(-5,15)*{}="1";
(5,15)*{}="2";
(0,5)*{}="3";
(-3,11)*{}="6";
(-7,11)*{m}="b";
(0,15)*{}="5";
(4,5)*{m}="a";
(0,-5)*{}="4";
"1";"3" **\dir{-};
"2";"3" **\dir{-};
"3";"4" **\dir{-};
"6";"5" **\dir{-};
\endxy
\qquad \qquad \qquad
\xy (0,15)*{}="1";
(0,5)*{}="2";
(-4,5)*{\Delta}="a";
(3,-1)*{}="5";
(7,-1)*{\Delta}="b";
(0,-5)*{}="6";
(-5,-5)*{}="3";
(5,-5)*{}="4";
"1";"2" **\dir{-};
"2";"3" **\dir{-};
"2";"4" **\dir{-};
"5";"6" **\dir{-};
\endxy
\quad = \quad
\xy
(0,15)*{}="1";
(0,5)*{}="2";
(4,5)*{\Delta}="a";
(-3,-1)*{}="5";
(-7,-1)*{\Delta}="b";
(0,-5)*{}="6";
(-5,-5)*{}="3";
(5,-5)*{}="4";
"1";"2" **\dir{-};
"2";"3" **\dir{-};
"2";"4" **\dir{-};
"5";"6" **\dir{-};
\endxy
\]
while commutativity and cocommutativity can be depicted as:
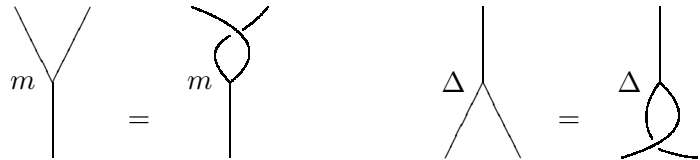
\[
\xy (0,-5)*{}="1";
(0,5)*{}="2";
(-4,5)*{m}="a";
(-5,15)*{}="3";
(5,15)*{}="4";
"1";"2" **\dir{-};
"2";"3" **\dir{-};
"2";"4" **\dir{-};
\endxy
\quad = \quad \xy (0,-5)*{}="1"; (0,5)*{}="2"; (-4,5)*{m}="a";
(-5,15)*{}="3"; (0,11.2)*{}="2a"; (1.7,12)*{}="4a"; (5,15)*{}="4";
"1";"2" **\dir{-}; "2";"3" **\crv{(7, 10.5)}; "2";"2a"
**\crv{(-4,9)}; "4a";"4" **\crv{(3.8,13.4)}; \endxy \qquad \qquad
\qquad \xy (0,15)*{}="1"; (0,5)*{}="2"; (-4,5)*{\Delta}="a";
(-5,-5)*{}="3"; (5,-5)*{}="4"; "1";"2" **\dir{-}; "2";"3"
**\dir{-}; "2";"4" **\dir{-};
\endxy
\quad = \quad \xy (0,15)*{}="1"; (0,5)*{}="2";
(-4,5)*{\Delta}="a"; (-5,-5)*{}="3"; (-1,-2.5)*{}="2a";
(-.1,-3.9)*{}="4a"; (5,-5)*{}="4"; "1";"2" **\dir{-}; "2";"3"
**\crv{(7, -1.5)}; "2";"2a" **\crv{(-3,1)}; "4a";"4"
**\crv{(2.5,-4.9)};
\endxy \]

With these diagrams, it becomes a simple exercise to show that the
tensor product two commutative algebras is again a commutative
algebra, and the tensor product of two cocommutative coalgebras is
again a cocommutative coalgebra.

We have shown that the space of cojets is isomorphic as a vector
space to the symmetric algebra, but, as already noted, both
$J_{\infty}(M,p)$ and $S(T_{p}M)$ have much richer structures!
They are both cocommutative coalgebras, with a compatible
comultiplication. The symmetric algebra becomes a coalgebra with
comultiplication $\Delta \maps S(T_p M) \to S(T_p M) \otimes S(T_p
M)$ defined as $\Delta(v) = v \otimes 1 + 1 \otimes v$ for $v \in
T_p M$. Notice that this determines $\Delta$ on all of $S(T_p M)$
since $\Delta$ is a morphism of algebras and elements $v \in T_p
M$ generate all of $S(T_p M)$ as an algebra.  Thus, knowing what
$\Delta$ does to elements in $T_p M$ allows us to determine
$\Delta$ on all of $S(T_p M)$. For example, we must have
$\Delta(1) = 1 \otimes 1$ since $\Delta$ is a algebra homomorphism
and $1 \otimes 1$ is the multiplicative identity in the algebra
$S(T_p M) \otimes S(T_p M)$. We must also have $\Delta(vw) =
\Delta(v) \Delta(w)$ for $v,w \in T_p M$.  Similarly, the counit
$\epsilon \maps S(T_pM) \to k$ is defined by $\epsilon(v) = 0$ for
$v \in T_p M$, and computations show that this determines
$\epsilon$ on all of $S(T_p M)$.  Moreover, $S(T_p M)$ is a
cocommutative coalgebra, since
$$\Delta(v) = v \otimes 1 + 1 \otimes v = 1 \otimes v + v \otimes
1 = S(v \otimes 1 + 1 \otimes v) = S(\Delta(v)),$$ where $S \maps
T_p M \otimes T_p M \to T_p M \otimes T_p M$ is the map that
switches its inputs. Of course, none of this is particular to the
tangent space, $T_p M$, and holds for an arbitrary vector space
$V$.

Furthermore, $J_{\infty}(M,p)$ becomes a coalgebra, as it is the
union over $k$ of the coalgebras $J_{k}(M,p)$, with
comultiplication and counit defined in terms of these operations
on the appropriate $J_{k}(M,p)$.  Given $[f] \in J_{\infty}(M,p) =
\bigcup _{k} J_{k}(M,p)$, $[f] \in J_{k}(M,p)$ for some $k$, so we
define $\Delta [f] = \Delta_{k}[f]$ and $\epsilon [f] = \epsilon
_{k} [f]$.  Recall that $\Delta_k$ and $\epsilon_k$ are defined as
the adjoints of the multiplication and unit maps on $J^{k}(M,p)$
since $J_{k}(M,p) = J^{k}(M,p)^{\ast}$.  Like $S(T_p M)$,
$J_{\infty}(M,p)$ is a cocommutative coalgebra, which can be seen
from a computation similar to the one given above for $S(T_p M)$.

Thus, any of the coordinate-dependent isomorphisms in Lemma
\ref{cojetsisotoSV} give rise to coalgebra isomorphisms:

\begin{lem} \et The space of cojets $J_{\infty}(M,p)$ is isomorphic as a
cocommutative coalgebra to the symmetric algebra $S(T_{p}M)$.
\end{lem}

\noindent {\bf Proof.}
\begin{center}$\displaystyle{J_{\infty}(M,p) = \bigcup _{k} J_{k}(M,p)
\cong \bigcup_{k} \bigoplus _{i=0} ^{k} S^{i}(T_{p}M) \cong
\bigoplus _{i=0} ^{\infty} S^{i}(T_{p}M) = S(T_{p}M)}$ \qed
\end{center}

In fact, we remark that both $SV$ and $J_{\infty}(M,p)$ are {\it
bialgebras}, meaning that they are algebras and coalgebras in a
compatible way. That is, the comultiplication and counit are {\it
algebra} homomorphisms.

Now that we know the coalgebra of cojets is isomorphic to the
symmetric algebra, and since the symmetric algebra functor is
product-preserving, we can use this to construct the target
category, $\mathcal{C}$, of our functor $J_{\infty}$.

\begin{defn} \et \label{defnspecialcoalgs} We define {\bf $\Co$},
the category of \textbf{special coalgebras}, to be the category
whose objects are triples $(C, V, \alpha)$ where $C$ is a
cocommutative coalgebra, $V$ is a vector space, and $\alpha \maps
C \to SV$ is a cocommutative coalgebra isomorphism. The morphisms
of $\mathcal{C}$ are coalgebra homomorphisms $f \maps C \to C'$.
We will often denote $(C, V, \alpha)$ by $C \cong SV$.
\end{defn}

Recall that since coalgebra homomorphisms send primitive elements
to primitive elements, $f \maps SV \to SV'$ automatically maps $V$
to $V'$. Notice that this category has finite products, since it
has binary products and a terminal object.  That is, given $(C, V,
\alpha)$ and $(C', V', \alpha ')$, their product is $(C \otimes
C', V \oplus V', (\alpha \otimes \alpha ') \phi)$ where we have
used the fact that the tensor product of cocommutative coalgebras
is a cocommutative coalgebra and $\phi \maps SV \otimes SV' \to
S(V \oplus V')$ is the standard isomorphism.  We remark that the
cocommutativity of our special coalgebras is necessary for
$\mathcal{C}$ to have products.  The terminal object in
$\mathcal{C}$ is $(k, \{0\}, 1)$.

We are now in a position to describe the functor $J_{\infty} \maps
\Diff_{\ast} \to \mathcal{C}$, since $J_{\infty}(M,p)$ is a
coalgebra that is isomorphic to $S(T_{p}M)$.  This functor sends
any object $(M,p) \in \Diff_{\ast}$ to $(J_{\infty}(M,p),
T_{p}(M), \psi)$ and any smooth map $f \maps (M,p) \to (N,q)$ to
\newline $J_{\infty}(f) \maps J_{\infty}(M,p) \to J_{\infty}(N,q)$.  The
functorality and covariance of $J_{\infty}$ follow from routine
computations.  Moreover, $J_{\infty}$ preserves groups, and hence
quandles in $\mathcal{C}$, since it preserves products:

\begin{prop} \et \label{finiteproducts} The functor $J_{\infty} \maps \Diff_{\ast} \to
\mathcal{C}$ preserves finite products.
\end{prop}

\noindent{\bf Proof. } Given $(M,p)$ and $(N,q)$ in
$\Diff_{\ast}$, we must show that $J_{\infty}(M \times N, (p,q))$
is the product of $J_{\infty}(M,p)$ and $J_{\infty}(N,q)$ in
$\mathcal{C}$. We have
\begin{eqnarray*}
J_{\infty}(M,p) &:=& (J_{\infty}(M,p), T_{p}M, \alpha) \\
J_{\infty}(N,q) &:=& (J_{\infty}(N,q), T_{q}N, \alpha ') \\
J_{\infty}(M \times N,(p,q)) &:=& (J_{\infty}(M \times N,(p,q)),
T_{(p,q)}(M \times N), \beta)  \\
J_{\infty}(M,p) \times J_{\infty}(N,q) &:=& (J_{\infty}(M,p)
\otimes J_{\infty}(N,q), T_{p}M \oplus T_{q}N, (\alpha \otimes
\alpha ')\phi)
\end{eqnarray*}
We need only to demonstrate that $J_{\infty}(M,p) \otimes
J_{\infty}(N,q) \cong J_{\infty}(M \times N,(p,q))$. By lemma
\ref{cojetsisotoSV}, $J_{\infty}(M,p) \cong S(T_{p}M)$, so that we have
\begin{eqnarray*}
J_{\infty}(M,p) \otimes J_{\infty}(N,q) & \cong & S(T_{p}M)
\otimes S(T_{q}N) \\
& \cong & S(T_{p}M \oplus T_{q}N) \\
& = & S(T_{(p,q)}(M \times N)) \\
& \cong & J_{\infty}(M \times N, (p,q)).
\end{eqnarray*}
It follows immediately that $\beta \cong (\alpha \otimes
\alpha)\phi$ since the following diagram commutes:
\[ \vcenter{
\xymatrix{
 J_{\infty}(M,p) \otimes J_{\infty}(N,q) \ar[rr]^{\alpha \otimes \alpha '}
 \ar[dd]_{\wr}
&& S(T_{p}M) \otimes S(T_{q}N) \ar[dd]^{\phi} \\ \\
 J_{\infty}(M \times N, (p,q)) \ar[rr]^<<<<<<<<<<<<{\beta}
&& S(T_{(p,q)}(M \times N)) = S(T_{p}M \oplus T_{q}N)}}
\]
where the isomorphism one the left side of this square is the one
described above. \qed Since $J_{\infty}$ preserves products, we
obtain the following corollary:

\begin{cor} \et The functor $J_{\infty} \maps \Diff_{\ast} \to
\mathcal{C}$ sends groups in $\mathcal{C}$  to groups in
$\mathcal{C}$.
\end{cor}
Thus, we have completely described the first piece of our process:
$$\xymatrix{
   {\rm Lie \; groups}
   \ar[dd] \\ \\
   {\rm Groups \; in \; \Diff_\ast}
   \ar[rr]^<<<<<<<<<<<{U}
   \ar[dd]
   && \Diff_{\ast}
   \ar[dd]^{J_{\infty}} \\ \\
   {\rm Groups \; in \; \mathcal{C}}
   \ar[rr]^<<<<<<<<<<<<<<<{U}
   && \mathcal{C}}$$

Note that a group in the category of coalgebras is a Hopf algebra,
and we are getting $U\mathfrak{g} = J_{\infty}(G,1)$. Now, given a
Lie group, $(G,1)$, it remains to show how to peel off
$\mathfrak{g}$ from $J_{\infty}(G,1)$.

\subsection{Unital Spindles} \label{final}
We begin this section by reminding the reader of the connection
between Lie algebras and spindles.  Both Lie algebras and spindles
are gadgets that give rise to Yang--Baxter operators.  In Section
\ref{topology2} we saw that the Yang--Baxter equation is none
other than the distributive law for the spindle operation, and
that it is the Jacobi identity in disguise when we define our
Yang--Baxter operator on the space $k \oplus L$, where $L$ is a
Lie algebra over $k$. Moreover, and most relevant to our task, we
noticed how conjugation plays a prominent role in both the theory
of spindles and that of Lie algebras.  The operation of
conjugation in groups satisfies the spindle axioms; indeed, groups
were our primordial examples of spindles.  Conjugation also
appears at the heart of the theory of Lie algebras since the
bracket in a Lie algebra arises from differentiating conjugation
twice. Furthermore, the Jacobi identity in a Lie algebra arises
from differentiating the self-distributivity axiom in a spindle
when the spindle operation is conjugation.  That is, given curves
$e^{sx}, e^{ty}$ and $e^{uz}$ in a matrix Lie group $G$ where $x,
y, z \in \mathfrak{g}$, we have:
$$[x,[y,z]] = \frac{d}{dsdtdu} (e^{sx} e^{ty} e^{uz} e^{-ty} e^{-sx}) |_{s=t=u=0}$$
while
$$[y,[x,z]] + [[x,y],z] = \frac{d}{dsdtdu} (e^{sx} e^{ty} e^{-sx} e^{sx} e^{uz} e^{-sx}
e^{sx} e^{-ty} e^{-sx}) |_{s=t=u=0}.$$ These observations taken
together inspired the novel passage from Lie groups to Lie
algebras that we continue describing.  In this section, we
describe the following aspect of our diagram:
$$\xymatrix{
   {\rm Groups \; in \; \mathcal{C}}
   \ar[rr]^<<<<<<<<<<<<<<<{U}
   \ar[dd]
   && \mathcal{C}
   \ar[dd]^{1} \\ \\
   {\rm Unital \; Spindles \; in \; \mathcal{C}}
   \ar[rr]^<<<<<<<<<<<{U}
   && \mathcal{C}}$$
Given a Lie group $G$, we will define an operation of conjugation
on the coalgebra $J_{\infty}(G,1)$ to make it into a spindle. More
precisely, $J_{\infty}(G,1)$ will be a spindle in $\mathcal{C}.$
Recall that in Section \ref{internalization} we demonstrated how a
group in $K$, where $K$ is a category with finite products, gives
a quandle in $K$.  This amounted to an internalization of the fact
that usual groups give quandles. Actually we do not need the whole
structure of a quandle to obtain a Lie algebra: it suffices to use
left conjugation to define the bracket, so we can work with a
spindle.  But not any spindle object in $\mathcal{C}$ gives a Lie
algebra. This internalized spindle must posses two special
properties, namely the internalized versions of:
$$x1x^{-1} = 1 \quad {\rm and} \quad 1x1^{-1} = x,$$ where
$x \in G$ for some group $G$ and $1 \in G$ is the identity.  We
call a spindle with these two properties a `unital spindle'.

\begin{defn} \et \label{wellpoint} Let $K$ be a category with finite products.
A {\bf unital spindle in $K$} is a spindle $Q$ in $K$ equipped
with a special point $\id \maps I \to Q$, where $I$ is the
terminal object in $K$, such that the following diagrams commute:
\begin{itemize}
\item[(a)]
\[
\xymatrix{
 Q
 \ar[rr]
 \ar@/^2pc/[rrrrrr]_{\quad}^{1}="1"
&& I \times Q
  \ar[rr]^{\id \times 1}
&& Q \times Q
   \ar[rr]^{\rhd}
&& Q }
\]
\item[(b)]
\[
\xy (-23,0)*{}; (23,0)*{}; (-30,0)*+{Q}="L"; (0,0)*+{Q \times
I}="M"; (30,0)*+{Q \times Q}="R"; (15,10)*+{I}="LT";
(60,0)*+{Q}="RR";
 (-11,7)="L1";
 (-11,-7)="L2";
 (11,7)="R1";
 (11,-7)="R2";
 {\ar^{}  "L";"M"};
 {\ar^{1 \times \id}  "M";"R"};
 {\ar^{\rhd}     "R";"RR"};
 "L";"LT" **\crv{(-20,7.5)};
             ?(.98)*\dir{>};
 "LT";"RR" **\crv{(35,7.5)};
             ?(.96)*\dir{>};
     (35,9)*{\id};
\endxy
\]
\end{itemize}
\end{defn}

Since the two conditions given in the definition of a unital
spindle are just the internalizations of the conjugation
properties involving the identity element of a group, it should
not come as a surprise that we can obtain a unital spindle in $K$
from a group in $K$.

\begin{prop} \et \label{ptdquandobj} Let $K$ be a category with
finite products.  A group $Q$ in $K$ gives a unital spindle in
$K$.
\end{prop}
\noindent{\bf Proof. } In Section \ref{internalization}, we began
to illustrate how a group in $K$ gives a quandle in $K$, and hence
spindle in $K$, by defining the quandle operations in terms of
conjugation. Since we are only concerned with spindles, we recall
that we defined left conjugation $\rhd \maps Q \times Q \to Q$ as
the composite:
\[
 \xy
     (-70,0)*+{(x,y)}="1";
     (-50,0)*+{(x,x,y)}="2";
     (-30,0)*+{(x,y,x)}="3a";
     (-0,0)*+{(x,y,x^{-1})}="3";
     (25,0)*+{(x,yx^{-1})}="4";
     (50,0)*+{xyx^{-1}}="5";
            {\ar^{\Delta \times 1} "1";"2"};
            {\ar^{1 \times S}      "2";"3a"};
            {\ar^<<<<<<<{1 \times 1 \times \inv} "3a";"3"};
            {\ar^{1 \times m} "3";"4"};
            {\ar^{m} "4";"5"};
 \endxy
\]
We then demonstrated how to obtain the internalized left
idempotence law. Showing that the internalized left distributive
law can be obtained by replacing $\rhd$ in diagram $(i^{ \prime})$
of Definition \ref{lshelfobj} amounts to a rather large, though
mostly trivial, diagram.

Obtaining the two additional properties of a unital quandle
require showing that the identity in our group object, $\id \maps
I \to Q$, makes the diagrams $(a)$ and $(b)$ above commute, where
$\rhd$ is replaced by the composite given above.

Property $(a)$ is the result of the diagram:
$$ \vcenter{
   \xy
   (0,35)*+{Q}="1";
   (-40,20)*+{I \times Q}="2";
   (40,20)*+{Q}="3";
   (-40,0)*+{Q \times Q}="4";
   (40,0)*+{Q \times Q}="5";
   (-40,-20)*+{Q \times Q \times Q}="6";
   (40,-20)*+{Q \times Q \times Q}="7";
   (0,-40)*+{Q \times Q \times Q}="8";
   (10, 10)*+{I \times Q}="9";
   (15, 20)*+{A}="12";
   (-20,0)*+{B}="13";
   (10,-10)*+{Q \times Q \times I}="10";
   (15,-25)*+{F}="14";
   (30,-10)*+{E}="15";
   (20,0)*+{D}="16";
   (33,12)*+{C}="17";
        {\ar_{}                   "1";"2"};
        {\ar^{1}                               "1";"3"};
        {\ar_{\id \times 1}   "2";"4"};
        {\ar_{m}                   "5";"3"};
        {\ar_{\Delta \times 1}                   "4";"6"};
        {\ar_{1 \times m}                   "7";"5"};
        {\ar_{1 \times S}                   "6";"8"};
        {\ar_{1 \times 1 \times \inv}                                "8";"7"};
        {\ar_{} "9";"3"};
        {\ar^{1 \times \id} "9";"5"};
        {\ar_{} "10";"5"};
        {\ar_{\id \times 1 \times 1} "10";"7"};
        {\ar_{} "1";"9"};
        {\ar_{\id \times 1} "9";"10"};
        {\ar_{1 \times 1 \times \id} "10";"8"};
\endxy }
\\ \\
$$
where region $C$ is the left unit law, region $E$ is the right
unit law, region $F$ is the internalization of the fact that the
inverse of the identity of a group is the identity.  Regions $A$,
$B$ and $D$ clearly commute.

Property $(b)$ is given by the diagram:
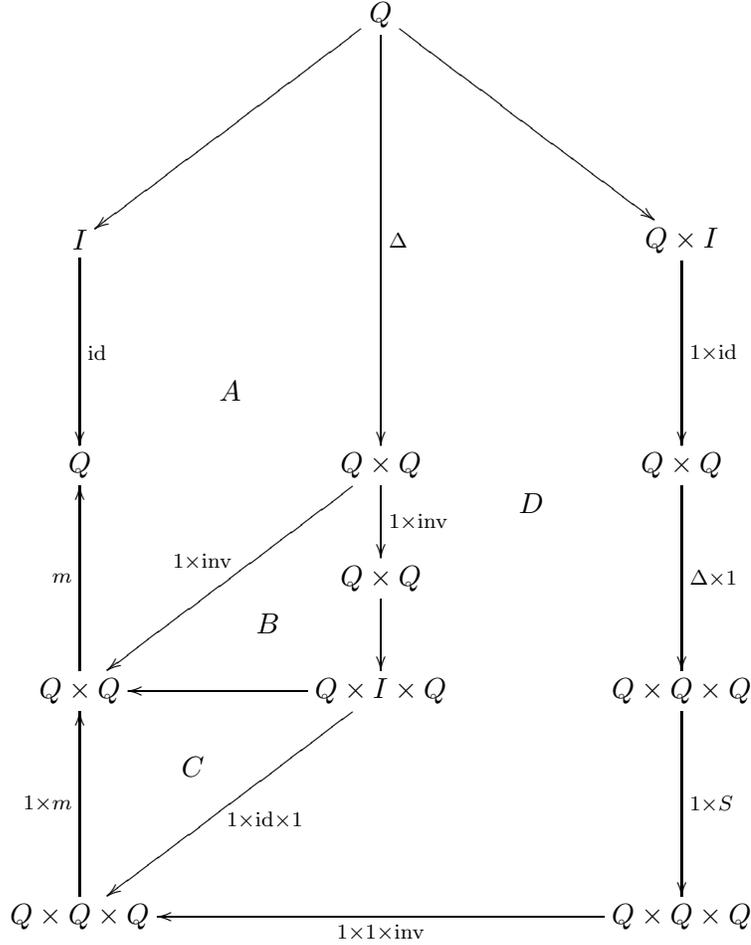
\[ \vcenter{
\begin{xy}
    (0, 60)*+{Q}="1";
    (-40,30)*+{I}="2";
    (40, 30)*+{Q \times I}="3";
    (-40,0)*+{Q}="4";
    (40,0)*+{Q \times Q}="5";
    (-40,-30)*+{Q \times Q}="6";
    (40, -30)*+{Q \times Q \times Q}="7";
   (-40,-60)*+{Q \times Q \times Q}="8";
    (40, -60)*+{Q \times Q \times Q}="9";
    (0,0)*+{Q \times Q}="10";
    (0,-15)*+{Q \times Q}="11";
    (0,-30)*+{Q \times I \times Q}="12";
    (-20,10)*{A}="A";
    (20,-5)*{D}="D";
    (-25,-40)*{C}="C";
    (-15,-21)*{B}="B";
        {\ar^{} "1";"2"};
        {\ar_{} "1";"3"};
        {\ar^{\id} "2";"4"};
        {\ar^{1 \times \id} "3";"5"};
        {\ar^{m} "6";"4"};
        {\ar^{\Delta \times 1} "5";"7"};
        {\ar^{1 \times m} "8";"6"};
        {\ar^{1 \times S} "7";"9"};
        {\ar^{1 \times 1 \times \inv} "9";"8"};
        {\ar^{\Delta} "1";"10"};
        {\ar_{1 \times \inv} "10";"6"};
        {\ar^{1 \times \inv} "10";"11"};
        {\ar^{} "11";"12"};
        {\ar^{} "12";"6"};
        {\ar^{1 \times \id \times 1} "12";"8"};
\end{xy}
} \] where region $A$ is the right inverse law and region $C$ is
the left unit law.  Regions $B$ and $D$ clearly commute. Thus we
have that a group in $K$ gives a unital spindle in $K$. \qed

We will now use the result of Proposition \ref{ptdquandobj} to
turn our groups in $\mathcal{C}$ into unital spindles in
$\mathcal{C}$.  We remark that while we have chosen to use the
less familiar notion of spindle for this process, it is also true
that a unital quandle in $\mathcal{C}$ will give a Lie algebra.
Moreover, Lie groups actually give unital quandles in
$\mathcal{C}$, where a unital quandle is merely a unital spindle
which is a quandle.

It remains, then, to describe the final aspect of our process: how
to obtain a Lie algebra from a unital spindle in $\mathcal{C}$.
Thus, we need only to define the functor $F \maps \mathcal{C} \to
\Vect$ and show that when we start with a unital spindle $Q \cong
SV$ for some vector space $V$, that we can define a bracket
operation on $V$ making it into a Lie algebra.

\subsection{From Unital Spindles to Lie Algebras} \label{punchline}
In this section, we finish describing the process of obtaining a
Lie algebra from a Lie group using the language of spindles.  That
is, we explain this final aspect of our diagram:
$$\xymatrix{
  {\rm Unital \; Spindles \; in \; \mathcal{C}}
   \ar[rr]^<<<<<<<<<<<{U}
   \ar[dd]
   && \mathcal{C}
   \ar[dd]^{F} \\ \\
   {\rm Lie \; algebras}
   \ar[rr]^<<<<<<<<<<<<<<{U}
   && \Vect
}$$  Given a special coalgebra $C \cong SV$ for some vector space
$V$, we define $F \maps \mathcal{C} \to \Vect$ to send any object
$(C, V, \alpha)$ to $V$ and any morphism $f \maps C \to C'$ to
$\alpha ^{-1} \circ f \circ \alpha ' |_{V} \maps V \to V'$.

Now, in order to show that a unital spindle $Q$ in $\mathcal{C}$
gives a Lie algebra, it is beneficial to know how to interpret
conditions $(a)$ and $(b)$ in Definition \ref{wellpoint} for an
element $x \in Q$.  To do so, we need the following lemmas:

\begin{lem} \et \label{terminal} If $C$ is a special coalgebra, the
counit $\epsilon_C \maps C \to k$ is the unique homomorphism to
the terminal object $I$ in $\mathcal{C}$.
\end{lem}

\noindent{\bf Proof.}  Recall that the terminal object, $I$, in
$\mathcal{C}$ is given by $(k, \{ 0 \}, 1)$.  Thus, by definition
of a morphism in $\mathcal{C}$, a morphism from the special
coalgebra $(C, V, \alpha)$ to the terminal object $(k, \{ 0 \},
1)$ is a coalgebra homomorphism $f \maps C \to k$.  Therefore, we
must show $\epsilon_C \maps C \to k$ is a coalgebra homomorphism,
meaning we must show that these two diagrams:
\[
  \xymatrix{
   C
     \ar[rr]^{\Delta_C}
     \ar[dd]_{\epsilon _C}
     && C \otimes C
     \ar[dd]^{\epsilon_C \otimes \epsilon_C} \\ \\
   k
     \ar[rr]^{\Delta_k}
     && k \otimes k }
\]
\[
\xymatrix{
     C
     \ar[r]^{\epsilon_C}
      \ar[dr]_{\epsilon_C}
      & k \\
     & k \ar[u]_{\epsilon_k}}
\]
commute, which is clear. \qed

\begin{lem}
\et \label{diagonal} If $C$ is a special coalgebra, the
comultiplication $\Delta_C \maps C \to C \otimes C$ is the
diagonal map in $\mathcal{C}$.
\end{lem}

\noindent{\bf Proof.}  Since $C \otimes C$ is the product in
$\mathcal{C}$, there is a unique morphism $\phi \maps C \to C
\otimes C$ which makes the following diagram commute:
\[
  \vcenter{
  \xymatrix@!C{
   &   C
    \ar[ddl]_{1_C}
    \ar[dd]^{\phi}
    \ar[ddr]^{1_C}
    \\ \\
  C &
   C \otimes C
   \ar[l]^{\pi_{1}}
   \ar[r]_{\pi_{2}} &
   C}}
\]
where the $\pi_1$ and $\pi_2$ are projection maps.  This map
$\phi$ is called the diagonal.   Since $C$ is a coalgebra, we have
the comultiplication map $\Delta_C \maps C \to C \otimes C$.  This
makes the same sort of diagram commute:
\[
  \vcenter{
  \xymatrix@!C{
   &   C
    \ar[ddl]_{1_C}
    \ar[dd]^{\Delta_C}
    \ar[ddr]^{1_C}
    \\ \\
  C &
   C \otimes C
   \ar[l]^{\pi_{1}}
   \ar[r]_{\pi_{2}} &
   C}}
\]
since the projection maps are the following composites:
$$\pi_1 : = \xymatrix{
    A \otimes B \ar[rr]^{1_A \otimes \epsilon_B}
    &&
    A \otimes k \ar[rr]^{\sim}
    &&
    A}$$
$$\pi_2 : = \xymatrix{
    A \otimes B \ar[rr]^{\epsilon_A \otimes 1_B}
    &&
    k \otimes B \ar[rr]^{\sim}
    &&
    B}$$
where $\epsilon_A$ and $\epsilon_B$ are the counit maps for
coalgebras $A$ and $B$.  Since the diagonal map $\phi$ is unique,
it must equal $\Delta$.  \qed

\begin{lem} \et \label{multid} If $(C, V, \alpha)$ is a special
coalgebra with the structure of a unital spindle in $\mathcal{C}$,
then the map $\id \maps k \to SV$, obtained by composing $\id
\maps k \to C$ with the isomorphism $\alpha \maps C \to SV$, maps
the multiplicative identity $1$ in $k$ to the multiplicative
identity $1$ in $SV$. That is, $\id(1) = 1$.
\end{lem}

\noindent{\bf Proof.}  Notice that in the statement of this
result, we took advantage of Lemma \ref{terminal} so that we could
refer to the terminal object of $\mathcal{C}$ as $k$.

To prove this result, note first that since $1 \in k$ satisfies
$\Delta(1) = 1 \tensor 1$ and \newline $\id \maps k \to SV$ is a
coalgebra homomorphism, we have $\Delta(\id(1)) = \id(1) \tensor
\id(1)$. Thus, it suffices to show that the only element $x \in
SV$ with $\Delta(x) = x \tensor x$ is $x = 1$.

To do this, note that we can think of $SV$ as the algebra of
polynomial functions on $V^\ast$, and the comultiplication $\Delta
\maps SV \to SV \tensor SV \iso S(V \oplus V)$ as the following
map from polynomials on $V$ to polynomials on $V^\ast \oplus
V^\ast$:
$$   \Delta(x)(a,b) = x(a + b)  $$
for all $a,b \in V^\ast$ and all polynomials $x$ on $V^\ast$. The
condition that $\Delta(x) = x \tensor x$ thus says that
$$    x(a + b) = x(a) x(b)  $$
for all $a,b \in V^\ast$.  If $x$ is a polynomial of degree $n$ on
$V^\ast$, the left side of the above equation describes a
polynomial of degree $n$ on $V^\ast \oplus V^\ast$, while the
right side describes a polynomial of degree $2n$.  This is only
possible if $n = 0$, so $x$ must be a constant polynomial, and the
only constant polynomial satisfying the above equation is $x = 1$.
\qed

Lemma \ref{multid} together with Lemma \ref{terminal} allow us to
interpret conditions $(a)$ and $(b)$ in Definition \ref{wellpoint}
for an element $x \in Q$, where $Q$ is a unital spindle in
$\mathcal{C}$.  Axiom $(a)$ becomes
$$1 \rhd x = x,$$
while axiom $(b)$ becomes
$$x \rhd 1 = \epsilon(x) 1. $$
We will make use of both of these identities below.

We now have all the ingredients to construct a Lie algebra from a
unital spindle in $\mathcal{C}$. Proposition \ref{ptdquandobj}
demonstrated that groups in $\mathcal{C}$ give unital spindles in
$\mathcal{C}$.  Thus, given a group in $\mathcal{C}$, it becomes a
unital spindle in $\mathcal{C}$, say $(Q, V, \alpha)$.  Now we
shall use the spindle operation $\rhd \maps Q \otimes Q \to Q$ to
define a bracket on $V$. Since we want a map $[\cdot, \cdot] \maps
V \otimes V \to V,$ we first use the isomorphism $\alpha$ to
transfer $\rhd$ from $Q$ to $SV$:
\[
\xy (-15,0)*+{SV \otimes SV}="X1";
    (-15,-3.5)*{}="a";
    (15,0)*+{SV}="X2";
    (15,-3.5)*{}="b";
    (-15,-30)*+{Q \otimes Q}="X3";
    (15,-30)*+{Q}="X4";
             {\ar^{\rhd}  "X1";"X2"};
             {\ar_{\alpha ^{-1} \otimes \alpha ^{-1}} "X1";"X3"};
             {\ar^{\rhd }  "X3";"X4"};
             {\ar_{\alpha}  "X4";"X2"};
\endxy
\]
We then define the bracket as the following composite:
\[
\xy (-15,0)*+{V \otimes V}="X1";
    (-15,-3.5)*{}="a";
    (15,0)*+{V}="X2";
    (15,-3.5)*{}="b";
    (-15,-30)*+{SV \otimes SV}="X3";
    (15,-30)*+{SV}="X4";
             {\ar^{[\cdot, \cdot]}  "X1";"X2"};
             {\ar@{_{(}->} "a";"X3"};
             {\ar^{\rhd }  "X3";"X4"};
             {\ar_{\pi}  "X4";"X2"};
\endxy
\]
Using this, we finally obtain our desired result:

\begin{thm} \et \label{punchlinethm} Given a unital spindle
$Q \iso SV$ in the category $\mathcal{C}$ of special coalgebras,
and defining $[\cdot, \cdot] \maps V \otimes V \to V$ as
\[
\xy (-15,0)*+{V \otimes V}="X1";
    (-15,-3.5)*{}="a";
    (15,0)*+{V}="X2";
    (15,-3.5)*{}="b";
    (-15,-30)*+{SV \otimes SV}="X3";
    (15,-30)*+{SV}="X4";
             {\ar^{[\cdot, \cdot]}  "X1";"X2"};
             {\ar@{_{(}->} "a";"X3"};
             {\ar^{\rhd}  "X3";"X4"};
             {\ar^{\pi}  "X4";"X2"};
\endxy
\]
then $V$ becomes a Lie algebra.
\end{thm}

\noindent{\bf Proof.} To prove the antisymmetry of the bracket it
suffices to show that $[v,v] = 0$ for all $v \in V$, and by the
above definition of the bracket it suffices to show that $v \rhd v
= 0$. We will prove this using the internalized version of the
left idempotence law of a spindle. Recall that the left
idempotence law says that
\[
  \vcenter{
  \xymatrix@!C{
  & SV \times SV \ar[d]^{\rhd} \\
    SV
    \ar[ur]^{\Delta}
    \ar[r]_{1}
  & SV
     }}
\]
commutes, where $SV \times SV$ is computed using the product in
$\mathcal{C}$, namely the tensor product of special coalgebras.

Let $v \in V$.  Then $v^2 \in S(V) \cong Q$. Applying the
internalized idempotence law to $v^2$, we have
$$v^2 = \rhd(\Delta(v^2)) = \rhd(\Delta(v) \Delta(v)),$$
since $\Delta$ is an algebra homomorphism, due to the fact that
$SV$, and hence $Q$, is a bialgebra.  But by Lemma \ref{diagonal},
the diagonal $\Delta$ is the same as the comultiplication map on
$Q$, so we have:
$$
\begin{array}{ccl}
\rhd(\Delta(v) \Delta(v)) &=& \rhd[(v \otimes
1 + 1 \otimes v)(v \otimes 1 + 1 \otimes v)] \\
&=& \rhd(v^2 \otimes 1 + v \otimes v + v \otimes v + 1 \otimes v^2) \\
&=& v^2 \rhd 1 + 2v \rhd v + 1 \rhd v^2.
\end{array}
$$
and thus
$$ v^2 = v^2 \rhd 1 + 2v \rhd v + 1 \rhd v^2.$$
Now we recall that conditions $(a)$ and $(b)$ of a unital spindle
say, respectively, that
$$1 \rhd x = x \; \; \;  \textrm{and} \; \; \; x \rhd 1 = \epsilon(x) 1$$
for $x \in Q$.  We thus have
$$v^2 = \epsilon(v^2) + 2 v \rhd v + v^2 $$
but the counit in $SV$ satisfies $\epsilon(v^2) = 0$, so
$$v \rhd v = 0$$
as desired.

Similarly, the Jacobi identity follows from the fact that
\[       u \rhd (v \rhd w) = (u \rhd v) \rhd w + v \rhd (u \rhd w)  \]
for all $u,v,w \in V$.  This in turn follows from the internalized
left distributive law of our spindle in $\mathcal{C}$:
\[
\begin{xy} 0;/r.28pc/:
    (-2, 20)*+{SV \times SV \times SV}="1";
    (-30,10)*+{SV \times SV \times SV \times SV}="2";
    (-30,-10)*+{SV \times SV \times SV \times SV}="3";
    (-15,-20)*+{SV \times SV \times SV}="4";
    (15,-20)*+{SV \times SV}="4a";
    (30,-10)*+{SV}="5";
    (30, 10)*+{SV \times SV}="6";
        {\ar^{1 \times \rhd} "1";"6"};
        {\ar_{\Delta \times 1 \times 1} "1";"2"};
        {\ar_{1 \times S \times 1} "2";"3"};
        {\ar_{\rhd \times 1 \times 1} "3";"4"};
        {\ar_{1 \times \rhd} "4";"4a"};
        {\ar_{\rhd} "4a";"5"};
        {\ar^{\rhd} "6";"5"};
\end{xy}
\]
which says
\[
\def\objectstyle{\scriptstyle}
  \def\labelstyle{\scriptstyle}
\begin{xy} 0;/r.30pc/:
    (-2, 20)*+{u \otimes v \otimes w}="1";
    (-30,10)*+{(u \otimes 1 + 1 \otimes u) \otimes v \otimes w}="2";
    (-30,6)*+{= u \otimes 1 \otimes v \otimes w + 1 \otimes u
    \otimes v \otimes w}="2a";
    (-30,-10)*+{u \otimes v \otimes 1 \otimes w + 1 \otimes v \otimes u \otimes w}="3";
    (-15,-20)*+{u \rhd v \otimes 1 \otimes w + v \otimes u \otimes w}="4";
    (15,-20)*+{u \rhd v \otimes w + v \otimes u \rhd w}="4a";
    (30,-10)*+{(u \rhd v) \rhd w + v \rhd (u \rhd w)}="5";
    (30,-6)*+{u \rhd (v \rhd w) = }="5a";
    (30, 10)*+{u \otimes [v,w]}="6";
        {\ar^{1 \otimes \rhd} "1";"6"};
        {\ar_{\Delta \otimes 1 \otimes 1} "1";"2"};
        {\ar_{1 \otimes S \otimes 1} "2a";"3"};
        {\ar_{\rhd \otimes 1 \otimes 1} "3";"4"};
        {\ar_{1 \otimes \rhd} "4";"4a"};
        {\ar_{\rhd} "4a";"5"};
        {\ar^{\rhd} "6";"5a"};
\end{xy}
\]
Thus, unital spindles in $\mathcal{C}$ give Lie algebras. \qed

But, our original goal was to obtain the Lie algebra of a given
Lie group, not of a unital spindle in $\mathcal{C}$!  Applying the
previous theorem to the special case when the unital spindle in
$\mathcal{C}$ came from a group in $\mathcal{C}$, which, in turn,
came from a group in $\Diff_{\ast}$, we obtain the result we
desire.

To summarize, we carry out the following procedure:  We begin by
using the language of internalization to say that a Lie group $G$
is a group in $\Diff_{\ast}$.  We next apply the cojet functor
$J_{\infty}$, obtaining a group in $\mathcal{C}$ since this
functor preserves products. A group in $\mathcal{C}$ automatically
gives a unital spindle in $\mathcal{C}$. Then, applying the
functor $F$ to this result, we extract a vector space which we
call $\mathfrak{g}$. Finally, as above, we define the bracket on
this vector space $\mathfrak{g}$ as:
\[
\xy (-15,0)*+{\mathfrak{g} \otimes \mathfrak{g}}="X1";
    (-15,-3.5)*{}="a";
    (15,0)*+{\mathfrak{g}}="X2";
    (15,-3.5)*{}="b";
    (-15,-30)*+{S\mathfrak{g} \otimes S\mathfrak{g}}="X3";
    (15,-30)*+{S\mathfrak{g}}="X4";
             {\ar^{[\cdot, \cdot]}  "X1";"X2"};
             {\ar@{_{(}->} "a";"X3"};
             {\ar^{\rhd}  "X3";"X4"};
             {\ar^{\pi}  "X4";"X2"};
\endxy
\]
and use the axioms of a unital spindle in $\mathcal{C}$ to show
that this bracket operation is antisymmetric and satisfies the
Jacobi identity, so that $\mathfrak{g}$ is the Lie algebra of $G$,
as desired.

\begin{cor} \et If $G$ is a Lie group and $Q$ is the unital
spindle in $\mathcal{C}$ defined as described in Proposition
\ref{ptdquandobj}, the Lie algebra $\mathfrak{g}$ constructed in
Theorem \ref{punchlinethm} is isomorphic to the Lie algebra of
$G$.
\end{cor}

Thus, we have described a novel method of obtaining the Lie
algebra of a Lie group using spindles.  In the next section, we
indicate how we anticipate categorifying this procedure in order
to acquire the Lie $2$-algebra of a Lie $2$-group.

\section{Lie 2-algebras, 2-Quandles and 2-Braids} \label{conclusions}

We conclude by outlining how we suspect the categorification of
the previous three sections will occur.  We start by categorifying
the notions of shelf, rack, spindle and quandle. Recall that
categorifying a mathematical concept involves finding
category-theoretic analogs of set-theoretic concepts.  Therefore,
to obtain a categorified shelf, or `$2$-shelf', we will replace
the set $Q$ in Definition \ref{lshelf} by a category, the function
$\rhd \maps Q \times Q \to Q$ by a functor, called `left
conjugation' and the equation $(i)$ by a natural isomorphism,
called the `left distributor', which will take us from one side of
the left distributive law to the other. In addition, we will
require that the left distributor satisfy and equation of its own,
known as a coherence law.  While we only present the definition of
a `$2$-shelf' here, we hint at the correct formulations of the
notions of `$2$-rack', `$2$-spindle' and `$2$-quandle'.

\begin{defn} \et \label{l2-shelf}
A {\bf left 2-shelf}, $(Q, \rhd)$, consists of:
\begin{itemize}
\item a category $Q$,
\end{itemize}
equipped with:
\begin{itemize}
\item a functor, {\bf left conjugation}, $\rhd \maps Q \times Q \to Q$,
\item a natural isomorphism, the {\bf left distributor},
$$LD_{x,y,z} \maps x \rhd (y \rhd z) \to (x \rhd y) \rhd (x \rhd z)$$
\end{itemize}
that is required to satisfy:
\begin{itemize}
\item the {\bf distributor identity:}
$$(z \rhd LD_{y,x,w}) LD_{z, y \rhd x, y \rhd w}
((z \rhd ( y \rhd x)) \rhd LD_{z,y,w})  (LD_{z,y,x} \rhd ((z \rhd
y) \rhd (z \rhd w))) =$$
$$LD_{z,y, x \rhd w} ((z \rhd y) \rhd LD_{z,x,w})  LD_{z \rhd y, z \rhd
x, z \rhd w}$$
\end{itemize}
for all $w,x,y,z \in Q$.
\end{defn}

As in the case of the Jacobiator identity, the left distributor
identity becomes more manageable if we draw it as a commutative
diagram. Doing so enables us to see that this law relates two ways
of using the left distributor to reparenthesize the expression $z
\rhd (y \rhd (x \rhd w))$:
$$ \def\objectstyle{\scriptstyle}
  \def\labelstyle{\scriptstyle}
   \xy
   (0,35)*+{z \rhd (y \rhd (x \rhd w))}="1";
   (-40,20)*+{z \rhd ((y \rhd x) \rhd (y \rhd w))}="2";
   (40,20)*+{z \rhd (y \rhd (x \rhd w))}="3";
   (-40,0)*+{(z \rhd (y \rhd x)) \rhd (z \rhd (y \rhd w))}="4'";
   (-40,-4)*+{}="4";
   (40,0)*+{(z \rhd y) \rhd (z \rhd (x \rhd w))}="5'";
   (-40,-20)*+{}="6'";
   (-40,-24)*+{(z \rhd (y \rhd x)) \rhd ((z \rhd y) \rhd (z \rhd w))}="6";
   (40,-20)*+{(z \rhd y) \rhd ((z \rhd x) \rhd (z \rhd w))}="7'";
   (40,-24)*+{}="7";
   (0,-40)*+{((z \rhd y) \rhd (z \rhd x)) \rhd ((z \rhd y) \rhd (z \rhd w))}="8'";
   (0,-44)*+{}="8";
            (32,-32)*{LD_{z \rhd y, z \rhd x, z \rhd w} }; 
            (-40, -34)*+{LD_{z,y,x} \rhd ((z \rhd y) \rhd (z \rhd
            w))};
        {\ar_{z \rhd LD_{y,x,w}}               "1";"2"};
        {\ar^{1}                               "1";"3"};
        {\ar_{LD_{z, y \rhd x, y \rhd w}}      "2";"4'"};
        {\ar_{(z \rhd (y \rhd x)) \rhd LD_{z,y,w}}                      "4";"6'"};
        {\ar^{LD_{z,y, x \rhd w}}              "3";"5'"};
        {\ar^{(z \rhd y) \rhd LD_{z,x,w}}                      "5'";"7'"};
        {\ar_{}                      "6";"8'"};
        {\ar^{}                                "7'";"8'"};
\endxy
$$
\\
\noindent At this point, we remark that we can take a cue from
topology to assist us in our categorification!  We saw in Section
\ref{topology2} that the Yang--Baxter equation, or third
Reidemeister move, is equivalent to the left distributive law for
the shelf operation.  Therefore, rather than thinking of the left
distributor solely as an algebraic object, we can think of it as
the process of performing the third Reidemeister move. With this
in mind, the coherence law for the left distributor will be the
higher dimensional analogue of the third Reidemeister move, the
Zamolodchikov tetrahedron equation, familiar from the theory of
$2$-knots and braided monoidal $2$-categories
\cite{BLan,BN,CS,Crans,KV}.   This equation plays a role in the
theory of knotted surfaces in 4-space which is closely analogous
to that played by the Yang--Baxter equation, or third Reidemeister
move, in the theory of ordinary knots in 3-space.

The analogue of the Yang--Baxter equation is called the
`Zamolodchikov tetrahedron equation':

\begin{defn} \et Given a category $C$ and an
invertible functor $B \maps C \times C \to C \times C$, a natural
isomorphism
$$Y \maps (B \times 1)(1 \times B)(B \times 1) \To
          (1 \times B)(B \times 1)(1 \times B) $$
satisfies the {\bf Zamolodchikov tetrahedron equation} if
$$[Y \circ (1 \ttimes 1 \ttimes B)(1 \ttimes B \ttimes 1)(B \ttimes 1
\ttimes 1)] [(1 \ttimes B \ttimes 1)(B \ttimes 1 \ttimes 1) \circ
Y \circ (B \ttimes 1 \ttimes 1)]$$ $$[(1 \ttimes B \ttimes 1)(1
\ttimes 1 \ttimes B) \circ Y \circ (1 \ttimes 1 \ttimes B)] [Y
\circ (B \ttimes 1 \ttimes 1)(1 \ttimes B \ttimes 1)(1 \ttimes 1
\ttimes B)]$$
$$ = $$
$$[(B \ttimes 1 \ttimes 1)(1 \ttimes B \ttimes 1)(1 \ttimes 1 \ttimes B)
\circ Y ] [(B \ttimes 1 \ttimes 1) \circ Y \circ (B \ttimes 1
\ttimes 1)(1 \ttimes B \ttimes 1)]$$
$$[(1 \ttimes 1 \ttimes B) \circ Y \circ (1 \ttimes 1 \ttimes
B)(1 \ttimes B \ttimes 1)] [(1 \ttimes 1 \ttimes B)(1 \ttimes B
\ttimes 1)(B \ttimes 1 \ttimes 1) \circ Y],$$ where $\circ$
represents the whiskering of a functor by a natural
transformation.  We will often refer to $Y$ as a {\bf
Yang--Baxterator}.
\end{defn}

To see the significance of this complex but beautifully
symmetrical equation, one should think of $Y$ as the surface in
4-space traced out by the process of performing the third
Reidemeister move:

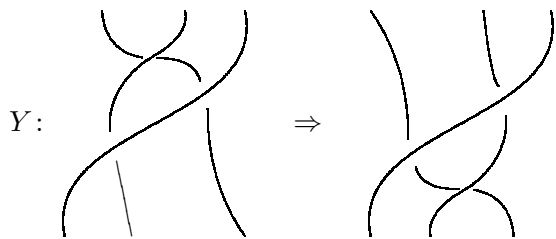
\[   Y \maps
\def\objectstyle{\scriptstyle}
\def\labelstyle{\scriptstyle}
  \xy
   (12,15)*{}="C";
   (4,15)*{}="B";
   (-7,15)*{}="A";
   (12,-15)*{}="3";
   (-3,-15)*{}="2";
   (-12,-15)*{}="1";
       "C";"1" **\crv{(15,0)& (-15,0)};
       (-5,-5)*{}="2'";
       (7,2)*{}="3'";
     "2'";"2" **\crv{};
     "3'";"3" **\crv{(7,-8)};
       \vtwist~{"A"}{"B"}{(-6,-1)}{(6,5.7)};
\endxy
 \quad \To \quad \xy
   (-12,-15)*{}="C";
   (-4,-15)*{}="B";
   (7,-15)*{}="A";
   (-12,15)*{}="3";
   (3,15)*{}="2";
   (12,15)*{}="1";
       "C";"1" **\crv{(-15,0)& (15,0)};
       (5,5)*{}="2'";
       (-7,-2)*{}="3'";
     "2'";"2" **\crv{(4,6)};
     "3'";"3" **\crv{(-7,8)};
       \vtwist~{"A"}{"B"}{(6,1)}{(-6,-5.7)};
\endxy
\]

\noindent Then the Zamolodchikov tetrahedron equation says the
surface traced out by first performing the third Reidemeister move
on a threefold crossing and then sliding the result under a fourth
strand:
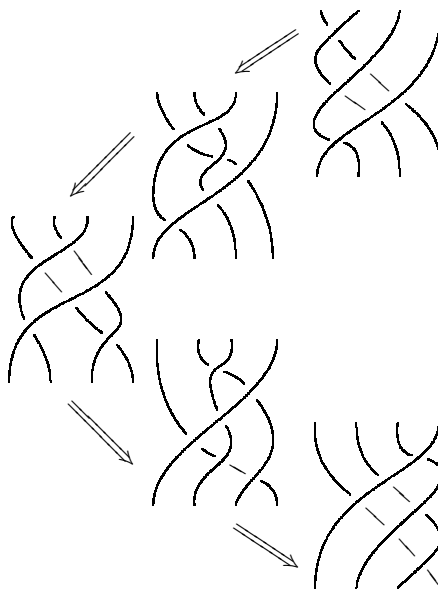
\[
  \xy 0;/r.13pc/:   
    (0,50)*+{       
 \xy 
 (-15,-20)*{}="T1";
 (-5,-20)*{}="T2";
 (5,-20)*{}="T3";
 (15,-20)*{}="T4";
 (-14,20)*{}="B1";
 (-5,20)*{}="B2";
 (5,20)*{}="B3";
 (15,20)*{}="B4";
    "T1"; "B4" **\crv{(-15,-7) & (15,-5)}
        \POS?(.25)*{\hole}="2x" \POS?(.47)*{\hole}="2y" \POS?(.6)*{\hole}="2z";
    "T2";"2x" **\crv{(-4,-12)};
    "T3";"2y" **\crv{(5,-10)};
    "T4";"2z" **\crv{(16,-9)};
 (-15,-5)*{}="3x";
    "2x"; "3x" **\crv{(-18,-10)};
    "3x"; "B3" **\crv{(-13,0) & (4,10)}
        \POS?(.3)*{\hole}="4x" \POS?(.53)*{\hole}="4y";
    "2y"; "4x" **\crv{};
    "2z"; "4y" **\crv{};
 (-15,10)*{}="5x";
    "4x";"5x" **\crv{(-17,6)};
    "5x";"B2" **\crv{(-14,12)}
        \POS?(.6)*{\hole}="6x";
    "6x";"B1" **\crv{(-14,18)};
    "4y";"6x" **\crv{(-8,10)};
 \endxy
    }="T";
    (-40,30)*+{
 \xy 
 (-15,-20)*{}="b1";
 (-5,-20)*{}="b2";
 (5,-20)*{}="b3";
 (14,-20)*{}="b4";
 (-14,20)*{}="T1";
 (-5,20)*{}="T2";
 (5,20)*{}="T3";
 (15,20)*{}="T4";
    "b1"; "T4" **\crv{(-15,-7) & (15,-5)}
        \POS?(.25)*{\hole}="2x" \POS?(.47)*{\hole}="2y" \POS?(.65)*{\hole}="2z";
    "b2";"2x" **\crv{(-5,-15)};
    "b3";"2y" **\crv{(5,-10)};
    "b4";"2z" **\crv{(14,-9)};
 (-15,-5)*{}="3x";
    "2x"; "3x" **\crv{(-15,-10)};
    "3x"; "T3" **\crv{(-15,15) & (5,10)}
        \POS?(.38)*{\hole}="4y" \POS?(.65)*{\hole}="4z";
    "T1";"4y" **\crv{(-14,16)};
    "T2";"4z" **\crv{(-5,16)};
    "2y";"4z" **\crv{(-10,3) & (10,2)} \POS?(.6)*{\hole}="5z";
    "4y";"5z" **\crv{(-5,5)};
    "5z";"2z" **\crv{(5,4)};
 \endxy
    }="TL";
    (-75,0)*+{
 \xy 
 (-14,20)*{}="T1";
 (-4,20)*{}="T2";
 (4,20)*{}="T3";
 (15,20)*{}="T4";
 (-15,-20)*{}="B1";
 (-5,-20)*{}="B2";
 (5,-20)*{}="B3";
 (15,-20)*{}="B4";
    "B1";"T4" **\crv{(-15,5) & (15,-5)}
        \POS?(.25)*{\hole}="2x" \POS?(.49)*{\hole}="2y" \POS?(.65)*{\hole}="2z";
    "2x";"T3" **\crv{(-20,10) & (5,10) }
        \POS?(.45)*{\hole}="3y" \POS?(.7)*{\hole}="3z";
    "2x";"B2" **\crv{(-5,-14)};
        "T1";"3y" **\crv{(-16,17)};
        "T2";"3z" **\crv{(-5,17)};
        "3z";"2z" **\crv{};
        "3y";"2y" **\crv{};
        "B3";"2z" **\crv{ (5,-5) &(20,-10)}
            \POS?(.4)*{\hole}="4z";
        "2y";"4z" **\crv{(6,-8)};
        "4z";"B4" **\crv{(15,-15)};
 \endxy
    }="ML";
    (-40,-30)*+{
 \xy 
 (-14,20)*{}="T1";
 (-4,20)*{}="T2";
 (4,20)*{}="T3";
 (15,20)*{}="T4";
 (-15,-20)*{}="B1";
 (-5,-20)*{}="B2";
 (5,-20)*{}="B3";
 (15,-20)*{}="B4";
    "B1";"T4" **\crv{(-15,-5) & (15,5)}
        \POS?(.38)*{\hole}="2x" \POS?(.53)*{\hole}="2y" \POS?(.7)*{\hole}="2z";
    "T1";"2x" **\crv{(-15,5)};
    "2y";"B2" **\crv{(10,-10) & (-6,-10)}
        \POS?(.45)*{\hole}="4x";
    "2z";"B3" **\crv{ (15,0)&(15,-10) & (6,-16)}
        \POS?(.7)*{\hole}="5x";
    "T3";"2y" **\crv{(5,10)& (-6,18) }
        \POS?(.5)*{\hole}="3x";
    "T2";"3x" **\crv{(-5,15)};
    "3x";"2z" **\crv{(7,11)};
    "2x";"4x" **\crv{(-3,-7)};
    "4x";"5x" **\crv{};
    "5x";"B4" **\crv{(15,-15)};
 \endxy
    }="BL";
(0,-50)*+{
 \xy 
 (15,20)*{}="T1";
 (5,20)*{}="T2";
 (-5,20)*{}="T3";
 (-15,20)*{}="T4";
 (15,-20)*{}="B1";
 (5,-20)*{}="B2";
 (-5,-20)*{}="B3";
 (-15,-20)*{}="B4";
    "T1"; "B4" **\crv{(15,7) & (-15,5)}
        \POS?(.25)*{\hole}="2x" \POS?(.45)*{\hole}="2y" \POS?(.6)*{\hole}="2z";
    "T2";"2x" **\crv{(4,12)};
    "T3";"2y" **\crv{(-5,10)};
    "T4";"2z" **\crv{(-16,9)};
 (15,5)*{}="3x";
    "2x"; "3x" **\crv{(18,10)};
    "3x"; "B3" **\crv{(13,0) & (-4,-10)}
        \POS?(.3)*{\hole}="4x" \POS?(.53)*{\hole}="4y";
    "2y"; "4x" **\crv{};
    "2z"; "4y" **\crv{};
 (15,-10)*{}="5x";
    "4x";"5x" **\crv{(17,-6)};
    "5x";"B2" **\crv{(14,-12)}
        \POS?(.6)*{\hole}="6x";
    "6x";"B1" **\crv{};
    "4y";"6x" **\crv{};
 \endxy
    }="B";
            (-20,65)*{}="X1";
            (-35,55)*{}="X2";               
                {\ar@{=>} "X1";"X2"};       
            (-60,40)*{}="X1";               
            (-75,25)*{}="X2";               
                {\ar@{=>} "X1";"X2"};
            (-60,-40)*{}="X2";
            (-75,-25)*{}="X1";
                {\ar@{=>} "X1";"X2"};
            (-20,-65)*{}="X2";
            (-35,-55)*{}="X1";
                {\ar@{=>} "X1";"X2"};
  \endxy
\]
is isotopic to that traced out by first sliding the threefold
crossing under the fourth strand and then performing the third
Reidemeister move:
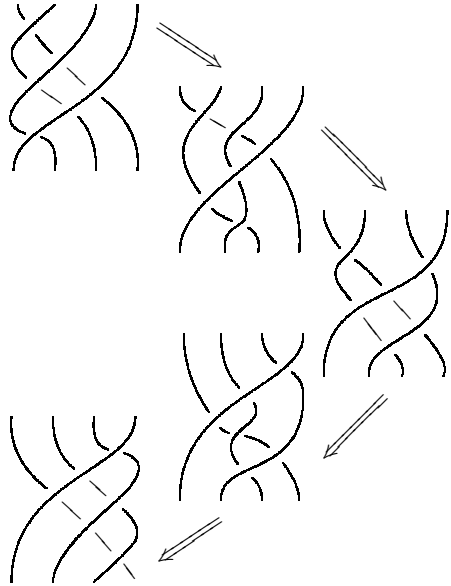
\[
  \xy 0;/r.13pc/:   
    (0,50)*+{       
 \xy 
 (-15,-20)*{}="T1";
 (-5,-20)*{}="T2";
 (5,-20)*{}="T3";
 (15,-20)*{}="T4";
 (-14,20)*{}="B1";
 (-5,20)*{}="B2";
 (5,20)*{}="B3";
 (15,20)*{}="B4";
    "T1"; "B4" **\crv{(-15,-7) & (15,-5)}
        \POS?(.25)*{\hole}="2x" \POS?(.47)*{\hole}="2y" \POS?(.6)*{\hole}="2z";
    "T2";"2x" **\crv{(-4,-12)};
    "T3";"2y" **\crv{(5,-10)};
    "T4";"2z" **\crv{(16,-9)};
 (-15,-5)*{}="3x";
    "2x"; "3x" **\crv{(-18,-10)};
    "3x"; "B3" **\crv{(-13,0) & (4,10)}
        \POS?(.3)*{\hole}="4x" \POS?(.53)*{\hole}="4y";
    "2y"; "4x" **\crv{};
    "2z"; "4y" **\crv{};
 (-15,10)*{}="5x";
    "4x";"5x" **\crv{(-17,6)};
    "5x";"B2" **\crv{(-14,12)}
        \POS?(.6)*{\hole}="6x";
    "6x";"B1" **\crv{(-14,18)};
    "4y";"6x" **\crv{(-8,10)};
 \endxy
    }="T";
(40,30)*+{
 \xy 
 (14,-20)*{}="T1";
 (4,-20)*{}="T2";
 (-4,-20)*{}="T3";
 (-15,-20)*{}="T4";
 (15,20)*{}="B1";
 (5,20)*{}="B2";
 (-5,20)*{}="B3";
 (-15,20)*{}="B4";
    "B1";"T4" **\crv{(15,5) & (-15,-5)}
        \POS?(.38)*{\hole}="2x" \POS?(.53)*{\hole}="2y" \POS?(.7)*{\hole}="2z";
    "T1";"2x" **\crv{(15,-5)};
    "2y";"B2" **\crv{(-10,10) & (6,10)}
        \POS?(.45)*{\hole}="4x";
    "2z";"B3" **\crv{ (-15,0)&(-15,10) & (-6,16)}
        \POS?(.7)*{\hole}="5x";
    "T3";"2y" **\crv{(-5,-10)& (6,-18) }
        \POS?(.5)*{\hole}="3x";
    "T2";"3x" **\crv{(5,-15)};
    "3x";"2z" **\crv{(-7,-11)};
    "2x";"4x" **\crv{(3,7)};
    "4x";"5x" **\crv{};
    "5x";"B4" **\crv{(-15,15)};
 \endxy
    }="TR";
    (75,0)*+{
 \xy 
 (14,-20)*{}="T1";
 (4,-20)*{}="T2";
 (-4,-20)*{}="T3";
 (-15,-20)*{}="T4";
 (15,20)*{}="B1";
 (5,20)*{}="B2";
 (-5,20)*{}="B3";
 (-15,20)*{}="B4";
    "B1";"T4" **\crv{(15,-5) & (-15,5)}
        \POS?(.25)*{\hole}="2x" \POS?(.49)*{\hole}="2y" \POS?(.65)*{\hole}="2z";
    "2x";"T3" **\crv{(20,-10) & (-5,-10) }
        \POS?(.45)*{\hole}="3y" \POS?(.7)*{\hole}="3z";
    "2x";"B2" **\crv{(5,14)};
        "T1";"3y" **\crv{(16,-17)};
        "T2";"3z" **\crv{(5,-17)};
        "3z";"2z" **\crv{};
        "3y";"2y" **\crv{};
        "B3";"2z" **\crv{ (-5,5) &(-20,10)}
            \POS?(.4)*{\hole}="4z";
        "2y";"4z" **\crv{(-6,8)};
        "4z";"B4" **\crv{(-15,15)};
 \endxy
    }="MR";
    (40,-30)*+{
 \xy 
 (15,20)*{}="b1";
 (5,20)*{}="b2";
 (-5,20)*{}="b3";
 (-14,20)*{}="b4";
 (14,-20)*{}="T1";
 (5,-20)*{}="T2";
 (-5,-20)*{}="T3";
 (-15,-20)*{}="T4";
    "b1"; "T4" **\crv{(15,7) & (-15,5)}
        \POS?(.25)*{\hole}="2x" \POS?(.47)*{\hole}="2y" \POS?(.65)*{\hole}="2z";
    "b2";"2x" **\crv{(5,15)};
    "b3";"2y" **\crv{(-5,10)};
    "b4";"2z" **\crv{(-14,9)};
 (15,5)*{}="3x";
    "2x"; "3x" **\crv{(15,10)};
    "3x"; "T3" **\crv{(15,-15) & (-5,-10)}
        \POS?(.38)*{\hole}="4y" \POS?(.65)*{\hole}="4z";
    "T1";"4y" **\crv{(14,-16)};
    "T2";"4z" **\crv{(5,-16)};
    "2y";"4z" **\crv{(10,-3) & (-10,-2)} \POS?(.6)*{\hole}="5z";
    "4y";"5z" **\crv{(5,-5)};
    "5z";"2z" **\crv{(-5,-4)};
 \endxy
    }="BR";
    (0,-50)*+{
 \xy 
 (15,20)*{}="T1";
 (5,20)*{}="T2";
 (-5,20)*{}="T3";
 (-15,20)*{}="T4";
 (15,-20)*{}="B1";
 (5,-20)*{}="B2";
 (-5,-20)*{}="B3";
 (-15,-20)*{}="B4";
    "T1"; "B4" **\crv{(15,7) & (-15,5)}
        \POS?(.25)*{\hole}="2x" \POS?(.45)*{\hole}="2y" \POS?(.6)*{\hole}="2z";
    "T2";"2x" **\crv{(4,12)};
    "T3";"2y" **\crv{(-5,10)};
    "T4";"2z" **\crv{(-16,9)};
 (15,5)*{}="3x";
    "2x"; "3x" **\crv{(18,10)};
    "3x"; "B3" **\crv{(13,0) & (-4,-10)}
        \POS?(.3)*{\hole}="4x" \POS?(.53)*{\hole}="4y";
    "2y"; "4x" **\crv{};
    "2z"; "4y" **\crv{};
 (15,-10)*{}="5x";
    "4x";"5x" **\crv{(17,-6)};
    "5x";"B2" **\crv{(14,-12)}
        \POS?(.6)*{\hole}="6x";
    "6x";"B1" **\crv{};
    "4y";"6x" **\crv{};
 \endxy
    }="B";
            (20,65)*{}="X1";                
            (35,55)*{}="X2";                
                {\ar@{=>} "X1";"X2"};       
            (60,40)*{}="X1";                
            (75,25)*{}="X2";                
                {\ar@{=>} "X1";"X2"};       
            (60,-40)*{}="X2";
            (75,-25)*{}="X1";
                {\ar@{=>} "X1";"X2"};
            (20,-65)*{}="X2";
            (35,-55)*{}="X1";
                {\ar@{=>} "X1";"X2"};
  \endxy
\]
In short, the Zamolodchikov tetrahedron equation is a
formalization of this commutative octagon:
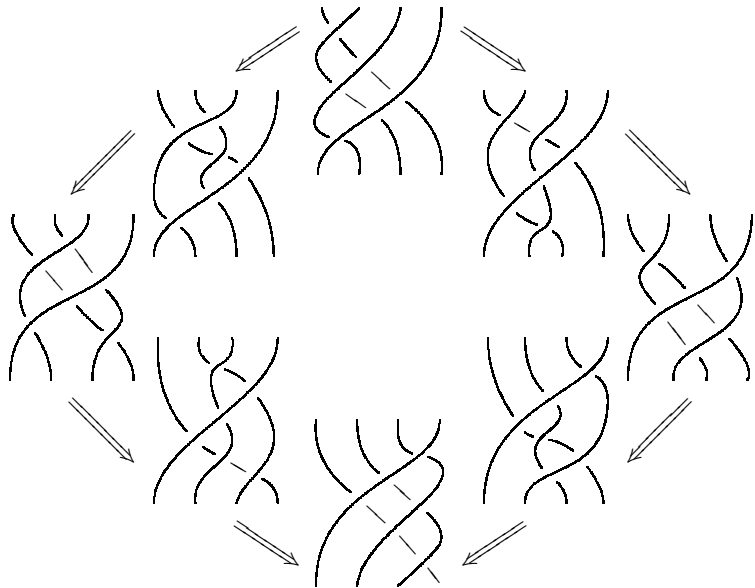
\[
  \xy 0;/r.13pc/:   
    (0,50)*+{       
 \xy 
 (-15,-20)*{}="T1";
 (-5,-20)*{}="T2";
 (5,-20)*{}="T3";
 (15,-20)*{}="T4";
 (-14,20)*{}="B1";
 (-5,20)*{}="B2";
 (5,20)*{}="B3";
 (15,20)*{}="B4";
    "T1"; "B4" **\crv{(-15,-7) & (15,-5)}
        \POS?(.25)*{\hole}="2x" \POS?(.47)*{\hole}="2y" \POS?(.6)*{\hole}="2z";
    "T2";"2x" **\crv{(-4,-12)};
    "T3";"2y" **\crv{(5,-10)};
    "T4";"2z" **\crv{(16,-9)};
 (-15,-5)*{}="3x";
    "2x"; "3x" **\crv{(-18,-10)};
    "3x"; "B3" **\crv{(-13,0) & (4,10)}
        \POS?(.3)*{\hole}="4x" \POS?(.53)*{\hole}="4y";
    "2y"; "4x" **\crv{};
    "2z"; "4y" **\crv{};
 (-15,10)*{}="5x";
    "4x";"5x" **\crv{(-17,6)};
    "5x";"B2" **\crv{(-14,12)}
        \POS?(.6)*{\hole}="6x";
    "6x";"B1" **\crv{(-14,18)};
    "4y";"6x" **\crv{(-8,10)};
 \endxy
    }="T";
    (-40,30)*+{
 \xy 
 (-15,-20)*{}="b1";
 (-5,-20)*{}="b2";
 (5,-20)*{}="b3";
 (14,-20)*{}="b4";
 (-14,20)*{}="T1";
 (-5,20)*{}="T2";
 (5,20)*{}="T3";
 (15,20)*{}="T4";
    "b1"; "T4" **\crv{(-15,-7) & (15,-5)}
        \POS?(.25)*{\hole}="2x" \POS?(.47)*{\hole}="2y" \POS?(.65)*{\hole}="2z";
    "b2";"2x" **\crv{(-5,-15)};
    "b3";"2y" **\crv{(5,-10)};
    "b4";"2z" **\crv{(14,-9)};
 (-15,-5)*{}="3x";
    "2x"; "3x" **\crv{(-15,-10)};
    "3x"; "T3" **\crv{(-15,15) & (5,10)}
        \POS?(.38)*{\hole}="4y" \POS?(.65)*{\hole}="4z";
    "T1";"4y" **\crv{(-14,16)};
    "T2";"4z" **\crv{(-5,16)};
    "2y";"4z" **\crv{(-10,3) & (10,2)} \POS?(.6)*{\hole}="5z";
    "4y";"5z" **\crv{(-5,5)};
    "5z";"2z" **\crv{(5,4)};
 \endxy
    }="TL";
    (-75,0)*+{
 \xy 
 (-14,20)*{}="T1";
 (-4,20)*{}="T2";
 (4,20)*{}="T3";
 (15,20)*{}="T4";
 (-15,-20)*{}="B1";
 (-5,-20)*{}="B2";
 (5,-20)*{}="B3";
 (15,-20)*{}="B4";
    "B1";"T4" **\crv{(-15,5) & (15,-5)}
        \POS?(.25)*{\hole}="2x" \POS?(.49)*{\hole}="2y" \POS?(.65)*{\hole}="2z";
    "2x";"T3" **\crv{(-20,10) & (5,10) }
        \POS?(.45)*{\hole}="3y" \POS?(.7)*{\hole}="3z";
    "2x";"B2" **\crv{(-5,-14)};
        "T1";"3y" **\crv{(-16,17)};
        "T2";"3z" **\crv{(-5,17)};
        "3z";"2z" **\crv{};
        "3y";"2y" **\crv{};
        "B3";"2z" **\crv{ (5,-5) &(20,-10)}
            \POS?(.4)*{\hole}="4z";
        "2y";"4z" **\crv{(6,-8)};
        "4z";"B4" **\crv{(15,-15)};
 \endxy
    }="ML";
    (-40,-30)*+{
 \xy 
 (-14,20)*{}="T1";
 (-4,20)*{}="T2";
 (4,20)*{}="T3";
 (15,20)*{}="T4";
 (-15,-20)*{}="B1";
 (-5,-20)*{}="B2";
 (5,-20)*{}="B3";
 (15,-20)*{}="B4";
    "B1";"T4" **\crv{(-15,-5) & (15,5)}
        \POS?(.38)*{\hole}="2x" \POS?(.53)*{\hole}="2y" \POS?(.7)*{\hole}="2z";
    "T1";"2x" **\crv{(-15,5)};
    "2y";"B2" **\crv{(10,-10) & (-6,-10)}
        \POS?(.45)*{\hole}="4x";
    "2z";"B3" **\crv{ (15,0)&(15,-10) & (6,-16)}
        \POS?(.7)*{\hole}="5x";
    "T3";"2y" **\crv{(5,10)& (-6,18) }
        \POS?(.5)*{\hole}="3x";
    "T2";"3x" **\crv{(-5,15)};
    "3x";"2z" **\crv{(7,11)};
    "2x";"4x" **\crv{(-3,-7)};
    "4x";"5x" **\crv{};
    "5x";"B4" **\crv{(15,-15)};
 \endxy
    }="BL";
    (40,30)*+{
 \xy 
 (14,-20)*{}="T1";
 (4,-20)*{}="T2";
 (-4,-20)*{}="T3";
 (-15,-20)*{}="T4";
 (15,20)*{}="B1";
 (5,20)*{}="B2";
 (-5,20)*{}="B3";
 (-15,20)*{}="B4";
    "B1";"T4" **\crv{(15,5) & (-15,-5)}
        \POS?(.38)*{\hole}="2x" \POS?(.53)*{\hole}="2y" \POS?(.7)*{\hole}="2z";
    "T1";"2x" **\crv{(15,-5)};
    "2y";"B2" **\crv{(-10,10) & (6,10)}
        \POS?(.45)*{\hole}="4x";
    "2z";"B3" **\crv{ (-15,0)&(-15,10) & (-6,16)}
        \POS?(.7)*{\hole}="5x";
    "T3";"2y" **\crv{(-5,-10)& (6,-18) }
        \POS?(.5)*{\hole}="3x";
    "T2";"3x" **\crv{(5,-15)};
    "3x";"2z" **\crv{(-7,-11)};
    "2x";"4x" **\crv{(3,7)};
    "4x";"5x" **\crv{};
    "5x";"B4" **\crv{(-15,15)};
 \endxy
    }="TR";
    (75,0)*+{
 \xy 
 (14,-20)*{}="T1";
 (4,-20)*{}="T2";
 (-4,-20)*{}="T3";
 (-15,-20)*{}="T4";
 (15,20)*{}="B1";
 (5,20)*{}="B2";
 (-5,20)*{}="B3";
 (-15,20)*{}="B4";
    "B1";"T4" **\crv{(15,-5) & (-15,5)}
        \POS?(.25)*{\hole}="2x" \POS?(.49)*{\hole}="2y" \POS?(.65)*{\hole}="2z";
    "2x";"T3" **\crv{(20,-10) & (-5,-10) }
        \POS?(.45)*{\hole}="3y" \POS?(.7)*{\hole}="3z";
    "2x";"B2" **\crv{(5,14)};
        "T1";"3y" **\crv{(16,-17)};
        "T2";"3z" **\crv{(5,-17)};
        "3z";"2z" **\crv{};
        "3y";"2y" **\crv{};
        "B3";"2z" **\crv{ (-5,5) &(-20,10)}
            \POS?(.4)*{\hole}="4z";
        "2y";"4z" **\crv{(-6,8)};
        "4z";"B4" **\crv{(-15,15)};
 \endxy
    }="MR";
    (40,-30)*+{
 \xy 
 (15,20)*{}="b1";
 (5,20)*{}="b2";
 (-5,20)*{}="b3";
 (-14,20)*{}="b4";
 (14,-20)*{}="T1";
 (5,-20)*{}="T2";
 (-5,-20)*{}="T3";
 (-15,-20)*{}="T4";
    "b1"; "T4" **\crv{(15,7) & (-15,5)}
        \POS?(.25)*{\hole}="2x" \POS?(.47)*{\hole}="2y" \POS?(.65)*{\hole}="2z";
    "b2";"2x" **\crv{(5,15)};
    "b3";"2y" **\crv{(-5,10)};
    "b4";"2z" **\crv{(-14,9)};
 (15,5)*{}="3x";
    "2x"; "3x" **\crv{(15,10)};
    "3x"; "T3" **\crv{(15,-15) & (-5,-10)}
        \POS?(.38)*{\hole}="4y" \POS?(.65)*{\hole}="4z";
    "T1";"4y" **\crv{(14,-16)};
    "T2";"4z" **\crv{(5,-16)};
    "2y";"4z" **\crv{(10,-3) & (-10,-2)} \POS?(.6)*{\hole}="5z";
    "4y";"5z" **\crv{(5,-5)};
    "5z";"2z" **\crv{(-5,-4)};
 \endxy
    }="BR";
    (0,-50)*+{
 \xy 
 (15,20)*{}="T1";
 (5,20)*{}="T2";
 (-5,20)*{}="T3";
 (-15,20)*{}="T4";
 (15,-20)*{}="B1";
 (5,-20)*{}="B2";
 (-5,-20)*{}="B3";
 (-15,-20)*{}="B4";
    "T1"; "B4" **\crv{(15,7) & (-15,5)}
        \POS?(.25)*{\hole}="2x" \POS?(.45)*{\hole}="2y" \POS?(.6)*{\hole}="2z";
    "T2";"2x" **\crv{(4,12)};
    "T3";"2y" **\crv{(-5,10)};
    "T4";"2z" **\crv{(-16,9)};
 (15,5)*{}="3x";
    "2x"; "3x" **\crv{(18,10)};
    "3x"; "B3" **\crv{(13,0) & (-4,-10)}
        \POS?(.3)*{\hole}="4x" \POS?(.53)*{\hole}="4y";
    "2y"; "4x" **\crv{};
    "2z"; "4y" **\crv{};
 (15,-10)*{}="5x";
    "4x";"5x" **\crv{(17,-6)};
    "5x";"B2" **\crv{(14,-12)}
        \POS?(.6)*{\hole}="6x";
    "6x";"B1" **\crv{};
    "4y";"6x" **\crv{};
 \endxy
    }="B";
            (-20,65)*{}="X1";
            (-35,55)*{}="X2";               
                {\ar@{=>} "X1";"X2"};       
            (20,65)*{}="X1";                
            (35,55)*{}="X2";                
                {\ar@{=>} "X1";"X2"};       
            (60,40)*{}="X1";                
            (75,25)*{}="X2";                
                {\ar@{=>} "X1";"X2"};       
            (-60,40)*{}="X1";               
            (-75,25)*{}="X2";               
                {\ar@{=>} "X1";"X2"};
            (-60,-40)*{}="X2";
            (-75,-25)*{}="X1";
                {\ar@{=>} "X1";"X2"};
            (60,-40)*{}="X2";
            (75,-25)*{}="X1";
                {\ar@{=>} "X1";"X2"};
            (-20,-65)*{}="X2";
            (-35,-55)*{}="X1";
                {\ar@{=>} "X1";"X2"};
            (20,-65)*{}="X2";
            (35,-55)*{}="X1";
                {\ar@{=>} "X1";"X2"};
  \endxy
\]
\noindent in a $2$-category whose $2$-morphisms are isotopies of
surfaces in 4-space --- or more precisely, `$2$-braids'. Details
can be found in HDA1, HDA4 and a number of other references, going
back to the work of Kapranov and Voevodsky
\cite{BLan,BN,CaS,Crans,KV}.

As suggested above, we now see a relationship between the
algebraic version of the coherence law for the left distributor
and the Zamolodchikov tetrahedron equation!  They are both
octagons, which turn out to be equivalent in the following
context:

\begin{thm} \et  Let $Q$ be a category equipped with a functor
$\rhd \maps Q \times Q \to Q$ and a natural isomorphism
$$LD_{x,y,z} \maps x \rhd (y \rhd z) \to (x \rhd y) \rhd (x \rhd
z).$$ Define $B \maps Q \times Q \to Q \times Q$ by
$$B(x,y) = (y,y \rhd x)$$ whenever $x$ and $y$ are
both either objects or morphisms in $Q$. Let
$$Y \maps (B \times 1)(1 \times B)(B \times 1) \To (1 \times
B)(B \times 1)(1 \times B)$$ be defined as:
$$Y_{x,y,z} = 1_{z} \times 1_{z \rhd y} \times LD_{z,y,x}.$$  Then
$Y$ is a solution of the Zamolodchikov tetrahedron equation if and
only if $(Q, \rhd, LD)$ is a left $2$-shelf.
\end{thm}

\noindent {\bf Proof.}  Applying the left-hand side of the
Zamolodchikov tetrahedron equation to an object $(w,x,y,z)$ of $Q
\times Q \times Q \times Q$ produces an expression consisting of
various uninteresting terms together with one involving
$$(z \rhd LD_{y,x,w}) LD_{z, y \rhd x, y \rhd w}
((z \rhd ( y \rhd x)) \rhd LD_{z,y,w})  (LD_{z,y,x} \rhd ((z \rhd
y) \rhd (z \rhd w)))$$ while applying the right-hand side gives an
expression with the same uninteresting terms, but also one with
$$LD_{z,y, x \rhd w} ((z \rhd y) \rhd LD_{z,x,w})  LD_{z \rhd y, z \rhd
x, z \rhd w}$$ in exactly the same way.  Thus, the two sides are
equal if and only if the left distributor identity holds.  \qed

Therefore, we can think of the left distributor and its coherence
law from a completely topological viewpoint!

The definition of a right $2$-shelf is similar, but uses the
versions of the third Reidemeister move and Zamolodchikov equation
having left-handed crossings.   Thus far, this topological
connection has proven to be a useful guiding light in formulating
the correct definitions of $2$-shelves.  As we move to racks and
quandles, however, the situation becomes significantly more
complex.  The presence of the \emph{two} conjugation operations
complicates matters since each operation corresponds to its own
type of crossing.  We have seen that left conjugation corresponds
to a positive, or right-handed, crossings, while right conjugation
corresponds to a negative, or left-handed, crossings.

Thus, a $2$-rack will consist of a left and right $2$-shelf
equipped with two inverse natural isomorphisms
$$L_{x,y}
\maps x \rhd (y \lhd x) \to x \qquad R_{x,y} \maps (y \rhd x) \lhd
y \to x$$ which we can draw as:
\[
 \begin{xy}
\vtwist~{(-5,5)}{(5,5)}{(-5,-5)}{(5,-5)};
\vcross~{(-5,-5)}{(5,-5)}{(-5,-15)}{(5,-15)};
    (8,-5)*{}="3";
    (25,-5)*{}="4";
      {\ar^{~} "3";"4"};
    (30,5)*{}="A";
    (30,-15)*{}="A'";
    (35,5)*{}="B";
    (35,-15)*{}="B'";
      {\ar@{-} "A";"A'"};
      {\ar@{-} "B";"B'"};
    (60,-5)*{R_{x,y} \maps (y \rhd x) \lhd y \Rightarrow x};
 \end{xy}
\]

\[
 \begin{xy}
\vcross~{(-5,5)}{(5,5)}{(-5,-5)}{(5,-5)};
\vtwist~{(-5,-5)}{(5,-5)}{(-5,-15)}{(5,-15)};
    (8,-5)*{}="3";
    (25,-5)*{}="4";
      {\ar^{~} "3";"4"};
    (30,5)*{}="A";
    (30,-15)*{}="A'";
    (35,5)*{}="B";
    (35,-15)*{}="B'";
      {\ar@{-} "A";"A'"};
      {\ar@{-} "B";"B'"};
    (60,-5)*{L_{x,y} \maps x \rhd (y \lhd x) \Rightarrow y};
 \end{xy}
\]
that are required to satisfy additional coherence laws. Due to the
presence of the two conjugation functors corresponding to two
different crossing changes, describing all of these coherence laws
is no easy task.  This is a consequence of the fact that these two
conjugation functors allow for $2^3 =8$ possible ways of
performing the third Reidemeister move, corresponding to the two
choices of crossing for each of the three crossings of one side of
the third Reidemeister move.  Thus, a priori, we should expect $8$
different versions of $Y$.  However, drawing pictures of each of
these higher Yang--Baxter operators reveals that $2$ are
impossible, since the strands become too tangled up.  That is,
there is no version of the Yang--Baxter equation with either $(B
\times 1)(1 \times B^{-1})(B \times 1)$ or $(B^{1-} \times 1)(1
\times B)(B^{-1} \times 1)$ on the left-hand side, regardless of
what you try for the right-hand side. These six versions of $Y$,
then suggest at most $2^6 = 64$ versions of the Zamolodchikov
tetrahedron equation corresponding to the two choices of crossing
for each of the six crossings in the top braid of the equation.
Fortunately, we have shown that, of these six, there is only one
independent $Y$.  That is, we are able to derive any of the other
five from the $Y$ related to a left $2$-shelf described above.
However, this does not mean that we have now reduced the number of
independent Zamolodchikov tetrahedron equations down to one.  We
certainly can express each of the Zamolodchikov tetrahedron
equations in terms of our favorite $Y$, but it is not immediately
obvious how many remain.

In addition to the Zamolodchikov tetrahedron equation, we
anticipate the natural isomorphisms of a $2$-rack to satisfy
coherence laws relating the untanglers to the distributors.  That
is, we expect to require the commutativity of various versions of
the following two diagrams:

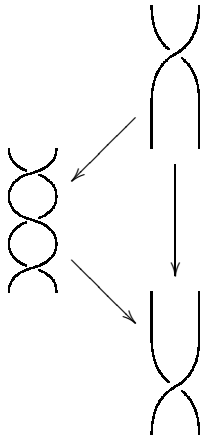
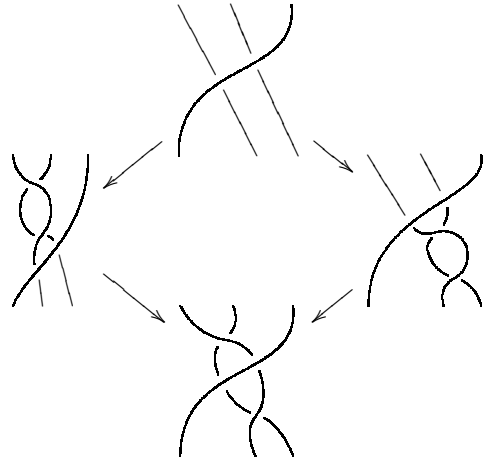
\begin{itemize}
\item[(i)] \[
 \xy 0;/r.15pc/:
 (0,30)*++{ \xy
  \vtwist~{(-5,10)}{(5,10)}{(-5,-10)}{(5,-10)};
  (-5,-10)*{}; (-5,-20)*{} **\dir{-};
  (5,-10)*{}; (5,-20)*{} **\dir{-};
 \endxy}="t";
 (-30,0)*++{\xy
  \vtwist~{(-5,15)}{(5,15)}{(-5,5)}{(5,5)};
   \vcross~{(-5,5)}{(5,5)}{(-5,-5)}{(5,-5)};
   \vtwist~{(-5,-5)}{(5,-5)}{(-5,-15)}{(5,-15)};
 \endxy}="l";
 (0,-30)*++{\xy
  \vtwist~{(-5,0)}{(5,0)}{(-5,-20)}{(5,-20)};
  (-5,10)*{}; (-5,0)*{} **\dir{-};
  (5,10)*{}; (5,0)*{} **\dir{-};
 \endxy}="b";
  {\ar "t";"l"};
  {\ar "t";"b"};
  {\ar "l";"b"};
 \endxy
\]
\item[(ii)] \[
 \xy
 (0,20)*++{ \xy
  (5,10)*{}; (-10,-10)*{} **\crv{(6,0)&(-10,3)}
   \POS?(.4)*+{}="x" \POS?(.6)*+{}="y";
   (-3,10)*{}; "x" **\crv{};
   (-10,10)*{}; "y" **\crv{};
    "y"+(.5,-1); (.5,-10)*{}  **\dir{-};
    "x"+(.5,-1); (6,-10)*{}  **\dir{-};
 \endxy}="t";
 (-25,0)*++{ \xy
    (0,10)*{}; (-10,-10)*{}
     **\crv{(0,-1)&(-8,-5)}
     \POS?(.5)*+{}="x" \POS?(.7)*+{}="y";
     \vcross~{(-10,10)}{(-5,10)}{(-9,2)}{(-5,3)};
     \vtwist~{(-9,2)}{(-5,3)}{"y"}{"x"};
     "y"+(.5,-1); (-6,-10)*{}  **\dir{-};
    "x"+(.5,-1); (-2,-10)*{}  **\dir{-};
 \endxy}="l";
 (25,0)*++{ \xy
  (5,10)*{}; (-10,-10)*{} **\crv{(6,4)&(-10,3)}
   \POS?(.4)*+{}="x" \POS?(.6)*+{}="y";
   (-3,10)*{}; "x" **\crv{};
   (-10,10)*{}; "y" **\crv{};
   \vcross~{"y"+(.5,-.5)}{"x"}{(-2,-3)}{(3,-2)};
   \vtwist~{(-2,-3)}{(3,-2)}{(0,-10)}{(5,-10)};
 \endxy}="r";
 (0,-20)*++{ \xy
  (5,10)*{}; (-10,-10)*{} **\crv{(6,0)&(-10,3)}
   \POS?(.4)*{\hole}="x" \POS?(.6)*{\hole}="y";
   \vcross~{(-10,10)}{(-3,10)}{"y"}{"x"};
   \vtwist~{"y"}{"x"}{(0,-10)}{(5,-10)};
 \endxy}="b";
 {\ar "t";"l"};
  {\ar "t";"r"};
  {\ar "r";"b"};
  {\ar "l";"b"};
 \endxy
\]
\end{itemize}
We have not yet determined the number of independent coherence
laws which arise from using both left and right-handed crossings
in the two diagrams above.

Since we have already begun to wade into a sea of uncertainty, and
since $2$-quandles are equipped with still more natural
isomorphisms, we will end our discussion of the categorifications
of these four algebraic objects.  Before we move on, however, we
remark that we expect that the coherent $2$-groups described in
HDA5 \cite{BLau} will give examples of $2$-quandles just as groups
provide examples of quandles. Moreover, just as Lie algebra
cohomology was used to classify semistrict Lie $2$-algebras in
Section \ref{skeletallie2algs}, we anticipate that quandle and
rack cohomology \cite{CJKLS, EG} will serve to classify
$2$-quandles and $2$-racks. Finally, we hope to prove a
categorified version of Theorem \ref{squares} from Section
\ref{braids}, since the categorified versions of the braid and
framed braid groups and monoids are known.  This categorified
result will give an invariant of $2$-braids.

Ultimately, we desire to demonstrate that every Lie $2$-group has
a Lie $2$-algebra, using a categorified version of the process
described in Section \ref{liealgofliegrp}.  Recall that this
process was described by the following diagram:
$$
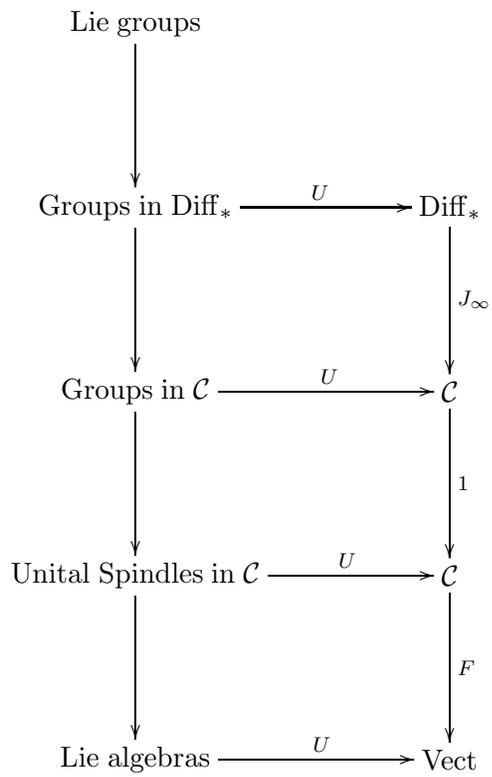
\xymatrix{
   {\rm Lie \; groups}
   \ar[dd] \\ \\
   {\rm Groups \; in \; \Diff_\ast}
   \ar[rr]^<<<<<<<<<<<{U}
   \ar[dd]
   && \Diff_{\ast}
   \ar[dd]^{J_{\infty}} \\ \\
   {\rm Groups \; in \; \mathcal{C}}
   \ar[rr]^<<<<<<<<<<<<<<<{U}
   \ar[dd]
   && \mathcal{C}
   \ar[dd]^{1} \\ \\
   {\rm Unital \; Spindles \; in \; \mathcal{C}}
   \ar[rr]^<<<<<<<<<<<{U}
   \ar[dd]
   && \mathcal{C}
   \ar[dd]^{F} \\ \\
   {\rm Lie \; algebras}
   \ar[rr]^<<<<<<<<<<<<<<{U}
   && \Vect
}$$ One might hope that the passage from Lie $2$-group to Lie
$2$-algebra will be described by a diagram in which we have
categorified everything in sight:
$$
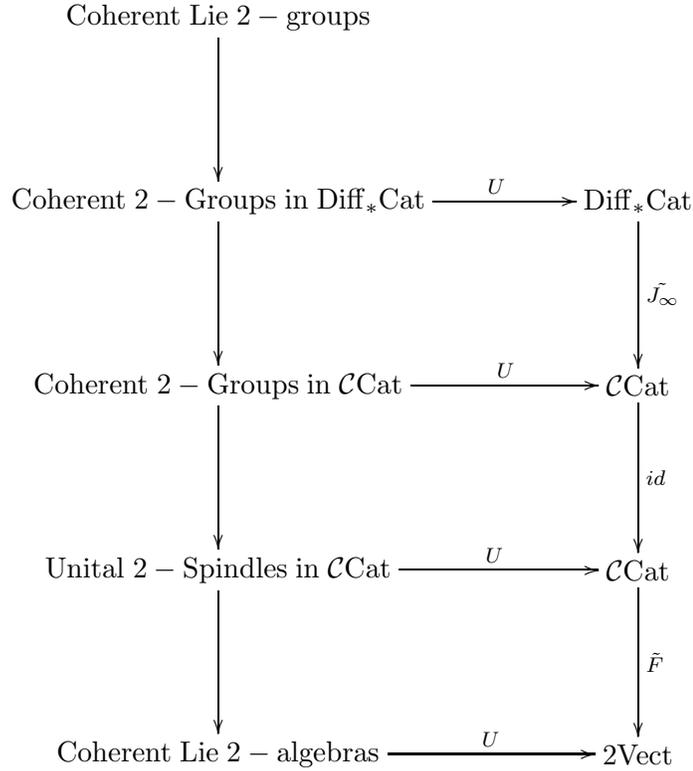
\xymatrix{
   {\rm Coherent \; Lie \; 2-groups}
   \ar[dd] \\ \\
   {\rm Coherent \; 2-Groups \; in \; \Diff_{\ast} \Cat}
   \ar[rr]^<<<<<<<<<{U}
   \ar[dd]
   && \Diff_{\ast}\Cat
   \ar[dd]^{\tilde{J_{\infty}}} \\ \\
   {\rm Coherent \; 2-Groups \; in \; \mathcal{C}\Cat}
   \ar[rr]^<<<<<<<<<<<<<{U}
   \ar[dd]
   && \mathcal{C}\Cat
   \ar[dd]^{id} \\ \\
   {\rm Unital \; 2-Spindles \; in \; \mathcal{C}\Cat }
   \ar[rr]^<<<<<<<<<<<<<{U}
   \ar[dd]
   && \mathcal{C}\Cat
   \ar[dd]^{\tilde{F}} \\ \\
   {\rm Coherent \; Lie \; 2-algebras}
   \ar[rr]^<<<<<<<<<<<<<<{U}
   && 2\Vect
}$$ The good news is that we already have a definition of a
coherent $2$-group in $K$, where $K$ is a $2$-category with finite
products, which then gives us the notion of coherent Lie
$2$-group.  Moreover, we know what a $2$-vector space is, and that
$\Diff \Cat$ and $K \Cat$ have finite products since $\Diff$ and
$K$ do.  This is promising because it means that the notion of a
$2$-spindle in one of these $2$-categories makes sense.

Unfortunately, it seems the bad news outweighs the good.  As
mentioned at the end of Section \ref{strictlie2algs} of Chapter
\ref{ch1}, this diagram will probably require what we call
`coherent' Lie $2$-algebras, in which the antisymmetry law of the
bracket is replaced by a natural isomorphism which we call the
`antisymmetrizer', but where bilinearity continues to hold on the
nose. We suspect that such Lie $2$-algebras are necessary because
there are various laws that hold only up to isomorphism in a
coherent Lie $2$-group, most likely exactly the laws that we need
to get skew-symmetry of the bracket via differentiation.

We have already begun to investigate the result of weakening the
skew-symmetry condition by introducing an additional natural
transformation into the definition of a semistrict Lie $2$-algebra
called the {\bf antisymmetrizer:}
$$A_{x,y} \maps [x,y] \to -[y,x],$$
which, together with its coherence laws will give the definition
of a `coherent' Lie $2$-algebra.   Thus far we have found four
coherence laws for the antisymmetrizer, three of which relate it
to the Jacobiator.   Though we do not, yet, have a complete
definition of a $2$-quandle, we know in principle what it should
be like, and we have determined some of its coherence laws. At
this point, it is interesting to remark that we were only able to
derive two of the coherence laws for a coherent Lie $2$-algebra
from two of the known coherence laws of a $2$-quandle. Therefore,
the other two laws, including
$$ \xy
   (0,20)*+{[[x,y],z]}="1";
   (-30,0)*+{[x,[y,z]] + [[x,z],y]}="2";
   (30,0)*+{[x,[y,z]] + [x,[z,y]] + [[x,y],z]}="3";
   {\ar_{J_{x,y,z}}     "1";"2"}
   {\ar_<<<<<<<<<{J_{x,z,y}}     "2";"3"}
   {\ar_{[x,A_{y,z}]}   "3";"1"}
\endxy
\\ \\
$$
appear to be unrelated to topology, in particular, to the
$2$-braid diagrams.  If so, this suggests that our strategy for
passing from coherent Lie $2$-groups to coherent Lie $2$-algebras
via $2$-spindles may be too naive.  We hope that after this
categorified version is understood, we will be able to show that a
coherent Lie $2$-algebra gives a $2$-spindle in $\mathcal{C} \Cat$
in analogy to Theorem \ref{punchlinethm}.

We conclude with evidence supporting our conjectures relating
categorified Lie theory to higher-dimensional topology, as a way
of showing that our guesses are not too far afield. We show that,
in a suitable context, the Jacobiator identity is equivalent to
the Zamolodchikov tetrahedron equation!  However, we need a
different version of the Zamolodchikov tetrahedron equation than
the one given earlier in this section.  This is because we defined
a $2$-shelf in the $2$-category $\Cat$, which is a monoidal
$2$-category where the tensor product is the Cartesian product,
whereas the definition of a Lie $2$-algebra lives in the
$2$-category $2\Vect$, which is a monoidal $2$-category with the
usual tensor product.  Thus, we need the following:

\begin{defn} \et Given a $2$-vector space $V$ and an
invertible linear functor $B \maps V \otimes V \to V \otimes V$, a
linear natural isomorphism
$$Y \maps (B \otimes 1)(1 \otimes B)(B \otimes 1) \To
          (1 \otimes B)(B \otimes 1)(1 \otimes B) $$
satisfies the {\bf Zamolodchikov tetrahedron equation} if
$$[Y \circ (1 \ltimes 1 \ltimes B)(1 \ltimes B \ltimes 1)(B \ltimes 1
\ltimes 1)] [(1 \ltimes B \ltimes 1)(B \ltimes 1 \ltimes 1) \circ
Y \circ (B \ltimes 1 \ltimes 1)]$$ $$[(1 \ltimes B \ltimes 1)(1
\ltimes 1 \ltimes B) \circ Y \circ (1 \ltimes 1 \ltimes B)] [Y
\circ (B \ltimes 1 \ltimes 1)(1 \ltimes B \ltimes 1)(1 \ltimes 1
\ltimes B)]$$
$$ = $$
$$[(B \ltimes 1 \ltimes 1)(1 \ltimes B \ltimes 1)(1 \ltimes 1 \ltimes B)
\circ Y ] [(B \ltimes 1 \ltimes 1) \circ Y \circ (B \ltimes 1
\ltimes 1)(1 \ltimes B \ltimes 1)]$$
$$[(1 \ltimes 1 \ltimes B) \circ Y \circ (1 \ltimes 1 \ltimes
B)(1 \ltimes B \ltimes 1)] [(1 \ltimes 1 \ltimes B)(1 \ltimes B
\ltimes 1)(B \ltimes 1 \ltimes 1) \circ Y],$$ where $\circ$
represents the whiskering of a linear functor by a linear natural
transformation.
\end{defn}

Though we have seen that coherent Lie $2$-algebras may not be as
closely related to topology as we would like, it is the case,
however, that semistrict Lie $2$-algebras have a connection to
higher-braid theory.  Recall from Section \ref{definitions} that
the coherence law for the Jacobiator, the Jacobiator identity in
Definition \ref{defnlie2alg}, seems rather arcane.  It turns out
to be related to the Zamolodchikov tetrahedron equation in
$2\Vect$ in the same way that the Jacobi identity is related to
the Yang--Baxter equation.  That is, Proposition \ref{YBEJacobi},
which roughly stated that a Lie algebra gives a solution of the
Yang--Baxter equation, has a higher-dimensional analog, obtained
by categorifying everything in sight!   Recall in the definition
of a semistrict Lie $2$-algebra, Definition \ref{defnlie2alg}, we
clarified the Jacobiator identity by drawing it as a commutative
octagon.  In fact, that commutative octagon becomes {\it
equivalent} to the octagon for the Zamolodchikov tetrahedron
equation given earlier in this section in the following context:

\begin{thm} \et \label{Zameqn} Let $L$ be a $2$-vector space,
let $[\cdot, \cdot] \maps L \times L \to L$ be a skew-symmetric
bilinear functor, and let $J$ be a completely antisymmetric
trilinear natural transformation with $J_{x,y,z}\maps [[x,y],z]
\to [x,[y,z]] + [[x,z],y]$.  Let $L' = K \oplus L$, where $K$ is
the categorified ground field. Let $B \maps L' \otimes L' \to L'
\otimes L'$ be defined as follows:
$$B((a,x) \otimes (b,y))
= (b,y) \otimes (a,x) + (1,0) \otimes (0, [x,y])$$ whenever
$(a,x)$ and $(b,y)$ are both either objects or morphisms in $L'$.
Finally, let
$$Y \maps (B \otimes 1)(1 \otimes B)(B \otimes 1) \To
          (1 \otimes B)(B \otimes 1)(1 \otimes B) $$
be defined as follows:
$$Y = (p \otimes p \otimes p) \circ J \circ j$$
where $p \maps L' \to L$ is the projection functor given by the
fact that $L' = K \oplus L$ and
$$j \maps L \rightarrow L' \otimes L' \otimes L'$$
is the linear functor defined by
$$j(x) = (1,0) \otimes (1,0) \otimes (0,x),$$
where $x$ is either an object or morphism of $L$. Then $Y$ is a
solution of the Zamolodchikov tetrahedron equation if and only if
$J$ satisfies the Jacobiator identity.
\end{thm}

\noindent{\bf Proof. }  Equivalently, we must show that $Y$
satisfies the Zamolodchikov tetrahedron equation if and only if
$J$ satisfies the Jacobiator identity.  Applying the left-hand
side of the Zamolodchikov tetrahedron equation to an object $(a,w)
\otimes (b,x) \otimes (c,y) \otimes (d,z)$ of $L' \otimes L'
\otimes L' \otimes L'$ yields an expression consisting of various
uninteresting terms together with one involving
$$J_{[w,x],y,z} ([J_{w,x,z},y]+1) (J_{w, [x,z], y} +
       J_{[w,z],x,y} + J_{w,x, [y,z]}), $$
while applying the right-hand side produces an expression with the
same uninteresting terms, but also one involving
$$[J_{w,x,y},z] (J_{[w,y],x,z} + J_{w, [x,y],z}) ([J_{w,y,z},x]+1)
([w, J_{x,y,z}]+1) $$ in precisely the same way.  Thus, the two
sides are equal if and only if the Jacobiator identity holds.  The
detailed calculation resembles that of the proof of Proposition
\ref{YBEJacobi} and is quite lengthy. \qed

\begin{cor} \et  If $L$ is a semistrict Lie $2$-algebra,
then $Y$ defined as in Theorem \ref{Zameqn} is a solution of the
Zamolodchikov tetrahedron equation.
\end{cor}

This result could be just the beginning of a beautiful
relationship between categorified Lie theory and
higher-dimensional braids --- a relationship that extends to
higher and higher dimensions.


\end{document}